\documentclass{article}
\ifdefined\pdfminorversion\pdfminorversion=7\fi

\usepackage[a4paper,top=2.5cm,bottom=2.5cm,left=2.5cm,right=2.5cm]{geometry}
\usepackage{amsmath,amsfonts,amssymb,amsthm}
\usepackage[svgnames]{xcolor}
\usepackage{natbib}
\usepackage{enumerate}
\usepackage{array}
\usepackage{booktabs,longtable}
\usepackage{graphicx,placeins}
\usepackage{tikz}
\usepackage[normalem]{ulem}
\newif\ifinlinesupp
\inlinesupptrue
\ifinlinesupp\else
\usepackage{zref-xr}
\zxrsetup{toltxlabel}
\zexternaldocument*{supplement}
\fi

\usepackage{hyperref}
\hypersetup{
    colorlinks,
    citecolor=SeaGreen,
    filecolor=black,
    linkcolor=Indigo,
    urlcolor=black
}

\newcommand{\E}{\mathbb{E}}
\renewcommand{\P}{\mathbb{P}}
\renewcommand{\d}{\mathrm{d}}
\newcommand{\var}{\mathrm{Var}}
\newcommand{\cov}{\mathrm{Cov}}

\newcommand{\one}{\mathbf{1}}
\newcommand{\Pcal}{\mathcal{P}}

\newcommand{\Nor}{\mathcal{N}}
\newcommand{\dto}{\xrightarrow{d}}
\newcommand{\Vrcj}{\widehat V_{\mathrm{RCJ}}}
\newcommand{\Vu}{\widehat V_{\mathrm{u}}}
\newcommand{\Thetah}{\widehat\Theta}
\newcommand{\wbar}{\bar w}
\newcommand{\disc}{\mathrm{disc}}
\newcommand{\conn}{\mathrm{conn}}

\newtheorem{theorem}{Theorem}
\newtheorem{proposition}[theorem]{Proposition}
\newtheorem{lemma}[theorem]{Lemma}
\newtheorem{corollary}[theorem]{Corollary}
\newtheorem{lemmaL}{Lemma}
\renewcommand{\thelemmaL}{\Alph{lemmaL}}
\theoremstyle{definition}
\newtheorem{definition}{Definition}

\newtheorem{remark}{Remark}
\newtheorem{example}{Example}

\newcommand{\regimegraph}[1]{\tikz[baseline=-.5ex,x=3mm,y=3mm,
row/.style={circle,draw,fill=white,inner sep=0pt,minimum size=2mm},
col/.style={rectangle,draw,fill=gray!25,inner sep=0pt,minimum size=2mm}]{
\ifcase#1
\node[row] at (0,0) {};
\or\node[col] at (0,0) {};
\or\draw (0,0)--(1,0);\node[row] at (0,0) {};\node[col] at (1,0) {};
\or\node[row] at (0,0) {};\node[row] at (1,0) {};
\or\node[col] at (0,0) {};\node[col] at (1,0) {};
\or\node[row] at (0,0) {};\node[col] at (1,0) {};
\or\draw (0,0)--(1,0)--(1,1)--(0,1)--cycle;
\node[row] at (0,0) {};\node[col] at (1,0) {};\node[row] at (1,1) {};\node[col] at (0,1) {};
\fi}}

\title{Uniform inference for $U$-statistics\\ of exchangeable arrays under degeneracy}
\author{T\^am \textsc{Le Minh}}
\date{}
\begin{document}
\maketitle
\begin{abstract}
Many $U$-statistics of row-column exchangeable arrays become degenerate at natural null hypotheses: their convergence rates change and their limiting laws may or may not remain Gaussian. We establish conditions under which inference based on the classical unbiased variance estimator is uniformly valid across these changes, without knowledge of the degeneracy order. It is known that $U$-statistics can be decomposed into orthogonal projections indexed by bipartite graphs. Under uniform fourth-moment assumptions, we establish that the standardized statistic is asymptotically Gaussian and the classical unbiased variance estimator is ratio-consistent, as long as a certain measure of concentration for projections indexed by disconnected graphs tends to zero. Under compactness of the normalized projections, both properties hold if and only if the proportion of variance contributed by disconnected graphs tends to zero. The necessity of this condition follows from a general limit theorem for these $U$-statistics, in which each graph contributes to a Wiener chaos of order equal to its number of connected components. We note that the row-column jackknife variance estimator, although nonnegative, does not generally give a studentized statistic with a standard normal limit in degenerate regimes. We illustrate our theory with an application to estimation and testing based on a distance from product form, exhibiting Gaussian and non-Gaussian limits at various rates depending on the model. We also propose a conservative inference procedure that does not require a Gaussian limit, valid under additional moment and projection conditions.

\end{abstract}
\section{Introduction}\label{sec:intro}

A $U$-statistic on a two-way array $(Y_{ij})_{i \ge 1, j \ge 1}$ is a statistic of the form
\begin{equation}\label{eq:ustat}
  U_N=\binom mp^{-1}\binom nq^{-1}\sum_WH_W, \qquad H_W:=h(Y_W),
\end{equation}
where the sum is over the pairs $W=(W_1,W_2)$ of a row set $W_1\subseteq[m]$ and a column set $W_2\subseteq[n]$, with $|W_1|=p$ and $|W_2|=q$, $Y_W$ is the subarray induced by the indices of $W$, and the kernel $h$ is a function of a $p\times q$ subarray, symmetric under separate permutations of the rows and columns. We assume that the array is row-column exchangeable (RCE), that means, its law is invariant under separate permutations of its rows and its columns, and dissociated, that means, subarrays sharing neither rows nor columns are independent. The study of RCE matrices is motivated by bipartite network analysis (for example, recommender systems \citep{zhou2007bipartite}, scientific authorship networks \citep{newman2001structure}, ecological pollination networks \citep{dormann2009indices}), where the RCE property implies vertex-exchangeability.

$U$-statistics are a large class of statistics on these arrays, yet, for a wide range of null models, U-statistics are degenerate in that their first-order Hoeffding projection is zero \citep{leminh2025hoeffdingtype}. In this case, first-order variance estimators are no longer suitable to describe the actual variance of the $U$-statistics, the rate of convergence might be different, and the limiting distribution might be Gaussian or non-Gaussian depending on the situation \citep{leminh2026characterization}. Thus, a test that is asymptotically valid at each fixed law does not necessarily remain valid when the law approaches such a null hypothesis. This paper characterizes the region under which \textit{Gaussian inference}, that is, tests and confidence intervals based on Gaussian critical values, is uniformly valid, and shows that there is a boundary beyond which it fails.

\subsection{Contributions}

\paragraph{Main result: Uniform inference and its boundary.}
We establish conditions under which Gaussian confidence intervals based on the classical unbiased variance estimator is uniformly valid across degeneracy regimes. More precisely, under uniform fourth-moment bounds and a condition on the graph projections, Theorem~\ref{thm:contraction} states that the standardized U-statistic is asymptotically Gaussian, and the classic unbiased estimator of variance is ratio-consistent. Assuming compactness of the normalized projections, a necessary-and-sufficient condition is derived: the two properties are satisfied if and only if the fraction of variance due to disconnected graphs is asymptotically vanishing (Thm.~\ref{thm:sharp} and Cor.~\ref{cor:S-consistency}). 

\paragraph{Limit distributions underlying the boundary.}
To determine when the standardized statistic can have a Gaussian limit, we establish a general limit theorem under square integrability: each principal support graph contributes to a Wiener chaos of order equal to its number of connected components (Thm.~\ref{thm:chaos}). This is the RCE analogue of \citet[Thm.~3]{janson1991asymptotic} and resolves the multi-component fixed-law case left open in \cite{leminh2026characterization}. Under compactness, a nonvanishing contribution from disconnected graphs then prevents the standardized statistic from being asymptotically Gaussian (Thm.~\ref{thm:sharp}). 

\paragraph{Jackknife bias and correction.}
Even when the statistic is asymptotically Gaussian, the row-column jackknife does not generally yield valid Gaussian inference because it is asymptotically biased. We show that its expectation reweights each variance contribution by a factor tending to its number of vertices, irrespective of its edges. Then, we establish an identity showing that the classical unbiased variance estimator can be interpreted as a bias-corrected jackknife that removes these excess weights exactly.

\paragraph{Conservative inference without a Gaussian limit.}
At a null hypothesis with a known structure, the conditions for Gaussian inference can be verified analytically, but they do not hold in general and cannot be verified from a single observed array. We therefore also construct a conservative test using independent diagonal blocks. Under uniform fourth-moment bounds, this test remains valid even when these conditions fail, although its confidence intervals are wider. Under additional moment and projection conditions, a statistic computed from the data determines whether to use this conservative test or the Gaussian test. For any fixed \(\alpha\in(0,1)\), the resulting procedure has a rejection probability asymptotically at most \(\alpha\), uniformly over the null class, and selects the Gaussian test when the quantity used for the diagnostic tends to zero sufficiently quickly.

\paragraph{Applications across degeneracy regimes.}
We apply the theory to estimation and testing based on a measure of departure from the product-form model, where the expected array entries can be factorized into a row effect and a column effect. Assuming conditionally Poisson entries, we use a ratio statistic that is degenerate under the product-form null hypothesis but has a Gaussian limit. Inference near the null hypothesis requires both the uniform theorem and a nonlinear remainder negligible compared with the standard deviation. We also consider two variations of this problem: a departure from product form that remains degenerate but has a non-Gaussian limit, and a statistic allowing Gaussian inference without the Poisson assumption, at a different rate. Table~\ref{tab:regimes} summarizes these regimes: neither the convergence rate nor the degeneracy order determines whether the limit is Gaussian. With our theoretical results, Gaussian inference based on the unbiased variance estimator applies in all Gaussian regimes using the corresponding kernel or linearization. In simulation studies, we analyze the coverage and rejection probabilities of the Gaussian test and the conservative procedures across degeneracy regimes. In an application to UN General Assembly roll calls, we assess whether additive effects can account for the departures from product form detected by our test, by fitting the classical Rasch model \citep{rasch1960probabilistic,clinton2004statistical}.

\begin{table}[!ht]
\centering\small
\caption{Product-form regimes at fixed laws with $m,n\asymp N$, under the conditions of Section~\ref{sec:product}. Circles denote rows and squares denote columns.}\label{tab:regimes}
\setlength{\tabcolsep}{7pt}\renewcommand{\arraystretch}{1.1}
\begin{tabular}{@{}lclc@{}}
\toprule
Regime & Principal support graphs & Limit & SD order\\
\midrule
Nondegenerate departure & \regimegraph{0}\; /\;\regimegraph{1} & Gaussian & $N^{-1/2}$\\
Product form & \regimegraph{2} & Gaussian & $N^{-1}$\\
Two-group model & \regimegraph{2}\; /\;\regimegraph{3}\; /\;\regimegraph{4}\; /\;\regimegraph{5} & Non-Gaussian & $N^{-1}$\\
Product form, general cells ($\Psi$) & \regimegraph{6} & Gaussian & $N^{-2}$\\
\bottomrule

\end{tabular}
\end{table}

\subsection{Relation to earlier work}

Our starting point is the theory of RCE $U$-statistics developed by \citet{leminh2023ustatisticsbipartite}, which studied nondegenerate array $U$-statistics, \citet{leminh2025hoeffdingtype}, which developed a coarse decomposition and first-order variance estimation, and \citet{leminh2026characterization}, which supplied the bipartite graph decomposition and Gaussian fixed-law limits. We address uniform inference over changing laws, including changes in degeneracy, and establish the general fixed-law limit theorem to determine its boundary. 

On i.i.d. observations, degenerate $U$-statistics have a non-Gaussian limiting law which lives in a higher Wiener chaos \citep{rubin1980asymptotic}. For generalized $U$-statistics, \citet{janson1991asymptotic} showed that their Hoeffding decomposition can be indexed by graphs, and that, under degeneracy, connected principal support graphs can still produce Gaussian limits. They also obtained Wick-product limits for disconnected principal support graphs. In jointly exchangeable arrays, similar limits arise for subgraph counts \citep{kaur2021higher,bhattacharya2023fluctuations}. 

Unbiased variance estimators have been studied for $U$-statistics of i.i.d. observations \citep{schucany1989small,wang2014variance} and in multi-reader multi-case inference \citep{gallas2006oneshot,gallas2009framework}. \cite{fuchs2020asymptotic} and \cite{xu2024variance} stated consistency results that cover nondegenerate regimes. Our results establish the ratio-consistency of the unbiased variance estimator for general kernels on two-way arrays, even under degeneracy.

\citet{davezies2021empirical} developed empirical-process results for RCE arrays, notably with a chaos limiting law for degenerate dyadic processes. For dyadic density estimation, \citet{graham2024kernel} established Gaussian inference across degenerate and nondegenerate cases. \citet{cattaneo2024uniform} constructed uniform confidence bands under degeneracy. Here, we study RCE $U$-statistics with general kernels, and establish conditions for the standardized statistic to be asymptotically Gaussian and the variance estimator to be ratio-consistent, uniformly over changing laws, which includes potential changes in degeneracy regimes.

For two-way means, \citet{menzel2021bootstrap} established uniform bootstrap guarantees and a conservative inference procedure. \citet{chiang2025multiway}, \citet{davezies2025analytic} and \citet{hounyo2026two} developed similar results for means, GMMs, and quantile regression. However, general array kernels introduce higher-order connected and disconnected graph contributions that do not arise in the two-way mean. \citet{chatterjee2026higherorder} developed a multiplier bootstrap and joint confidence sets for network moments in jointly exchangeable graphs, allowing both Gaussian and non-Gaussian limits.

Our analysis of the behavior of the row-column jackknife also relates to earlier findings that resampling can overcount variance components in crossed data \citep{mccullagh2000resampling,owen2007pigeonhole,owen2012bootstrapping}, and to work on two-way cluster jackknives \citep{mackinnon2024jackknife}. We give the exact component weights responsible for this behavior, with a graph-theoretic interpretation.

\subsection{Outline}
Section~\ref{sec:framework} introduces the general framework. Section~\ref{sec:results} establishes uniform Gaussian inference and its boundary. Section~\ref{sec:jackknife} discusses the jackknife. Section~\ref{sec:hybrid} gives conservative procedures, combining them with Gaussian inference. Section~\ref{sec:product} develops the product-form application, and Section~\ref{sec:numerics} presents simulations and the UN roll-call illustration. The supplementary material contains the proofs and additional numerical results.

\subsection{Notation and convention}\label{sec:notation} $\xrightarrow{d}$ denotes convergence in distribution, and $O_P$, $o_P$ are stochastic orders. When a statement is made uniformly over a class of laws, both are understood along every sequence of laws in the class. For sequences, $a_N\asymp b_N$ means that $a_N/b_N$ is bounded away from $0$ and $\infty$. We write $\Phi$ for the standard normal distribution function and $z_u=\Phi^{-1}(u)$ for its $u$-quantile. For graphs, $\sqcup$ denotes the disjoint union, $\mathrm{Aut}(G)$ the automorphism group of $G$ preserving the row and column vertex sets, and $\sum_{G\,\disc}$ and $\sum_{G\,\conn}$ are sums over disconnected and connected classes (definition in Sec.~\ref{sec:graph-decomposition}). By convention, we set \(\sup\varnothing:=0\). We call a test with rejection event $R_N$ \emph{uniformly valid} at level $\alpha$ on a null class $\Pcal_{0,N}$ if $\limsup_N\sup_{P\in\Pcal_{0,N}}\P_P(R_N)\le\alpha$.
Analogously, a confidence interval $I_N$ for $\theta_{N,P}$ is uniformly valid at level $1-\alpha$ over $\Pcal_N$ if $\liminf_N\inf_{P\in\Pcal_N}\P_P\{\theta_{N,P}\in I_N\}\ge1-\alpha$.

\section{Framework and preliminaries}\label{sec:framework}

\subsection{Row-column exchangeable arrays with drifting laws}\label{sec:framework-arrays}

Let $Y=(Y_{ij})_{i\le m,\,j\le n}$ be the upper-left block of a dissociated RCE array. By the Aldous--Hoover--Kallenberg (AHK) representation \citep{aldous1981representations, hoover1979relations,kallenberg2005probabilistic},
\[
Y_{ij}\overset{a.s.}{=}\phi(\xi_i,\eta_j,\zeta_{ij}),
\]
with $(\xi_i)_{i\ge 1}$, $(\eta_j)_{j\ge 1}$, $(\zeta_{ij})_{i \ge 1, j \ge 1}$ independent families of i.i.d.\ $\mathcal U[0,1]$ variables (the \textit{AHK variables}). Given the row and column variables, the cells are independent. When $\E|Y_{11}|<\infty$, we call $w(x,y)=\E[Y_{ij}\mid\xi_i=x,\eta_j=y]$ the \emph{mean graphon} of the array. Throughout, $m=m_N$, $n=n_N$ with $m_N+n_N=N$ and $m_N/N\to\rho\in(0,1)$, so $m,n\asymp N$.  We also define $L_N=\min(m,n)$.

We allow the law of the array to change with $N$ by replacing $\phi$ with $\phi_N$, while the AHK variables remain independent and uniformly distributed. Classes of laws can also be fixed or depend on $N$. In notation, we omit the indices $N$ or $P$ indicating dependence when there is no ambiguity.

\subsection{The graph-indexed Hoeffding decomposition}\label{sec:graph-decomposition}

Write $\theta_N=\E[U_N]$ and $\sigma_N^2=\var(U_N)$, where $U_N$ is defined in~\eqref{eq:ustat}. Since $h$ is symmetric, the same statistic is obtained by averaging over ordered tuples of distinct rows and columns. We use this when counting graph copies, with $(x)_k=x(x-1)\cdots(x-k+1)$.

\paragraph{Graphs.} 
A bipartite graph $G$ has row vertices $V_1(G)$, column vertices $V_2(G)$ and edges $E(G)\subseteq V_1(G)\times V_2(G)$. Isolated vertices are allowed. We write $r(G)$ and $c(G)$ for its numbers of row and column vertices, $d(G)=r(G)+c(G)$ for its total number of vertices, and $b(G)$ for its number of connected components. A \emph{subgraph} of $G$, written $F\subseteq G$, is a bipartite graph $F$ with $V_1(F)\subseteq V_1(G)$, $V_2(F)\subseteq V_2(G)$ and $E(F)\subseteq E(G)$. It is \emph{proper}, written $F\subsetneq G$, when $F\ne G$.

The AHK variables associated with $G$ are
\begin{equation}\label{eq:graph-variables}
\mathcal X(G)=\big((\xi_a)_{a\in V_1(G)},(\eta_b)_{b\in V_2(G)},(\zeta_{ab})_{(a,b)\in E(G)}\big).
\end{equation}
The complete bipartite graph $K_{p,q}$ has $p$ row vertices, $q$ column vertices and all $pq$ edges. It represents one subarray, so the kernel $h$ is a function of $\mathcal X(K_{p,q})$ through the AHK representation. For each subgraph $G$, let
\begin{equation}\label{eq:graph-space}
\mathcal M_G=\{U\in L^2(\mathcal X(G)):\ \E[UV]=0\ \text{for all } V\in L^2(\mathcal X(F)),\ F\subsetneq G\}.
\end{equation}
Thus, $\mathcal M_G$ consists of the square-integrable functions of $\mathcal X(G)$ that are orthogonal to every function depending only on a proper subgraph.

Let $\Gamma$ contain one representative of each isomorphism class of nonempty subgraphs of $K_{p,q}$. We call its elements \emph{classes} and identify each class with its representative. A graph is \emph{connected} when $b(G)=1$ and \emph{disconnected} otherwise. A single isolated vertex, as in $K_{1,0}$ or $K_{0,1}$, is therefore connected. A \emph{copy} of $G$ is a labeled subgraph $\Lambda\subseteq K_{m,n}$ isomorphic to $G$. For a $G$-symmetric function $p_G\in \mathcal M_G$, we write $p_G(\Lambda)$ for its value on the variables of that copy. 

\paragraph{Hoeffding decomposition and fixed-law limits.} 
We recall the RCE version of the graph decomposition of \citet[Lem.~1]{janson1991asymptotic}, established in \citet{leminh2026characterization}.

\begin{lemmaL}[\citealp{leminh2026characterization}, Thm.~2 and eq.~(3.3)]\label{lem:decomp}
\[
U_N-\theta_N=\sum_{G\in\Gamma}T_G,\qquad T_G:=\frac1{(m)_{r(G)}(n)_{c(G)}}\sum_{\Lambda\text{ copy of }G}p_G(\Lambda),
\]
for every $m\ge p$, $n\ge q$,
where
\[
p_G=(p)_{r(G)}(q)_{c(G)}q_G,
\]
and $q_G$ is the orthogonal projection of the kernel, expressed in the AHK variables, onto $\mathcal M_G$. These functions do not depend on $(m,n)$.
\end{lemmaL}

The $T_G$ are pairwise orthogonal, and
\begin{equation}\label{eq:graph-variance-weight}
\var(T_G)=w_G\,\E[p_G^2],\qquad w_G(m,n):=\frac1{|\mathrm{Aut}(G)|\,(m)_{r(G)}(n)_{c(G)}}\asymp N^{-d(G)},
\end{equation}
and $\sigma_N^2=\var(U_N)=\sum_{G}\var(T_G)$. The \emph{principal degree} is the smallest $d(G)$ with $p_G\ne0$, and the corresponding $G$ are the \emph{principal support graphs}.

The numbers of vertices and connected components have different roles. The number of vertices $d(G)$ determines the factor $w_G\asymp N^{-d(G)}$. At a fixed law with $p_G\ne0$, this is also the order of $\var(T_G)$. When the law changes, however, $\|p_{N,G}\|_2$ may tend to zero. Contributions with different numbers of vertices can then have comparable variances.

The number of connected components $b(G)$ determines the form of the fixed-law limit. The connected case is already known.

\begin{lemmaL}[\citealp{leminh2026characterization}, Thm.~8]\label{lem:fixedclt}
Fix the law and let $d$ be the principal degree. If all principal support graphs are connected, then $N^{d/2}(U_N-\theta_N)$ is asymptotically Gaussian.
\end{lemmaL}

\noindent In our paper, Theorem~\ref{thm:chaos} will extend this result to disconnected graphs: each principal support graph with $b$ connected components contributes a term in the $b$-th Wiener chaos (defined in Sec.~\ref{sec:chaos}).

\subsection{The disconnected share}\label{sec:variance-proportions}

When the law changes with $N$, the principal degree does not completely describe the contributions that determine the variance. For instance, a low-degree projection can be nonzero for every $N$ but contribute a vanishing proportion of the variance (Ex.~\ref{ex:k12}). We therefore use the \emph{variance share} of each $G\in\Gamma$,
\begin{equation}\label{eq:variance-proportions}
s_{N,G}^2:=\var(T_G)/\sigma_N^2,\qquad \sum_{G\in\Gamma}s_{N,G}^2=1,
\end{equation}
and the sum of these shares over disconnected graphs, namely the \emph{disconnected share}
\[
D_N:=\sum_{G\ \disc}s_{N,G}^2.
\]

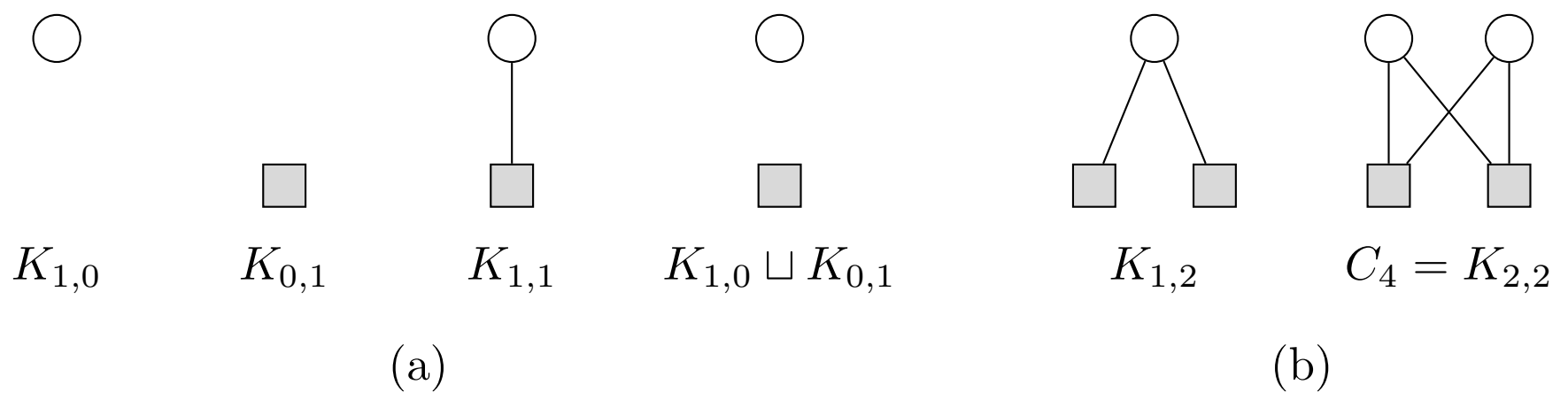
\begin{figure}[!ht]
\centering
\begin{tikzpicture}[
  row/.style={circle, draw, fill=white, minimum size=10pt, inner sep=0pt},
  col/.style={rectangle, draw, fill=black!15, minimum size=9pt, inner sep=0pt},
  lab/.style={anchor=base}]
\node[row] at (0,1.1) {};
\node[lab] at (0,-0.7) {$K_{1,0}$};
\node[col] at (1.7,0) {};
\node[lab] at (1.7,-0.7) {$K_{0,1}$};
\node[row] (a1) at (3.4,1.1) {}; \node[col] (b1) at (3.4,0) {};
\draw (a1) -- (b1);
\node[lab] at (3.4,-0.7) {$K_{1,1}$};
\node[row] at (5.4,1.1) {}; \node[col] at (5.4,0) {};
\node[lab] at (5.4,-0.7) {$K_{1,0}\sqcup K_{0,1}$};
\node[row] (r) at (8.2,1.1) {};
\node[col] (c1) at (7.75,0) {}; \node[col] (c2) at (8.65,0) {};
\draw (r) -- (c1) (r) -- (c2);
\node[lab] at (8.2,-0.7) {$K_{1,2}$};
\node[row] (r1) at (9.95,1.1) {}; \node[row] (r2) at (10.85,1.1) {};
\node[col] (s1) at (9.95,0) {}; \node[col] (s2) at (10.85,0) {};
\draw (r1) -- (s1) (r1) -- (s2) (r2) -- (s1) (r2) -- (s2);
\node[lab] at (10.4,-0.7) {$C_4=K_{2,2}$};
\node[lab] at (2.7,-1.45) {(a)};
\node[lab] at (9.3,-1.45) {(b)};
\end{tikzpicture}
\caption{Support graphs, with row vertices as circles and column vertices as squares. (a) The four classes of the two-way mean (Ex.~\ref{ex:mean}). (b) The connected graphs $K_{1,2}$, of degree three (Ex.~\ref{ex:k12}), and the $4$-cycle $C_4=K_{2,2}$, of degree four (Sec.~\ref{sec:d3}).}\label{fig:support}
\end{figure}

\begin{example}[The two-way mean]\label{ex:mean}
Take $p=q=1$ and $h(y)=y$, so that $U_N=\bar Y$. The class set $\Gamma$ has four elements: $K_{1,0}$ and $K_{0,1}$ (one row vertex or one column vertex), the single-edge graph $K_{1,1}$, and the edgeless graph $K_{1,0}\sqcup K_{0,1}$ with one row and one column vertex (Fig.~\ref{fig:support}(a)). Writing $a(\xi_i)=\E[Y_{ij}\mid\xi_i]-\theta_N$, $b(\eta_j)=\E[Y_{ij}\mid\eta_j]-\theta_N$, $c(\xi_i,\eta_j)=\E[Y_{ij}\mid\xi_i,\eta_j]-\theta_N-a(\xi_i)-b(\eta_j)$ and $e_{ij}=Y_{ij}-\E[Y_{ij}\mid\xi_i,\eta_j]$, the graph projections are the row effect $a$ (due to $K_{1,0}$), the column effect $b$ (due to $K_{0,1}$), the cell noise $e$ (due to $K_{1,1}$), and the latent interaction $c$ (due to $K_{1,0}\sqcup K_{0,1}$), with
\[
\sigma_N^2=\frac{\E[a^2]}{m}+\frac{\E[b^2]}{n}+\frac{\E[e^2]+\E[c^2]}{mn},\qquad D_N=\frac{\E[c^2]/(mn)}{\sigma_N^2}.
\]
Thus, $D_N$ is the proportion of variance due to the latent interaction $c$, denoted by $q_v$ in \citet{menzel2021bootstrap}. At a fixed law, when the row and column effects $a$ and $b$ vanish and the remaining variance is positive, the standard deviation $\sigma_N$ has order $(mn)^{-1/2} \asymp N^{-1}$ instead of $N^{-1/2}$. 

The projections $e$ and $c$ have the same degree. If $\E[e^2]$ stays bounded away from zero and $\E[c^2]\to0$, then $D_N\to0$ and the limit is Gaussian (Thm.~\ref{thm:contraction}). At a fixed law, the singular value decomposition of $c$ is $c(x,y)=\sum_k\lambda_ku_k(x)v_k(y)$. Then
\[
(mn)^{-1/2}\sum_{i,j}c(\xi_i,\eta_j)\dto\sum_k\lambda_kZ_kW_k,
\]
where the $Z_k,W_k$ are independent standard normal variables. In particular, for a rank-one interaction, this gives a single product $\lambda ZW$. Thus, $K_{1,1}$ and $K_{1,0}\sqcup K_{0,1}$ can give different limits despite having the same degree.
\end{example}

\begin{example}[A connected graph of degree three]\label{ex:k12}
Let $Y_{ij}=tA_i+e_{ij}$ with $t=t_N\ge0$, where the row effects $A_i=a(\xi_i)$ are i.i.d. and bounded, with $\E[A_i]=0$ and $\var(A_i^2)>0$, and the cell noises $e_{ij}=\varepsilon(\zeta_{ij})$ are i.i.d. with finite moments of all orders, with $\E[e_{ij}]=0$ and $\E[e_{ij}^2]=1$. The two families are assumed independent. Take the $1\times2$ kernel $h(y_1,y_2)=y_1y_2$. The nonzero classes are $K_{1,0}$, $K_{1,1}$ and $K_{1,2}$ (Fig.~\ref{fig:support}), all connected, so $D_N=0$, and
\[
\sigma_N^2=\frac{t^4\var(A^2)}{m}+\frac{4t^2\E[A^2]}{mn}+\frac{2}{mn(n-1)}.
\]
With $m=n=L_N$, the three terms are of orders $t^4L_N^{-1}$, $t^2L_N^{-2}$ and $L_N^{-3}$. The row effect dominates when $t\gg L_N^{-1/2}$, the three terms are comparable when $t\asymp L_N^{-1/2}$, and the class $K_{1,2}$ dominates when $t\ll L_N^{-1/2}$. In particular, for $t=L_N^{-1}$ the principal degree is $1$ for every $N$, yet $K_{1,2}$ carries asymptotically all the variance.
\end{example}

\subsection{Variance estimators}\label{sec:variance-estimators}

For $0\le a\le p$ and $0\le b\le q$, let $\Thetah(a,b)$ be the average of $H_WH_{W'}$ over all ordered pairs of subarrays sharing exactly $a$ rows and $b$ columns. It is unbiased for $\Theta(a,b):=\E[H_WH_{W'}]$ for any such pair. Throughout, $m\ge2p$ and $n\ge2q$.

\begin{definition}\label{def:ov}
The \emph{unbiased variance estimator} is
\[
\Vu:=U_N^2-\Thetah(0,0).
\]
\end{definition}

\noindent Subarrays that share neither rows nor columns are independent by dissociation. Hence, $\E[\Thetah(0,0)]=\theta_N^2$ and $\E[\Vu]=\sigma_N^2$. This is a classical variance estimator for ordinary $U$-statistics \citep{wang2014variance} and the two-way mean \citep{mccullagh2000resampling}, and has been discussed for general $U$-statistics on RCE arrays in \citet[Sec.~5]{leminh2025hoeffdingtype}. Section~\ref{sec:S-split} of the supplementary material gives a computationally convenient approximation based on random splits of the array and bounds its error.

$\Vu$ can be negative, so for inference, we use its positive part 
\[
\Vu^{+}:=\max(\Vu,0).
\]
In this paper, if the variance estimate used to studentize is nonpositive, we set the studentized statistic to zero, do not reject the null hypothesis, and report the whole parameter space as the confidence set.

The row-column jackknife is computed by deleting one row or one column at a time. Write $U^{(-i)}$ for the statistic without row $i$, and $U^{(-j)}$ for the statistic without column $j$. Let $\bar U_{\mathrm{row}}$ and $\bar U_{\mathrm{col}}$ be the averages of the row-deleted and column-deleted statistics, respectively.

\begin{definition}\label{def:rcj}
The \emph{row-column jackknife} is the sum of the delete-one-row and delete-one-column jackknife variance estimates,
\begin{equation}\label{eq:rcj}
\Vrcj:=\frac{m-1}{m}\sum_{i=1}^m\big(U^{(-i)}-\bar U_{\mathrm{row}}\big)^2+\frac{n-1}{n}\sum_{j=1}^n\big(U^{(-j)}-\bar U_{\mathrm{col}}\big)^2 .
\end{equation}
\end{definition}

\noindent For the two-way mean, it corresponds to the two-term cluster jackknife of \citet{mackinnon2024jackknife}. We compare it to the unbiased variance estimator in Section~\ref{sec:jackknife}.

\section{Uniform Gaussian inference and its boundary}\label{sec:results}

\subsection{Uniform inference}\label{sec:res-main}
First, we give conditions under which the standardized statistic $(U_N-\theta_N)/\sigma_N$ is asymptotically Gaussian and the unbiased variance estimator $\Vu$ is ratio-consistent as the law changes. Then, we show that, under compactness, both properties require the disconnected share $D_N$ to vanish.

\paragraph{The contraction condition.}
As the law changes, the variance contributions can move between graph classes and the principal degree can change. At a fixed law, connected classes have Gaussian limits (Lem.~\ref{lem:fixedclt}), whereas disconnected classes have limits in higher Wiener chaoses (Thm.~\ref{thm:chaos}). When the law changes, however, disconnected classes can also have Gaussian limits if their contributions spread over increasingly many independent directions.

To include this possibility, we measure how concentrated the disconnected contributions are. For a disconnected graph $G$, the projection $p_G$ is a function of the independent groups of variables associated with its connected components. For a nonempty proper set $B$ of components, we write
\[
(p_G\otimes_Bp_G)(x,x')=\E[p_G(x,\mathbf{X}_B)\,p_G(x',\mathbf{X}_B)],
\]
where $\mathbf{X}_B$ contains the variables of $B$, and $x,x'$ contain those of the remaining components. This operation, called a contraction, averages the product of two copies of $p_G$ over the variables they share. We define
\begin{equation}\label{eq:gamma}
\Upsilon_N:=\sigma_N^{-4}\sum_{G\ \disc}w_G^2\sum_{\substack{B\subsetneq\{\text{components of }G\}\\B \ne \varnothing}}\|p_G\otimes_Bp_G\|_2^2,
\end{equation}
where the contraction is computed on any copy of $G$.

The condition $\Upsilon_N\to0$ therefore allows disconnected contributions either to become negligible or to spread over increasingly many directions. For the two-way mean of Example~\ref{ex:mean}, with $c=\sum_k\lambda_ku_k\otimes v_k$,
\begin{equation}\label{eq:two-way-mean}
\Upsilon_N=\frac{2\sum_k\lambda_k^4}{(mn\sigma_N^2)^2},
\qquad
D_N=\frac{\sum_k\lambda_k^2}{mn\sigma_N^2}.
\end{equation}
Thus, $D_N$ measures the total disconnected contribution, while $\Upsilon_N$ also accounts for how it is distributed across directions. 

\paragraph{Main theorem: Uniform inference.} For a class of laws $\Pcal$, we use the following fourth-moment assumption:
\begin{description}
\item[(M4)] For all $G\in\Gamma$, $\E[p_{N,G}^4]\le C_{\mathrm{M4}}(\E[p_{N,G}^2])^2$, uniformly over $\Pcal$.
\end{description}
Under (M4), the following theorem, proved in Section~\ref{app:core} of the supplementary material, shows that $\Upsilon_N\to0$ implies that the standardized statistic is asymptotically Gaussian and the unbiased variance estimator is ratio-consistent. It relies on a multidimensional version of de Jong's central limit theorem stated in \citet{dobler2017quantitative}.

\begin{theorem}\label{thm:contraction}
Let $m\ge2p$, $n\ge2q$, and assume (M4) on a class $\Pcal$ with constant $C_{\mathrm{M4}}$ and $\sigma_{N,P}>0$ for every $P\in\Pcal$. Then
\[
\E\left[\left(\frac{\Vu}{\sigma_N^2}-1\right)^2\right]\le C\big(L_N^{-1}+L_N^{-1/2}D_N^{1/2}+\Upsilon_N^{1/4}\big),
\]
where $C$ depends only on $(p,q,C_{\mathrm{M4}})$. If moreover $\Upsilon_N\to0$, then $(U_N-\theta_N)/\sigma_N\dto\Nor(0,1)$, $\Vu/\sigma_N^2\to1$, and the studentized statistic $\widehat Z_N=(U_N-\theta_N)/(\Vu^+)^{1/2}$ satisfies $\widehat Z_N\dto\Nor(0,1)$. If $\Upsilon_N\to0$ uniformly over $\Pcal$, then these convergences hold uniformly over $\Pcal$.
\end{theorem}

\paragraph{When disconnected contributions are negligible.}\label{sec:res-gauss}\label{sec:res-attain}
The next proposition shows that $D_N\to0$ implies $\Upsilon_N\to0$, with equivalence under the following compactness assumption on a class $\Pcal$:
\begin{description}
\item[(K)] For each $G$, $\{p_{P,G}/\|p_{P,G}\|_2:P\in\Pcal,\ p_{P,G}\ne0\}$ is precompact in $L^2(\mathcal X(K_{p,q}))$.
\end{description}
Assumption (K) means that every sequence of normalized projections $p_{P,G}/\|p_{P,G}\|_2$ has a subsequence converging in $L^2$. It holds, for example, when these projections belong to a fixed finite-dimensional subspace.

\begin{proposition}\label{prop:gammaD}
Let $\Pcal$ be a class of laws with $\sigma_{N,P}>0$ for every $P\in\Pcal$. Then
\[
\Upsilon_N\le CD_N^2
\]
for every $P\in\Pcal$, where $C$ depends only on $(p,q)$. If (K) holds on $\Pcal$, then
\[
\Upsilon_N\ge cD_N^2
\]
for every $P\in\Pcal$, where $c>0$ depends only on $(p,q)$ and $\Pcal$.
\end{proposition}

Proposition~\ref{prop:gammaD} is proved in Section~\ref{sec:S-gammaD} of the supplementary material. Combining it with Theorem~\ref{thm:contraction}, under the same moment and sample-size assumptions, gives
\begin{equation}\label{eq:variance-D-bound}
\E\left[\left(\frac{\Vu}{\sigma_N^2}-1\right)^2\right]\le C\big(L_N^{-1}+D_N^{1/2}\big).
\end{equation}
In particular, when $D_N=0$, $\Vu/\sigma_N^2=1+O_P(L_N^{-1/2})$, which tends to $1$ as $L_N \to \infty$.

\begin{remark}\label{rem:diffuse}
The condition $\Upsilon_N\to0$ also covers cases where disconnected contributions are not negligible. For the two-way mean, let 
\[
Y_{ij}=\nu_N^{-1/2}\sum_{k=1}^{\nu_N}A_{ik}B_{jk},
\]
where all $A_{ik},B_{jk}$ are independent Rademacher variables and $\nu_N\to\infty$.
In this model, $A_{ik}$ and $B_{jk}$ can be expressed in terms of the AHK variables as
\[
A_{ik}=1-2b_k(\xi_i),\qquad B_{jk}=1-2b_k(\eta_j),
\]
where $b_k(x)$ is the $k$-th binary digit of $x$.
The only nonzero projection is the one on $K_{1,0}\sqcup K_{0,1}$. Its squared norm is $\E[Y_{11}^2]=1$, so $D_N=1$ and $mn\sigma_N^2=1$. Since the products $A_{1k}B_{1k}$ are independent Rademacher variables, $\E[Y_{11}^4]=3-2/\nu_N\le3$, which verifies (M4). Moreover, the functions $1-2b_k$ are orthonormal in $L^2([0,1])$, so the nonzero coefficients in \eqref{eq:two-way-mean} are $\lambda_k=\nu_N^{-1/2}$ for $1\le k\le\nu_N$. Because $mn\sigma_N^2=1$, we obtain
\[
\Upsilon_N=2\sum_{k=1}^{\nu_N}(\nu_N^{-1/2})^4=\frac{2}{\nu_N}\to0.
\]
Theorem~\ref{thm:contraction} therefore gives $(mn)^{1/2}\bar Y\dto\Nor(0,1)$. Assumption (K) fails on any class containing these laws: since $D_N=1$, Proposition~\ref{prop:gammaD} would otherwise bound $\Upsilon_N$ away from zero. Thus, this confirms that connectedness is not necessary for a Gaussian limit.
\end{remark}

\subsection{The Gaussian boundary}\label{sec:res-sharp}

Remark~\ref{rem:diffuse} shows that disconnected graphs can contribute all the variance without preventing a Gaussian limit. However, under the compactness assumption (K), this possibility is excluded: the distribution of the standardized statistic $(U_N-\theta_N)/\sigma_N$ converges uniformly to the standard normal distribution exactly when the disconnected share $D_N$ tends to zero uniformly. This is what Theorem~\ref{thm:sharp}, proved in Section~\ref{app:sharp} of the supplementary material, states.

\begin{theorem}\label{thm:sharp}
Let $\Pcal$ be a class of laws with $\sigma_{N,P}>0$ for every $P\in\Pcal$. Assume (K) on $\Pcal$, and consider:
\begin{enumerate}[(a)]
\item $(U_N-\theta_N)/\sigma_N\dto\Nor(0,1)$ uniformly over $P\in\Pcal$, \item $D_N\to0$ uniformly over $P\in\Pcal$.
\end{enumerate}
Then (a) and (b) are equivalent. For the implication (b)$\Rightarrow$(a), it suffices that (K) holds for connected graphs $G$.
\end{theorem}

\paragraph{Consequence for Gaussian inference.}\label{cor:nopivot}
Under (K), a uniformly ratio-consistent variance estimator can give a uniformly Gaussian studentized statistic only if $D_N\to0$ uniformly, by Slutsky's lemma. Moreover, under (M4), this is also exactly the condition for ratio-consistency of the unbiased variance estimator $\Vu$.

\begin{corollary}\label{cor:S-consistency}
Let $\Pcal$ be a class of laws with $\sigma_{N,P}>0$ for every $P\in\Pcal$. Assume (M4) and (K) on $\Pcal$. Then $\Vu/\sigma_N^2\to1$ in probability uniformly over $\Pcal$ if and only if $D_N\to0$ uniformly over $\Pcal$.
\end{corollary}

Corollary~\ref{cor:S-consistency} is proved in Section~\ref{sec:S-cor-consistency} of the supplementary material. Its proof and that of Theorem~\ref{thm:sharp} use the joint limits described below.

\paragraph{The projection joint limits.}\label{sec:chaos}
Under the compactness assumption (K), any sequence of normalized projections has a convergent subsequence. The following lemma, proved in Section~\ref{app:wick-construction} of the supplementary material, gives the joint limits associated with these projections (two-way analogue of \citealp[Lem. 9]{janson1991asymptotic}). Let $T_{N,G}(f)$ be the sum defining $T_G$ in Lemma~\ref{lem:decomp}, with $p_G$ replaced by $f$. An \emph{isonormal Gaussian process} on a Hilbert space $\mathfrak h$ is a centered Gaussian process $(X(f))_{f\in\mathfrak h}$ with $\E[X(f)X(g)]=\langle f,g\rangle$. Its $b$-th \emph{Wiener chaos} is the closed linear span of the variables $H_b(X(f))$, $\|f\|=1$, where $H_b$ is the $b$-th Hermite polynomial \citep[Ch.~2]{janson1997gaussian}: the first chaos consists of centered Gaussian variables, and different chaoses are orthogonal.

\begin{lemma}\label{lem:wick}
Let $\Gamma_0$ be a finite set of nonempty classes and, for $G\in\Gamma_0$, let $f_G\in \mathcal M_G$ be fixed, nonzero and $G$-symmetric. Then the standardized sums $W_{N,G}=T_{N,G}(f_G)/\var(T_{N,G}(f_G))^{1/2}$ converge jointly in distribution to limits $W_G$ built from one isonormal Gaussian process. Each $W_G$ lies in the $b(G)$-th Wiener chaos, with $\E[W_G^2]=1$ and $\E[W_GW_{G'}]=0$ for $G\ne G'$.
\end{lemma}

For connected graphs, the limits are independent standard normal variables. In contrast, a rank-one projection $u\otimes u$ on $K_{2,0}$ gives $(Z^2-1)/\sqrt2$, with $Z$ standard normal.

Now, suppose that $D_N$ stays bounded away from zero. Since there are finitely many graph classes, the variance shares $s_{N,G}^2$ converge along a subsequence, with a positive limit for at least one disconnected graph. By (K), the normalized projections of graphs with positive limiting shares converge in $L^2$ to unit-norm functions along a further subsequence. Moreover, limits from distinct graph classes are orthogonal, so their contributions cannot cancel. The resulting limit therefore has a nonzero term in a Wiener chaos of order at least two. A finite sum of chaoses with such a term cannot be Gaussian (Lem.~\ref{lem:chaosnotgauss} in Sec.~\ref{app:wick} of the supplementary material). This explains the necessity of $D_N\to0$ under (K) in Theorem~\ref{thm:sharp}.

\paragraph{Fixed-law limits of $U$-statistics.}
At a fixed law, the same construction gives the limit of $U_N$ by grouping the principal support graph contributions according to their number of connected components. The following theorem, proved in Section~\ref{sec:S-thm-chaos} of the supplementary material, is the two-way analogue of \citet[Thm.~3]{janson1991asymptotic}. It completes the multi-component case left open by \citet{leminh2026characterization}.

\begin{theorem}\label{thm:chaos}
Fix the law and let $d$ be the principal degree. Then $N^{d/2}(U_N-\theta_N)\dto\sum_{b=1}^dL_b$, where 
\begin{equation}\label{eq:fixed-law-coefficients}
L_b=\sum_{G:d(G)=d,\,b(G)=b}\tau_GW_G, \qquad \tau_G^2=\frac{\|p_G\|_2^2}{\rho^{r(G)}(1-\rho)^{c(G)}|\operatorname{Aut}(G)|}.
\end{equation}
Each $L_b$ lies in the $b$-th Wiener chaos of one common
isonormal Gaussian process, and $L_b\ne0$ if and only if
some principal support graph has exactly $b$ connected components.
\end{theorem}

The proof first treats finite sums of products of bounded functions, then uses the variance identity~\eqref{eq:graph-variance-weight} to pass to the square-integrable projections.

\section{The row-column jackknife}\label{sec:jackknife}\label{sec:counts}

The preceding results use the unbiased variance estimator $\Vu$. A natural alternative is the row-column jackknife $\Vrcj$ in~\eqref{eq:rcj}, which is always nonnegative. 

\subsection{The jackknife bias}\label{sec:jackknife-bias} 
The proposition below, proved in Section~\ref{app:estimators} of the supplementary material, shows that in expectation the row-column jackknife weights each variance contribution $\var(T_G)$ by a factor tending to its number of vertices. It explains why the jackknife does not generally yield valid Gaussian inference in degenerate regimes.

\begin{proposition}\label{prop:count}
Within the dissociated RCE framework, for every law and all $m>p$, $n>q$,
\[
\E[\Vrcj]=\sum_{G\in\Gamma}\kappa_G(m,n)\,\var(T_G),\qquad
\kappa_G(m,n)=r(G)\frac{m-1}{m-r(G)}+c(G)\frac{n-1}{n-c(G)} .
\]
\end{proposition}

The weights satisfy $\kappa_G\to d(G)$ and $\kappa_G\ge1$, with exact equality only for $K_{1,0}$ and $K_{0,1}$. Indeed, deleting a row or a column removes every graph copy containing that vertex. The row jackknife therefore counts a contribution with \(r\) row vertices \(r\) times asymptotically, and the column jackknife counts it \(c\) times. This explains why the weights count vertices rather than edges.

Since $\sigma_N^2=\sum_{G\in\Gamma}\var(T_G)$, the bias of the jackknife is
\[
\E[\Vrcj]-\sigma_N^2=\sum_{G\in\Gamma}(\kappa_G-1)\var(T_G)=\sum_{r+c\ge2}(\kappa_{r,c}-1)V_{r,c},
\]
where $V_{r,c}:=\sum_{G:\,r(G)=r,\,c(G)=c}\var(T_G)$ and $\kappa_{r,c}=r(m-1)/(m-r)+c(n-1)/(n-c)$ is the common value of $\kappa_G$ for graphs with $r$ row and $c$ column vertices, and the terms with $r+c=1$ vanish because $\kappa_{1,0}=\kappa_{0,1}=1$. Relative to the variance, the bias is $\sum_G(\kappa_G-1)s_{N,G}^2=\sum_G\{d(G)-1\}s_{N,G}^2+O(L_N^{-1})$.

\subsection{Removing the jackknife bias}\label{sec:overlaps}
The bias involves the variances only through the sums $V_{r,c}$, which can be estimated without knowing the graph projections. The projection on a copy of a graph contributes to $\E[H_WH_{W'}]$ exactly when both subarrays contain all vertices of that copy. Lemma~\ref{lem:overlap} (Sec.~\ref{sec:S-overlap} in the supplementary material) counts these copies and gives, for $0\le a\le p$ and $0\le b\le q$,
\[
\Theta(a,b)=\theta_N^2+\sum_{\substack{r\le a,\ c\le b\\ r+c\ge1}}\binom ar\binom bc\frac{\binom mr\binom nc}{\binom pr^2\binom qc^2}\,V_{r,c}.
\]
Inverting these relations and replacing $\Theta$ by $\Thetah$ gives the unbiased variance-component estimators
\begin{equation}\label{eq:variance-component-estimator}
\widehat V_{r,c}:=
\frac{\binom pr^2\binom qc^2}{\binom mr\binom nc}
\sum_{a=0}^{r}\sum_{b=0}^{c}
(-1)^{r-a+c-b}\binom ra\binom cb\,\Thetah(a,b),
\qquad r+c\ge1,
\end{equation}
for $0\le r\le p$ and $0\le c\le q$, where $\Thetah(a,b)$ is the average product of kernel values over pairs of subarrays sharing exactly $a$ rows and $b$ columns (Sec.~\ref{sec:variance-estimators}).

We have the following identity, proved in Section~\ref{sec:S-pathwise} of the supplementary material.
\begin{proposition}\label{prop:pathwise}
For every array and every kernel with $m\ge2p$ and $n\ge2q$,
\begin{enumerate}[(i)]
\item $\Vrcj=\sum_{r+c\ge1}\kappa_{r,c}\,\widehat V_{r,c}$;
\item $\Vu=\sum_{r+c\ge1}\widehat V_{r,c}=U_N^2-\Thetah(0,0)$;
\item consequently,
\begin{equation}\label{eq:pathwise}
\Vu=\Vrcj-\sum_{r+c\ge2}(\kappa_{r,c}-1)\,\widehat V_{r,c}.
\end{equation}
\end{enumerate}
\end{proposition}
Thus, subtracting the estimated bias from the jackknife gives exactly the unbiased variance estimator, for every observed array and not only in expectation.

\section{Inference without a Gaussian limit}\label{sec:hybrid}
Gaussian inference can fail when disconnected graphs contribute to the limit. The following example shows why a conservative test cannot always rely on a uniformly ratio-consistent variance estimate, although (M4) and (K) hold on its class of laws.

\begin{example}\label{ex:menzel-product}
Let $Y_{ij}=X_iZ_j$, where $X_i\sim\Nor(0,1)$ and $Z_j\sim\Nor(\mu,1)$ are independent i.i.d.\ sequences, and $|\mu|\le M$ for a fixed $M>0$. The mean $U_N=\bar X\bar Z$ has variance $\sigma_N^2=(1+n\mu^2)/(mn)$. This is precisely the example used in the proof of \citet[Lem.~A.1(b)]{menzel2021bootstrap}, which shows that no estimator is ratio-consistent for $\sigma_N^2$ uniformly over this class, even if all latent variables are observed. 
\end{example}

To retain uniformity even in non-Gaussian regimes, we propose an inference procedure with two ingredients: a conservative test that does not require $\Upsilon_N\to0$ (Sec.~\ref{sec:S-block-hybrid}), and a diagnostic to choose between this test and the Gaussian test (Sec.~\ref{sec:diagnostic-main}). However, this diagnostic is computationally impractical except for small arrays. We therefore propose an alternative procedure, valid under additional assumptions, whose diagnostic only requires quantities already obtained when computing $\Vu$ (Sec.~\ref{sec:alternative-procedure}).

\subsection{A conservative test from independent blocks}\label{sec:S-block-hybrid}

We can use independent blocks to obtain a conservative bound for the variance of the complete statistic $U_N$, used to estimate $\theta_N$.

We start with a deterministic integer $K_N\ge2$ with $K_N\to\infty$ and $K_N=o(L_N)$. We partition the rows into $K_N$ disjoint groups of size $m_b=\lfloor m/K_N\rfloor$ and the columns into $K_N$ disjoint groups of size $n_b=\lfloor n/K_N\rfloor$. Then, we pair the $j$-th row group with the $j$-th column group to form block $j$. Any remaining rows and columns are omitted from the block calculation, but remain in $U_N$. The partition may be deterministic or chosen independently of the observations.

Now, we compute the complete statistic $U_{N,j}$ on each block, using the same kernel as for $U_N$. Set
\begin{equation}\label{eq:block-statistics}
 \overline U_{\rm block}=K_N^{-1}\sum_{j=1}^{K_N}U_{N,j},\qquad
 S_{K_N}^2=\frac1{K_N-1}\sum_{j=1}^{K_N}(U_{N,j}-\overline U_{\rm block})^2.
\end{equation}
Conditional on the partition, the block statistics are i.i.d. by dissociation and exchangeability, with a law independent of the chosen partition.

The following proposition, proved in Section~\ref{app:protected} of the supplementary material, uses their sample variance $S_{K_N}^2$ to construct a conservative test. 

\begin{proposition}[Conservative block inference]\label{prop:block-safe}
Let $0<\alpha<1$, and let $\Pcal_0$ be a null class of laws with $\theta_N=\theta_0$ under the standing array and sample-growth assumptions. Assume (M4) on $\Pcal_0$, with $\sigma_{N,P}>0$ for every $P\in\Pcal_0$. Use the preceding block construction. For any nonnegative statistic $A_N$, let $B_N=\max(A_N,S_{K_N}^2/K_N)$. Then
\[
\limsup_N\sup_{P\in\Pcal_0}\P_P\{|U_N-\theta_0|\ge\sqrt{B_N/\alpha},\ B_N>0\}\le\alpha.
\]
The same conclusion holds for the upper one-sided test, with rejection event
\[
\left\{U_N-\theta_0\ge\sqrt{\frac{(1-\alpha)B_N}{\alpha}},\ B_N>0\right\}.
\]
\end{proposition}

A typical choice for $A_N$ is $A_N:=\Vu^+$. This construction has a cost. At a fixed law with principal degree $d$, the variance of the block average $\overline U_{\rm block}$ is of order $K_N^{d-1}\sigma_N^2$. Therefore, the block term $S_{K_N}^2/K_N$, which estimates this variance, alone gives an interval half-width of order $K_N^{(d-1)/2}\sigma_N$. When the law changes, the cost depends on the contributions of all graph classes, as detailed in Section~\ref{sec:S-block-cost} of the supplementary material.

\begin{remark}\label{rem:block-average}
  One can also use the block average $\overline U_{\rm block}$ itself as an incomplete $U$-statistic to estimate $\theta_N$. Under the same assumptions,
\[
\frac{\sqrt{K_N}(\overline U_{\rm block}-\theta_N)}{S_{K_N}}
\]
is uniformly asymptotically standard normal. This gives Gaussian inference without a diagnostic, but increases the variance compared with $U_N$. Corollary~\ref{cor:S-block-average} (Sec.~\ref{sec:S-block-average} in the supplementary material) proves this result.
\end{remark}

\subsection{Diagnostic and procedure}\label{sec:diagnostic-main}
Theorem~\ref{thm:contraction} requires $\Upsilon_N\to0$, but $\Upsilon_N$ depends on the unknown graph projections. To use the Gaussian test when possible, we therefore need to estimate a quantity comparable to $\Upsilon_N$. However, the expected products \(\Theta(a,b)\) of kernel values, defined in Section~\ref{sec:variance-estimators}, do not determine \(\Upsilon_N\). For instance, in the two-way mean of Example~\ref{ex:mean} without row and column effects, an interaction with $\E[c^2]=s$ and no cell noise, and cell noise with $\E[e^2]=s$ and no interaction, give the same values of \(\Theta(a,b)\). Yet, $D_N=1$ and $\Upsilon_N>0$ for the interaction, whereas $D_N=\Upsilon_N=0$ for the cell noise. Instead, Proposition~\ref{prop:ident} shows that products of four suitably chosen differences of kernel values identify such a quantity. We use them to construct the diagnostic below.

\paragraph{Computing the diagnostic.} 
The diagnostic index $\tau=(r_A,c_A\,|\,r_B,c_B)$ records the row and column counts of two vertex sets, with $r_A+c_A\ge1$, $r_B+c_B\ge1$, $r=r_A+r_B\le p$ and $c=c_A+c_B\le q$. We choose disjoint sets $A$ and $B$ containing these numbers of rows and columns.

For each vertex $u\in A\cup B$, we choose a distinct replacement $u''$ of the same kind: a row replaces a row and a column replaces a column. We add enough rows and columns to complete $A\cup B$ to a $p\times q$ subarray, with all additional and replacement vertices distinct from each other and from $A\cup B$. For $T\subseteq A\cup B$, let $W_T$ be the subarray obtained by replacing each $u\in T$ by $u''$. Let us define the alternating sum
\begin{equation}\label{eq:diagnostic-differences}
\mathcal D(A,B)=\sum_{T\subseteq A\cup B}(-1)^{|T|}H_{W_T}.
\end{equation}
Next, choose four mutually disjoint sets $A_1,A_2,B_1,B_2$, with each $A_i$ containing $r_A$ rows and $c_A$ columns, and each $B_j$ containing $r_B$ rows and $c_B$ columns. Form the four sums $\mathcal D(A_i,B_j)$ with separate replacement and additional vertices, all distinct across the four sums and from $A_1,A_2,B_1,B_2$. Put
\begin{equation}\label{eq:diagnostic-product}
\Pi_\tau=\mathcal D(A_1,B_1)\,\mathcal D(A_2,B_1)\,\mathcal D(A_1,B_2)\,\mathcal D(A_2,B_2).
\end{equation}

Let $\widehat Z_\tau$ be the average of $\Pi_\tau$ over all these choices, and define the diagnostic statistic $\widehat\Upsilon_N$ by
\begin{equation}\label{eq:diagnostic-statistic}
\widehat\Upsilon_N=\max_\tau\frac{c_\tau\widehat Z_\tau^+}{(\Vu^+)^2},\qquad c_\tau=\binom r{r_A}\binom c{c_A}\{(m)_r(n)_c\}^{-2},
\end{equation}
with $\widehat\Upsilon_N=+\infty$ if $\Vu\le0$.

\paragraph{Diagnostic assumptions and validity.}
The following proposition is proved in Section~\ref{sec:S-ident} of the supplementary material.

\begin{proposition}\label{prop:ident}
Assume $\E[|h(Y_{[p],[q]})|^4]<\infty$ and $\sigma_N>0$.
There exist constants $c,C>0$, depending only on $(p,q)$, such that
\[
c\,\Upsilon_N
\le
\frac{1}{\sigma_N^4}\sum_\tau c_\tau\,\E[\widehat Z_\tau]
\le
C\,\Upsilon_N.
\]
\end{proposition}
The statistic $\widehat\Upsilon_N$ replaces \(\mathbb E[\widehat Z_\tau]\) by \(\widehat Z_\tau^+\) and \(\sigma_N^4\) by \((\Vu^+)^2\). Theorem~\ref{thm:S-check} bounds the estimation errors under the additional assumptions below. We use the ordinary graph projections $q_G$ defined in Lemma~\ref{lem:decomp}. For a set of row and column vertices $A$, let $x_A$ collect the values of their vertex variables and the cell variables joining the rows and the columns in $A$.
For disjoint $A,B$ and a square-integrable function $F(x_A,x_B)$, define
\[
\|F\|_{S_4}^4
:=\E\!\left[F(\mathbf{X}_A,\mathbf{X}_B)F(\mathbf{X}_A',\mathbf{X}_B)F(\mathbf{X}_A,\mathbf{X}_B')F(\mathbf{X}_A',\mathbf{X}_B')\right],
\]
where $\mathbf{X}_A,\mathbf{X}_A',\mathbf{X}_B,\mathbf{X}_B'$ are independent, $\mathbf{X}_A,\mathbf{X}_A'$ have the distribution of the AHK variables in $x_A$, and $\mathbf{X}_B,\mathbf{X}_B'$ have that of the variables in $x_B$.
For a stated class $\Pcal$:
\begin{description}
\item[(M8)] For all $G\in\Gamma$, $\E[|p_{N,G}|^8]\le C_{\mathrm{M8}}(\E[p_{N,G}^2])^4$, uniformly over $\Pcal$.
\item[(S2)] Fix an index $\tau$, sets $A,B$ and a nonempty set $U$ of additional vertices such that $A\cup B\cup U$ fits in the subarray. Let $F=\sum_\Lambda q_\Lambda$, where the sum is over graphs $\Lambda$ with vertex set exactly $A\cup B\cup U$ and no edge joining different sets among $A,B,U$. For each value \(x_U\), let \(F_U(x_A,x_B)=F(x_A,x_B,x_U)\). Write $\mathcal A(F)$ for the splits $U=U_1\sqcup U_2$ such that no nonzero graph term of $F$ has an edge between $U_1$ and $U_2$. For each split, group the arguments of $F$ as $(x_{A\cup U_1},x_{B\cup U_2})$. 
Uniformly over $\Pcal$ and over all such $\tau$, $A$, $B$ and $U$, we assume
\[
\left(\E_U[\|F_U\|_{S_4}^2]\right)^2
\le C_{\mathrm{S2}}\max_{(U_1,U_2)\in\mathcal A(F)}
\big\|F(x_{A\cup U_1},x_{B\cup U_2})\big\|_{S_4}^4.
\]

\end{description}

Assumption (M8) is equivalent to the same bound for \(q_G\) and implies (M4). In (S2), fixing \(x_U\) gives a function \(F_U\) of \(x_A,x_B\). The condition requires that its squared $S_4$ norm, averaged over \(x_U\), cannot be arbitrarily large relative to the squared $S_4$ norm of \(F\) when the additional variables remain among its arguments, for the best split in $\mathcal A(F)$. The following lemma ensures this when the dependence on \(x_U\) involves a bounded number of functions.

\begin{lemma}\label{lem:S-filler-span}
Let $\Pcal$ be a class of laws. Suppose there is a fixed integer $M$ such that, for every $N$ and $P\in\Pcal$, each function $F$ defined in (S2) can be written as
\[
F(x_A,x_B,x_U)=\sum_{k=1}^M f_k(x_A,x_B)g_k(x_U),
\]
with square-integrable $f_k$ and $g_k$. Then (S2) holds on $\Pcal$ with $C_{\mathrm{S2}}=M^3$.
\end{lemma}
The proof is in Section~\ref{sec:S-ident} of the supplementary material.
Under these assumptions, the following theorem, proved in Section~\ref{sec:S-ident} of the supplementary material, shows that the diagnostic distinguishes when $\Upsilon_N$ is bounded away from zero from when $\Upsilon_N$ is small relative to a chosen threshold $\varepsilon_N$.

\begin{theorem}\label{thm:S-check}
Let $\Pcal$ be a fixed class of laws, and assume (M8) and (S2) on $\Pcal$, with $\sigma_{N,P}>0$ for every $P\in\Pcal$, and with $m\ge4p+2r$ and $n\ge4q+2c$ for every index $\tau$. Let $\varepsilon_N>0$ be deterministic with $\varepsilon_N\to0$ and $\varepsilon_NL_N^{1/4}\to\infty$. Then
\begin{enumerate}
\item[(i)] for each fixed $\delta>0$,
\[
\sup_{P\in\Pcal:\,\Upsilon_N(P)\ge\delta}
\P_P(\widehat\Upsilon_N\le\varepsilon_N)\longrightarrow0;
\]
\item[(ii)] for every deterministic $a_N\to0$ with $a_N\ge0$,
\[
\sup_{P\in\Pcal:\,\Upsilon_N(P)\le a_N\varepsilon_N}
\P_P(\widehat\Upsilon_N>\varepsilon_N)\longrightarrow0.
\]
\end{enumerate}
\end{theorem}

\paragraph{Combining the conservative and Gaussian tests.}\label{sec:block-switch}
Let $C_N$ be the event on which we use the Gaussian test, and let $R_{G,N}$ and $R_{S,N}$ be the Gaussian and conservative rejection events. The procedure rejects on
\begin{equation}\label{eq:hybrid-definition}
R_{H,N}=(C_N\cap R_{G,N})\cup(C_N^c\cap R_{S,N}).
\end{equation}
Even if each of two tests rejects with probability at most $\alpha$ under the null hypothesis, choosing between them using the same data can increase that probability above $\alpha$. We therefore require that $R_{S,N} \subseteq R_{G,N}$. Using Theorem~\ref{thm:S-check}, this inclusion ensures the uniform validity of this procedure, as formalized in Lemma~\ref{lem:hybrid-switch} (Sec.~\ref{sec:S-procedure} in the supplementary material).

The following proposition, also proved in Section~\ref{sec:S-procedure} of the supplementary material, applies this argument to the conservative block test. 

\begin{proposition}\label{prop:S-block-hybrid}
Let $\alpha$ and the null class $\Pcal_0$ be as in Proposition~\ref{prop:block-safe}, with its assumptions and block construction, and set $B_N=\max(A_N,S_{K_N}^2/K_N)$.
Take $A_N=\Vu^+$ and use the diagnostic defined in \eqref{eq:diagnostic-statistic}, with (M8) and (S2) on $\Pcal_0$, the sample-size requirements of Theorem~\ref{thm:S-check} for all sufficiently large $N$, and deterministic $\varepsilon_N\to0$ satisfying $\varepsilon_NL_N^{1/4}\to\infty$. Define
\[
\begin{aligned}
 C_N&=\{\widehat\Upsilon_N\le\varepsilon_N\},\\
 R_{G,N}&=\{A_N>0,\ |U_N-\theta_0|\ge z_{1-\alpha/2}\sqrt{A_N}\},\\
 R_{S,N}&=\{A_N>0,\ |U_N-\theta_0|\ge\sqrt{B_N/\alpha}\}.
\end{aligned}
\]
Then $R_{H,N}=(C_N\cap R_{G,N})\cup(C_N^c\cap R_{S,N})$ satisfies
\[
 \limsup_N\sup_{P\in\Pcal_0}\P_P(R_{H,N})\le\alpha.
\]
For every specified deterministic $a_N\ge0$ tending to zero,
\[
 \sup_{\substack{P\in\Pcal_0:\Upsilon_N(P)\le a_N\varepsilon_N}}
 \big|\P_P(R_{H,N})-\alpha\big|\longrightarrow0.
\]
Both conclusions also hold for the upper one-sided procedure, with
\[
\begin{aligned}
 R_{G,N}&=\{A_N>0,\ U_N-\theta_0\ge z_{1-\alpha}\sqrt{A_N} \},\\
 R_{S,N}&=\{A_N>0,\ U_N-\theta_0\ge\sqrt{(1-\alpha)B_N/\alpha} \}.
\end{aligned}
\]
\end{proposition}

\subsection{An alternative diagnostic under additional assumptions}\label{sec:alternative-procedure}

In general, the diagnostic and the conservative test can be chosen separately. They do not have to be those constructed above, as long as each satisfies its stated guarantees and every rejection by the conservative test would also be a rejection by the Gaussian test (Lem.~\ref{lem:hybrid-switch} in Sec.~\ref{sec:S-procedure} of the supplementary material).

Notably, the diagnostic of Section~\ref{sec:diagnostic-main} averages products over placements of up to $4p+2r$ rows and $4q+2c$ columns, so its computation is impractical except for small arrays. Under additional assumptions, we can use a simpler diagnostic for the block procedure.

Let $\Pcal$ be a class of laws with $\sigma_{N,P}>0$ for every $P\in\Pcal$, and let $\Gamma_{\mathrm{disc}}$ be a set of disconnected graphs. We assume (M4) on $\Pcal$ and the following condition:
\begin{description}
\item[(L)] (a) $s_{N,G}^2(P)\to0$ uniformly over $P\in\Pcal$ for every disconnected $G\notin\Gamma_{\mathrm{disc}}$, and \\
(b) for every $G\in\Gamma_{\mathrm{disc}}$, no union of some, but not all, of the connected components of $G$ has the same total numbers of rows and columns as a graph of $\Gamma_{\mathrm{disc}}$.
\end{description}
In particular, if $\Gamma_{\mathrm{disc}}$ contains a single graph, then (b) holds because every union of some, but not all, of its connected components has fewer vertices than that graph.

Let $F_N$ be the sum of the positive parts $\widehat V_{r,c}^+$ over the row and column counts $(r,c)$ of the graphs in $\Gamma_{\mathrm{disc}}$. Since $\Vu=\sum_{r+c\ge1}\widehat V_{r,c}$, computing $\Vu$ also gives every $\widehat V_{r,c}$, so $F_N$ requires no computation beyond that of $\Vu$. We choose the Gaussian test when $\Vu>0$ and $F_N\le\varepsilon_N\Vu$. 

These assumptions can also be exploited to improve the power of the block test in the procedure. Since $F_N$ sums the positive parts of some variance components, it estimates only part of the variance. Under (M4) and (L), $F_N/\sigma_N^2\le1+o_P(1)$ uniformly over the class (Lem.~\ref{thm:S-test} in Sec.~\ref{sec:S-chaos-procedure} of the supplementary material). Thus, $\max\{\Vu,F_N\}/\sigma_N^2\to1$ when $\Upsilon_N\to0$, as for $\Vu$. Therefore, we can replace $\Vu^+$ by $A_N=\max\{\Vu,F_N\}$ in the procedure. Since both tests reject only when $A_N>0$, the block test can then reject whenever $F_N>0$, whereas with $\Vu^+$, it never rejects when $\Vu\le0$. The resulting procedure is given by~\eqref{eq:S-floor-gate} in the supplementary material, and Proposition~\ref{prop:S-alt-procedure} proves its validity.

\section{Product form across Gaussian and non-Gaussian regimes}\label{sec:product}

We apply our results to estimation and testing based on measures of departure from product form. First, under a Poisson graphon model, we consider a ratio statistic whose limit remains Gaussian as the law approaches product form, although its convergence rate changes. Then, we give an alternative departure from product form with a non-Gaussian limit at the same rate. Finally, we consider a different measure of departure that allows other conditional cell distributions, and has a Gaussian limit at a faster rate. In the supplementary material, Section~\ref{sec:S-product} gives the additional results needed and all the proofs.

\subsection{A product-form measure of departure under a Poisson graphon model}\label{sec:dw}

Consider the Poisson graphon model
\[
Y_{ij}\mid\xi_i,\eta_j\sim\mathcal P(\lambda\wbar(\xi_i,\eta_j)),
\qquad \iint\wbar=1.
\]
The mean graphon is $w=\lambda\wbar$, where $\wbar$ is a normalized graphon, and $\lambda>0$ is the overall intensity. We write $f(x)=\int\wbar(x,y)\,\d y$ and $g(y)=\int\wbar(x,y)\,\d x$ for its row and column margins, and set $F_2=\int f^2$ and $G_2=\int g^2$. The model has product form if $\wbar(x,y)=f(x)g(y)$. Following \citet[Sec.~6.3]{leminh2025hoeffdingtype}, we measure departure from product form by
\[
\gamma(\wbar)=\|\wbar-f\otimes g\|_2^2=\|\Delta\|_2^2,
\qquad \Delta:=\wbar-f\otimes g.
\]
Thus, $\gamma(\wbar)=0$ exactly at product form. By construction, $\int\Delta(x,y)\,\d y=\int\Delta(x,y)\,\d x=0$. For a law $P$ of this model, we write $\wbar_P$ for its normalized graphon and $\gamma_P=\gamma(\wbar_P)$.

Let $U$ be the $U$-statistic obtained by symmetrizing the following $3\times3$ kernel over the rows and the columns:
\begin{equation}\label{eq:gamma-kernel}
h_\gamma=Y_{11}(Y_{11}-1)Y_{22}Y_{33}-2\,Y_{11}Y_{12}Y_{22}Y_{33}+Y_{11}Y_{12}Y_{23}Y_{33}.
\end{equation}
Its expectation is
\[
\E[U]=\lambda^4\gamma(\wbar).
\]
At product form, all projections indexed by graphs without edges vanish, whereas the projection on $K_{1,1}$ remains nonzero. We have the following Gaussian limit, proved in Section~\ref{sec:S-product} of the supplementary material.

\begin{corollary}\label{cor:dwnull}
Let $0<\lambda<\infty$. Then under product form ($\gamma=0$), $U$ has principal degree $2$, a unique principal support graph $K_{1,1}$, and
\[
(mn)^{1/2}\,U\dto\Nor\big(0,\,2\lambda^6F_2G_2\big).
\]
\end{corollary}

\subsection{Testing and estimation}\label{sec:dw-uniform}

\paragraph{Testing.}
We would like to test $\mathcal H_0:\gamma=0$ against $\mathcal H_1:\gamma>0$. Since $\lambda>0$ and $\E[U]=\lambda^4\gamma$, these hypotheses are equivalent to $\E[U]=0$ and $\E[U]>0$, respectively. Let $\Pcal_0$ be a class of product-form Poisson graphon laws with fixed $\lambda>0$, $\sup_{P\in\Pcal_0}\|\wbar_P\|_\infty<\infty$, and (M4) holding for $h_\gamma$. Let $\Vu$ be the unbiased variance estimator applied to $h_\gamma$, and $\Vu^+=\max(\Vu,0)$. Since $\E[U]=0$ under product form, the test rejects when
\[
\Vu>0
\qquad\text{and}\qquad
\frac{U}{\sqrt{\Vu^+}}>z_{1-\alpha}.
\]

The variance contributed by $K_{1,1}$ satisfies, by \eqref{eq:graph-variance-weight},
\[
\var(T_{K_{1,1}})=\frac{\E[p_{K_{1,1}}^2]}{mn}
\ge\frac{2\lambda^6}{mn}\ge cN^{-2},
\]
uniformly over $\Pcal_0$, since $\lambda>0$ is fixed and $m,n\asymp N$ (Prop.~\ref{prop:edge} in the supplementary material). On the other hand, the disconnected variance $\sigma_{\mathrm{disc}}^2=\sum_{G\,\disc}\var(T_G)$ is $O(N^{-3})$ uniformly over $\Pcal_0$ (proof of Lem.~\ref{lem:S-poisson-disconnected} in the supplementary material).
Thus, $D_N=O(N^{-1})$ uniformly, and Proposition~\ref{prop:gammaD} and Theorem~\ref{thm:contraction} therefore give both consistency of $\Vu$ relative to $\var_P(U)$ and $U/\sqrt{\var_P(U)}\dto\Nor(0,1)$, uniformly over $\Pcal_0$. By Slutsky's lemma, the rejection probability tends to $\alpha$ uniformly over $\Pcal_0$.

\paragraph{Estimation.} 
The statistic $U$ estimates $\lambda^4\gamma$, rather than $\gamma$ itself.
Since $\E[\bar Y]=\lambda$, we estimate the normalized departure $\gamma$ by
\[
\widehat\gamma_U=\frac{U}{\bar Y^4}.
\]

To construct confidence intervals for $\gamma$ that remain valid as
the law approaches product form, we must show that the error from
linearizing the ratio $\widehat\gamma_U$ is negligible compared with the standard deviation of the linear approximation. The linear approximation is
\begin{equation}\label{eq:gamma-linearization}
\ell_{N,P}
=\lambda^{-4}
\{(U-\lambda^4\gamma_P)
-4\gamma_P\lambda^3(\bar Y-\lambda)\},
\qquad
\sigma_{N,P}^2=\var_P(\ell_{N,P}).
\end{equation}
This is the centered $U$-statistic associated with the population
kernel $\lambda^{-4}\widetilde h_{\ell,P}$, where
\[
\widetilde h_{\ell,P}
=h_\gamma-4\gamma_P\lambda^3Y_{11}.
\]

To estimate $\sigma_{N,P}^2$, we first replace the unknown quantities
$\lambda$ and $\gamma_P$ in $\lambda^{-4}\widetilde h_{\ell,P}$
by $\bar Y$ and $\widehat\gamma_U=U/\bar Y^4$. This gives the estimated
linearization kernel
\[
\bar Y^{-4}
\{h_\gamma-4\widehat\gamma_U\bar Y^3Y_{11}\}
=
\bar Y^{-4}h_\gamma-4U\bar Y^{-5}Y_{11}.
\]
Then, we apply the unbiased variance estimator to this kernel,
where the coefficients $\bar Y^{-4}$ and $4U\bar Y^{-5}$ are only computed once and fixed throughout the calculation.
Denote the resulting variance estimate by
$\Vu^{(\mathrm{lin})}$ and its positive part by
$\Vu^{(\mathrm{lin}),+}$.

The following result, proved in Section~\ref{sec:S-product} of the supplementary material, gives uniform inference with $\Vu^{(\mathrm{lin}),+}$, over classes approaching product form and classes where the first-order projections satisfy a lower bound.

\begin{corollary}\label{cor:dwuniform}
For each $N$, let $\Pcal_N$ be a class of laws generated by the Poisson graphon model, with fixed $\lambda>0$ and $\sup_N\sup_{P\in\Pcal_N}\|\wbar_P\|_\infty<\infty$. Assume that the graph projections $p_G^{\widetilde\ell}$ of $\widetilde h_{\ell,P}$ satisfy (M4) on $\{(N,P):P\in\Pcal_N\}$, and that one of the following assumptions holds:
\begin{description}
\item[(a)] $\sup_{P\in\Pcal_N}\gamma_P\to0$,
\item[(b)] there is $c>0$ such that $\E[(p^{\widetilde\ell}_{K_{1,0}})^2]+\E[(p^{\widetilde\ell}_{K_{0,1}})^2]\ge c\,\gamma_P^2$ for every $N$ and every $P\in\Pcal_N$.
\end{description}
Then $\Vu^{(\mathrm{lin}),+}/\sigma_{N,P}^2\to1$ in probability, and the studentized statistic $\widehat Z_N=(\widehat\gamma_U-\gamma_P)/(\Vu^{(\mathrm{lin}),+})^{1/2}$ satisfies $\widehat Z_N\dto\Nor(0,1)$ uniformly over $P\in\Pcal_N$.
\end{corollary}

For $\Vu^{(\mathrm{lin})}>0$, the corollary gives the confidence
interval
\[
\left[
\widehat\gamma_U-z_{1-\alpha/2}\sqrt{\Vu^{(\mathrm{lin}),+}},
\quad
\widehat\gamma_U+z_{1-\alpha/2}\sqrt{\Vu^{(\mathrm{lin}),+}}
\right],
\]
with coverage tending to $1-\alpha$ uniformly over the stated classes. Assumption (a) covers laws approaching product form, whereas (b) also allows laws that remain away from product form, provided that their first-order projections satisfy the lower bound. In contrast, Section~\ref{sec:planted} gives an alternative path that satistifies neither assumption, with a non-Gaussian limiting law.

\subsection{An alternative departure with a non-Gaussian limit}\label{sec:planted}

Let $s(x)=1$ on $[0,1/2)$ and $s(x)=-1$ on $[1/2,1]$, and set
\[
\wbar_\varepsilon(x,y)=1+\varepsilon s(x)s(y),\qquad 0<\varepsilon\le1.
\]
This two-group model corresponds to a block model with two equally sized row groups and two equally sized column groups, uniform margins and $\gamma(\wbar_\varepsilon)=\varepsilon^2$. We develop this example in Proposition~\ref{prop:S-two-group} (Sec.~\ref{sec:S-two-group} in the supplementary material), which shows that its first-order projections vanish, but not the ones indexed by $K_{1,1}$, $K_{2,0}$, $K_{0,2}$ and $K_{1,0}\sqcup K_{0,1}$. Since the last three graphs are disconnected, Theorem~\ref{thm:chaos} implies a non-Gaussian limiting law for the linear approximation $\ell_{N,P}$ at every fixed $\varepsilon>0$. Then, the ratio $\widehat\gamma_U$ has the same limiting law because its linearization error is negligible compared with $N^{-1}$. This is an example showcasing that neither the convergence rate nor the degeneracy order determines whether the limiting law is Gaussian.

\subsection{Models beyond Poisson}\label{sec:d3}

Previously, the kernel $h_\gamma$ worked because of the Poisson identity $\E[Y(Y-1)\mid\cdot]=\lambda^2\wbar^2$. If the cell distribution is not Poisson, we can instead use products of distinct cells to define a new measure of departure. Given the row and column variables, the cells remain independent, so the expectations of these products depend only on the conditional means.

Let the mean graphon $w$ of Section~\ref{sec:framework-arrays} be bounded, with singular values $s_k\ge0$. The cycle densities $t(C_4,w)=\sum_ks_k^4$ and $t(C_8,w)=\sum_ks_k^8$ can be used to define the departure
\[
\Psi(w):=t(C_4,w)^2-t(C_8,w)\ge0,
\]
which is zero if and only if $w$ has product form. An unbiased estimator for $\Psi(w)$ is the $U$-statistic $\widehat\Psi$ with the $4\times4$ kernel
\[
h_\Psi=\frac1{18}\sum_{P}\prod_{e\in P}Y_e-\frac1{72}\sum_{Q}\prod_{e\in Q}Y_e,
\]
where $P$ runs over the edge sets of the $18$ copies of $C_4\sqcup C_4$, and $Q$ over those of the $72$ copies of $C_8$ in the $4 \times 4$ subarray. Figure~\ref{fig:psi-patterns} of the supplementary material illustrates the two families.

At product form, all projections of degree at most four vanish except possibly $p_{C_4}$, as proved in Proposition~\ref{prop:d3} (Sec.~\ref{sec:S-degree-four} in the supplementary material). Since $C_4$ is connected, our uniform inference theorem (Thm.~\ref{thm:contraction}) gives the following corollary, also proved there.

\begin{corollary}\label{cor:d3}
Let $\Pcal_0$ be a class of laws with a product-form mean graphon $w=\lambda f\otimes g$. Assume that (M4) holds for $h_\Psi$ on $\Pcal_0$, and that $\sup_{\Pcal_0}\E[h_\Psi^2]<\infty$ and $\inf_{\Pcal_0}\|p_{C_4}\|_2>0$. Then at every law in $\Pcal_0$, the $4$-cycle $C_4$ is the unique principal support graph of $\widehat\Psi$, and the studentized statistic $\widehat Z_N=\widehat\Psi/(\Vu^+)^{1/2}$ satisfies $\widehat Z_N\dto\Nor(0,1)$, uniformly over $\Pcal_0$.
\end{corollary}

Thus, the test that rejects when $\widehat Z_N>z_{1-\alpha}$ has a rejection probability tending to $\alpha$ uniformly over $\Pcal_0$. The standard deviation $\sigma_N$ of $\widehat\Psi$ is now of order $N^{-2}$, but inference still uses the same unbiased variance estimator.

\section{Simulation studies and illustration}\label{sec:numerics}
The main results of this paper are asymptotic. In this section, we examine their finite-sample behavior under different models and array sizes. First, we compare variance estimates and the coverage of Gaussian confidence intervals as the degeneracy changes. Then, we examine the rejection probabilities and thresholds of the conservative procedures. Finally, we apply the product-form test to UN roll-call data. 

\subsection{Confidence intervals across degeneracy regimes}\label{sec:transition-evidence}\label{sec:conc-numerics}
Each simulation replicate consists of the analysis of one independently generated array. All simulations use square arrays of size $L=m=n=N/2$. 
\paragraph{Changing degeneracy order.}
The following simulation examines the coverage of Gaussian confidence intervals along the drifting laws of Example~\ref{ex:k12}. We consider $Y_{ij}=t_NA_i+e_{ij}$, where $A_i$ is uniform on $[-\sqrt3,\sqrt3]$ and the cell noise is centered standard exponential. With kernel $y_1y_2$, the parameter of interest is $t_N^2$. Along $t_N=cL^{-1/2}$, the principal degree is one for every $c>0$ and three at $c=0$. Since $D_N=0$, the Gaussian confidence intervals remain uniformly valid. For every fixed $c>0$, the variance contributions of degrees one, two and three all have order $L^{-3}$, with shares depending on $c$.

For each $L\in\{16,32,64,128\}$ and $c\in\{0,0.5,1,2,4,8\}$, we use $500$ independent simulation replicates. We compare the positive part $\Vu^+$ of the unbiased variance estimator, the row-column jackknife $\Vrcj$, the pigeonhole bootstrap variance $\widehat V_{\mathrm{pig}}$ \citep{owen2007pigeonhole}, and the positive part of the first-order estimator $\widehat V_{\mathrm{fo}}^+:=\max(\widehat V_{\mathrm{fo}},0)$, where $\widehat V_{\mathrm{fo}} := \widehat V_{1,0}+\widehat V_{0,1}$, and $\widehat V_{1,0}$ and $\widehat V_{0,1}$ are the unbiased estimators~\eqref{eq:variance-component-estimator} of the variance components of the projections indexed by $K_{1,0}$ and $K_{0,1}$. The pigeonhole bootstrap uses $500$ resamples. We construct $95\%$ confidence intervals using Gaussian critical values.

\begin{figure}[!ht]
\centering
\includegraphics[width=\textwidth]{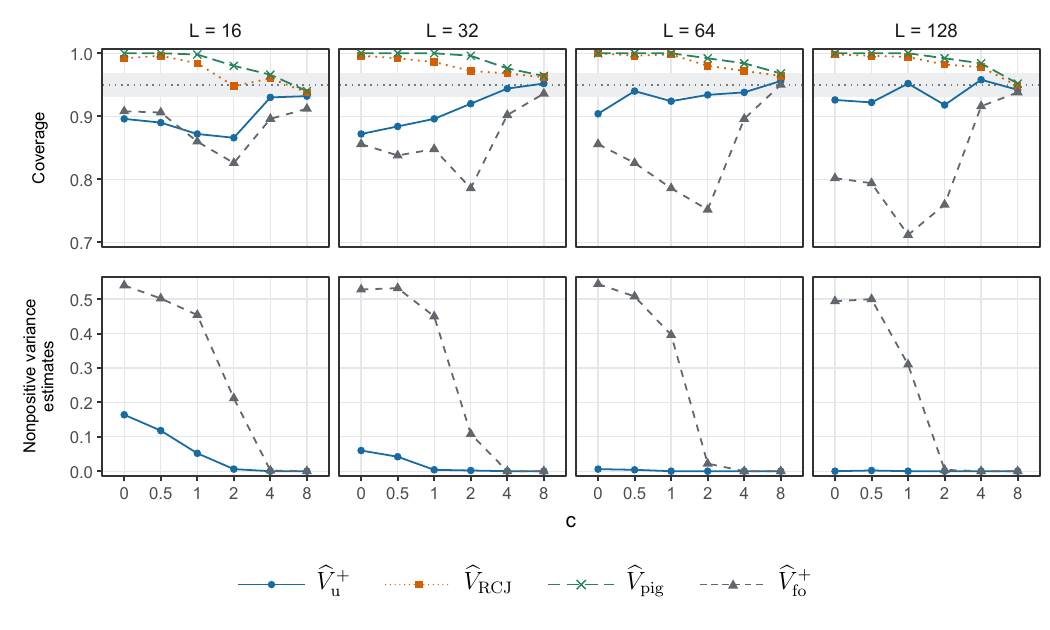}
\caption{Gaussian confidence intervals at $95\%$ along $t_N=cL^{-1/2}$, with centered exponential noise at $L\in\{16,32,64,128\}$, for different values of $c$ (equally spaced). Top: coverage for the four studentized statistics. The dotted grey line is the target $95\%$ coverage. The shaded band is the central $95\%$ range of $\operatorname{Bin}(500,0.95)/500$. Bottom: frequency of nonpositive $\Vu$ and $\widehat V_{\mathrm{fo}}$ estimates, before taking positive parts.}\label{fig:S-drift}
\end{figure}

Figure~\ref{fig:S-drift} gives the coverage in the top row and the frequency of nonpositive unbiased and first-order variance estimates in the bottom row. At $c=8$, where the first-order contribution to the variance dominates, coverage remains near $95\%$ for all four estimators as the size increases. Near $c=0$, where the third-order contribution becomes dominant, $\Vu^+$ undercovers at the smaller sizes and moves closer to $95\%$ at $L=128$, whereas $\Vrcj$ and $\widehat V_{\mathrm{pig}}$ lead to overcoverage at all sizes. The bottom row shows that nonpositive $\Vu$ estimates become rare as the array grows. In contrast, the first-order estimator severely undercovers for $c\le 2$, despite nonpositive estimates in about half the replicates at $c=0$ and $0.5$, for which the confidence interval is the whole parameter space. 

\paragraph{Estimating the product-form departure.}\label{sec:gamma-evidence}
The following simulation illustrates the estimation of the departure $\gamma$ as the law approaches product form, as considered in Corollary~\ref{cor:dwuniform}. We use the Poisson stochastic blockmodel with two equally sized row groups, two equally sized column groups, and group means proportional to $\left(\begin{smallmatrix}4+2\varepsilon&2\\2&1+2\varepsilon\end{smallmatrix}\right)$, normalized to intensity $\lambda=9/4$. The model has product form at $\varepsilon=0$.

For each $L\in\{16,32,64\}$ and $\varepsilon\in\{0,0.1,0.25,0.5,1,2\}$, we use $500$ independent simulation replicates and compute the ratio $\widehat\gamma_U$. We compute $\Vu^{(\mathrm{lin}),+}$, the positive part of the unbiased variance estimator applied to the estimated linearization kernel (Sec.~\ref{sec:dw-uniform}), approximating it with $2000$ random half-splits (Sec.~\ref{sec:S-split}). We compare it with the estimates of the row-column jackknife $\Vrcj$ and the pigeonhole bootstrap $\widehat V_{\mathrm{pig}}$, using Gaussian critical values for the $95\%$ confidence intervals. We apply the row-column jackknife to the ratio itself. The pigeonhole bootstrap uses $200$ resamples. 

\begin{figure}[!ht]
\centering
\includegraphics[width=\textwidth]{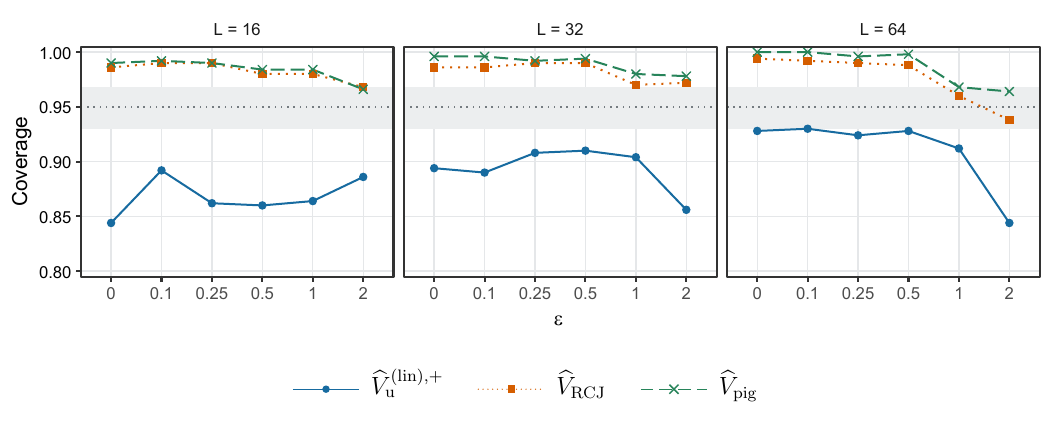}
\caption{Coverage of $95\%$ Gaussian confidence intervals for $\gamma$ for the model of Section~\ref{sec:gamma-evidence} at $L\in\{16,32,64\}$, for different values of $\varepsilon$ (equally spaced). The dotted grey line is the target $95\%$ coverage. The shaded band is the central $95\%$ range of $\operatorname{Bin}(500,0.95)/500$.}\label{fig:transition}
\end{figure}

Figure~\ref{fig:transition} gives the coverage of the confidence intervals obtained with the three variance estimates. Near product form ($\gamma \approx 0$), we observe undercoverage with $\Vu^{(\mathrm{lin}),+}$ but this reduces as the size increases, whereas $\Vrcj$ and $\widehat V_{\mathrm{pig}}$ lead to overcoverage. However, away from product form, undercoverage with $\Vu^{(\mathrm{lin}),+}$ remains visible at $L=64$.

\paragraph{Testing for product form.}\label{sec:psi-numerics}
The following simulation illustrates the product-form test of Corollary~\ref{cor:d3}, for which the principal support graph of $\widehat\Psi$ under the null hypothesis is $C_4$. We draw row and column propensities $f_i$ and $g_j$ independently from the uniform distribution on $[0.5,1.5]$, and independent interaction factors $a_i$ and $b_j$ from the uniform distribution on $[-1,1]$. Conditional on these variables, the cells are independent, with
\[
Y_{ij}\mid \xi_i,\eta_j\sim\operatorname{Poisson}\!\left(3\max\{f_i g_j+\delta a_i b_j,0.05\}\right).
\]
At $\delta=0$, $f_ig_j\ge0.25>0.05$, so the truncation is inactive and the model has product form. 

For each $L\in\{16,32,64\}$, we use $500$ independent simulation replicates at each $\delta\in\{0,0.25,\allowbreak 0.5,0.75,1\}$. We compute the test statistics based on the unbiased $\Vu^+$, the row-column jackknife $\Vrcj$ and the pigeonhole $\widehat V_{\mathrm{pig}}$ estimates. The unbiased variance estimate uses $10000$ random half-splits at $L=16,32$ and $4000$ at $L=64$. The pigeonhole variance uses $200$ resamples. All three tests use a one-sided Gaussian critical value at level $5\%$. 

\begin{figure}[!ht]
\centering
\includegraphics[width=\textwidth]{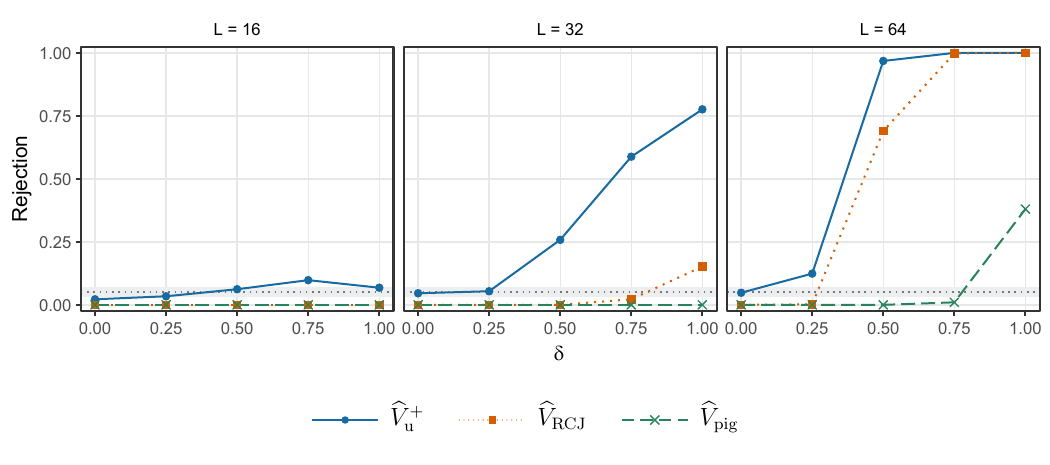}
\caption{Rejection frequencies of the one-sided product-form tests at $5\%$, with $L\in\{16,32,64\}$, for different values of $\delta$. Each setting uses $500$ replicates. The dotted grey line is the target $5\%$ rejection probability. The shaded band is the central $95\%$ range of $\operatorname{Bin}(500,0.05)/500$.}\label{fig:psi-power}
\end{figure}

Figure~\ref{fig:psi-power} gives the rejection frequencies under the null hypothesis and the alternative laws. At $\delta=0$, the rejection frequency using $\Vu^+$ increases from $2.2\%$ at $L=16$ to $4.6\%$ and $4.8\%$ at $L=32$ and $64$. In contrast, the jackknife and pigeonhole tests never reject. At $\delta \ge 0.25$, all three tests have little power at $L=16$, but a clear separation appears at larger sizes, with the $\Vu^+$ test being the most powerful. 

\subsection{Finite-sample behavior of conservative procedures}\label{sec:fallback-experiment}
The following simulation compares the coverage, power and interval widths of Gaussian inference (without a diagnostic), the block-average test (Rem.~\ref{rem:block-average}), and the block procedure using the diagnostic based on $F_N$ (Sec.~\ref{sec:alternative-procedure}).

We consider the factor model:
\[
Y_{ij}=aA_i+ \frac{t_N}{\sqrt r}\sum_{k=1}^r U_{ik}V_{jk}+e_{ij},
\]
where $a, t_N\ge0$, and $A_i$, $U_{ik}$, $V_{jk}$ and $e_{ij}$ are independent standard normal random variables. Let the null hypothesis be $\mathcal{H}_0: a = 0$ and the alternative $\mathcal{H}_1: a \neq 0$. To test $\mathcal{H}_0$, we consider the $U$-statistic averaging $(Y_{ij}-Y_{kj})(Y_{il}-Y_{kl})$ over row and column pairs and targeting $\theta=2a^2$. The exact variance $\sigma_0^2$ under $\mathcal{H}_0$ is given in \eqref{eq:S-model-variance} in Section~\ref{sec:S-model-hybrid} of the supplementary material, alongside a verification of the assumptions needed to use the procedure.

We parameterize the model with $\vartheta \ge 0$, such that $a = \sqrt{\vartheta \sigma_0/2}$, and therefore $\theta = \vartheta \sigma_0$. We consider three settings: (i) no interaction, $t_N=0$, (ii) decreasing rank-four interaction, $r=4$, $t_N=2L^{-1/4}$, (iii) fixed rank-one interaction, $r=1$, $t_N=1$. For every fixed $\vartheta$, $D_N$ tends to $0$, $4/5$ and $1$, respectively (Lem.~\ref{lem:S-model-alternative} in Sec.~\ref{sec:S-model-limits} of the supplementary material). Thus, the limit is Gaussian in setting (i), but Gaussian inference is not justified in settings (ii) and (iii).

For each setting, $L\in\{64,128,256,512,1024\}$ and $\vartheta\in\{0,2,4,8,16,32\}$, we use $3000$ independent simulation replicates. We compute the $95\%$ confidence-interval coverage frequency at the true parameter $\theta$ and the rejection frequency of the null hypothesis $\theta_0=0$ for the block procedure, the block-average test, and the Gaussian test. The block procedure uses the alternative diagnostic of Section~\ref{sec:alternative-procedure}, with \(\varepsilon_N=1/\lceil L^{1/4}\rceil\) and \(F_N\) given by \eqref{eq:S-model-floor} in the supplementary material. Both Gaussian tests use the cutoff $1.96$. The block procedure and the block-average test both use $K_N=\lfloor L^{1/3}\rfloor$ consecutive disjoint diagonal blocks, with side length $b=\lfloor L/K_N\rfloor$. The block-average test uses the $0.975$ quantile of Student's $t_{K_N-1}$ distribution without a diagnostic.

\begin{figure}[!ht]
\centering
\includegraphics[width=\textwidth]{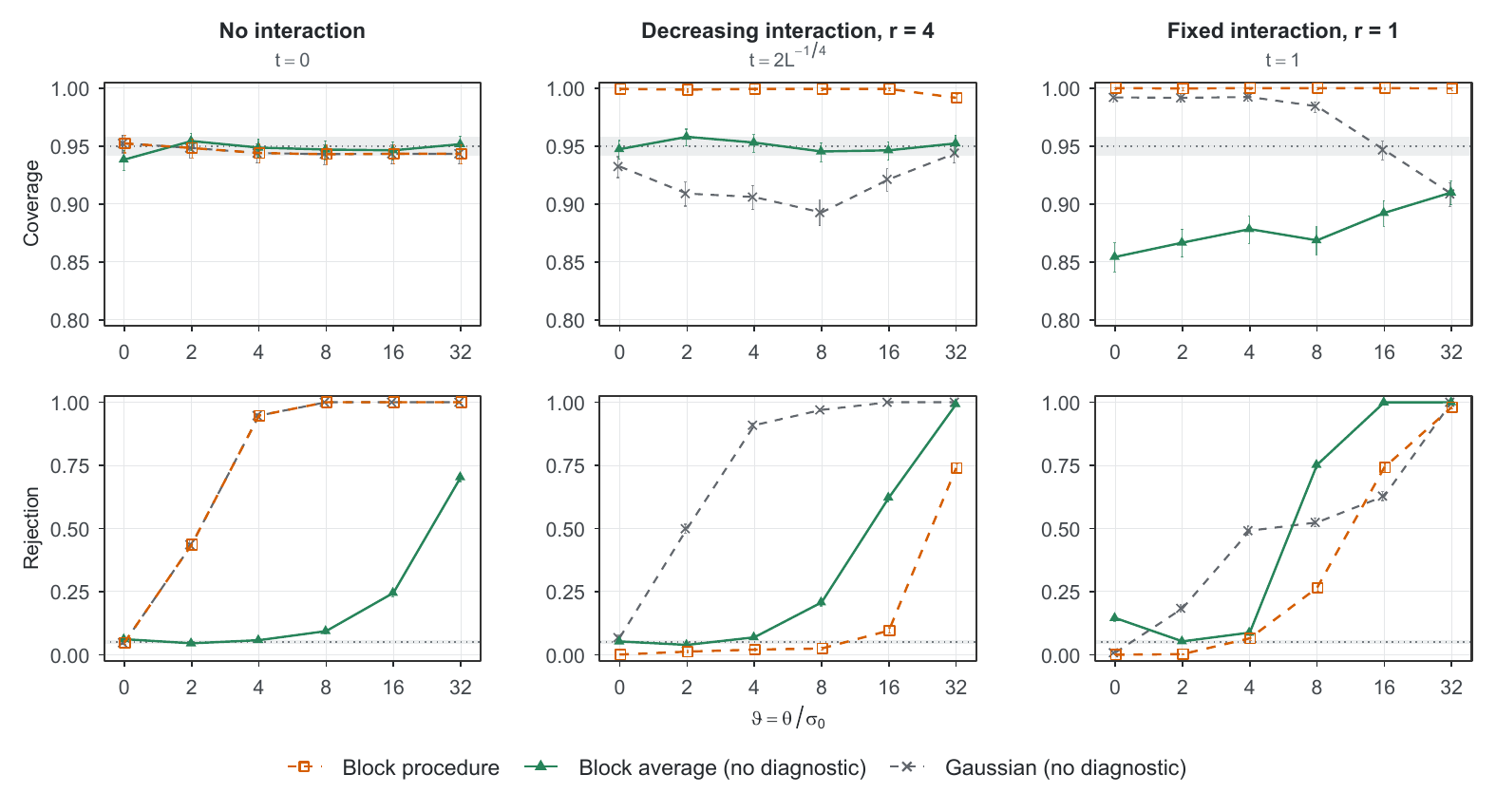}
\caption{Coverage of $95\%$ confidence intervals (top)
and rejection frequencies for tests of zero (bottom), at $L=1024$, for different values of $\vartheta=\theta/\sigma_0$ (equally spaced). The bars are pointwise $95\%$ Wilson intervals. The dotted grey line is the target probability $p$ and the shaded bands are the central $95\%$ ranges of $\operatorname{Bin}(3000,p)/3000$, with $p=0.95$ for coverage and $p=0.05$ for rejection under the null hypothesis.}\label{fig:block-chaos}
\end{figure}

Figure~\ref{fig:block-chaos} shows all six levels considered for $\vartheta$ at $L=1024$. The top and bottom rows respectively give the coverage and rejection frequencies.
Without interaction, the block procedure retains the power of the Gaussian test, whereas the block-average test has low power. Under decreasing rank-four interaction, the Gaussian confidence intervals undercover. By contrast, the block procedure overcovers, at the cost of lower power. In this setting, the block-average test has coverage close to $95\%$ and is more powerful than the block procedure at every $\vartheta > 0$. Under fixed rank-one interaction, the block procedure intervals contain the true parameter in almost every replicate. The block-average test rejects the true null hypothesis in $15\%$ of replicates and its intervals undercover heavily. The Gaussian confidence intervals also undercover for $\vartheta = 32$. The block procedure is approximately as powerful as the Gaussian test for $\vartheta \ge 16$.

Figures~\ref{fig:S-block-coverage}, \ref{fig:S-block-power} and~\ref{fig:S-block-widths} in Section~\ref{sec:S-figures} of the supplementary material report every size, with, additionally, power and interval half-width comparisons. We see that no method maintains the target coverage across all cases. However, at \(L=1024\), the block procedure is the only method that avoids substantial undercoverage across all three settings, although it is conservative and less powerful in some situations.

\subsection{Illustration: roll-call structure beyond additive effects}\label{sec:data}
To illustrate our results, we consider United Nations General Assembly roll-call votes \citep{voeten2013data,michalskasmith2019telling}. First, we would like to test whether the votes exhibit a product-form structure. After that, we would like to investigate whether the classical Rasch model \citep{rasch1960probabilistic,clinton2004statistical}, modelling additive state and roll-call effects, accounts for the observed departure from product form. In this dataset, member states and roll calls are the rows and the columns respectively. The arrays consist of binary entries, where one means a vote in favor, zero means a vote against, an abstentions or an absence. 

\citet{michalskasmith2019telling} supplies $68$ yearly arrays. For each, we fit the Rasch model by maximum likelihood, with probabilities $\operatorname{logit}^{-1}(a_i+b_j)$, where $a_i$ is the state effect and $b_j$ is the roll-call effect. We analyze these arrays and simulate $200$ arrays of the same dimensions from the fitted model. We apply the one-sided test of Corollary~\ref{cor:d3}, which rejects when $\widehat\Psi/(\Vu^+)^{1/2}>z_{0.95}$, using $200$ random half-splits per array to estimate the variance. We also compute the normalized departure $\widehat\Psi/\widehat t(C_4)^2$, where $\widehat t(C_4)$ is the $U$-statistic averaging $Y_{ij}Y_{il}Y_{i'j}Y_{i'l}$ over distinct rows $i\ne i'$ and distinct columns $j\ne l$, an unbiased estimator of $t(C_4,w)$.

\begin{figure}[!ht]
\centering
\includegraphics[width=0.75\textwidth]{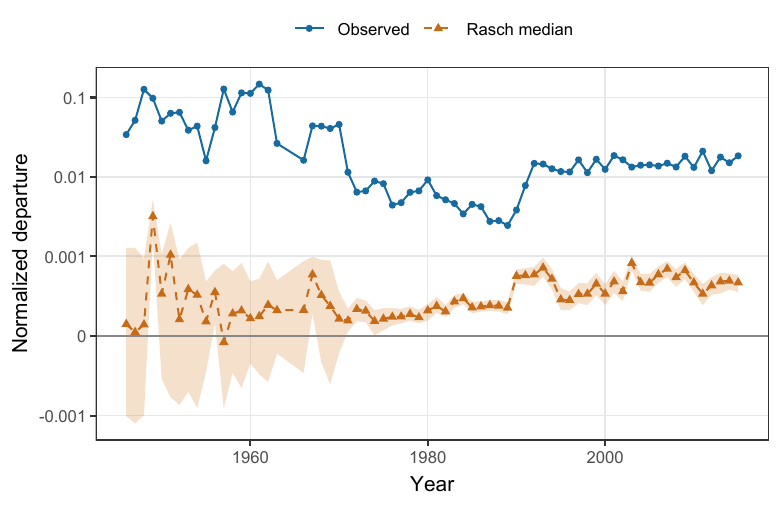}
\caption{Estimated normalized departure from product form, $\widehat\Psi/\widehat t(C_4)^2$, for the observed UN roll-call arrays and the unclipped fitted Rasch model, using all states and roll calls in the $68$ supplied yearly arrays. The solid blue curve with circles gives the observed values. The dashed orange curve with triangles and the shaded band give the median of the $200$ Rasch simulations per year. The shaded band is their central $95\%$ range. The vertical axis uses the transformation \(\operatorname{asinh}(y/0.0002)/\log 10\), which is approximately linear near zero and logarithmic farther away.}\label{fig:un}
\end{figure}

The test rejects product form in all $68$ observed yearly arrays, compared with $70.7\%$ of the $13600$ Rasch simulations. Thus, the Rasch model generates departures from product form that are detectable by this test. However, it does not reproduce the departures found in the observed arrays. Figure~\ref{fig:un} gives the normalized departure $\widehat\Psi/\widehat t(C_4)^2$ for the observed yearly arrays and the Rasch simulations. These results suggest that there may be structure, other interactions between states and roll calls, beyond the fitted additive effects. 

\FloatBarrier

\section*{Acknowledgments}
The author would like to thank Sophie Donnet, François Massol and Stéphane Robin for many fruitful discussions
and insights. The author was supported by the Australian Research Council (grant number DP250100860). The author used ChatGPT and Claude: (1) to explore proof strategies, (2) to review the proofs for correctness, (3) to assist with coding, and (4) to copyedit the paper. The author takes full responsibility for the content of the paper. More specifically: 
\begin{enumerate}
    \item AI suggested the use of a de Jong-type central limit theorem and replacing the sufficient condition \(D_N\to0\) initially considered with the weaker condition \(\Upsilon_N\to0\) in Theorem~\ref{thm:contraction}. AI helped adapt Janson and Nowicki's proof strategy for the general fixed-law limit of generalized statistics to the two-way setting of the paper. AI proposed the additional assumption (L) and the simpler diagnostic that applies under it.
    \item For every proof, AI was prompted to check its correctness and improve the exposition. Most of the time, it identified steps hidden behind a terse exposition, and on a few occasions, it identified gaps in the arguments. In these cases, it suggested solutions, sometimes overcomplicated. The author evaluated these suggestions and decided how to revise the proofs.
    \item In the simulation studies and the UN roll-call illustration, the code was partially written with AI assistance. 
    \item After the initial draft, AI was used to refine the language and to streamline the notation, in order to improve exposition.
\end{enumerate}
The author takes full responsibility for the content of the paper. 

\section*{Supplementary material}
The supplement contains the proofs of the results of the paper, auxiliary results and additional experimental results. 

\bibliographystyle{apalike}
\bibliography{biblio}

\ifinlinesupp
\clearpage
\section*{Supplementary material}

\setcounter{section}{0}
\setcounter{equation}{0}
\setcounter{figure}{0}
\setcounter{table}{0}
\setcounter{theorem}{0}
\setcounter{lemmaL}{0}
\setcounter{remark}{0}
\setcounter{assumption}{0}
\renewcommand{\thesection}{S\arabic{section}}
\renewcommand{\theequation}{S\arabic{equation}}
\renewcommand{\thefigure}{S\arabic{figure}}
\renewcommand{\thetable}{S\arabic{table}}
\renewcommand{\thetheorem}{\thesection.\arabic{theorem}}
\makeatletter
\@addtoreset{theorem}{section}
\let\c@lemmaL\c@theorem
\let\c@remark\c@theorem
\let\c@assumption\c@theorem
\makeatother
\renewcommand{\thelemmaL}{\thetheorem}
\renewcommand{\theremark}{\thetheorem}
\renewcommand{\theassumption}{\thetheorem}
\renewcommand{\theHsection}{supp.\arabic{section}}
\renewcommand{\theHequation}{supp.\arabic{equation}}
\renewcommand{\theHfigure}{supp.\arabic{figure}}
\renewcommand{\theHtable}{supp.\arabic{table}}
\renewcommand{\theHtheorem}{supp.\arabic{section}.\arabic{theorem}}
\renewcommand{\theHlemmaL}{supp.\arabic{section}.\arabic{lemmaL}}
\renewcommand{\theHremark}{supp.\arabic{section}.\arabic{remark}}
\renewcommand{\theHassumption}{supp.\arabic{section}.\arabic{assumption}}
This supplement contains the proofs of the results of the paper, auxiliary results and additional experimental results. Section~\ref{sec:S-part-projections} sets up the graph projections and the counting lemma used throughout. It also proves the identities for the row-column jackknife and the unbiased variance estimator. Sections~\ref{sec:S-part-uniform} and~\ref{sec:S-part-limits} give the proofs for Section~\ref{sec:results}. Section~\ref{sec:S-part-procedures} gives the proofs and alternative procedures for Section~\ref{sec:hybrid}. Section~\ref{sec:S-applications} gives the proofs for Section~\ref{sec:product}, and Section~\ref{sec:S-numerics} gives the derivations and some additional figures for Section~\ref{sec:fallback-experiment}. Table~\ref{tab:S-notation} provides an index of notation, with pointers to the definitions.

\begingroup
\small
\setlength{\tabcolsep}{4pt}
\renewcommand{\arraystretch}{1.12}
\setlength{\LTpre}{\medskipamount}
\setlength{\LTpost}{\medskipamount}
\newcommand{\ntrow}[3]{#1 & #2 & #3\\}
\newcommand{\ntgroup}[1]{\midrule\multicolumn{3}{@{}l}{\textbf{#1}}\\*\addlinespace[2pt]}
\begin{longtable}{@{}>{\raggedright\arraybackslash}p{.20\dimexpr\textwidth-4\tabcolsep\relax}>{\raggedright\arraybackslash}p{.66\dimexpr\textwidth-4\tabcolsep\relax}>{\raggedright\arraybackslash}p{.14\dimexpr\textwidth-4\tabcolsep\relax}@{}}
\caption{Notation used in the supplementary material. References are to the main text unless their number begins with S.}\label{tab:S-notation}\\
\toprule
Symbol & Meaning & Definition\\
\midrule
\endfirsthead
\multicolumn{3}{@{}l}{\tablename~\thetable. Notation (continued)}\\
\toprule
Symbol & Meaning & Definition\\
\midrule
\endhead
\midrule
\multicolumn{3}{r@{}}{Continued on the next page}\\
\endfoot
\bottomrule
\endlastfoot
\multicolumn{3}{@{}l}{\textbf{Arrays and subarrays}}\\*
\ntrow{$m,n$; $p,q$}{Array dimensions; kernel dimensions.}{\eqref{eq:ustat}}
\ntrow{$h$ (array kernel)}{Kernel symmetric separately in its $p$ rows and $q$ columns.}{\eqref{eq:ustat}}
\ntrow{$N,L_N,\rho$}{$N=m+n$, $L_N=\min(m,n)$, and $m/N\to\rho\in(0,1)$.}{Sec.~\ref{sec:framework-arrays}}
\ntrow{$\xi_i,\eta_j,\zeta_{ij}$; $\phi_N$}{Independent row, column and cell AHK variables; $Y_{ij}=\phi_N(\xi_i,\eta_j,\zeta_{ij})$.}{Sec.~\ref{sec:framework-arrays}}
\ntrow{$W=(W_1,W_2)$; $H_W$}{Row and column sets of a subarray; its kernel value $H_W=h(Y_W)$.}{Sec.~\ref{sec:S-subarrays}}
\ntrow{$K$}{Number of subarrays, $K=\binom mp\binom nq$.}{Sec.~\ref{sec:S-subarrays}}
\ntrow{$o(W,W')$}{Pair $(|W_1\cap W'_1|,|W_2\cap W'_2|)$ of numbers of shared rows and columns.}{Sec.~\ref{sec:S-subarrays}}
\ntrow{$K(a,b)$}{Number of ordered subarray pairs with $o(W,W')=(a,b)$.}{Sec.~\ref{sec:S-subarrays}}
\ntrow{$U_N,\theta_N,\sigma_N^2$}{$U$-statistic, expectation and variance.}{Sec.~\ref{sec:graph-decomposition}}
\ntrow{$\Theta(a,b)$; $\Thetah(a,b)$}{Expected product of two kernel values sharing $a$ rows and $b$ columns; its average over such pairs.}{\eqref{eq:kernel-product-averages}}
\ntrow{$\Pcal,\Pcal_0,\Pcal_N$}{Class of laws, null class, and a class that may depend on $N$.}{Sec.~\ref{sec:framework-arrays}}
\ntgroup{Graphs and vertex sets}
\ntrow{$\Gamma$; $G$; $\Lambda$}{Nonempty graph classes within $K_{p,q}$; a representative graph; a labeled copy in the array.}{Sec.~\ref{sec:graph-decomposition}}
\ntrow{$V_1(G),V_2(G),E(G)$; $V(G)$}{Row vertices, column vertices, edges; all vertices.}{Sec.~\ref{sec:graph-decomposition}; Sec.~\ref{app:wick-products}}
\ntrow{$r(G),c(G),d(G),b(G)$}{Numbers of rows, columns, vertices, and connected components.}{Sec.~\ref{sec:graph-decomposition}}
\ntrow{$d$; principal support graphs}{The smallest $d(G)$ for which $p_G\ne0$; the graphs attaining it.}{Sec.~\ref{sec:graph-decomposition}}
\ntrow{$K_{r,c},C_{2k},\sqcup$}{Complete bipartite graph, bipartite cycle of length $2k$, and disjoint union.}{Sec.~\ref{sec:graph-decomposition}; Sec.~\ref{sec:S-degree-four}}
\ntrow{$(x)_k$; $\operatorname{Aut}(G)$}{Falling factorial $x(x-1)\cdots(x-k+1)$; graph automorphism group preserving the row and column vertex sets.}{Secs.~\ref{sec:graph-decomposition}, \ref{sec:notation}}
\ntrow{$N(\Lambda)$; $A=(A_1,A_2)$}{$N(\Lambda)=(V_1(\Lambda),V_2(\Lambda))$ is a pair of vertex sets; $A_1$ contains rows and $A_2$ columns.}{Sec.~\ref{app:projection-notation}}
\ntrow{$\ell(G),\ell(A)$}{$\ell(G)=(r(G),c(G))$ and $\ell(A)=(|A_1|,|A_2|)$.}{Sec.~\ref{sec:S-variance}}
\ntrow{$A\cap A',A\cup A',A\subseteq A'$}{Componentwise intersection, union and inclusion of vertex-set pairs.}{Sec.~\ref{sec:S-variance}}
\ntrow{$A\sim A'$}{The pairs share a row or a column.}{\eqref{eq:variance-decomposition}}
\ntgroup{Projections and their weights}
\ntrow{$\mathcal X(G)$; $\mathcal M_G$}{AHK variables associated with $G$; square-integrable functions of them orthogonal to functions of every proper subgraph.}{\eqref{eq:graph-variables}, \eqref{eq:graph-space}}
\ntrow{$\mathcal M_G^{\rm sym}$}{The subspace of $\mathcal M_G$ invariant under the automorphisms of $G$.}{Sec.~\ref{app:projection-notation}}
\ntrow{$q_G,q_\Lambda$; $p_G$}{Ordinary graph projection, evaluated on a representative or a labeled copy; $p_G=(p)_{r(G)}(q)_{c(G)}q_G$.}{Lem.~\ref{lem:decomp}; Sec.~\ref{app:projection-notation}}
\ntrow{$f_S$; $f_S(\Lambda)$}{Hoeffding component indexed by $S$, centered in each variable; the labeled projection of $p_G(\Lambda)$.}{\eqref{eq:S-efron-stein}; Sec.~\ref{sec:S-ahk}}
\ntrow{$R,G(R),p_R,T_R$}{Isomorphism type of $(\Lambda,S)$, its graph class, the sum of labeled projections of that type, and its contribution to the statistic.}{Sec.~\ref{sec:S-ahk}}
\ntrow{$X_A$; $c_A$}{$X_A=\sum_{N(\Lambda)=A}q_\Lambda$; $c_A=(p)_{|A_1|}(q)_{|A_2|}/\{(m)_{|A_1|}(n)_{|A_2|}\}$ is the fraction of subarrays containing $A$.}{\eqref{eq:grouped-projections}}
\ntrow{$\beta_G=\beta_G(m,n)$}{Coefficient $1/\{(m)_{r(G)}(n)_{c(G)}\}$ in the graph-copy sum.}{\eqref{eq:copy-normalization}}
\ntrow{$M_N(G)$}{Number of copies: $(m)_{r(G)}(n)_{c(G)}/|\operatorname{Aut}(G)|$.}{Sec.~\ref{app:projection-notation}}
\ntrow{$w_G$}{Variance weight $w_G=\beta_G/|\operatorname{Aut}(G)|$.}{\eqref{eq:graph-variance-weight}}
\ntrow{$T_G$; $\varsigma_G^2,\varsigma_R^2$}{Graph contribution $T_G=\beta_G\sum_\Lambda p_G(\Lambda)$; variances of $T_G$ and $T_R$.}{\eqref{eq:copy-normalization}; Sec.~\ref{sec:S-ahk}}
\ntrow{$s_{N,G}^2$; $\sigma_{\rm disc}^2,D_N$}{$s_{N,G}^2=\varsigma_G^2/\sigma_N^2$; $\sigma_{\rm disc}^2=\sum_{G\,\disc}\varsigma_G^2$ and $D_N=\sigma_{\rm disc}^2/\sigma_N^2$.}{\eqref{eq:variance-proportions}; Sec.~\ref{app:projection-notation}}
\ntgroup{Variance estimators and their decomposition}
\ntrow{$V_{r,c},V_\ell$}{Variance summed over graph classes with $r$ rows and $c$ columns; $\ell=(r,c)$.}{Sec.~\ref{sec:jackknife-bias}; Sec.~\ref{sec:S-variance}}
\ntrow{$\widehat V_{r,c},\widehat V_\ell$}{Unbiased estimator of $V_{r,c}=V_\ell$.}{\eqref{eq:variance-component-estimator}}
\ntrow{$\omega(r,c)$}{Sum of $\E[q_\Lambda^2]$ over graphs with vertex sets $[r]$ and $[c]$.}{Lem.~\ref{lem:overlap}}
\ntrow{$\Vu$; $\Vu^+$}{$\Vu=U_N^2-\Thetah(0,0)$; $\Vu^+=\max(\Vu,0)$.}{Def.~\ref{def:ov}, Prop.~\ref{prop:pathwise}}
\ntrow{$\widehat Z_N$}{Studentized statistic.}{Sec.~\ref{sec:variance-estimators}, Thm.~\ref{thm:contraction}}
\ntrow{$\Vrcj$; $U^{(-i)},U^{(-j)}$; $\kappa_G,\kappa_{r,c}$}{Row-column jackknife; statistic without row $i$ or without column $j$; jackknife weights on variance components.}{\eqref{eq:rcj}, Prop.~\ref{prop:count}}
\ntrow{$J$; $E_N$}{$J=\sum_{A\sim A'}c_Ac_{A'}X_AX_{A'}$ over ordered pairs; $E_N=\Vu-J$.}{\eqref{eq:variance-decomposition}}
\ntrow{$J_\ell$; $E_{N,\ell}$}{$J_\ell$ restricts the sum for $J$ to $\ell(A\cap A')=\ell$; $E_{N,\ell}=\widehat V_\ell-J_\ell$.}{\eqref{eq:component-J}}
\ntrow{$b_\ell(A,A')$; $N_{a,b}(A,A')$}{Coefficient of $X_AX_{A'}$ in $\widehat V_\ell$; number of ordered subarray pairs containing $A,A'$ and sharing $a$ rows and $b$ columns.}{\eqref{eq:S-level-counting}}
\ntrow{$\Vu(f,k)$}{Symmetric bilinear form $\{\Vu(f+k)-\Vu(f)-\Vu(k)\}/2$ for kernels $f,k$.}{Sec.~\ref{sec:S-smooth}}
\ntgroup{Configurations and contractions}
\ntrow{$\pi_k=f_{S_k}(\Lambda_k)$; configuration}{A labeled projection; a configuration consists of four such projections.}{Sec.~\ref{sec:S-ahk}}
\ntrow{Pattern; labeling}{A pattern records shared vertices and selected-variable types; a labeling assigns actual row and column indices to its vertices.}{Sec.~\ref{sec:S-ahk}}
\ntrow{$v,v_r,v_c$; $r_{\rm tot},c_{\rm tot}$}{Numbers of distinct vertices, rows and columns in a configuration; totals of row and column occurrences, counted with multiplicity.}{Sec.~\ref{sec:S-ahk}; Lem.~\ref{lem:auto-pattern}}
\ntrow{$d_{\rm exc}$}{$d_{\rm exc}=\sum_kd(\Lambda_k)-2v$.}{Sec.~\ref{sec:S-ahk}}
\ntrow{$g\otimes_Bg$}{Contraction over variables in $B$: $(g\otimes_Bg)(x,x')=\E[g(x,\mathbf X_B)g(x',\mathbf X_B)]$.}{Sec.~\ref{sec:res-main}; \eqref{eq:S-contr}}
\ntrow{$\Upsilon_N$}{Sum of $w_G^2\|p_G\otimes_Bp_G\|_2^2/\sigma_N^4$ over disconnected $G$ and nonempty proper component sets $B$.}{\eqref{eq:gamma}}
\ntrow{$x_A,\mathbf X_A$ (AHK coordinates)}{AHK coordinate values and corresponding random vectors for the vertices in $A$ and cells with both endpoints in $A$.}{Sec.~\ref{sec:diagnostic-main}}
\ntrow{$T_g,T_g^*$; $\|g\|_{S_4},\|T_g\|_{S_4}$}{Integral operator with kernel $g$ and its adjoint; $S_4$ norms of the kernel and operator.}{\eqref{eq:S-contraction-identity}}
\ntgroup{Gaussian space and joint limits}
\ntrow{$\Gamma_0$}{Finite family of graph classes considered jointly in the Wick limits.}{Lem.~\ref{lem:wick}}
\ntrow{$\mathfrak h_H$; $\mathfrak h$; $X,X_H$}{$\mathfrak h_H=\mathcal M_H^{\rm sym}$ for connected $H$; $\mathfrak h=\bigoplus_H\mathfrak h_H$; a common isonormal Gaussian process and its restrictions.}{\eqref{eq:wick-gaussian-space}}
\ntrow{$\widehat\otimes,\mathfrak h^{\otimes b},\mathfrak h^{\odot b}$; $\operatorname{Sym}_b$}{Completed Hilbert tensor product and tensor power, its symmetric subspace; averaging over permutations of the $b$ arguments.}{Sec.~\ref{app:wick}}
\ntrow{$I_b$}{Multiple Wiener integral of order $b$.}{\eqref{eq:multiple-integral-isometry}}
\ntrow{$m_H,\mu_G$; $\mathcal W_G(f)$}{Component multiplicities and $\mu_G=\prod_Hm_H!$; $\mathcal W_G(f)=I_{b(G)}(\operatorname{Sym}_{b(G)}f)$.}{Sec.~\ref{app:wick-construction}; \eqref{eq:wick-isometry}}
\ntrow{$T_{N,G}(f),W_{N,G},W_G$}{Graph contribution with $p_G$ replaced by $f$, its standardized version, and its unit-variance limit.}{Lem.~\ref{lem:wick}}
\ntrow{$\mathcal J\subseteq[b]$}{Set of positions of connected components.}{Sec.~\ref{app:wick-products}}
\ntrow{$Z_{N,j}$; $\Xi_N(\mathcal J)$}{Normalized sum over individual copies of a connected graph; normalized product sum over ordered vertex-disjoint component copies indexed by $\mathcal J$.}{Sec.~\ref{app:wick-products}}
\ntrow{$\chi_{rs},\widehat\chi_{N,rs}$}{Covariance of the connected-graph sums and its average over copies.}{Sec.~\ref{app:wick-products}}
\ntrow{$\mathfrak M(\mathcal J)$; $\mathfrak m,V(\mathfrak m)$}{Set of matchings of $\mathcal J$; one matching and its covered positions.}{Sec.~\ref{app:wick-products}}
\ntrow{$S_N(f)$}{Sum of $f$ over ordered vertex-disjoint copies of the connected components.}{Sec.~\ref{app:wick-construction}}
\ntrow{$a_{N,G}$}{$a_{N,G}^2=\prod_{j=1}^bM_N(H_j)$ for the connected components $H_1,\ldots,H_b$ of $G$.}{Sec.~\ref{app:wick-construction}}
\ntrow{$\tau_G,L_b$}{Coefficients and graph contributions grouped by chaos order: $L_b=\sum_{d(G)=d,\,b(G)=b}\tau_GW_G$.}{\eqref{eq:fixed-law-coefficients}}
\ntrow{$s_G^*,f_G^*$}{Along the chosen subsequence, $s_{k,G}\to s_G^*$ and, when $s_G^*>0$, $p_{P_k,G}/\|p_{P_k,G}\|_2\xrightarrow{L^2}f_G^*$.}{Lem.~\ref{lem:S-replacement}}
\ntrow{$F$; $\mathcal V(F)$; $D(F)$}{$F=\sum_bI_b(f_b)$ is the statistic's limit; $\mathcal V(F)$ is the accompanying random variance-estimator limit; $D(F)=\sum_{b\ge2}b!\|f_b\|^2$.}{\eqref{eq:joint-limit-variances}; Thm.~\ref{thm:S-joint}}
\ntrow{$f\widetilde\otimes_tg$}{Contraction of $t$ tensor coordinates followed by symmetrization.}{\eqref{eq:joint-limit-variances}}
\ntrow{$X_A(a),\sigma(a),J(a,b)$}{Grouped projections for a family $a$, the standard deviation of their weighted sum, and the bilinear extension of $J$.}{Sec.~\ref{sec:S-beyond}}
\ntrow{$o_{L^k},O_{L^k}$, $k=1,2$}{Asymptotic orders in the $L^k$ norm.}{\hyperref[def:S-Lk-orders]{Sec.~\ref*{app:wick-products}}}
\ntgroup{Block tests and diagnostics}
\ntrow{$K_N,m_b,n_b,U_{N,j}$}{Number of diagonal blocks, their dimensions, and the statistic on block $j$.}{Sec.~\ref{sec:S-block-hybrid}}
\ntrow{$\overline U_{\rm block},S_{K_N}^2,\tau_N^2$}{Block average, sample variance, and variance $\tau_N^2=\var(U_{N,1})$ of one block statistic.}{\eqref{eq:block-statistics}; Sec.~\ref{app:protected}}
\ntrow{$A_N,B_N$ (procedures)}{Nonnegative variance statistic and $B_N=\max\{A_N,S_{K_N}^2/K_N\}$.}{Prop.~\ref{prop:block-safe}; Sec.~\ref{sec:S-chaos-procedure}}
\ntrow{$C_N$; $R_{G,N}$, $R_{S,N}$, $R_{H,N}$}{Event selecting the Gaussian test; Gaussian, conservative and combined rejection events.}{\eqref{eq:hybrid-definition}}
\ntrow{$\tau=(r_A,c_A\mid r_B,c_B)$}{Diagnostic index giving the row and column counts of two disjoint vertex sets.}{Sec.~\ref{sec:diagnostic-main}}
\ntrow{$W_T$; $\mathcal D(A,B)$; $\Pi_\tau$}{Subarray with the vertices in $T\subseteq A\cup B$ replaced; alternating kernel difference; product of four prescribed differences.}{\eqref{eq:diagnostic-differences}, \eqref{eq:diagnostic-product}}
\ntrow{$g_\tau$; $Z_\tau,\widehat Z_\tau$}{Sum of ordinary projections on exactly $A\cup B$ with no edge between $A$ and $B$; $Z_\tau=\|g_\tau\|_{S_4}^4$ and its estimator.}{Sec.~\ref{sec:S-ident}}
\ntrow{$c_\tau,\widehat\Upsilon_N$}{Weight and diagnostic statistic.}{\eqref{eq:diagnostic-statistic}}
\ntrow{$x_\tau,\widehat x_\tau,X_\tau$}{Normalized quantities $x_\tau=c_\tau Z_\tau/\sigma_N^4$, $\widehat x_\tau=c_\tau\widehat Z_\tau/\sigma_N^4$, and $X_\tau=\max_{\tau'\succ\tau}x_{\tau'}$.}{\hyperref[def:S-diagnostic-quantities]{Sec.~\ref*{sec:S-ident}}}
\ntrow{$\Omega_\tau,\omega,\Pi_\omega$}{Set of diagnostic placements, one placement, and its four-factor product.}{\eqref{eq:S-diagnostic-covariance}}
\ntrow{$\mathcal C_{\tau,1}$; $\mathcal C_{\tau,2}$}{Pairs sharing a vertex that is original or a replacement in at least one placement; pairs sharing at least one vertex, with every shared vertex additional in both placements.}{\eqref{eq:S-diagnostic-covariance}}
\ntrow{$\Gamma_{\rm disc}$; $\mathcal L_{\rm disc}$}{Selected disconnected graph classes under (L); $\mathcal L_{\rm disc}=\{\ell(G):G\in\Gamma_{\rm disc}\}$.}{Sec.~\ref{sec:alternative-procedure}; Sec.~\ref{sec:S-chaos-procedure}}
\ntrow{$F_N$}{$F_N=\sum_{\ell\in\mathcal L_{\rm disc}}\widehat V_\ell^+$.}{\eqref{eq:chaos-variance-adjustment}}
\ntrow{$\rho_N$; $v_N$}{$\rho_N$ is the disconnected share outside $\Gamma_{\rm disc}$; $v_N=\sum_{\ell\in\mathcal L_{\rm disc}}V_\ell/\sigma_N^2$.}{Lem.~\ref{thm:S-test}}
\ntrow{$\varepsilon_N$; $C_N^{\rm ctr},C_N^{\rm fl}$}{Threshold sequence; selection by $\widehat\Upsilon_N\le\varepsilon_N$ or by $\Vu>0$ and $F_N\le\varepsilon_N\Vu$.}{Sec.~\ref{sec:S-chaos-procedure}; \eqref{eq:S-floor-gate}}
\ntrow{$a_N$ (diagnostic rate)}{Rate sequence in the subclasses $\Upsilon_N\le a_N\varepsilon_N$.}{Thm.~\ref{thm:S-check}(ii)}
\ntrow{$\alpha,\theta_0$}{Test significance level and hypothesized value of $\theta_N$.}{Prop.~\ref{prop:block-safe}}
\ntgroup{Assumptions}
\ntrow{(M4), $C_{\rm M4}$; (K)}{Uniform fourth-moment bound relative to squared variance; compactness of normalized graph projections.}{Sec.~\ref{sec:res-main}}
\ntrow{(M8), $C_{\rm M8}$; (S2), $C_{\rm S2}$; $\mathcal A(F)$}{Uniform eighth-moment bound; comparison of kernel norms before and after fixing additional variables; allowed splits of the additional vertices.}{Sec.~\ref{sec:diagnostic-main}}
\ntrow{(L)}{Conditions on the variance outside $\Gamma_{\rm disc}$ and the row and column counts of unions of its graph components.}{Sec.~\ref{sec:alternative-procedure}}
\ntgroup{Product-form applications}
\ntrow{$w,\wbar,\lambda$}{Mean graphon $w=\lambda\wbar$, with $\iint\wbar=1$ and intensity $\lambda>0$.}{Sec.~\ref{sec:dw}}
\ntrow{$f,g$; $F_2,G_2$}{Row and column margins of $\wbar$; $F_2=\int f^2$ and $G_2=\int g^2$.}{Sec.~\ref{sec:dw}}
\ntrow{$\Delta,\gamma$; $\gamma_P$}{$\Delta=\wbar-f\otimes g$ and $\gamma=\|\Delta\|_2^2$; $\gamma_P=\gamma(\wbar_P)$.}{Sec.~\ref{sec:dw}}
\ntrow{$h_\gamma,U$; $\widehat\gamma_U$}{Kernel and its symmetrized $U$-statistic; ratio $\widehat\gamma_U=U/\bar Y^4$.}{\eqref{eq:gamma-kernel}, Sec.~\ref{sec:dw-uniform}}
\ntrow{$\ell_{N,P},\sigma_{N,P}^2$}{Linear approximation to $\widehat\gamma_U-\gamma_P$ and its variance.}{\eqref{eq:gamma-linearization}}
\ntrow{$\widetilde h_{\ell,P},p_G^{\widetilde\ell}$}{$\widetilde h_{\ell,P}=h_\gamma-4\gamma_P\lambda^3Y_{11}$ and its projections.}{Sec.~\ref{sec:dw-uniform}}
\ntrow{$A_N,B_N$; $a_N,b_N$; $\tau_A,\tau_B$ in Sec.~\ref{sec:S-smooth}}{Two $U$-statistics, their means and their standard deviations in the general linearization lemma.}{Sec.~\ref{sec:S-smooth}}
\ntrow{$\psi,c_1,c_2,\ell_N,\sigma_N^2$ in Sec.~\ref{sec:S-smooth}}{Scalar transformation, coefficients $c_1=\psi(b_N)$ and $c_2=a_N\psi'(b_N)$, the resulting centered linear term, and its variance $\sigma_N^2=\var(\ell_N)$.}{Sec.~\ref{sec:S-smooth}}
\ntrow{$\mu,\nu$; $\mathcal E(\mu,\nu)$}{Row and column laws; normalized kernel expectation $\mathcal E(\mu,\nu)=\lambda^{-4}\E_{\mu,\nu}[h_\gamma]$.}{Sec.~\ref{sec:S-ratio-projections}}
\ntrow{$s(x),\wbar_\varepsilon,\varepsilon$}{Two-group sign function, graphon $\wbar_\varepsilon(x,y)=1+\varepsilon s(x)s(y)$, and interaction parameter.}{Sec.~\ref{sec:planted}}
\ntrow{$t(C_{2k},w),\Psi(w)$}{Cycle density and the distance $\Psi(w)=t(C_4,w)^2-t(C_8,w)$.}{Sec.~\ref{sec:d3}}
\ntrow{$h_\Psi,\widehat\Psi$}{Kernel formed from products on $C_4\sqcup C_4$ and $C_8$, and its unbiased $U$-statistic.}{Sec.~\ref{sec:d3}; Fig.~\ref{fig:psi-patterns}}
\ntgroup{Simulation model}
\ntrow{$a,t_N\;(t),r$}{Row-effect strength, interaction strength and rank.}{Sec.~\ref{sec:fallback-experiment}; \eqref{eq:S-model-array}}
\ntrow{$A_i,U_{ik},V_{jk},e_{ij}$}{Independent standard normal coordinates of the simulation model.}{Sec.~\ref{sec:fallback-experiment}}
\ntrow{$\theta$; $\sigma_0^2$; $\vartheta$}{$\theta=2a^2$; variance of $U_N$ under $a=0$; standardized parameter $\vartheta=\theta/\sigma_0$.}{Sec.~\ref{sec:fallback-experiment}; \eqref{eq:S-model-variance}}
\end{longtable}
\endgroup

\section{Projections, counting and variance identities}\label{sec:S-part-projections}

\subsection{Preliminaries on projections}\label{app:projection-notation}\label{sec:S-normalization}\label{sec:mechanisms}\label{sec:S-ahk}\label{sec:S-subarrays}

\paragraph{Subarrays.}
Recall from \eqref{eq:ustat} that $U_N=K^{-1}\sum_WH_W$ with $H_W=h(Y_W)$, where $K=\binom mp\binom nq$ and $W=(W_1,W_2)$ runs over the subarrays with $|W_1|=p$ rows and $|W_2|=q$ columns. For two subarrays $W$ and $W'=(W'_1,W'_2)$, put $o(W,W')=(|W_1\cap W'_1|,|W_2\cap W'_2|)$, and let $K(a,b)$ be the number of ordered pairs with $o(W,W')=(a,b)$. The average product of kernel values in Section~\ref{sec:variance-estimators} is
\begin{equation}\label{eq:kernel-product-averages}
\Thetah(a,b)=K(a,b)^{-1}\sum_{o(W,W')=(a,b)}H_WH_{W'},
\end{equation}
an unbiased estimator of $\Theta(a,b)=\E[H_WH_{W'}]$ for two subarrays sharing $a$ rows and $b$ columns.

\paragraph{Graph projections and normalization.}
Recall the AHK variables $\mathcal X(G)$ and the space $\mathcal M_G$ from \eqref{eq:graph-variables} and~\eqref{eq:graph-space}. Write $\mathcal M_G^{\rm sym}$ for the subspace of $\mathcal M_G$ consisting of functions invariant under the automorphisms of $G$.

The graph decomposition gives
\[
L^2(\mathcal X(K_{p,q}))=\bigoplus_{G\subseteq K_{p,q}}\mathcal M_G.
\]
Applied to the union $\Lambda\cup \Lambda'$, it shows that projections on distinct copies are orthogonal, that means, if $\Lambda=\Lambda'$, then $\E[p(\Lambda)p(\Lambda')]=0$.

For a labeled subgraph $\Lambda$ of a subarray, with $r$ rows and $c$ columns, let $q_\Lambda$ be the ordinary graph projection, that means, the orthogonal projection of the subarray kernel onto $\mathcal M_\Lambda$. By the symmetry of $h$, it is a function of $\mathcal X(\Lambda)$, identical for every subarray containing $\Lambda$. If $\Lambda$ is a copy of $G$, the normalization in Lemma~\ref{lem:decomp} gives $q_\Lambda=p_G(\Lambda)/\{(p)_r(q)_c\}$. Equivalently, $p_G$ is the sum of the projections of $h\circ\phi_N$ over the $(p)_r(q)_c$ placements of $G$ in one subarray. Expanding each centered kernel as $H_W-\theta_N=\sum_{\Lambda\subseteq W}q_\Lambda$ and averaging over subarrays, each copy appears in a fraction $(p)_r(q)_c/\{(m)_r(n)_c\}$ of them. 

Lemma~\ref{lem:decomp} therefore reads
\begin{equation}\label{eq:copy-normalization}
\begin{gathered}U_N-\theta_N=\sum_{G\in\Gamma}T_G,\qquad T_G=\beta_G(m,n)\sum_{\Lambda\text{ copy of }G}p_G(\Lambda),\\ \beta_G(m,n)=\frac{(m-r)!\,(n-c)!}{m!\,n!}=\frac1{(m)_r(n)_c},\end{gathered}
\end{equation}
for every $m\ge p$, $n\ge q$, with the same functions $p_G$. Since projections on distinct copies are orthogonal, the $T_G$ are pairwise orthogonal. Write $M_N(G)=(m)_{r(G)}(n)_{c(G)}/|\mathrm{Aut}(G)|$ for the number of copies of $G$ in $K_{m,n}$. Recall that $w_G=1/\{|\mathrm{Aut}(G)|(m)_{r(G)}(n)_{c(G)}\}$ is the variance weight in \eqref{eq:graph-variance-weight}. Summing the variances of the copies gives
\[
\var(T_G)=\frac{\E[p_G^2]}{|\mathrm{Aut}(G)|\,(m)_{r(G)}(n)_{c(G)}}=w_G\E[p_G^2].
\]

We write $\varsigma_G^2=\var(T_G)$, so that the variance share in~\eqref{eq:variance-proportions} satisfies $s_{N,G}=\varsigma_G/\sigma_N$. We define the \emph{disconnected variance} by
\[
\sigma_{\mathrm{disc}}^2:=\sum_{G\ \disc}\varsigma_G^2,
\qquad D_N=\sigma_{\mathrm{disc}}^2/\sigma_N^2.
\]

\paragraph{Grouped projections.}
Write $N(\Lambda):=(V_1(\Lambda),V_2(\Lambda))$ for the row and column vertex sets of a labeled subgraph $\Lambda$. For a pair of vertex sets $A=(A_1,A_2)$, put
\begin{equation}\label{eq:grouped-projections}
X_A:=\sum_{\Lambda:\,N(\Lambda)=A}q_\Lambda,
\qquad
c_A:=\frac{(p)_{|A_1|}(q)_{|A_2|}}{(m)_{|A_1|}(n)_{|A_2|}}.
\end{equation}
The sum defining $X_A$ is over the graphs with exactly these vertices, and $c_A$ is the fraction of subarrays containing them. In all sums below, $A$ is nonempty, with $|A_1|\le p$ and $|A_2|\le q$. The centered kernel on a subarray $W$ is
\[
H_W-\theta_N
 =\sum_{\substack{A=(A_1,A_2),\,A\ne\varnothing\\ A_1\subseteq W_1,\,A_2\subseteq W_2}}X_A.
\]
Averaging over subarrays gives
\begin{equation}\label{eq:auto-U-expansion}
U_N-\theta_N=\sum_Ac_AX_A.
\end{equation}
If $A=N(\Lambda)$ for a copy $\Lambda$ of $G$, then $c_Aq_\Lambda=\beta_Gp_G(\Lambda)$ by the definitions of $q_\Lambda$ and $\beta_G$. Since $T_G$ sums the $M_N(G)$ orthogonal terms $\beta_Gp_G(\Lambda)$, which have the same second moment,
\begin{equation}\label{eq:copy-scale}
c_A\|q_\Lambda\|_2=\beta_G\|p_G\|_2=\varsigma_G/M_N(G)^{1/2}.
\end{equation}

\paragraph{Counting products of graph projections.}
The following lemma describes which products of graph projections can have a nonzero expectation.

\begin{lemmaL}[\citealp{leminh2026characterization}, Lem.~9]\label{lem:config}
Let $\Lambda_1,\dots,\Lambda_K$ be copies with $p_{G_k}\in \mathcal M_{G_k}$. If $\E[\prod_kp_{G_k}(\Lambda_k)]\ne0$, then every vertex and edge of each $\Lambda_k$ belongs to another $\Lambda_\ell$. If, in addition, all $\Lambda_k$ are connected and nonempty, then either they coincide in pairs, or some vertex belongs to at least three of them.
\end{lemmaL}

\paragraph{Decomposition within a graph projection (Efron--Stein).}
To handle the uniform guarantees of this paper, we need a finer Hoeffding decomposition, similar to \citet{efron1981jackknife}. We further decompose each space $\mathcal M_G$ over the independent AHK variables. For a copy $\Lambda$ of $G$ and $f\in L^2(\mathcal X(\Lambda))$, write 
\begin{equation}\label{eq:S-efron-stein}
f=\sum_{S\subseteq \mathcal X(\Lambda)}f_S, \qquad f_S=\sum_{S'\subseteq S}(-1)^{|S|-|S'|}\E[f\mid S']. 
\end{equation}
Each $f_S$ depends on the variables in $S$ and satisfies $\E[f_S\mid\mathcal F_K]=0$ unless $S\subseteq K$, where $\mathcal F_K$ is the sigma-field generated by the variables in $K$. The functions $f_S$ are orthogonal.

If $f\in \mathcal M_G$, then $f_S=0$ unless $S$ contains every edge variable of $\Lambda$ and the variable of every isolated vertex of $\Lambda$. Indeed, removing a missing edge or an isolated vertex gives a proper subgraph $F$ with $S\subseteq \mathcal X(F)$. Then, $\|f_S\|_2^2=\E[ff_S]=0$, since $f\perp L^2(\mathcal X(F))$. We call sets $S$ satisfying this requirement \emph{admissible}: variables of non-isolated vertices may be included or omitted. An admissible $S$ is contained in no $\mathcal X(F)$ with $F\subsetneq \Lambda$, so $f_S\in \mathcal M_G$. 

For a copy $\Lambda$ and an admissible $S\subseteq \mathcal X(\Lambda)$, call $f_S(\Lambda):=(p_G(\Lambda))_S$ the \emph{labeled Hoeffding projection} indexed by $(\Lambda,S)$, or \emph{labeled projection} for short. Its variables determine the copy: the edge variables specify edges and their endpoints, and the remaining vertex variables specify isolated vertices. Thus, distinct labeled projections have distinct variable sets and are orthogonal. We use the same terminology for the corresponding Hoeffding projections of an ordinary graph projection.

Write $R$ for the isomorphism type of $(\Lambda,S)$, with $G(R)=G$. We call $R$ disconnected when $G(R)$ is disconnected. Define $p_R(\Lambda):=\sum_Sf_S(\Lambda)$, where the sum is over admissible sets $S\subseteq \mathcal X(\Lambda)$ of type $R$. Then, $p_G=\sum_Rp_R$ is an orthogonal decomposition, and 
\[
T_G=\sum_RT_R, \qquad T_R=\beta_G\sum_\Lambda p_R(\Lambda).
\]
The $T_R$ are pairwise orthogonal, with $\var(T_R)=w_G\E[p_R^2]=:\varsigma_R^2\le\varsigma_G^2$.

\paragraph{Patterns and configurations.} We call a \emph{configuration} a set of four labeled projections $\pi_k=f_{S_k}(\Lambda_k)$, $k=1,\dots,4$, on copies $\Lambda_k\subseteq K_{m,n}$ of graphs $G_k$. A configuration is \emph{paired} if the four labeled projections form two identical pairs, with the same copy and variable set in each pair. The \emph{pattern} of copies $\Lambda_1,\dots,\Lambda_4$ of $G_1,\dots,G_4$ records which of their vertices coincide, but not which rows and columns they occupy. A \emph{labeling} of a pattern is a choice of distinct rows and columns realizing it. The pattern of a configuration is the pattern of its copies, together with the type of each $S_k$ in its copy. Define $d_{\mathrm{exc}}=\sum_kd(\Lambda_k)-2v$, where $v$ counts the distinct vertices in the union of the copies.

\subsection{Variance identities}\label{sec:S-decomp}\label{sec:S-overlap}

Throughout, we assume $m\ge2p$ and $n\ge2q$. Recall the variance sums $V_{r,c}$ of Section~\ref{sec:jackknife}.

\begin{lemma}\label{lem:overlap}
Let $\omega(r,c)=\sum_\Lambda\E[q_\Lambda^2]$, the sum running over the labeled $\Lambda$ with row set $[r]$ and column set $[c]$, with $\omega(0,0)=\theta_N^2$. Then
\[
\Theta(a,b)=\sum_{r\le a,\,c\le b}\binom ar\binom bc\,\omega(r,c),\qquad
V_{r,c}=\frac{\{(p)_r(q)_c\}^2}{r!\,c!\,(m)_r(n)_c}\,\omega(r,c)\quad(r+c\ge1).
\]
Hence, $\omega(r,c)=\sum_{a\le r,b\le c}(-1)^{r-a+c-b}\binom ra\binom cb\Theta(a,b)$.
\end{lemma}

\begin{proof}[Proof of Lemma~\ref{lem:overlap}]
First, we express $\Theta(a,b)$ in terms of the projections supported on the shared rows and columns. Fix two subarrays \(W,W'\) sharing exactly \(a\) rows and \(b\) columns, and let \(\mathcal H_{a,b}\) be the sigma-field generated by all AHK variables belonging to their shared \(a\times b\) block. Conditional on \(\mathcal H_{a,b}\), the AHK variables used only by \(W\) and only by \(W'\) are independent. Moreover, $\E[H_W\mid\mathcal H_{a,b}]=\E[H_{W'}\mid\mathcal H_{a,b}]$, since $h$ is symmetric and both subarrays contain the shared block. Hence
\[
\Theta(a,b)=\E[H_WH_{W'}]
 =\E\left[\E[H_W\mid\mathcal H_{a,b}]^2\right].
\]
Write the orthogonal decomposition of the subarray kernel as
\[
H_W=\theta_N+\sum_{\Lambda\subseteq W,\,\Lambda\neq\varnothing}q_\Lambda.
\]
If a projection involves a variable outside the shared block, its conditional expectation given $\mathcal H_{a,b}$ is zero. In contrast, a projection supported entirely on the shared block is measurable with respect to $\mathcal H_{a,b}$, so its conditional expectation is equal to the projection itself. Because the spaces \(\mathcal M_\Lambda\) are orthogonal, we have
\[
\Theta(a,b)=\theta_N^2+
\sum_{\Lambda\subseteq K_{a,b},\,\Lambda\neq\varnothing}\E[q_\Lambda^2].
\]
If the vertex set of \(\Lambda\) contains \(r\) shared rows and \(c\) shared columns, there are \(\binom ar\binom bc\) ways to choose those vertices. By exchangeability and the symmetry of $h$, the sum of $\E[q_\Lambda^2]$ over all graphs on any one of these vertex sets is equal to $\omega(r,c)$. Thus,
\[
\Theta(a,b)=\sum_{r\le a,\,c\le b}\binom ar\binom bc\,\omega(r,c).
\]
Inverting these relations for every $(r,c)$ gives the formula for $\omega(r,c)$.

To relate \(\omega(r,c)\) to $V_{r,c}$, consider a projection \(q_\Lambda\) on a vertex set of size \((r,c)\). It occurs in every subarray containing those vertices. The fraction of subarrays containing them is
\[
c_{r,c}=\frac{\binom{m-r}{p-r}\binom{n-c}{q-c}}{\binom mp\binom nq}=\frac{(p)_r(q)_c}{(m)_r(n)_c}.
\]
This fraction is $c_A$ in \eqref{eq:grouped-projections} for the vertex set $A$ of $\Lambda$, so \eqref{eq:auto-U-expansion} reads
\[
U_N-\theta_N=\sum_{\Lambda\neq\varnothing}c_{r(\Lambda),c(\Lambda)}q_\Lambda.
\]
Since projections on distinct labeled graphs are orthogonal, the variance contributed by vertex sets of size \((r,c)\) is
\[
V_{r,c}
 =c_{r,c}^2\binom mr\binom nc\,\omega(r,c).
\]
Using \(\binom mr=(m)_r/r!\) and \(\binom nc=(n)_c/c!\) gives
\[
V_{r,c}
 =\frac{\{(p)_r(q)_c\}^2}{r!\,c!\,(m)_r(n)_c}\,\omega(r,c),
\]
as claimed.
\end{proof}

\subsection{The jackknife bias: proof of Proposition~\ref{prop:count}}\label{app:estimators}

\begin{proof}[Proof of Proposition~\ref{prop:count}]
We compute the expectation of each graph contribution to the row and column jackknives. Deleting a row removes every copy containing that row. Since the projections of distinct copies are orthogonal, we then obtain the second moments, while the change from \(\beta_G(m,n)\) to \(\beta_G(m-1,n)\) gives the finite-sample factor $(m-1)/(m-r)$.

\textbf{(i)} First, we express each deletion difference as a sum of graph contributions.
Fix a graph class \(G\), abbreviate \(r=r(G)\), and recall from \eqref{eq:copy-normalization} that \(\beta_G(m,n)=1/\{(m)_r(n)_{c(G)}\}\). Write
\[
\Sigma_G:=\sum_{\Lambda\text{ copy of }G}p_G(\Lambda),
\qquad
A_{i,G}:=\sum_{\Lambda:\,i\in V_1(\Lambda)}p_G(\Lambda).
\]
After deleting row \(i\), only copies not containing \(i\) remain. Thus, Lemma~\ref{lem:decomp} on the \((m-1)\times n\) array gives
\[
U^{(-i)}-\theta_N
 =\sum_G\beta_G(m-1,n)\{\Sigma_G-A_{i,G}\}.
\]
Every copy of \(G\) contains exactly \(r\) rows, so each summand \(p_G(\Lambda)\) occurs in exactly \(r\) of the sums \(A_{i,G}\). Therefore,
\[
\sum_{i=1}^m A_{i,G}=r\Sigma_G.
\]
Averaging these statistics over \(i\) gives
\[
\bar U_{\mathrm{row}}-\theta_N
 =\sum_G\beta_G(m-1,n)\left(1-\frac rm\right)\Sigma_G,
\]
and subtraction gives the exact identity
\begin{equation}\label{eq:J-row-difference}
U^{(-i)}-\bar U_{\mathrm{row}}
 =-\sum_G\beta_G(m-1,n)\left(A_{i,G}-\frac r m\Sigma_G\right).
\end{equation}

\textbf{(ii)} Next, we compute the expected squared deletion differences in \eqref{eq:J-row-difference}.
Squaring \eqref{eq:J-row-difference} and taking expectations removes cross terms between distinct graph classes, by Lemma~\ref{lem:decomp}. For a fixed class \(G\), expand
\[
\sum_i A_{i,G}^2
 =\sum_{\Lambda,\Lambda'}p_G(\Lambda)p_G(\Lambda')\,|V_1(\Lambda)\cap V_1(\Lambda')|.
\]
If \(\Lambda\neq \Lambda'\), then \(\E[p_G(\Lambda)p_G(\Lambda')]=0\). If  instead \(\Lambda=\Lambda'\), the intersection contains exactly \(r\) rows. Thus, summing over the \(M_N(G)\) copies of \(G\),
\[
\E[\sum_iA_{i,G}^2]=rM_N(G)\E[p_G^2].
\]
Similarly,
\[
\E[\Sigma_G^2]=M_N(G)\E[p_G^2].
\]
Since \(\sum_iA_{i,G}=r\Sigma_G\), centering by \(r\Sigma_G/m\) gives
\begin{align*}
\E[\sum_i\left(A_{i,G}-\frac r m\Sigma_G\right)^2]
&=\E[\sum_iA_{i,G}^2]-\frac{r^2}{m}\E[\Sigma_G^2]\\
&=r\left(1-\frac r m\right)M_N(G)\E[p_G^2].
\end{align*}

\textbf{(iii)} Now, we insert the second moments from (ii) into the row jackknife and compare its normalization with that of $T_G$.
The row jackknife term is
\[
\widehat V_{\mathrm{row}}:=\frac{m-1}{m}\sum_i(U^{(-i)}-\bar U_{\mathrm{row}})^2.
\]
Using \eqref{eq:J-row-difference} from (i) and the second moments in (ii),
\[
\E[\widehat V_{\mathrm{row}}]
 =\sum_G\frac{m-1}{m}\beta_G(m-1,n)^2
 r\left(1-\frac r m\right)M_N(G)\E[p_G^2].
\]
Now
\[
\frac{\beta_G(m-1,n)}{\beta_G(m,n)}
 =\frac{m}{m-r},
\qquad
\var(T_G)=\beta_G(m,n)^2M_N(G)\E[p_G^2].
\]
Substituting these identities gives
\[
\E[\widehat V_{\mathrm{row}}]
 =\sum_G r\frac{m-1}{m-r}\,\var(T_G).
\]
The same argument with rows and columns interchanged gives
\[
\E[\widehat V_{\mathrm{col}}]
 =\sum_G c(G)\frac{n-1}{n-c(G)}\,\var(T_G).
\]
Adding the row and column terms proves the formula for \(\E[\Vrcj]\). Since \(1\le r+c\le p+q\) for every nonempty class, the stated properties of \(\kappa_G\) follow immediately.
\end{proof}

\subsection{Removing the jackknife bias: proof of Proposition~\ref{prop:pathwise}}\label{sec:S-pathwise}

To prove Proposition~\ref{prop:pathwise}, we first show that a quadratic statistic whose coefficients depend only on the numbers of shared rows and columns is determined by its expectations under all dissociated RCE laws and kernels. Let $\mathcal Q$ be the set of statistics of the form
\[
B=\sum_{W,W'}\beta\bigl(o(W,W')\bigr)H_WH_{W'},
\]
where, for each pair $(a,b)$ with $0\le a\le p$ and $0\le b\le q$, $\beta(a,b)$ is a real coefficient, which may depend on $m,n,p,q$. Grouping the products by their numbers of shared rows and columns gives
\[
B=\sum_{a=0}^{p}\sum_{b=0}^{q}\beta(a,b)K(a,b)\Thetah(a,b).
\]

\begin{lemma}\label{lem:unique}
If $B\in\mathcal Q$ and $\E[B]=0$ for every dissociated RCE law and every kernel, then $B\equiv0$ identically in $(H_W)$.
\end{lemma}

\begin{proof}[Proof of Lemma~\ref{lem:unique}]
By unbiasedness of $\Thetah(a,b)$, the assumption gives
\begin{equation}\label{eq:expectation_identity}
       0=\E[B]=\sum_{a,b}\beta(a,b)K(a,b)\Theta(a,b)
\end{equation}
for every law and every kernel. We choose laws and kernels whose vectors of expected products $(\Theta(a,b))_{0\le a\le p,\,0\le b\le q}$ form a basis, so that this equality forces every fixed coefficient $\beta(a,b)K(a,b)$ to vanish.

For each $(r,c)\ne(0,0)$, take $Y_{ij}=(\xi_i,\eta_j,\zeta_{ij})$ and let $h$ be the row-column symmetrization of
\[
\prod_{i\in I}u(\xi_i)\prod_{j\in J}v(\eta_j),
\qquad I\subseteq[p],\quad J\subseteq[q],\quad |I|=r,\quad |J|=c,
\]
where $u,v$ are centered functions with positive, finite second moments, and, by convention, we take an empty product to equal one. Since $u$ and $v$ are centered, integrating out $\xi_i$ for some $i\in I$, or $\eta_j$ for some $j\in J$, makes the product equal to zero. Thus, its symmetrization has labeled projections only on copies of $K_{r,0}\sqcup K_{0,c}$. Lemma~\ref{lem:overlap} gives
\[
\Theta(a,b)=\omega(r,c)\binom ar\binom bc,
\qquad \omega(r,c)>0.
\]
For $(r,c)=(0,0)$, the kernel $h\equiv1$ gives the same formula with $\omega(0,0)=1$.

Now, we place the vectors
\[
e^{(r,c)}=\left(\binom ar\binom bc\right)_{0\le a\le p,\,0\le b\le q}
\]
as the columns of a square matrix, using the same lexicographic order for the row indices $(a,b)$ and the column indices $(r,c)$. Above the diagonal, either $a<r$, or $a=r$ and $b<c$, so the entry is zero. Moreover, every diagonal entry is $\binom aa\binom bb=1$. The matrix is therefore invertible, and its columns form a basis.
Finally, applying~\eqref{eq:expectation_identity} to the law and kernel constructed for $(r,c)$ gives
\[
0=\sum_{a,b}\beta(a,b)K(a,b)\Theta(a,b)
=\omega(r,c)\sum_{a,b}\beta(a,b)K(a,b)\binom ar\binom bc.
\]
Since $\omega(r,c)>0$, the vector of coefficients $\beta(a,b)K(a,b)$ is orthogonal to $e^{(r,c)}$. As $(r,c)$ is arbitrary, this vector is orthogonal to every vector in this basis.
Hence, $\beta(a,b)K(a,b)=0$ for every $(a,b)$. Since $K(a,b)>0$, we also have $\beta(a,b)=0$, and $B\equiv0$.
\end{proof}

Now, we use the uniqueness lemma (Lem.~\ref{lem:unique}) to prove the identities in Proposition~\ref{prop:pathwise}.

\begin{proof}[Proof of Proposition~\ref{prop:pathwise}]
\textbf{Proof of (i).} By Proposition~\ref{prop:count} and the unbiasedness of $\widehat V_{r,c}$, both sides of (i) have expectation $\sum_{r+c\ge1}\kappa_{r,c}V_{r,c}$. To apply Lemma~\ref{lem:unique}, we must also show that their coefficients depend only on the numbers of shared rows and columns. Let
\[
S:=\sum_WH_W,
\qquad
S_i:=\sum_{\substack{W=(W_1,W_2)\\ i\in W_1}}H_W.
\]
Since each subarray contains \(p\) rows, \(\sum_iS_i=pS\). The statistic based on subarrays avoiding row \(i\) is
\[
U^{(-i)}=\frac{S-S_i}{K(m-p)/m},
\]
because exactly a fraction \((m-p)/m\) of the \(p\)-row subsets avoid a fixed row. Therefore
\[
U^{(-i)}-U_N=\frac{pS-mS_i}{K(m-p)}.
\]
Averaging this identity over $i$ and using $\sum_iS_i=pS$ gives $\bar U_{\mathrm{row}}=U_N$. The row jackknife term is therefore
\begin{align*}
\frac{m-1}{m}\sum_i(U^{(-i)}-U_N)^2
&=\frac{m-1}{mK^2(m-p)^2}\sum_i(pS-mS_i)^2\\
&=\frac{m-1}{K^2(m-p)^2}\left\{m\sum_iS_i^2-p^2S^2\right\},
\end{align*}
where the second equality uses $\sum_iS_i=pS$. Expanding the two sums gives
\[
\sum_iS_i^2=\sum_{W,W'}|W_1\cap W'_1|H_WH_{W'},
\qquad S^2=\sum_{W,W'}H_WH_{W'}.
\]
The coefficients are the number of shared rows and the constant one, respectively, so both sums belong to $\mathcal Q$. The same argument applies to the column term, so \(\Vrcj\in\mathcal Q\).

Moreover, by definition~\eqref{eq:variance-component-estimator}, each $\widehat V_{r,c}$ is a fixed linear combination of the statistics $\Thetah(a,b)$, so \(\sum\kappa_{r,c}\widehat V_{r,c}\in\mathcal Q\). Their difference belongs to $\mathcal Q$ and has expectation zero for every law and every kernel, so Lemma~\ref{lem:unique} gives (i).

\textbf{Proof of (ii).}
Both \(U_N^2\) and \(\Thetah(0,0)\) are in \(\mathcal Q\). Also,
\[
\E[U_N^2-\Thetah(0,0)]
=\E[U_N^2]-\theta_N^2
=\sigma_N^2,
\]
because two disjoint subarrays are independent by dissociation, so \(\E[\Thetah(0,0)]=\theta_N^2\). On the other hand, Lemma~\ref{lem:overlap} gives
\[
\E[\sum_{r+c\ge1}\widehat V_{r,c}]
=\sum_{r+c\ge1}V_{r,c}
=\sigma_N^2.
\]
Again, the difference belongs to \(\mathcal Q\) and has expectation zero for every law and kernel. Lemma~\ref{lem:unique} therefore gives
\[
\Vu=\sum_{r+c\ge1}\widehat V_{r,c}=U_N^2-\Thetah(0,0)
\]
for every array.

\textbf{Proof of (iii).}
Subtracting the identity in (ii) from that in (i) gives the correction formula~\eqref{eq:pathwise}. The terms with $r+c=1$ vanish because $\kappa_{1,0}=\kappa_{0,1}=1$.
\end{proof}

\subsection{Computation by random splits}\label{sec:S-split}

Following the partition representation of \citet{wang2014variance}, split the rows uniformly at random into $R_1,R_2$ of sizes $\lfloor m/2\rfloor$ and $\lceil m/2\rceil$, and split the columns into $C_1,C_2$ in the same way. Assume $\lfloor m/2\rfloor\ge p$ and $\lfloor n/2\rfloor\ge q$. Let $U_{ab}$ be the $U$-statistic on $R_a\times C_b$, and put
\begin{equation}\label{eq:split-estimate}
\pi=-\frac12\{(U_{11}-U_N)(U_{22}-U_N)+(U_{12}-U_N)(U_{21}-U_N)\}.
\end{equation}
Let $\widehat V_{\mathrm{split}}$ be the average of $\pi$ over $S$ independently drawn splits. Write $\E_s$ for expectation over the splits, conditional on the data.

\begin{lemma}\label{rem:split}
The split estimator satisfies $\E_s[\widehat V_{\mathrm{split}}]=\Vu$. If (M4) holds on a class $\Pcal$, then for every $P\in\Pcal$ and all sufficiently large $N$,
\[
\E\big[\big(\widehat V_{\mathrm{split}}-\Vu\big)^2\big]\le C\,4^{p+q}\,\sigma_N^4/S,
\]
where $C$ depends only on the kernel dimensions and the constant $C_{\mathrm{M4}}$ in (M4).
\end{lemma}

\begin{proof}[Proof of Lemma~\ref{rem:split}]
First, we condition on the observed array and average over the random splits. Thus, $U_N$ and $\Vu$ are fixed. Since the split is uniformly random, every subarray $W$ has the same probability of lying in $R_a\times C_b$. Every ordered pair $(W,W')$ with $o(W,W')=(0,0)$ also has the same probability that $W\subseteq R_1\times C_1$ and $W'\subseteq R_2\times C_2$, and similarly for $R_1\times C_2$ and $R_2\times C_1$. Hence, by \eqref{eq:kernel-product-averages},
\[
\E_s[U_{ab}]=U_N,
\qquad \E_s[U_{11}U_{22}]=\E_s[U_{12}U_{21}]=\Thetah(0,0).
\]
Expanding the centered products in~\eqref{eq:split-estimate} and applying Proposition~\ref{prop:pathwise}(ii) gives
\[
\E_s[\pi]=U_N^2-\Thetah(0,0)=\Vu.
\]
Since the $S$ splits are independent conditional on the observed array, their average has conditional mean $\Vu$ and conditional variance $S^{-1}\var_s(\pi)$. Taking expectation over the observed array gives
\begin{equation}\label{eq:split-mse-reduction}
\E[(\widehat V_{\mathrm{split}}-\Vu)^2]
=\frac1S\E[\var_s(\pi)]
\le\frac1S\E[\pi^2],
\end{equation}
where $\var_s$ denotes variance over the splits conditional on the observed array.

Next, we bound $\E[\pi^2]$ using fourth moments. Let $V_{r,c}^{(a,b)}$ be the corresponding variance contribution to $U_{ab}$. The law and the kernel are unchanged, so $\omega(r,c)$ is unchanged. Applying Lemma~\ref{lem:overlap} at dimensions $(|R_a|,|C_b|)$ gives
\[
V_{r,c}^{(a,b)}
=\frac{\{(p)_r(q)_c\}^2}{r!\,c!\,(|R_a|)_r(|C_b|)_c}\,\omega(r,c)
=\frac{(m)_r(n)_c}{(|R_a|)_r(|C_b|)_c}\,V_{r,c}.
\]
The factor satisfies
\[
\frac{(m)_r(n)_c}{(|R_a|)_r(|C_b|)_c}
=2^{r+c}\{1+O(L_N^{-1})\},
\]
uniformly over the four blocks, because $|R_a|=m/2+O(1)$ and $|C_b|=n/2+O(1)$. Summing the contributions and using $r+c\le p+q$ gives
\[
\var(U_{ab})\le2^{p+q}\sigma_N^2\{1+o(1)\}.
\]
For all sufficiently large $N$, we have $|R_a|\ge2p$ and $|C_b|\ge2q$. Under (M4), Lemma~\ref{lem:S-four} at dimensions $(|R_a|,|C_b|)$ and $(m,n)$ then gives $\E[(U_{ab}-\theta_N)^4]\le C\var(U_{ab})^2$ and $\E[(U_N-\theta_N)^4]\le C\sigma_N^4$. By the preceding display, $\var(U_{ab})^2\le C\,4^{p+q}\sigma_N^4$. Using these bounds and $(x+y)^4\le8(x^4+y^4)$, we obtain
\begin{align*}
\E[(U_{ab}-U_N)^4]
&\le8\E[(U_{ab}-\theta_N)^4]+8\E[(U_N-\theta_N)^4]\\
&\le C\,4^{p+q}\sigma_N^4,
\end{align*}
where $C$ depends only on $(p,q,C_{\mathrm{M4}})$. Finally, $(x+y)^2\le2(x^2+y^2)$ and the Cauchy--Schwarz inequality, applied to the two products defining $\pi$, give
\[
\E[\pi^2]\le\max_{a,b\in\{1,2\}}\E[(U_{ab}-U_N)^4]
\le C\,4^{p+q}\sigma_N^4.
\]
Combining this bound with~\eqref{eq:split-mse-reduction} proves the result.
\end{proof}

\section{Uniform inference}\label{sec:S-part-uniform}

In this section, we prove Theorem~\ref{thm:contraction} and Proposition~\ref{prop:gammaD} (Sec.~\ref{sec:S-thm1} and~\ref{sec:S-gammaD}). Sections~\ref{sec:S-counting},~\ref{sec:S-cumulants} and~\ref{sec:S-variance} establish supporting lemmas.

We write $\beta_G$ for the coefficient $\beta_G(m,n)$ in \eqref{eq:copy-normalization}. Then $w_G=\beta_G/|\mathrm{Aut}(G)|$ and $\beta_G\asymp w_G\asymp M_N(G)^{-1}$. In Sections~\ref{sec:S-counting}, \ref{sec:S-cumulants}, \ref{sec:S-variance} and~\ref{sec:S-thm1}, we work on a class $\Pcal$ satisfying (M4) with constant $C_{\mathrm{M4}}$, with $m\ge2p$, $n\ge2q$ and $\sigma_{N,P}>0$ for every $P\in\Pcal$. There, $C$ denotes a constant depending only on $(p,q,C_{\mathrm{M4}})$, whose value may change from line to line.

\subsection{Counting and fourth-moment bounds}\label{sec:S-counting}\label{sec:concentration}

This section collects bounds used in several proofs. First, the following consequence of Lemma~\ref{lem:config} gives the sharper vertex count needed when all four graphs are connected.

\begin{lemma}\label{lem:vloss}
Let $\Lambda_1,\dots,\Lambda_4$ be connected copies such that every vertex and every edge belongs to at least two of them, $\E[\prod_kf_k(\Lambda_k)]\ne0$ for some $f_k\in \mathcal M_{G_k}$, and the list does not split into identical pairs. Set $d_{\mathrm{tot}}:=\sum_k d(\Lambda_k)$. Then $d_{\mathrm{exc}}=d_{\mathrm{tot}}-2v\ge2$, where $v$ is the number of distinct vertices.
\end{lemma}

\begin{proof}[Proof of Lemma~\ref{lem:vloss}]
For each distinct vertex $x$ in the union, let $\mathrm{mult}(x)$ be the number of copies containing it. Since every vertex belongs to at least two copies, we have
\[
d_{\mathrm{exc}}=\sum_x\{\mathrm{mult}(x)-2\}\ge0.
\]
By assumption, the copies do not split into identical pairs. Thus, Lemma~\ref{lem:config} ensures that there is a vertex occurring at least three times, and $d_{\mathrm{exc}}\ge1$.

Since $d_{\mathrm{exc}}$ is an integer, it remains to rule out $d_{\mathrm{exc}}=1$. Suppose this equality holds. Then, one vertex $x$ occurs three times and every other vertex occurs twice. Exactly one copy in the list, denoted by $Q$, omits $x$. We show that another copy is identical to $Q$, which would leave at most two copies containing $x$.

Write $V(Q)$ and $E(Q)$ for the vertex set and the edge set of $Q$, with $V(Q)$ containing both row and column vertices. For $u\in V(Q)$, let $P(u)$ be the unique other copy containing $u$. If $uw\in E(Q)$, the edge must occur in another copy. That copy contains both endpoints, so it is both $P(u)$ and $P(w)$. Thus, $P(u)=P(w)$ along every edge of $Q$. Because any two vertices of $Q$ are joined by a path, $P(u)$ is the same copy for every $u\in V(Q)$. Denote it by $P$. Then,
\[
V(Q)\subseteq V(P),\qquad E(Q)\subseteq E(P).
\]
Conversely, any edge of $P$ incident to a vertex $u\in V(Q)$ must occur in another copy. Since $u$ belongs only to $P$ and $Q$, that copy must be $Q$. Thus, $P$ has no additional edge between vertices of $Q$ and no edge joining a vertex of $Q$ to a vertex outside $Q$. Because $P$ is connected and contains $Q$, this implies $V(P)=V(Q)$ and then $E(P)=E(Q)$.

Therefore, $P$ and $Q$ are identical labeled graphs occurring at two different positions in the list. Both omit $x$, so $x$ belongs to at most the other two copies, contradicting $\mathrm{mult}(x)=3$. This proves $d_{\mathrm{exc}}\ge2$.
\end{proof}

The next lemma bounds sums of products of four graph projections by counting the possible vertex labels. It uses the patterns and labelings defined in Section~\ref{sec:S-ahk}.
\begin{lemma}\label{lem:auto-pattern}
Assume (M4). For a copy $\Lambda$ of $G$, write $q_G=q_\Lambda$ and $c_G=c_{N(\Lambda)}$, so that $c_G\|q_G\|_2=\varsigma_G/M_N(G)^{1/2}$ by \eqref{eq:copy-scale}.
\begin{enumerate}[(i)]
\item Consider copies \(\Lambda_1,\dots,\Lambda_4\) of \(G_1,\dots,G_4\), with a fixed pattern of shared vertices using \(v_r\) distinct rows and \(v_c\) distinct columns. Put \(r_{\mathrm{tot}}=\sum_kr(G_k)\) and \(c_{\mathrm{tot}}=\sum_kc(G_k)\). The sum of \(\prod_kc_{G_k}\,|\E[\prod_kq_{\Lambda_k}]|\) over all labelings of the pattern is at most \(C\,m^{v_r-r_{\mathrm{tot}}/2}n^{v_c-c_{\mathrm{tot}}/2}\prod_k\varsigma_{G_k}\). If the expectation is nonzero, then \(v_r\le r_{\mathrm{tot}}/2\) and \(v_c\le c_{\mathrm{tot}}/2\). Thus, a pattern with \((r_{\mathrm{tot}}-2v_r)+(c_{\mathrm{tot}}-2v_c)\ge j\) has a factor at most \(L_N^{-j/2}\).
\item \(\sum_G\varsigma_G\le C\sigma_N\) and \(\sum_{G\,\disc}\varsigma_G\le C\sigma_{\mathrm{disc}}\).
\end{enumerate}
\end{lemma}

\begin{proof}[Proof of Lemma~\ref{lem:auto-pattern}]

\textbf{Proof of (i).} Fix a pattern of shared vertices among copies $\Lambda_1,\dots,\Lambda_4$ of $G_1,\dots,G_4$. Since $q_G$ is a constant multiple of $p_G$, (M4) and H\"older's inequality give
\[
\left|\E\left[\prod_{k=1}^4q_{\Lambda_k}\right]\right|
\le\prod_{k=1}^4\|q_{G_k}\|_4
\le C\prod_{k=1}^4\|q_{G_k}\|_2.
\]
There are at most $m^{v_r}n^{v_c}$ labelings of the pattern. Thus, by \eqref{eq:copy-scale}, the sum in (i) is at most
\[
C m^{v_r}n^{v_c}\prod_{k=1}^4c_{G_k}\|q_{G_k}\|_2
=C m^{v_r}n^{v_c}\prod_{k=1}^4\frac{\varsigma_{G_k}}{\sqrt{M_N(G_k)}}.
\]
Since $M_N(G)\asymp m^{r(G)}n^{c(G)}$, this gives the bound
\[
C\,m^{v_r-r_{\mathrm{tot}}/2}n^{v_c-c_{\mathrm{tot}}/2}\prod_{k=1}^4\varsigma_{G_k}.
\]
If the expectation is nonzero, Lemma~\ref{lem:config} requires every vertex to occur at least twice. Thus, $v_r\le r_{\mathrm{tot}}/2$ and $v_c\le c_{\mathrm{tot}}/2$. Since $L_N=\min(m,n)$, it follows that
\[
m^{v_r-r_{\mathrm{tot}}/2}n^{v_c-c_{\mathrm{tot}}/2}
\le L_N^{-\{(r_{\mathrm{tot}}-2v_r)+(c_{\mathrm{tot}}-2v_c)\}/2}
\le L_N^{-j/2}
\]
whenever $(r_{\mathrm{tot}}-2v_r)+(c_{\mathrm{tot}}-2v_c)\ge j$.

\textbf{Proof of (ii).} Since $\Gamma$ is finite and $\sum_G\varsigma_G^2=\sigma_N^2$, the Cauchy--Schwarz inequality gives $\sum_G\varsigma_G\le C\sigma_N$. Restricting the sum to disconnected graphs gives the second inequality.
\end{proof}

Summing the bounds of Lemma~\ref{lem:auto-pattern} over the patterns of shared vertices gives the following bound on absolute fourth products.

\begin{lemma}\label{lem:S-four}
Assume (M4). For all $m\ge2p$ and $n\ge2q$, with $\sigma_N^2=\var(U_N)$ at these dimensions,
\begin{equation}\label{eq:auto-four-projection-bound}
\sum_{A,A',B,B'}c_Ac_{A'}c_Bc_{B'}\bigl|\E[X_AX_{A'}X_BX_{B'}]\bigr|\le C\sigma_N^4,
\end{equation}
and consequently $\E[(U_N-\theta_N)^4]\le C\sigma_N^4$, where $C$ depends only on $(p,q,C_{\mathrm{M4}})$.
\end{lemma}

\begin{proof}[Proof of Lemma~\ref{lem:S-four}]
Expand each $X_A$ into graph projections. Summing the bounds in Lemma~\ref{lem:auto-pattern}(i) over the finitely many patterns of shared vertices and applying part~(ii), we obtain
\[
\sum_{A,A',B,B'}c_Ac_{A'}c_Bc_{B'}\bigl|\E[X_AX_{A'}X_BX_{B'}]\bigr|\le C\Big(\sum_G\varsigma_G\Big)^4\le C\sigma_N^4.
\]
By \eqref{eq:auto-U-expansion}, $\E[(U_N-\theta_N)^4]$ is the same sum without absolute values, so it satisfies the same bound.
\end{proof}

\subsection{Contraction bounds and fourth cumulants}\label{sec:S-cumulants}\label{sec:S-contraction}\label{sec:S-support}
The central limit theorem requires a bound on the fourth cumulant of sums of the orthogonal contributions $T_R$ defined in Section~\ref{sec:S-ahk}. This section proves the following bound.

\begin{lemma}\label{lem:S-kappa}
Let $m\ge2p$, $n\ge2q$, $\sigma_N>0$, and assume (M4) with constant $C_{\mathrm{M4}}$. For any set $\mathcal R$ of indices $R$, put $T_{\mathcal R}=\sum_{R\in\mathcal R}T_R$. Then
\[
\big|\E[T_{\mathcal R}^4]-3(\E[T_{\mathcal R}^2])^2\big|
\le C\big(L_N^{-1}\sigma_N^4+L_N^{-1/2}\sigma_{\mathrm{disc}}\sigma_N^3
 +\Upsilon_N^{1/4}\sigma_N^4\big),
\]
where $C$ depends only on $(p,q,C_{\mathrm{M4}})$.
\end{lemma}

We use the configurations, patterns and $d_{\mathrm{exc}}$ defined in Section~\ref{sec:S-ahk}, and expand $\E[T_{\mathcal R}^4]$ over configurations. The paired configurations give $3(\E[T_{\mathcal R}^2])^2$, up to terms with $d_{\mathrm{exc}}\ge2$ that contribute $O(L_N^{-1}\sigma_N^4)$. It remains to bound the unpaired configurations (Lemma~\ref{lem:S-unpaired}). By Lemma~\ref{lem:S-cover}, a pattern with $d_{\mathrm{exc}}\ge1$ carries a factor $L_N^{-d_{\mathrm{exc}}/2}\le L_N^{-1/2}$. In contrast, patterns with $d_{\mathrm{exc}}=0$ carry no such factor, and we bound them with contractions, which give the term in $\Upsilon_N$.

\begin{lemma}\label{lem:S-cover}
Assume (M4).
\begin{enumerate}[(i)]
\item Consider a configuration with $\E[\prod_k\pi_k]\ne0$. Then every variable of $\bigcup_kS_k$ belongs to at least two sets $S_k$, and every vertex and every edge of $\bigcup_k\Lambda_k$ belongs to at least two copies. Thus, $d_{\mathrm{exc}}\ge0$.
\item Every labeled projection $f_S(\Lambda)$ of a copy $\Lambda$ of $G$ satisfies $\|f_S\|_2\le\|p_G\|_2$ and $\|f_S\|_4\le C\|p_G\|_2$.
\item Consider four labeled projections with a fixed pattern of shared vertices in which every vertex belongs to at least two copies. With the notation of Lemma~\ref{lem:auto-pattern}(i), the sum of $\prod_k\beta_{G_k}\|\pi_k\|_4$ over all labelings of the pattern is at most
\[
C\,m^{v_r-r_{\mathrm{tot}}/2}n^{v_c-c_{\mathrm{tot}}/2}\prod_k\varsigma_{G_k}
\le CL_N^{-d_{\mathrm{exc}}/2}\prod_k\varsigma_{G_k}.
\]
The same bound holds for the sum of $\prod_k\beta_{G_k}|\E[\prod_k\pi_k]|$.
\end{enumerate}
\end{lemma}

\begin{proof}[Proof of Lemma~\ref{lem:S-cover}]
\textbf{Proof of (i).} If a variable $x$ belongs to only one set $S_k$, integrating $\pi_k$ over $x$ gives zero by degeneracy. The other factors do not depend on $x$, so $\E[\prod_k\pi_k]=0$.

By Section~\ref{sec:S-ahk}, $S_k$ contains every edge variable of $\Lambda_k$ and the variable of every isolated vertex. Each occurs in another $S_l$, so every edge and every isolated vertex belongs to another copy. Every remaining vertex is an endpoint of an edge and therefore also belongs to another copy, even if its variable is absent from $S_k$. Hence, $\sum_kd(\Lambda_k)\ge2v$ and $d_{\mathrm{exc}}\ge0$.

\textbf{Proof of (ii).} Equation~\eqref{eq:S-efron-stein} expresses $f_S$ as a signed sum of $2^{|S|}$ conditional expectations of $p_G$. Conditional expectation does not increase the $L^4$ norm, so the triangle inequality and (M4) give
\[
\|f_S\|_4\le 2^{|S|}\|p_G\|_4\le C\|p_G\|_2,
\]
where $C$ depends only on $(p,q,C_{\mathrm{M4}})$, since $|S|\le p+q+pq$. Also, $\|f_S\|_2\le\|p_G\|_2$ by orthogonality.

\textbf{Proof of (iii).} We follow the proof of Lemma~\ref{lem:auto-pattern}(i), using this fourth-moment bound and the identity $\beta_G\|p_G\|_2=\varsigma_G/\sqrt{M_N(G)}$. Since every vertex belongs to at least two copies, $v_r\le r_{\mathrm{tot}}/2$ and $v_c\le c_{\mathrm{tot}}/2$, which gives the factor $L_N^{-d_{\mathrm{exc}}/2}$. Only the variables of non-isolated vertices can be included or omitted, so each copy has at most $2^{p+q}$ labeled projections. Summing over these choices only changes the constant. Finally, H\"older's inequality gives $|\E[\prod_k\pi_k]|\le\prod_k\|\pi_k\|_4$.
\end{proof}

It remains to treat the case $d_{\mathrm{exc}}=0$. We use the contractions defined before \eqref{eq:gamma} and the labeled projections of Section~\ref{sec:S-ahk}. For $f_S(\Lambda)$, we apply the same definition by grouping the variables in $S$ according to the connected components of $\Lambda$.

For square-integrable functions $g(x_A,x_B)$ and $h(x_A,x_{B'})$ on independent groups of variables $A,B,B'$, put $F(x_B,x_{B'})=\E_{\mathbf{X}_A}[g(\mathbf{X}_A,x_B)h(\mathbf{X}_A,x_{B'})]$. Expanding the square with an independent copy $\mathbf{X}_A'$ of $\mathbf{X}_A$ and integrating first over $B$ and $B'$ gives the equality below. The Cauchy--Schwarz inequality gives the first bound:
\begin{equation}\label{eq:S-contr}
\|F\|_2^2=\big\langle g\otimes_Bg,\ h\otimes_{B'}h\big\rangle\le\|g\otimes_Bg\|_2\,\|h\otimes_{B'}h\|_2,\qquad \|g\otimes_Bg\|_2\le\|g\|_2^2.
\end{equation}
The inner product is over $(\mathbf{X}_A,\mathbf{X}_A')$. For the second bound, apply the Cauchy--Schwarz inequality in $\mathbf{X}_B$ and integrate over these two remaining arguments.

For a square-integrable kernel $g(x_A,x_B)$, define the integral operator
\[
(T_gu)(x_A)=\int g(x_A,x_B)u(x_B)\,dP_B(x_B),
\]
where $P_B$ is the law of $x_B$. Its adjoint $T_g^*:L^2(P_A)\to L^2(P_B)$ is given by
\[
(T_g^*v)(x_B)=\int g(x_A,x_B)v(x_A)\,dP_A(x_A),
\]
where $P_A$ is the law of $x_A$. For a Hilbert--Schmidt operator $T$, we define its \emph{Schatten $4$-norm} by
\[
\|T\|_{S_4}:=\left(\operatorname{tr}[(TT^*)^2]\right)^{1/4}.
\]
Therefore, for the kernel norm defined in Section~\ref{sec:diagnostic-main} of the paper, integrating first over the variables in $B$ gives
\begin{equation}\label{eq:S-contraction-identity}
\|g\|_{S_4}^4=\|g\otimes_Bg\|_2^2
=\operatorname{tr}[(T_gT_g^*)^2] = \|T_g\|_{S_4}^4.
\end{equation}
The following lemma compares the contraction of a kernel with those of its labeled projections.

\begin{lemma}\label{lem:S-blocks}\label{lem:S-orbit}
Let $X$ be a Hilbert--Schmidt operator with $X=\sum_{k,l}X_{kl}$. Suppose $X_{kl}$ maps $\mathcal H_l$ into $\mathcal K_k$ and vanishes on $\mathcal H_l^\perp$, where $(\mathcal H_l)$ and $(\mathcal K_k)$ are orthogonal families of closed subspaces. If $K$ operator blocks are nonzero, then $\sum_{k,l}\|X_{kl}\|_{S_4}^4\le\|X\|_{S_4}^4\le K^3\sum_{k,l}\|X_{kl}\|_{S_4}^4$.

For finite disjoint groups of independent coordinates $A,B$, let $g(x_A,x_B)=\sum_S g_S$ be the Hoeffding decomposition of a square-integrable kernel. If $K$ projections are nonzero, then
\[
\sum_S\|g_S\otimes_Bg_S\|_2^2
\le\|g\otimes_Bg\|_2^2
\le K^3\sum_S\|g_S\otimes_Bg_S\|_2^2.
\]
\end{lemma}

\begin{proof}[Proof of Lemma~\ref{lem:S-blocks}]
Put $X_k=\sum_lX_{kl}$. Since the $\mathcal H_l$ are pairwise orthogonal, and so are the $\mathcal K_k$,
\[
X_kX_k^*=\sum_lX_{kl}X_{kl}^*,
\qquad
X^*X=\sum_kX_k^*X_k.
\]
For positive operators $P,Q$, we have $\operatorname{tr}(PQ)\ge0$. Expanding the squares therefore gives
\[
\begin{aligned}
\|X\|_{S_4}^4
&=\operatorname{tr}\!\left[\left(\sum_kX_k^*X_k\right)^2\right]
\ge\sum_k\|X_k\|_{S_4}^4\\
&=\sum_k\operatorname{tr}\!\left[\left(\sum_lX_{kl}X_{kl}^*\right)^2\right]
\ge\sum_{k,l}\|X_{kl}\|_{S_4}^4.
\end{aligned}
\]
For the upper bound, the triangle inequality and H\"older's inequality, applied to the $K$ nonzero blocks, give
\[
\|X\|_{S_4}^4
\le\left(\sum_{k,l}\|X_{kl}\|_{S_4}\right)^4
\le K^3\sum_{k,l}\|X_{kl}\|_{S_4}^4.
\]

For the kernel inequality, write $g=\sum_Sg_S$ and $S=S_A\cup S_B$. The operator $T_{g_S}$ maps the Hoeffding space indexed by $S_B$ into that indexed by $S_A$ and vanishes on the orthogonal complement of the former. These spaces are orthogonal for distinct variable sets, and distinct $S$ give distinct pairs $(S_A,S_B)$. Thus, the $T_{g_S}$ are the blocks of $T_g$. Applying the two operator bounds and the contraction identity \eqref{eq:S-contraction-identity} proves the claim.
\end{proof}

We then identify how the connected components must be shared between copies.

\begin{lemma}\label{lem:S-atoms}
If $\E[\prod_k\pi_k]\ne0$ and $d_{\mathrm{exc}}=0$, each connected component of a copy is also a connected component of exactly one other copy, with the same selected variables on that component.
\end{lemma}

\begin{proof}[Proof of Lemma~\ref{lem:S-atoms}]
By Lemma~\ref{lem:S-cover}(i), every vertex and every edge belongs to at least two copies. Since $d_{\mathrm{exc}}=0$, the total number of vertex occurrences is $2v$, so every vertex belongs to exactly two copies.

Let $Q$ be a connected component of $\Lambda_k$. As in the proof of Lemma~\ref{lem:vloss}, the endpoints of each edge of $Q$ have the same unique other copy. Since $Q$ is connected, a single copy $\Lambda_l$, with $l\ne k$, contains all vertices and edges of $Q$. Conversely, any edge of $\Lambda_l$ incident to a vertex of $Q$ must also belong to $\Lambda_k$, the only other copy containing that vertex. It therefore belongs to $Q$, so $Q$ is also a connected component of $\Lambda_l$.

Finally, any variable on $Q$ selected in $S_k$ or $S_l$ must occur in another set by Lemma~\ref{lem:S-cover}(i). Only these two copies contain the vertices of $Q$, so $S_k$ and $S_l$ select the same variables on $Q$.
\end{proof}

For such a configuration, the distinct matched components have disjoint vertex sets. Their selected variables therefore form independent groups, each occurring in exactly two factors. The following lemma supplies the contraction needed when $d_{\mathrm{exc}}=0$.

\begin{lemma}\label{lem:S-split}
Suppose the configuration is unpaired, with $\E[\prod_k\pi_k]\ne0$ and $d_{\mathrm{exc}}=0$. Then some copy $\Lambda_k$ shares connected components with at least two other copies, so it is disconnected. Let $\Lambda_l$ be one of these copies, and let $B$ be the set of components of $\Lambda_k$ not shared with $\Lambda_l$. Write $i,j$ for the other two indices. Then
\[
\Big|\E[\prod_{k'}\pi_{k'}]\Big|\le\|f_{S_k}\otimes_Bf_{S_k}\|_2^{1/2}\,\|f_{S_l}\|_2\,\|f_{S_i}\|_2\,\|f_{S_j}\|_2 ,
\]
with $\varnothing\ne B\subsetneq\{\text{components of }\Lambda_k\}$.
\end{lemma}

\begin{proof}[Proof of Lemma~\ref{lem:S-split}]
If every copy shared its components with only one other copy, Lemma~\ref{lem:S-atoms} would give two identical pairs of labeled projections, contrary to the hypothesis.

We integrate first over the components shared by $\Lambda_k$ and $\Lambda_l$. Let $A$ be the set of these components, and let $B'$ be the set of components of $\Lambda_l$ not shared with $\Lambda_k$. Both $A$ and $B$ are nonempty, since $\Lambda_k$ shares components with $\Lambda_l$ and with another copy. Let $E$ be the set of components shared by $\Lambda_i$ and $\Lambda_j$. Put
\[
F(x_B,x_{B'})=\E_{\mathbf{X}_A}[\pi_k\pi_l],
\qquad K(x_B,x_{B'})=\E_{\mathbf{X}_E}[\pi_i\pi_j].
\]
Since these groups are independent, integrating over $A$ and $E$ first gives
\begin{equation}\label{eq:S-split-product}
\E[\prod_{k'}\pi_{k'}]
=\E_{\mathbf{X}_{B\cup B'}}[FK]
=\langle F,K\rangle.
\end{equation}
By \eqref{eq:S-contr},
\[
\|F\|_2\le\|f_{S_k}\otimes_Bf_{S_k}\|_2^{1/2}\,\|f_{S_l}\|_2.
\]
Next, the Cauchy--Schwarz inequality in $\mathbf{X}_E$ gives
\[
|K|^2\le\E_{\mathbf{X}_E}[\pi_i^2]\,\E_{\mathbf{X}_E}[\pi_j^2].
\]
Each group in $B\cup B'$ occurs in exactly one of $\pi_i,\pi_j$, so the two factors on the right depend on disjoint groups of remaining variables. Integrating gives
\[
\|K\|_2^2\le\E[\pi_i^2]\,\E[\pi_j^2]
=\|f_{S_i}\|_2^2\,\|f_{S_j}\|_2^2.
\]
Finally, applying the Cauchy--Schwarz inequality to \eqref{eq:S-split-product} and substituting these bounds gives
\[
\Big|\E[\prod_{k'}\pi_{k'}]\Big| \le\|F\|_2\,\|K\|_2 \le\|f_{S_k}\otimes_Bf_{S_k}\|_2^{1/2}\,\|f_{S_l}\|_2\,\|f_{S_i}\|_2\,\|f_{S_j}\|_2.
\]
\end{proof}

Now, we can sum the bound in Lemma~\ref{lem:S-split} over all unpaired configurations with $d_{\mathrm{exc}}=0$.

\begin{lemma}\label{lem:S-loss0}
The sum of $\prod_k\beta_{G_k}|\E[\prod_k\pi_k]|$ over all unpaired configurations with $d_{\mathrm{exc}}=0$ is at most $C\sigma_N^4\Upsilon_N^{1/4}$.
\end{lemma}

\begin{proof}[Proof of Lemma~\ref{lem:S-loss0}]
Fix an unpaired pattern with $d_{\mathrm{exc}}=0$ and nonzero expectation, including the selected variables. Lemma~\ref{lem:S-cover}(i) and $d_{\mathrm{exc}}=0$ imply that every vertex occurs exactly twice. Thus, $v_r=r_{\mathrm{tot}}/2$ and $v_c=c_{\mathrm{tot}}/2$. The count in the proof of Lemma~\ref{lem:auto-pattern}(i) gives at most $C\prod_kM_N(G_k)^{1/2}$ labelings. By exchangeability, each labeling has the same expectation.

Choose $k,l,i,j$ and $B$ as in Lemma~\ref{lem:S-split}. Since $G_k$ is disconnected and $B$ is nonempty and proper, Lemma~\ref{lem:S-blocks} and \eqref{eq:gamma} give
\[
w_{G_k}^2\|f_{S_k}\otimes_Bf_{S_k}\|_2^2
\le w_{G_k}^2\|p_{G_k}\otimes_Bp_{G_k}\|_2^2
\le\sigma_N^4\Upsilon_N.
\]
By the normalization in Section~\ref{app:projection-notation}, $M_N(G)^{1/2}\beta_G=w_G^{1/2}$ and $w_G^{1/2}\|p_G\|_2=\varsigma_G$. Multiplying the bound in Lemma~\ref{lem:S-split} by the labeling count and using $\|f_S\|_2\le\|p_G\|_2$ from Lemma~\ref{lem:S-cover}(ii) bounds the contribution of the pattern by
\[
C\big(w_{G_k}\|f_{S_k}\otimes_Bf_{S_k}\|_2\big)^{1/2}
\varsigma_{G_l}\varsigma_{G_i}\varsigma_{G_j}
\le C\sigma_N^4\Upsilon_N^{1/4},
\]
where we used $\varsigma_G\le\sigma_N$. Summing over the finitely many graph types, choices of selected variables and patterns proves the result.
\end{proof}

Finally, the next lemma collects the bounds for unpaired configurations.

\begin{lemma}\label{lem:S-unpaired}
The sum of $\prod_k\beta_{G_k}|\E[\prod_k\pi_k]|$ over all unpaired configurations is at most
\[
C\big(L_N^{-1}\sigma_N^4+L_N^{-1/2}\sigma_{\mathrm{disc}}\sigma_N^3+\Upsilon_N^{1/4}\sigma_N^4\big).
\]
\end{lemma}

\begin{proof}[Proof of Lemma~\ref{lem:S-unpaired}]
First, we consider unpaired configurations whose graphs are all connected. For a nonzero expectation, Lemma~\ref{lem:S-cover}(i) ensures that every vertex and edge occurs at least twice. If the copies do not form identical pairs, Lemma~\ref{lem:vloss} gives $d_{\mathrm{exc}}\ge2$. Otherwise, relabel so that $\Lambda_1=\Lambda_2$ and $\Lambda_3=\Lambda_4$. A shared vertex belongs to all four copies, again giving $d_{\mathrm{exc}}\ge2$. If, instead, the two vertex sets are disjoint, the pairs $(\pi_1,\pi_2)$ and $(\pi_3,\pi_4)$ are independent, so
\[
\E[\pi_1\pi_2\pi_3\pi_4]
=\E[\pi_1\pi_2]\E[\pi_3\pi_4]=0.
\]
The last equality follows from orthogonality, since at least one pair contains distinct labeled projections. Lemmas~\ref{lem:S-cover}(iii) and~\ref{lem:auto-pattern}(ii) therefore bound the total absolute contribution by $CL_N^{-1}\sigma_N^4$.

Now, for unpaired configurations containing a disconnected graph and satisfying $d_{\mathrm{exc}}\ge1$, Lemma~\ref{lem:S-cover}(iii) bounds each pattern by $CL_N^{-1/2}\prod_k\varsigma_{G_k}$. By Lemma~\ref{lem:auto-pattern}(ii), the sum over one disconnected graph gives a factor at most $C\sigma_{\mathrm{disc}}$, and each of the other three graph sums gives at most $C\sigma_N$. Thus, the total absolute contribution is at most $CL_N^{-1/2}\sigma_{\mathrm{disc}}\sigma_N^3$.

Finally, we consider unpaired configurations with $d_{\mathrm{exc}}=0$. Lemma~\ref{lem:S-loss0} bounds their total contribution. Adding the three cases proves the lemma.
\end{proof}

Now, we can prove Lemma~\ref{lem:S-kappa}.

\begin{proof}[Proof of Lemma~\ref{lem:S-kappa}]
Write $T_{\mathcal R}=\sum_\pi\beta_G\pi$, where the sum runs over labeled projections $\pi=f_S(\Lambda)$ whose type belongs to $\mathcal R$, and $G$ is the graph type of $\Lambda$. By orthogonality, $\E[T_{\mathcal R}^2]=\sum_\pi\beta_G^2\E[\pi^2]$. We expand the fourth moment and subtract $3(\E[T_{\mathcal R}^2])^2$ to cancel independent pairs.

First, we consider paired configurations. Two distinct projections appearing twice have six orderings. The sum over ordered pairs $(\pi,\pi')$ counts each choice twice, giving the coefficient $3$. By contrast, four identical projections have one ordering. Thus, the paired contribution after subtraction is
\[
3\sum_{\pi\ne\pi'}\beta_G^2\beta_{G'}^2\cov(\pi^2,\pi'^2)+\sum_\pi\beta_G^4\big\{\E[\pi^4]-3(\E[\pi^2])^2\big\},
\]
where $G'$ is the graph type of $\pi'$. The covariance is zero when $\pi,\pi'$ depend on disjoint sets of AHK variables. Every remaining term, including the diagonal sum, has a vertex occurring four times, so $d_{\mathrm{exc}}\ge2$. By H\"older's inequality, each covariance and each diagonal term is at most $4\prod_{k=1}^4\|\pi_k\|_4$ for the corresponding configuration. Lemma~\ref{lem:S-cover}(iii) therefore gives a factor $L_N^{-d_{\mathrm{exc}}/2}\le L_N^{-1}$. Summing by Lemma~\ref{lem:auto-pattern}(ii) bounds their total absolute contribution by $CL_N^{-1}\sigma_N^4$.

Finally, Lemma~\ref{lem:S-unpaired} bounds the total absolute contribution of the unpaired configurations. Adding the two bounds proves the result.
\end{proof}

\subsection{Key decomposition for \texorpdfstring{$\Vu$}{Vu}}\label{sec:S-variance}\label{app:inference-bounds}

 To bound the error of $\Vu$, we write $\Vu=J+E_N$, where $J$ sums products of graph projections whose copies share a row or a column, and bound $\E[E_N^2]$.

We use the grouped projections $X_A$ and the fractions $c_A$ of \eqref{eq:grouped-projections}, so that $U_N-\theta_N=\sum_Ac_AX_A$ by \eqref{eq:auto-U-expansion}.
Write $A\sim A'$ when the two pairs of vertex sets share a row or a column, including $A=A'$. Define
\begin{equation}\label{eq:variance-decomposition}
J:=\sum_{A\sim A'}c_Ac_{A'}X_AX_{A'},
\qquad E_N:=\Vu-J,
\end{equation}
where the sum is over ordered pairs.
For a graph $G$, write $\ell(G)=(r(G),c(G))$ for its numbers of rows and columns. For $\ell=(r,c)$ with $0\le r\le p$, $0\le c\le q$ and $r+c\ge1$, write $V_\ell=V_{r,c}$ and $\widehat V_\ell=\widehat V_{r,c}$ for the variance contributed by graphs with row and column counts $\ell$ and its unbiased estimator (Lem.~\ref{lem:overlap}). Then $\Vu=\sum_\ell\widehat V_\ell$ (Prop.~\ref{prop:pathwise}(ii)).

Write $\ell(A)=(|A_1|,|A_2|)$. Unions, intersections and inclusions of these pairs are taken separately for rows and columns. For every $\ell$, define
\begin{equation}\label{eq:component-J}
J_\ell:=\sum_{\ell(A\cap A')=\ell}c_Ac_{A'}X_AX_{A'},\qquad E_{N,\ell}:=\widehat V_\ell-J_\ell,
\end{equation}
where the sum is over the ordered pairs of vertex sets whose intersection has $r$ rows and $c$ columns, for $\ell=(r,c)$. Since $A\sim A'$ exactly when $\ell(A\cap A')\ne(0,0)$, we have $J=\sum_\ell J_\ell$. First, the following lemmas bound the coefficients in $\widehat V_\ell-J_\ell$.

\begin{lemma}\label{lem:S-coefficients}
Let $m\ge2p$ and $n\ge2q$. Then for every $\ell=(r,c)$ with $0\le r\le p$, $0\le c\le q$ and $r+c\ge1$, the estimator has the expansion
\[
\widehat V_\ell=\sum_{A,A'}b_\ell(A,A')X_AX_{A'}
\]
with deterministic coefficients satisfying
\begin{equation}\label{eq:S-level-coefficients}
\left|b_\ell(A,A')-c_Ac_{A'}\one\{\ell(A\cap A')=\ell\}\right|
\le CL_N^{-1}c_Ac_{A'},
\end{equation}
where $C$ depends only on $(p,q)$.
\end{lemma}

\begin{proof}[Proof of Lemma~\ref{lem:S-coefficients}]
Fix $\ell=(r,c)$. Recall from Section~\ref{sec:S-subarrays} that $K(a,b)$ counts the ordered pairs of subarrays sharing $a$ rows and $b$ columns. Every subarray $W$ has the same number $K(a,b)/K$ of such partners $W'$. Averaging $H_WH_{W'}$ over these pairs therefore gives
\[
\begin{aligned}
\Thetah(a,b)
&=\frac{1}{K(a,b)}
\sum_{\substack{W,W':\,|W_1\cap W'_1|=a\\|W_2\cap W'_2|=b}}
(H_W-\theta_N)(H_{W'}-\theta_N)+2\theta_NU_N-\theta_N^2\\
&=\frac{1}{K(a,b)}
\sum_{\substack{W,W':\,|W_1\cap W'_1|=a\\|W_2\cap W'_2|=b}}
\sum_{\substack{A\subseteq W\\A'\subseteq W'}}X_AX_{A'}
+2\theta_NU_N-\theta_N^2.
\end{aligned}
\]
Substituting into \eqref{eq:variance-component-estimator} gives
\[
\begin{aligned}
\widehat V_\ell
&=\sum_{A,A'}b_\ell(A,A')X_AX_{A'}\\
&\quad+\frac{\binom pr^2\binom qc^2}{\binom mr\binom nc}
(2\theta_NU_N-\theta_N^2)
\underbrace{\left[\sum_{a=0}^r(-1)^{r-a}\binom ra\right]
\left[\sum_{b=0}^c(-1)^{c-b}\binom cb\right]}_{=0\text{ since }r+c\ge1},
\end{aligned}
\]
where
\begin{equation}\label{eq:S-level-counting}
b_\ell(A,A')=\frac{\binom pr^2\binom qc^2}{\binom mr\binom nc}
\sum_{a=0}^r\sum_{b=0}^c(-1)^{r-a+c-b}\binom ra\binom cb
\frac{N_{a,b}(A,A')}{K(a,b)},
\end{equation}
and
\[
N_{a,b}(A,A')=
\#\left\{(W,W'):
\begin{array}{l}
|W_1\cap W'_1|=a,\quad |W_2\cap W'_2|=b,\\
A\subseteq W,\quad A'\subseteq W'
\end{array}
\right\}.
\]

Next, fix $A,A'$ and write
\[
\ell(A)=(x,y),\qquad \ell(A')=(x',y'),\qquad
\ell(A\cap A')=(s,t).
\]
If $a<s$ or $b<t$, no pair can contain $A$ and $A'$, so $N_{a,b}(A,A')=0$. Hence, if $s>r$ or $t>c$, both coefficients in \eqref{eq:S-level-coefficients} are zero.

Now, suppose $s\le r$, $t\le c$ and $(s,t)\ne(r,c)$. The coefficient in $J_\ell$ is zero, so we need to bound $|b_\ell(A,A')|$. For $a\ge s$ and $b\ge t$, the union of the two row sets has $2p-a$ elements and must contain the $x+x'-s$ rows of $A_1\cup A'_1$. Each possible union gives the same number of ordered pairs of row sets. Counting such unions for the rows and the columns gives
\[
\begin{aligned}
\frac{N_{a,b}(A,A')}{K(a,b)}
&\le
\frac{\binom{m-x-x'+s}{2p-a-x-x'+s}}{\binom m{2p-a}}
\frac{\binom{n-y-y'+t}{2q-b-y-y'+t}}{\binom n{2q-b}}\\
&=O\bigl(m^{-(x+x'-s)}n^{-(y+y'-t)}\bigr),
\end{aligned}
\]
where a binomial coefficient is zero when the prescribed choices are impossible. Substituting into \eqref{eq:S-level-counting} and using $c_Ac_{A'}\asymp m^{-x-x'}n^{-y-y'}$ gives
\[
\begin{aligned}
|b_\ell(A,A')|
&\le Cm^{-r}n^{-c}m^{-(x+x'-s)}n^{-(y+y'-t)}\\
&\le CL_N^{-1}c_Ac_{A'},
\end{aligned}
\]
because at least one of $r-s$ and $c-t$ is positive.

It remains to show that $b_\ell(A,A')=c_Ac_{A'}\{1+O(L_N^{-1})\}$ when $(s,t)=(r,c)$. Only $(a,b)=(r,c)$ contributes to \eqref{eq:S-level-counting}, and the common row set must be $A_1\cap A'_1$. First, choose $W_1$ containing $A_1$ and avoiding $A'_1\setminus A_1$, in $\binom{m-x-x'+r}{p-x}$ ways. Then, choose $W'_1$ containing $A'_1$ and avoiding $W_1\setminus A'_1$, in $\binom{m-p-x'+r}{p-x'}$ ways. There are $\binom mp\binom pr\binom{m-p}{p-r}$ ordered pairs of row sets sharing $r$ rows in total. Dividing the count for $A_1,A'_1$ by this total and multiplying by $\binom pr^2/\binom mr$ from \eqref{eq:S-level-counting} gives
\[
\frac{\binom pr^2}{\binom mr}
\frac{\binom{m-x-x'+r}{p-x}\binom{m-p-x'+r}{p-x'}}
{\binom mp\binom pr\binom{m-p}{p-r}}
=\frac{(p)_x(p)_{x'}}{(m)_r(m)_{x+x'-r}}.
\]
Row and column choices can be made separately. Applying the same count to the columns and using the definition of $c_A$ in \eqref{eq:grouped-projections} gives
\[
\begin{aligned}
\frac{b_\ell(A,A')}{c_Ac_{A'}}
&=\frac{(m)_x(m)_{x'}}{(m)_r(m)_{x+x'-r}}
\frac{(n)_y(n)_{y'}}{(n)_c(n)_{y+y'-c}}\\
&=1+O(m^{-1}+n^{-1})\\
&=1+O(L_N^{-1}),
\end{aligned}
\]
where we used $(m)_k=m^k\{1+O(m^{-1})\}$ and $(n)_k=n^k\{1+O(n^{-1})\}$ for bounded $k$. This proves \eqref{eq:S-level-coefficients} in the remaining case.
\end{proof}

\begin{lemma}\label{lem:S-level}
Let $m\ge2p$, $n\ge2q$, and assume (M4) on a class $\Pcal$ with constant $C_{\mathrm{M4}}$. Then for every $P\in\Pcal$ and every $\ell$,
\[
\qquad \E[E_{N,\ell}^2]\le CL_N^{-2}\sigma_N^4,
\]
where $C$ depends only on $(p,q,C_{\mathrm{M4}})$.
\end{lemma}

\begin{proof}[Proof of Lemma~\ref{lem:S-level}]
Fix $\ell$. Lemma~\ref{lem:S-coefficients} and \eqref{eq:component-J} give
\[
E_{N,\ell}
=\sum_{A,A'}
\bigl[b_\ell(A,A')-c_Ac_{A'}\one\{\ell(A\cap A')=\ell\}\bigr]X_AX_{A'}.
\]
Applying \eqref{eq:S-level-coefficients} and \eqref{eq:auto-four-projection-bound} gives
\[
\begin{aligned}
\E[E_{N,\ell}^2]
&\le CL_N^{-2}\sum_{A,A',B,B'}c_Ac_{A'}c_Bc_{B'}
\bigl|\E[X_AX_{A'}X_BX_{B'}]\bigr|\\
&\le CL_N^{-2}\sigma_N^4.
\end{aligned}
\]
\end{proof}

Recall that $L_N=\min(m,n)$ and that $\Upsilon_N$ is defined in \eqref{eq:gamma}.  The following lemma bounds $J$ and $E_N$.

\begin{lemma}\label{lem:auto-decomp}
Let $m\ge2p$, $n\ge2q$, and assume (M4) with constant $C_{\mathrm{M4}}$. Then
\[
\E[J]=\sigma_N^2,\qquad \E[E_N]=0,\qquad
\E[E_N^2]\le C L_N^{-2}\sigma_N^4.
\]
If $\sigma_N>0$, then
\begin{equation}\label{eq:app-J-bound}
\var(J)\le C\big(L_N^{-1}\sigma_N^4
 +L_N^{-1/2}\sigma_{\mathrm{disc}}\sigma_N^3
 +\Upsilon_N^{1/4}\sigma_N^4\big),
\end{equation}
where $C$ depends only on $(p,q,C_{\mathrm{M4}})$.
\end{lemma}

\begin{proof}[Proof of Lemma~\ref{lem:auto-decomp}]
\textbf{(i)}  First, we compute the expectations of $J$ and $E_N$, and then bound $\E[E_N^2]$. Since distinct graph projections are orthogonal, \(\E[X_AX_{A'}]=0\) unless \(A=A'\). Because the diagonal \(A=A'\) is included in \(J\) and $\E[\Vu]=\sigma_N^2$ (Sec.~\ref{sec:variance-estimators} of the paper),
\[
\E[J]=\sum_Ac_A^2\E[X_A^2]=\var(U_N)=\sigma_N^2,
\qquad
\E[E_N]=\E[\Vu]-\E[J]=0.
\]

Since $\Vu=\sum_\ell\widehat V_\ell$ and $J=\sum_\ell J_\ell$, Lemma~\ref{lem:S-level} gives
\[
E_N=\sum_\ell E_{N,\ell},\qquad
\E[E_N^2]\le C\sum_\ell\E[E_{N,\ell}^2]\le CL_N^{-2}\sigma_N^4.
\]

\textbf{(ii)} It remains to bound $\var(J)$ as stated in~\eqref{eq:app-J-bound}. Expanding the graph projections in $J$ into their labeled projections gives
\[
\var(J)=
\sum_{\substack{\pi_1,\dots,\pi_4\\N(\Lambda_1)\sim N(\Lambda_2)\\N(\Lambda_3)\sim N(\Lambda_4)}}
\left(\prod_{k=1}^4\beta_{G_k}\right)
\cov(\pi_1\pi_2,\pi_3\pi_4),
\]
where $\pi_k=f_{S_k}(\Lambda_k)$ and $N(\Lambda_k)$ is the pair of row and column vertex sets of $\Lambda_k$. The sum is over ordered quadruples. The normalization follows from $c_Aq_\Lambda=\beta_Gp_G(\Lambda)$ for a copy $\Lambda$ of $G$ with vertex sets $A$.

First, consider a paired configuration, so the four labeled projections form two identical pairs. For the pairing $\pi_1=\pi_2$, $\pi_3=\pi_4$, the covariance is zero if the two pairs depend on disjoint sets of AHK variables, because they are independent. Otherwise, a shared variable forces a vertex to belong to all four copies. For either crossed pairing $(\pi_1=\pi_3,\pi_2=\pi_4)$ or $(\pi_1=\pi_4,\pi_2=\pi_3)$, the copies $\Lambda_1$ and $\Lambda_2$ share a vertex, because $N(\Lambda_1)\sim N(\Lambda_2)$. The other two copies coincide with $\Lambda_1$ and $\Lambda_2$, so this vertex belongs to all four copies.

Every vertex of a paired configuration occurs at least twice. Thus, a nonzero covariance requires $d_{\mathrm{exc}}=\sum_kd(\Lambda_k)-2v\ge2$, where $v$ counts the distinct vertices. Moreover, H\"older's inequality gives
\[
|\cov(\pi_1\pi_2,\pi_3\pi_4)|
\le2\prod_{k=1}^4\|\pi_k\|_4.
\]
Lemma~\ref{lem:S-cover}(iii) therefore gives a factor $L_N^{-d_{\mathrm{exc}}/2}\le L_N^{-1}$. Summing over the finitely many patterns and using Lemma~\ref{lem:auto-pattern}(ii) bounds the total absolute contribution of paired configurations by $CL_N^{-1}\sigma_N^4$.

Next, for an unpaired configuration, at least one of $\E[\pi_1\pi_2]$, $\E[\pi_3\pi_4]$ is zero: if both were nonzero, orthogonality would force $\pi_1=\pi_2$ and $\pi_3=\pi_4$. Hence, the covariance is equal to $\E[\prod_k\pi_k]$. Lemma~\ref{lem:S-unpaired} therefore applies. It bounds sums of absolute contributions, so restricting to the configurations with $N(\Lambda_1)\sim N(\Lambda_2)$ and $N(\Lambda_3)\sim N(\Lambda_4)$ preserves the bound. Adding the paired contribution gives
\[
\var(J)\le C\bigl(L_N^{-1}\sigma_N^4
+L_N^{-1/2}\sigma_{\mathrm{disc}}\sigma_N^3
+\Upsilon_N^{1/4}\sigma_N^4\bigr).
\]

\end{proof}

\subsection{Proof of Theorem~\ref{thm:contraction}}\label{sec:S-thm1}\label{app:proofs}\label{app:core}

Lemmas~\ref{lem:auto-decomp} and~\ref{lem:S-kappa} give the two main parts of the proof: ratio-consistency of the unbiased variance estimator $\Vu$ and Gaussian convergence after division by the true standard deviation $\sigma_N$. We first prove these conclusions along arbitrary sequences of laws. We then deduce uniformity and replace the true variance $\sigma_N^2$ by its estimate $\Vu^+$.

\begin{proof}[Proof of Theorem~\ref{thm:contraction}]
\textbf{(i)} First, we prove the bound on the error of the variance estimator stated in the theorem. Lemma~\ref{lem:auto-decomp} gives $\Vu-\sigma_N^2=(J-\E[J])+E_N$, together with bounds on $\var(J)$ and $\E[E_N^2]$. Therefore,
\[
\E[(\Vu-\sigma_N^2)^2]\le2\var(J)+2\E[E_N^2].
\]
Divide by $\sigma_N^4$ and use $\sigma_{\mathrm{disc}}/\sigma_N=D_N^{1/2}$. The term $L_N^{-2}$ from the remainder $E_N$ is at most $L_N^{-1}$, which gives the asserted relative-error bound.

\textbf{(ii)} Next, we prove $(U_N-\theta_N)/\sigma_N\dto\Nor(0,1)$ along every sequence satisfying the assumptions and $\Upsilon_N\to0$. We use the de Jong-type central limit theorem of~\cite{dobler2017quantitative}. For a degenerate $U$-statistic of fixed order and unit variance, it requires a fourth moment tending to $3$ and a vanishing maximal contribution from any one independent variable.

The order in this theorem counts independent variables, so we use the decomposition~\eqref{eq:S-efron-stein}. For an index $R$ represented by $(\Lambda,S)$, put $k(R)=|S|$. Group the contributions as $U_N-\theta_N=\sum_{k=1}^{K}T^{[k]}$, where $T^{[k]}=\sum_{R:\,k(R)=k}T_R$ and $K=pq+p+q$. The $T^{[k]}$ are orthogonal. Each $T^{[k]}$ is a degenerate $U$-statistic of order $k$ in the sense of equation (1.6) of~\cite{dobler2017quantitative}: its summands are $\beta_Gf_S(\Lambda)$ with $|S|=k$, indexed by distinct sets of the independent variables $(\xi_a)_{a\le m}$, $(\eta_b)_{b\le n}$ and $(\zeta_{ab})$.

We only need Gaussian limits for the orders that contribute a nonvanishing proportion of the total variance $\sigma_N^2$. Indeed, the bound in Lemma~\ref{lem:S-kappa} is relative to $\sigma_N^4$. If $v_k:=\var(T^{[k]})$ is much smaller than $\sigma_N^2$, dividing that bound by $v_k^2$ need not give a vanishing error. Such a term can instead be omitted from $(U_N-\theta_N)/\sigma_N$ when $v_k/\sigma_N^2\to0$.

To separate these cases, take any subsequence and put $s_k=v_k^{1/2}/\sigma_N$. Since there are finitely many $s_k\in[0,1]$, along a further subsequence, $s_k\to s_{k,\infty}$ for every $k$, with $\sum_ks_{k,\infty}^2=1$. Let $\mathcal K=\{k:s_{k,\infty}>0\}$. By orthogonality,
\[
\E\left[\left(\frac{1}{\sigma_N}\sum_{k\notin\mathcal K}T^{[k]}\right)^2\right]
=\sum_{k\notin\mathcal K}s_k^2\longrightarrow0.
\]
For $k\in\mathcal K$, put $\widetilde T^{[k]}=T^{[k]}/v_k^{1/2}$.

Now, we check that no single independent variable has a nonnegligible contribution to $\widetilde T^{[k]}$. The squared maximal influence in the cited theorem is the largest, over these variables, of the sum of the variances of the summands involving that variable. By symmetry, a fraction $r(G)/m$ of the copies of $G$ contain row $a$. The summands involving $\xi_a$ must come from these copies, so their total variance in $\widetilde T^{[k]}$ is at most
\[
\frac{1}{v_k}\sum_{R:\,k(R)=k}\frac{r(G(R))}{m}\varsigma_R^2
\le\frac{p}{m v_k}\sum_{R:\,k(R)=k}\varsigma_R^2
=\frac{p}{m}.
\]
Similarly, for $\eta_b$, we count copies containing column $b$. For $\zeta_{ab}$, we count copies containing both row $a$ and column $b$. These fractions are $c(G)/n$ and $r(G)c(G)/(mn)$, so the same calculation gives bounds $q/n$ and $pq/(mn)$. Thus, the squared maximal influence is at most
\[
\max\left\{\frac{p}{m},\frac{q}{n},\frac{pq}{mn}\right\}
\le\frac{C}{L_N}\longrightarrow0.
\]

Next, for the fourth moment, Lemma~\ref{lem:S-kappa} gives
\[
\begin{aligned}
\left|\E[(\widetilde T^{[k]})^4]-3\right|
&=\frac{\left|\E[(T^{[k]})^4]-3v_k^2\right|}{v_k^2}\\
&\le\frac{C}{s_k^4}\left(L_N^{-1}+L_N^{-1/2}D_N^{1/2}+\Upsilon_N^{1/4}\right)
\longrightarrow0,
\end{aligned}
\]
since $s_k\to s_{k,\infty}>0$, $D_N\le1$ and $\Upsilon_N\to0$.

It remains to obtain joint convergence before adding the retained terms $\widetilde T^{[k]}$, $k\in\mathcal K$. Their covariance matrix is the identity by orthogonality. Moreover, their orders are distinct, so the additional condition (iv) of \citet[Thm.~1.7]{dobler2017quantitative}, concerning pairs of coordinates of the same order, is automatically satisfied. That theorem, or their Theorem~1.2 when $\mathcal K$ has one element, gives $(\widetilde T^{[k]})_{k\in\mathcal K}\dto\Nor(0,I)$. Since $\sum_{k\in\mathcal K}s_{k,\infty}^2=1$,
\[
\frac{U_N-\theta_N}{\sigma_N}
=\sum_{k\in\mathcal K}s_k\widetilde T^{[k]}+o_P(1)
\dto\Nor(0,1)
\]
along the further subsequence. Every subsequence has such a further subsequence with this same limit, which proves convergence along the whole sequence.

\textbf{(iii)} The argument in (ii) allows the law to change with $N$. Suppose that convergence to the standard normal distribution is not uniform over a class on which (M4) holds and $\Upsilon_N\to0$ uniformly. Then, there are $\varepsilon>0$ and laws $P_N$ in the class such that $\sup_x|\P_{P_N}((U_N-\theta_N)/\sigma_N\le x)-\Phi(x)|\ge\varepsilon$ along a subsequence. All bounds above depend only on $(p,q,C_{\mathrm{M4}})$, so the argument applies to this sequence of laws. It gives a further subsequence along which $(U_N-\theta_N)/\sigma_N$ converges to $\Nor(0,1)$. By P\'olya's theorem, its distribution function then converges uniformly in $x$, a contradiction. The relative-error bound is uniform because its constant $C$ also depends only on $(p,q,C_{\mathrm{M4}})$.

\textbf{(iv)} Finally, we replace the true variance $\sigma_N^2$ by $\Vu^+$. The bound in (i) gives $\Vu/\sigma_N^2\to1$ in probability. Taking the positive part preserves this convergence: the map $x\mapsto x^+$ is $1$-Lipschitz and leaves $\sigma_N^2$ unchanged. Thus, $\Vu^+/\sigma_N^2\to1$ and $\P(\Vu\le0)\to0$. The convention on $\{\Vu\le0\}$ therefore changes $\widehat Z_N$ only on an event whose probability tends to zero. Combining this with (ii), Slutsky's lemma gives $\widehat Z_N\dto\Nor(0,1)$. The argument in (iii) gives uniform convergence for $\widehat Z_N$ as well.
\end{proof}

\subsection{Proof of Proposition~\ref{prop:gammaD}}\label{sec:S-gammaD}

In this proof, $C$ denotes a constant depending only on $(p,q)$.

\begin{proof}[Proof of Proposition~\ref{prop:gammaD}]
By \eqref{eq:S-contr}, $w_G\|p_G\otimes_Bp_G\|_2\le w_G\|p_G\|_2^2=\varsigma_G^2$. There are at most $2^{p+q}$ choices of $B$ for each $G$, and $\sum_{G\ \disc}\varsigma_G^2=\sigma_N^2D_N$. Summing the squared bounds gives
\[
\Upsilon_N\le C\sigma_N^{-4}\sum_{G\ \disc}\varsigma_G^4
\le C\sigma_N^{-4}\Big(\sum_{G\ \disc}\varsigma_G^2\Big)^2
=CD_N^2.
\]

For the lower bound, fix a disconnected graph $G$ and a nonempty proper set $B$ of its connected components, with complement $A$.
Let $g(x_A,x_B)\ne0$ be square-integrable, and choose an orthonormal basis $(e_j)$ in the $A$ variables. Some coefficient $g_j(x_B)=\E_{\mathbf{X}_A}[g(\mathbf{X}_A,x_B)e_j(\mathbf{X}_A)]$ is nonzero, since otherwise $g=0$. By the definition of the contraction,
\[
\big\langle g\otimes_Bg,\ e_j(x_A)e_j(x_A')\big\rangle
=\|g_j\|_2^2>0.
\]
Thus, $\|g\otimes_Bg\|_2>0$. Moreover, by the bilinearity of the contraction and the Cauchy--Schwarz inequality, every square-integrable $h(x_A,x_B)$ satisfies
\[
\|g\otimes_Bg-h\otimes_Bh\|_2
\le(\|g\|_2+\|h\|_2)\|g-h\|_2,
\]
so the contraction norm $g\mapsto\|g\otimes_Bg\|_2$ is continuous.

Now, under (K), the closure of the normalized projections $p_{P,G}/\|p_{P,G}\|_2$ is compact, and every element has norm one. The contraction norm therefore has a positive lower bound $\eta_G$ on this closure. Rescaling gives
\[
w_G^2\|p_G\otimes_Bp_G\|_2^2
\ge\eta_G^2w_G^2\|p_G\|_2^4
=\eta_G^2\varsigma_G^4.
\]
By \eqref{eq:gamma}, summing over the finitely many disconnected graph classes gives
\[
\Upsilon_N \ge c\sigma_N^{-4}\sum_{G\ \disc}\varsigma_G^4 \ge\frac{c}{|\Gamma|}\sigma_N^{-4}\Big(\sum_{G\ \disc}\varsigma_G^2\Big)^2
=\frac{c}{|\Gamma|}D_N^2,
\]
where $c>0$ depends on $\Pcal$, and the second inequality follows from the Cauchy--Schwarz inequality.
\end{proof}

\section{Joint limits and the Gaussian boundary}\label{sec:S-part-limits}

In this section, we prove Lemma~\ref{lem:wick}, Theorems~\ref{thm:chaos} and~\ref{thm:sharp}, and Corollary~\ref{cor:S-consistency} (Sec.~\ref{app:wick-construction}, \ref{sec:S-thm-chaos}, \ref{app:sharp} and~\ref{sec:S-cor-consistency}). In particular, the proofs of Lemma~\ref{lem:wick} and Theorem~\ref{thm:chaos} follow the strategy of \citet[Lem.~9, Thm.~3]{janson1991asymptotic}, adapted to copies with row and column vertices. We develop it in detail.

\subsection{Gaussian space and Wiener chaos}\label{app:wick}\label{sec:S-wick}

Let $\Gamma_0$ be a finite set of nonempty classes, as in Lemma~\ref{lem:wick}. For each connected class $H$ occurring as a component of a graph in $\Gamma_0$, set $\mathfrak h_H=\mathcal M_H^{\rm sym}$, and let $X$ be an isonormal Gaussian process on
\begin{equation}\label{eq:wick-gaussian-space}
\mathfrak h=\bigoplus_H\mathfrak h_H.
\end{equation}
Let $X_H$ be its restriction to $\mathfrak h_H$.

We write $\widehat\bigotimes$ for the completed Hilbert tensor product. In particular, for $L^2$ spaces of independent variables, $L^2(\mu_1)\widehat\otimes L^2(\mu_2)=L^2(\mu_1\times\mu_2)$. Let $\mathfrak h^{\otimes b}$ be the $b$-fold completed tensor power of $\mathfrak h$, and $\mathfrak h^{\odot b}$ the symmetric subspace of this power. The operator $\operatorname{Sym}_b$ averages over the $b!$ permutations of the tensor arguments. For the multiple Wiener integrals $I_b$ associated with $X$, which take values in the $b$-th Wiener chaos \citep[Thm.~7.26]{janson1997gaussian}, we use the convention
\begin{equation}\label{eq:multiple-integral-isometry}
\E[I_b(u)I_b(v)]=b!\langle u,v\rangle,
\qquad u,v\in\mathfrak h^{\odot b}.
\end{equation}

To determine when a limit is Gaussian, we use the following consequence of the tail theorem of \citet{janson1997gaussian}.

\begin{lemmaL}[Consequence of \citealp{janson1997gaussian}, Thm.~6.12]\label{lem:chaosnotgauss}
Let $X$ be an isonormal Gaussian process over a real separable Hilbert space $\mathfrak h$, and let
\[
F=\sum_{b=0}^B I_b(f_b),\qquad f_b\in\mathfrak h^{\odot b},
\]
where $I_b(f_b)$ denotes the $b$-th multiple Wiener integral (with $I_0(f_0)$ constant). If $B\ge2$ and $f_B\ne0$, then $F$ is not Gaussian.
\end{lemmaL}

Indeed, the cited tail theorem gives $\P(|F|>t)\ge\exp(-c t^{2/B})$ for all sufficiently large $t$ and some $c>0$. When $B\ge2$, this lower bound is incompatible with Gaussian tails, whose logarithm has order $-t^2$.

\subsection{Gaussian limits for connected graphs}\label{sec:clt}\label{app:gauss}

\begin{lemma}[Consequence of \citealp{leminh2026characterization}, Lem.~10]\label{lem:clt}
On each of finitely many connected graphs $G\in\Gamma$, take a finite orthonormal family of fixed bounded functions $\phi_{G,\ell}\in \mathcal M_G^{\rm sym}$. Then the sums
\[
Z_{N,G,\ell}:=M_N(G)^{-1/2}\sum_{\Lambda\text{ copy of }G}\phi_{G,\ell}(\Lambda)
\]
converge jointly to independent standard normal random variables.
\end{lemma}

\begin{proof}[Proof of Lemma~\ref{lem:clt}]
For one function on each graph, the convergence in distribution is the cited result, with the normalization of \eqref{eq:copy-normalization}. Since distinct copies are orthogonal, $\var(Z_{N,G,\ell})=1$ exactly. To allow several functions on the same graph, apply the result to their linear combinations, separately on each graph. The squared $L^2$ norm of $\sum_\ell a_\ell\phi_{G,\ell}$ is $\sum_\ell a_\ell^2$. The joint limit for distinct graphs and the Cram\'er--Wold theorem therefore give the asserted joint convergence.

\end{proof}

\subsection{The Wick approximation}\label{app:wick-products}

Fix connected classes $H_1,\ldots,H_b$ and bounded $\phi_j\in \mathcal M_{H_j}^{\rm sym}$. Write $M_j=M_N(H_j)$ and
\[
Z_{N,j}:=M_j^{-1/2}\sum_{\Lambda\text{ copy of }H_j}\phi_j(\Lambda).
\]
For $\mathcal J\subseteq[b]$, let $\Xi_N(\mathcal J)$ be the following normalized sum over ordered tuples $(\Lambda_j)_{j\in \mathcal J}$ that are pairwise vertex-disjoint,
\[
\Xi_N(\mathcal J):=\prod_{j\in \mathcal J}M_j^{-1/2}
\sum_{\substack{\Lambda_j\text{ copy of }H_j,\ j\in \mathcal J\\ V(\Lambda_r)\cap V(\Lambda_s)=\varnothing,\ r\ne s}}
\prod_{j\in \mathcal J}\phi_j(\Lambda_j),
\qquad \Xi_N(\varnothing):=1.
\]
Here, $V(\Lambda) := V_1(\Lambda) \cup V_2(\Lambda)$ contains both row and column vertices. Put
\[
\chi_{rs}:=\E[Z_{N,r}Z_{N,s}]
=\one_{\{H_r=H_s\}}\langle\phi_r,\phi_s\rangle_2.
\]

For $\mathcal J\subseteq[b]$, let $\mathfrak M(\mathcal J)$ be the set of \textit{matchings} of $\mathcal J$, that is, sets $\mathfrak m$ of pairwise disjoint pairs $\{r,s\}\subseteq \mathcal J$, including the empty matching, and let $V(\mathfrak m) =\bigcup_{\{r,s\}\in\mathfrak m}\{r,s\}$ be the set of matched indices. Call a tuple \emph{regular} if each distinct copy occurs once or twice and distinct copies are vertex-disjoint. Let $B_N(\mathcal J)$ be the part of $\prod_{j\in \mathcal J}Z_{N,j}$ coming from all other tuples:
\[
B_N(\mathcal J)=\left(\prod_{j\in \mathcal J}M_j^{-1/2}\right)
\sum_{\substack{\Lambda_j\text{ copy of }H_j,\ j\in \mathcal J\\
(\Lambda_j)_{j\in \mathcal J}\text{ not regular}}}
\prod_{j\in \mathcal J}\phi_j(\Lambda_j).
\]
When two copies coincide, their product gives the following average. For $H_r=H_s$, define
\[
\widehat\chi_{N,rs}:=M_r^{-1}\sum_{\Lambda\text{ copy of }H_r}\phi_r(\Lambda)\phi_s(\Lambda),
\]
and put $\widehat\chi_{N,rs}=0$ when $H_r\ne H_s$.

For $\mathcal J\subseteq[b]$ and real numbers $(x_j)_{j\in \mathcal J}$, the \emph{Wick polynomial} associated with the covariance matrix $(\chi_{rs})$ is
\[
\sum_{\mathfrak m\in\mathfrak M(\mathcal J)}(-1)^{|\mathfrak m|}
\prod_{\{r,s\}\in\mathfrak m}\chi_{rs}
\prod_{j\in \mathcal J\setminus V(\mathfrak m)}x_j.
\]
When $(x_j)_{j\in \mathcal J}$ is a centered Gaussian vector with $\E[x_rx_s]=\chi_{rs}$, its value is the \emph{Wick product} of the $x_j$ \citep[Ch.~3]{janson1997gaussian}.

\phantomsection\label{def:S-Lk-orders}
For $u_N>0$ and $k=1,2$, we write $Z_N=o_{L^k}(u_N)$ when $\|Z_N\|_k=o(u_N)$, and similarly for $O_{L^k}$.

\begin{lemma}\label{lem:wick-products}
For every $\mathcal J\subseteq[b]$, the following hold.
\begin{enumerate}[(i)]
\item $\E[B_N(\mathcal J)^2]=O(N^{-1})$, and the same bound holds if $B_N(\mathcal J)$ is restricted to any subset of its tuples.
\item $\widehat\chi_{N,rs}=\chi_{rs}+o_{L^2}(1)$ for all $r,s\in \mathcal J$.
\item The product of the $Z_{N,j}$, $j\in \mathcal J$, satisfies
\begin{equation}\label{eq:matching-triangular}
\prod_{j\in \mathcal J}Z_{N,j}
=
\sum_{\mathfrak m\in\mathfrak M(\mathcal J)}
\left(\prod_{\{r,s\}\in\mathfrak m}\widehat\chi_{N,rs}\right)
\Xi_N(\mathcal J\setminus V(\mathfrak m))+o_{L^2}(1).
\end{equation}
\item Up to $o_P(1)$, the sum over vertex-disjoint copies is the Wick polynomial evaluated at $(Z_{N,j})_{j\in \mathcal J}$:
\begin{equation}\label{eq:wick-disjoint}
\Xi_N(\mathcal J)=
\sum_{\mathfrak m\in\mathfrak M(\mathcal J)}(-1)^{|\mathfrak m|}
\prod_{\{r,s\}\in\mathfrak m}\chi_{rs}
\prod_{j\in \mathcal J\setminus V(\mathfrak m)}Z_{N,j}
+o_P(1).
\end{equation}
\end{enumerate}
\end{lemma}

\begin{proof}[Proof of Lemma~\ref{lem:wick-products}]
We start from the product $\prod_{j\in \mathcal J}Z_{N,j}$. Its expansion includes all tuples of copies, whereas $\Xi_N(\mathcal J)$ keeps only those whose copies are vertex-disjoint. First, we show which tuples can be discarded. Then, we evaluate the identical pairs that remain. We group the remaining tuples by matchings to establish a system which, finally, can be solved by induction.

\textbf{Proof of (i).} Put $d_{\mathrm{tot}}(\mathcal J)=\sum_{j\in \mathcal J}d(H_j)$. Expanding $\E[B_N(\mathcal J)^2]$ gives
\[
\E[B_N(\mathcal J)^2]
=\left(\prod_{j\in \mathcal J}M_j^{-1}\right)
\sum_{\substack{\Lambda_j,\Lambda'_j\text{ copies of }H_j,\ j\in \mathcal J\\
(\Lambda_j)_{j\in \mathcal J},\ (\Lambda'_j)_{j\in \mathcal J}\text{ not regular}}}
\E\left[\prod_{j\in \mathcal J}\phi_j(\Lambda_j)\phi_j(\Lambda'_j)\right].
\]
The two tuples in this sum are chosen separately. For each $j\in \mathcal J$, they contain two copies of $H_j$. Counting each vertex once for every copy containing it, the total count is $2\sum_{j\in \mathcal J}d(H_j)=2d_{\mathrm{tot}}(\mathcal J)$.

Let $v$ be the number of distinct row and column vertices in these two tuples. If an expected product is nonzero, Lemma~\ref{lem:config} says that every vertex belongs to at least two copies. If the copies cannot be grouped into identical pairs, the same lemma gives a vertex belonging to at least three copies. Consequently,
\[
2v\le 2d_{\mathrm{tot}}(\mathcal J)-1,
\]
so $v\le d_{\mathrm{tot}}(\mathcal J)-1$.

Now, suppose that the copies can be grouped into identical pairs. Keep one copy from each pair. The sum of the numbers of vertices in these representatives is $d_{\mathrm{tot}}(\mathcal J)$. Since the original tuples are not regular, the representatives cannot all be distinct and vertex-disjoint. Either two distinct representatives share a vertex, or the same copy represents more than one pair. Thus, at least one vertex has been counted more than once, and again $v\le d_{\mathrm{tot}}(\mathcal J)-1$.

Therefore, every nonzero expected product involves at most $d_{\mathrm{tot}}(\mathcal J)-1$ distinct row and column indices. There are only finitely many ways for the fixed graphs to share vertices, and $m,n\asymp N$. Hence, there are $O(N^{d_{\mathrm{tot}}(\mathcal J)-1})$ contributing label choices. Moreover, each expected product is bounded independently of $N$, since the $\phi_j$ are bounded. The denominator is $\prod_{j\in \mathcal J}M_j\asymp N^{d_{\mathrm{tot}}(\mathcal J)}$, so
\[
\E[B_N(\mathcal J)^2]=O(N^{-1}).
\]
The same bound holds if $B_N(\mathcal J)$ is restricted to any subset of its tuples: the count above bounds the number and the absolute expectations of all contributing products. We will use this observation in (iii).

\textbf{Proof of (ii).} The expectation of $\widehat\chi_{N,rs}$ is $\chi_{rs}$. In its variance, two vertex-disjoint copies give independent products. There are only $O(N^{2d(H_r)-1})$ pairs sharing a vertex, and the $\phi_j$ are bounded, so dividing by $M_r^2\asymp N^{2d(H_r)}$ gives $\var(\widehat\chi_{N,rs})=O(N^{-1})$. Thus,
\begin{equation}\label{eq:contraction}
\widehat\chi_{N,rs}=\chi_{rs}+o_{L^2}(1).
\end{equation}

\textbf{Proof of (iii).} Now, we return to $\prod_{j\in \mathcal J}Z_{N,j}$. By (i), we can discard its nonregular tuples with an error $o_{L^2}(1)$. Each remaining tuple determines a unique $\mathfrak m\in\mathfrak M(\mathcal J)$: a pair $\{r,s\}$ belongs to $\mathfrak m$ exactly when $\Lambda_r=\Lambda_s$. The copies $\Lambda_j$, $j\in \mathcal J\setminus V(\mathfrak m)$, each occur once, and all distinct copies in the tuple are vertex-disjoint.

Conversely, fix $\mathfrak m\in\mathfrak M(\mathcal J)$. If $H_r\ne H_s$ for some $\{r,s\}\in\mathfrak m$, then $\widehat\chi_{N,rs}=0$ and no regular tuple determines $\mathfrak m$. Otherwise, we expand $\bigl(\prod_{\{r,s\}\in\mathfrak m}\widehat\chi_{N,rs}\bigr)\Xi_N(\mathcal J\setminus V(\mathfrak m))$ by choosing one copy $\Lambda_{rs}$ for each pair and pairwise vertex-disjoint copies for the unmatched indices. Set $\Lambda_r=\Lambda_s=\Lambda_{rs}$ for each pair. Since $M_r=M_s$, the normalization is
\[
\left(\prod_{\{r,s\}\in\mathfrak m}M_r^{-1}\right)
\left(\prod_{j\in \mathcal J\setminus V(\mathfrak m)}M_j^{-1/2}\right)
=\prod_{j\in \mathcal J}M_j^{-1/2}.
\]
In a regular tuple, only identical copies can be matched, and every identical pair must be matched because the unmatched copies are vertex-disjoint. Thus, each regular tuple appears exactly once, for its unique matching.

The remaining terms of the expansion are tuples in which some $\Lambda_{rs}$ equals or shares a vertex with another $\Lambda_{r's'}$ or with an unmatched copy. These tuples are nonregular. The last observation in (i) therefore bounds the second moment of their normalized sum by $O(N^{-1})$. Since $\mathfrak M(\mathcal J)$ is finite and does not depend on $N$, the triangle inequality bounds the total over $\mathfrak m$ by $O(N^{-1/2})$ in $L^2$. This proves (iii).

\textbf{Proof of (iv).} For $\mathfrak m=\varnothing$, the summand is $\Xi_N(\mathcal J)$. Every nonempty $\mathfrak m$ contains at least one pair, so its summand involves $\Xi_N(\mathcal J')$ with $|\mathcal J'|\le|\mathcal J|-2$. Therefore, solving \eqref{eq:matching-triangular} for $\Xi_N(\mathcal J)$ expresses it through $\prod_{j\in \mathcal J}Z_{N,j}$ and sums $\Xi_N(\mathcal J')$ over smaller sets. We prove \eqref{eq:wick-disjoint} by induction on $|\mathcal J|$. The case $\mathcal J=\varnothing$ is trivial, since $\Xi_N(\varnothing)=1$.

For the induction step, suppose that \eqref{eq:wick-disjoint} holds for every proper subset of $\mathcal J$. For every nonempty $\mathfrak m\in\mathfrak M(\mathcal J)$, the set $\mathcal J\setminus V(\mathfrak m)$ is a proper subset of $\mathcal J$. By the induction hypothesis, $\Xi_N(\mathcal J\setminus V(\mathfrak m))$ is then a polynomial in the $Z_{N,j}$, up to $o_P(1)$. Each $Z_{N,j}$ has second moment $\chi_{jj}$, so $\Xi_N(\mathcal J\setminus V(\mathfrak m))=O_P(1)$. Since $\widehat\chi_{N,rs}\to\chi_{rs}$ in probability by (ii), we can replace every $\widehat\chi_{N,rs}$ in \eqref{eq:matching-triangular} by $\chi_{rs}$, which gives
\[
\Xi_N(\mathcal J)=\prod_{j\in \mathcal J}Z_{N,j}
-\sum_{\mathfrak m\in\mathfrak M(\mathcal J) \setminus \{\varnothing\}}
\Bigl(\prod_{\{r,s\}\in\mathfrak m}\chi_{rs}\Bigr)\Xi_N(\mathcal J\setminus V(\mathfrak m))
+o_P(1).
\]
Now, we use \eqref{eq:wick-disjoint} to replace $\Xi_N(\mathcal J\setminus V(\mathfrak m))$ for each $\mathfrak m \in \mathfrak M(\mathcal J) \setminus \{ \varnothing \}$. The error remains $o_P(1)$ because only finitely many $o_P(1)$ terms are added. Each resulting term comes from some $\mathfrak m\ne\varnothing$ and some $\mathfrak m'\in\mathfrak M(\mathcal J\setminus V(\mathfrak m))$. Since $\mathfrak m'$ avoids $V(\mathfrak m)$, $\mathfrak m''=\mathfrak m\cup\mathfrak m'$ belongs to $\mathfrak M(\mathcal J)$, and the term equals
\[
-(-1)^{|\mathfrak m'|}\Bigl(\prod_{\{r,s\}\in\mathfrak m''}\chi_{rs}\Bigr)\prod_{j\in \mathcal J\setminus V(\mathfrak m'')}Z_{N,j}.
\]
Conversely, a fixed $\mathfrak m''\in\mathfrak M(\mathcal J)$ arises once from each nonempty $\mathfrak m\subseteq\mathfrak m''$, with $\mathfrak m'=\mathfrak m''\setminus\mathfrak m$. For $\mathfrak m''\ne\varnothing$, its total coefficient is therefore
\[
-\sum_{\substack{\mathfrak m\subseteq\mathfrak m'' \\\mathfrak m \ne \varnothing}}(-1)^{|\mathfrak m''|-|\mathfrak m|}
=(-1)^{|\mathfrak m''|}-\sum_{\mathfrak m\subseteq\mathfrak m''}(-1)^{|\mathfrak m''|-|\mathfrak m|}
=(-1)^{|\mathfrak m''|},
\]
because the last sum is zero for nonempty $\mathfrak m''$. The term $\mathfrak m''=\varnothing$ comes only from $\prod_{j\in \mathcal J}Z_{N,j}$, with coefficient $1=(-1)^{0}$. This is \eqref{eq:wick-disjoint} for $\mathcal J$, which concludes the proof.
\end{proof}

\subsection{Proof of Lemma~\ref{lem:wick}}\label{app:wick-construction}

Before proving Lemma~\ref{lem:wick}, we establish useful intermediate results. Let $G$ be a class with connected components $H_1,\dots,H_b$, so that $G=H_1\sqcup\cdots\sqcup H_b$. Let $m_H$ be the number of components of $G$ belonging to class $H$, and put $\mu_G=\prod_Hm_H!$. Let $f\in\widehat\bigotimes_{j=1}^b\mathcal M_{H_j}^{\rm sym}$. We view $f$ both as a function of $b$ component copies and, since $\mathcal M_{H_j}^{\rm sym}=\mathfrak h_{H_j}\subseteq\mathfrak h$, as a tensor in $\mathfrak h^{\otimes b}$. Set
\[
\mathcal W_G(f)=I_b(\operatorname{Sym}_b f),
\qquad
S_N(f)=\sum_{\substack{\Lambda_j\text{ copy of }H_j,\ 1\le j\le b\\
\Lambda_1,\ldots,\Lambda_b\text{ pairwise vertex-disjoint}}}f(\Lambda_1,\ldots,\Lambda_b),
\]
and $a_{N,G}^2=\prod_{j=1}^bM_N(H_j)$.

\begin{lemma}\label{lem:wick-tensor}
The following hold.
\begin{enumerate}[(i)]
\item $\mathcal M_G^{\rm sym}$ is the subspace of $\widehat\bigotimes_j\mathfrak h_{H_j}$ invariant under permutations of components of the same class. For each connected class $H$, the space $\mathfrak h_H$ has an orthonormal basis $(\phi_{H,\ell})_\ell$ of bounded functions. Averages over these permutations of finite linear combinations of the tensors $\phi_{H_1,\ell_1}\otimes\cdots\otimes\phi_{H_b,\ell_b}$ are dense in $\mathcal M_G^{\rm sym}$. We call them finite symmetrized tensor expansions.
\item For $f\in\mathcal M_G^{\rm sym}$, $S_N(f)=\mu_G\sum_{\Lambda\text{ copy of }G}f(\Lambda)$ and $\var(S_N(f))=\mu_G^2M_N(G)\|f\|_2^2$. Moreover, $a_{N,G}^2=\mu_GM_N(G)\{1+O(N^{-1})\}$.
\item For $f,g\in\mathcal M_G^{\rm sym}$,
\begin{equation}\label{eq:wick-isometry}
\E[\mathcal W_G(f)\mathcal W_G(g)]
=\mu_G\langle f,g\rangle_{L^2(\mathcal X(G))}.
\end{equation}
\end{enumerate}
\end{lemma}

\begin{proof}[Proof of Lemma~\ref{lem:wick-tensor}]
\textbf{Proof of (i).} First, we express a graph projection through functions on its connected components. Distinct components involve independent families of AHK variables, so
\[
L^2(\mathcal X(G))=\widehat\bigotimes_{j=1}^bL^2(\mathcal X(H_j))
=\bigoplus_{F_j\subseteq H_j}\widehat\bigotimes_{j=1}^b\mathcal M_{F_j}.
\]
Removing the spaces indexed by proper subgraphs of $G$, that is, with $F_j\subsetneq H_j$ for some $j$, gives
\[
\mathcal M_G=\widehat\bigotimes_{j=1}^b\mathcal M_{H_j}.
\]
Invariance under graph automorphisms first requires invariance within each connected component. It also requires invariance under permutations of components of the same class. Thus, $\mathcal M_G^{\rm sym}$ is the subspace of $\widehat\bigotimes_j\mathfrak h_{H_j}$ with this latter invariance.

Next, we can approximate its elements by bounded functions without leaving these spaces. For each connected $H$, the AHK variables are i.i.d.\ uniform, so polynomials in $\mathcal X(H)$ are bounded and dense in $L^2(\mathcal X(H))$. The orthogonal projection onto $\mathcal M_H$ is a finite alternating sum of conditional expectations given subfamilies $\mathcal X(F)$, $F\subseteq H$, and symmetrization averages over $\mathrm{Aut}(H)$. These operations preserve boundedness and their composition is the orthogonal projection onto $\mathcal M_H^{\rm sym}$. Applying it to the polynomials and the Gram--Schmidt procedure gives a bounded orthonormal basis $(\phi_{H,\ell})$ of $\mathcal M_H^{\rm sym}$. Tensors of these basis elements form an orthonormal basis of $\widehat\bigotimes_j\mathcal M_{H_j}^{\rm sym}$, which contains $\mathcal M_G^{\rm sym}$. Averaging finite combinations of them over permutations of components of the same class therefore gives a dense subset of $\mathcal M_G^{\rm sym}$.

\textbf{Proof of (ii).} Now, we relate $S_N$ to sums over copies of $G$. Each copy of $G$ has exactly $\mu_G$ ordered decompositions into copies of $H_1,\dots,H_b$, obtained by permuting components of the same class. Thus, for $f\in \mathcal M_G^{\rm sym}$, $S_N(f)=\mu_G\sum_{\Lambda\text{ copy of }G}f(\Lambda)$ and $T_{N,G}(f)=\beta_G\mu_G^{-1}S_N(f)$, with $\beta_G$ as in Section~\ref{app:projection-notation}. Since distinct graph copies are orthogonal, we have exactly
\[
\var(S_N(f))=\mu_G^2M_N(G)\|f\|_2^2.
\]
The product $a_{N,G}^2$ counts all ordered tuples of component copies, of which $\mu_GM_N(G)$ are vertex-disjoint. The proportion of the others is $O(N^{-1})$, so
\[
a_{N,G}^2=\mu_GM_N(G)\{1+O(N^{-1})\},\qquad
\frac{\var(S_N(f))}{a_{N,G}^2}\longrightarrow\mu_G\|f\|_2^2.
\]

\textbf{Proof of (iii).} The same factor $\mu_G$ appears in the Gaussian limit. Let $\mathfrak S_b$ be the set of permutations of $[b]$, and let $f^\pi$ be $f$ with its tensor arguments permuted by $\pi$. The $\mu_G$ permutations exchanging only components of the same class leave $f$ unchanged. In contrast, every other permutation puts $f^\pi$ in a tensor space orthogonal to that of $f$. Since $\operatorname{Sym}_b$ is an orthogonal projection,
\[
\|\operatorname{Sym}_b f\|_{\mathfrak h^{\otimes b}}^2
=\langle f,\operatorname{Sym}_b f\rangle_{\mathfrak h^{\otimes b}} =\frac{1}{b!}\sum_{\pi\in\mathfrak S_b}\langle f,f^\pi\rangle_{\mathfrak h^{\otimes b}}=\frac{\mu_G}{b!}\|f\|_{L^2(\mathcal X(G))}^2.
\]
By polarization and~\eqref{eq:multiple-integral-isometry}, this gives (iii).
\end{proof}

Now, we prove Lemma~\ref{lem:wick}.

\begin{proof}[Proof of Lemma~\ref{lem:wick}]
The counting argument in Lemma~\ref{lem:wick-products} requires bounded functions, since it bounds expectations of products of many factors. However, in the present lemma, the projections are only square-integrable: we first prove the result for bounded approximations, before passing to the original projections using exact $L^2$ identities. We use the Gaussian process $X$ and the multiple Wiener integrals of Section~\ref{app:wick}.

\textbf{(i)} First, we prove the limit for a product of bounded functions on the connected components. For basis elements $\phi_j=\phi_{H_j,\ell_j}$ from Lemma~\ref{lem:wick-tensor}(i), we use the sums $Z_{N,j}$ and $\Xi_N(\mathcal J)$ from Lemma~\ref{lem:wick-products}, with $\mathcal J=[b]$.

We apply Lemma~\ref{lem:clt} to the distinct selected basis elements on each graph, and repeat a coordinate whenever its basis element is repeated. This gives joint convergence $Z_{N,j}\dto Z_j:=X_{H_j}(\phi_j)$, with $\E[Z_rZ_s]=\chi_{rs}$. By~\eqref{eq:wick-disjoint}, the continuous mapping theorem and Slutsky's lemma,
\[
\Xi_N([b])\dto I_b\!\left(\operatorname{Sym}_b(\phi_1\otimes\cdots\otimes\phi_b)\right),
\]
where we use the multiple-integral representation of the Wick product \citep[Thms.~3.4 and~7.26]{janson1997gaussian}.
With the notation of Lemma~\ref{lem:wick-tensor}, $\Xi_N([b])=S_N(\phi_1\otimes\cdots\otimes\phi_b)/a_{N,G}$. Finitely many choices of the $\phi_j$ involve only finitely many $Z_{N,j}$, which converge jointly. By linearity, $S_N(f)/a_{N,G}\dto\mathcal W_G(f)$ for every finite tensor expansion $f$. Since the same Gaussian process is used for all graphs, this convergence holds jointly over $G\in\Gamma_0$.

\textbf{(ii)} Next, we standardize sums over copies of $G$. Let $f\in \mathcal M_G^{\rm sym}$ be a finite symmetrized tensor expansion, as in Lemma~\ref{lem:wick-tensor}(i). By (i), $S_N(f)/a_{N,G}\dto\mathcal W_G(f)$. Since $T_{N,G}(f)$ is a fixed multiple of $S_N(f)$, Lemma~\ref{lem:wick-tensor}(ii) gives
\[
\frac{T_{N,G}(f)}{\var(T_{N,G}(f))^{1/2}}
=\frac{S_N(f)/a_{N,G}}{\var(S_N(f))^{1/2}/a_{N,G}}
\dto\frac{\mathcal W_G(f)}{\sqrt{\mu_G}\,\|f\|_2}.
\]
The map $\mathcal W_G$ is defined and continuous on all of $\mathcal M_G^{\rm sym}$, by \eqref{eq:wick-isometry}. For the projection $f_G$, we therefore set
\[
W_G=\frac{\mathcal W_G(f_G)}{\sqrt{\mu_G}\,\|f_G\|_2},
\]
and (iii) shows that $W_{N,G}\dto W_G$.

\textbf{(iii)} To this end, we remove the bounded approximation. Let $f_G^{(M)}\in \mathcal M_G^{\rm sym}$ be finite symmetrized tensor expansions with $f_G^{(M)}\to f_G$ in $L^2$, and set
\[
\widetilde W_{N,G}^{(M)}
=\frac{T_{N,G}(f_G^{(M)})}{\var(T_{N,G}(f_G))^{1/2}},
\qquad
\widetilde W_G^{(M)}
=\frac{\mathcal W_G(f_G^{(M)})}{\sqrt{\mu_G}\,\|f_G\|_2}.
\]
By Lemma~\ref{lem:wick-tensor}(ii), $\var(T_{N,G}(f))$ is proportional to $\|f\|_2^2$ on $\mathcal M_G^{\rm sym}$. Thus, (ii) rescaled by $\|f_G^{(M)}\|_2/\|f_G\|_2$ gives $\widetilde W_{N,G}^{(M)}\dto\widetilde W_G^{(M)}$ jointly over $G\in\Gamma_0$ for each fixed $M$. Moreover, Lemma~\ref{lem:wick-tensor}(ii) and \eqref{eq:wick-isometry} give, for every $N$,
\[
\E[|W_{N,G}-\widetilde W_{N,G}^{(M)}|^2]
=\frac{\|f_G-f_G^{(M)}\|_2^2}{\|f_G\|_2^2}
=\E[|W_G-\widetilde W_G^{(M)}|^2],
\]
which tends to zero as $M\to\infty$, uniformly in $N$. The joint convergence of $(W_{N,G})_{G\in\Gamma_0}$ to $(W_G)_{G\in\Gamma_0}$ follows \citep[Thm.~3.2]{billingsley1999convergence}.

Finally, $W_G$ is a multiple of $I_{b(G)}(\operatorname{Sym}_{b(G)}f_G)$, so it lies in the $b(G)$-th Wiener chaos, and \eqref{eq:wick-isometry} gives $\E[W_G^2]=1$. Let $G\ne G'$. If $b(G)\ne b(G')$, the two limits lie in orthogonal chaoses. Otherwise, $G$ and $G'$ have different multisets of connected components, so their symmetrized tensors are orthogonal in $\mathfrak h^{\otimes b}$. In both cases, $\E[W_GW_{G'}]=0$.
\end{proof}

\subsection{Proof of Theorem~\ref{thm:chaos}}\label{sec:S-thm-chaos}

\begin{proof}[Proof of Theorem~\ref{thm:chaos}]
First, we show that the projections of degree greater than $d$ are negligible after multiplication by $N^{d/2}$. By Lemma~\ref{lem:decomp}, $U_N-\theta_N=\sum_GT_G$, and $T_G=0$ when $d(G)<d$ by the definition of the principal degree. Since the $T_G$ are pairwise orthogonal, \eqref{eq:graph-variance-weight} gives
\[
N^d\,\E\Bigl[\Bigl\{U_N-\theta_N-\sum_{G:\,d(G)=d}T_G\Bigr\}^2\Bigr]
=\sum_{G:\,d(G)>d}N^dw_G\,\E[p_G^2]
=O(N^{-1}),
\]
because the law is fixed, $w_G\asymp N^{-d(G)}$ and there are finitely many graphs. Thus,
\[
N^{d/2}\left\{U_N-\theta_N-\sum_{G:\,d(G)=d}T_G\right\}
\longrightarrow0\qquad\text{in }L^2.
\]

Next, for each principal support graph, the variance formula \eqref{eq:graph-variance-weight} and $m/N\to\rho$ give $N^d\var(T_G)\to\tau_G^2$, with $\tau_G$ defined in \eqref{eq:fixed-law-coefficients}. Since $N^{d/2}T_G=\{N^d\var(T_G)\}^{1/2}W_{N,G}$ with $f_G=p_G$, Lemma~\ref{lem:wick} and Slutsky's lemma then give joint convergence of $N^{d/2}T_G$ to $\tau_GW_G$.
Summing over the finitely many principal support graphs and using \eqref{eq:fixed-law-coefficients}, we obtain
\[
N^{d/2}(U_N-\theta_N)\dto
\sum_{G:\,d(G)=d}\tau_GW_G=\sum_{b=1}^dL_b.
\]

Now, each $W_G$ belongs to the chaos of order $b(G)$. By Lemma~\ref{lem:wick}, the $W_G$ have unit variance and are pairwise orthogonal. Thus,
\[
\E[L_b^2]=\sum_{\substack{G\in\Gamma:\ d(G)=d,\ b(G)=b}}\tau_G^2.
\]
Every principal support graph has $\tau_G>0$. Consequently, $L_b$ is nonzero exactly when a principal support graph has $b$ connected components. If all principal support graphs are connected, the limit lies in the first chaos and is Gaussian. Otherwise, let $B\ge2$ be the largest $b$ with $L_b\ne0$. Each $L_b$ lies in the $b$-th Wiener chaos, so $L_b=I_b(f_b)$ for some $f_b\in\mathfrak h^{\odot b}$, with $f_B\ne0$. Lemma~\ref{lem:chaosnotgauss} therefore shows that $\sum_bL_b$ is not Gaussian.
\end{proof}

\subsection{Proof of Theorem~\ref{thm:sharp}}\label{app:boundary}\label{app:sharp}

Lemma~\ref{lem:wick} identifies limits for fixed projection functions, whereas the laws considered here may vary with $N$. Assumption (K) lets us pass to a subsequence on which the normalized projections have fixed $L^2$ limits. By orthogonality, we can replace the projections by these limits, changing $(U_N-\theta_N)/\sigma_N$ only by a vanishing amount. 

Consider a sequence of laws $P_k$ and sample sizes $N_k$. We make this dependence explicit by writing $T_{N_k,G}$, $p_{P_k,G}$, $s_{k,G}$, $\theta_{P_k}$ and $\sigma_{N_k,P_k}$ for $T_G$, $p_G$, $s_{N,G}$, $\theta_N$ and $\sigma_N$ computed at size $N_k$ under $P_k$.

\begin{lemma}\label{lem:S-replacement}
Let $\Pcal$ be a class of laws, and let $N_k\to\infty$ and $P_k\in\Pcal$, with $\sigma_{N_k,P_k}>0$, be such that $s_{k,G}\to s_G^*$ for every $G\in\Gamma$. Assume (K) on $\Pcal$ for the graphs $G$ with $s_G^*>0$. Then, along a subsequence, $p_{P_k,G}/\|p_{P_k,G}\|_2$ converges in $L^2(\mathcal X(K_{p,q}))$ to a unit-norm $f_G^*\in\mathcal M_G^{\rm sym}$ for every such graph. Moreover, with $\widetilde Z_{k,G}=M_{N_k}(G)^{-1/2}\sum_{\Lambda\text{ copy of }G}f_G^*(\Lambda)$,
\begin{equation}\label{eq:fixed-shape-replacement}
\frac{U_{N_k}-\theta_{P_k}}{\sigma_{N_k,P_k}}
-
\sum_{G:\,s_G^*>0}s_G^*\widetilde Z_{k,G}
\longrightarrow0
\qquad\text{in }L^2.
\end{equation}
\end{lemma}

\begin{proof}[Proof of Lemma~\ref{lem:S-replacement}]
Since $\Gamma$ is finite, (K) lets us pass to one subsequence along which, for every $G$ with $s_G^*>0$,
\[
\frac{p_{P_k,G}}{\|p_{P_k,G}\|_2}\longrightarrow f_G^*
\qquad\text{in }L^2(\mathcal X(K_{p,q})),
\]
for some unit-norm $f_G^*\in \mathcal M_G^{\rm sym}$. The limit remains in this space because $\mathcal M_G^{\rm sym}$ is closed in $L^2$.

Now, we show that replacing the normalized projections by these limits changes $(U_{N_k}-\theta_{P_k})/\sigma_{N_k,P_k}$ by $o_{L^2}(1)$. Dividing each graph contribution by its standard deviation cancels its sampling coefficient $\beta_G$. Indeed, writing $M_{N_k}(G)$ for the number of copies of $G$,
\[
Z_{k,G}=M_{N_k}(G)^{-1/2}\sum_{\Lambda\text{ copy of }G}
\frac{p_{P_k,G}(\Lambda)}{\|p_{P_k,G}\|_2}.
\]
Since distinct copies are orthogonal, we then have exactly
\[
\E_{P_k}\big[\big(Z_{k,G}-\widetilde Z_{k,G}\big)^2\big]
=
\left\|
\frac{p_{P_k,G}}{\|p_{P_k,G}\|_2}-f_G^*
\right\|_2^2\longrightarrow0.
\]
By the normalization above, $T_{N_k,G}=\var_{P_k}(T_{N_k,G})^{1/2}Z_{k,G}$ whenever $p_{P_k,G}\ne0$, and by Lemma~\ref{lem:decomp},
\[
\frac{U_{N_k}-\theta_{P_k}}{\sigma_{N_k,P_k}}=\sum_{G\in\Gamma}s_{k,G}Z_{k,G},
\]
where terms with $p_{P_k,G}=0$ are omitted. The terms with $s_G^*=0$ have $L^2$ norm $s_{k,G}\to0$. For $s_G^*>0$, we have $\|Z_{k,G}\|_2=1$, $s_{k,G}\to s_G^*$ and $Z_{k,G}-\widetilde Z_{k,G}\to0$ in $L^2$, so $s_{k,G}Z_{k,G}-s_G^*\widetilde Z_{k,G}\to0$ in $L^2$. Since $\Gamma$ is finite, \eqref{eq:fixed-shape-replacement} follows.
\end{proof}

\begin{proof}[Proof of Theorem~\ref{thm:sharp}]
\textbf{Proof of (a)$\Rightarrow$(b).} We assume (a) and show that failure of (b) produces a non-Gaussian subsequential limit, a contradiction.
\begin{enumerate}
\item[\textbf{(i)}] First, choose a sequence on which the disconnected share stays positive. Failure of (b) gives $\delta>0$, $N_k\to\infty$ and $P_k\in\Pcal$ such that
\[
\sum_{G\,\disc}s_{k,G}^2\ge\delta
\]
for every $k$. Since $\Gamma$ is finite and $\sum_Gs_{k,G}^2=1$, pass to a subsequence on which $s_{k,G}\to s_G^*$ for every $G$. Then,
\begin{equation}\label{eq:disc-share-limit}
\sum_{G\,\disc}(s_G^*)^2\ge\delta.
\end{equation}

\item[\textbf{(ii)}] By Lemma~\ref{lem:S-replacement}, along a further subsequence, the normalized projections of the classes with $s_G^*>0$ converge in $L^2$ to unit-norm limits $f_G^*\in\mathcal M_G^{\rm sym}$, and \eqref{eq:fixed-shape-replacement} holds.

\item[\textbf{(iii)}] The variables $\widetilde Z_{k,G}$ are sums of the fixed functions $f_G^*$ of the AHK variables, which are i.i.d.\ uniform under every $P_k$. Their joint law therefore does not depend on $P_k$, and $\widetilde Z_{k,G}=W_{N_k,G}$ in the notation of Lemma~\ref{lem:wick}, with $f_G=f_G^*$. By Lemma~\ref{lem:wick}, jointly over the finitely many classes with $s_G^*>0$,
\[
(\widetilde Z_{k,G})_G\dto(W_G)_G,
\]
where $W_G$ belongs to the $b(G)$-th homogeneous Wiener chaos of one common isonormal Gaussian process. \phantomsection\label{def:S-subsequential-chaoses}Hence, by \eqref{eq:fixed-shape-replacement},
\[
\frac{U_{N_k}-\theta_{P_k}}{\sigma_{N_k,P_k}}
\dto
F:=\sum_{b=1}^{p+q}L_b^*,
\qquad
 L_b^*:=\sum_{G:\,b(G)=b}s_G^*W_G.
\]
Since the $W_G$ are orthonormal by Lemma~\ref{lem:wick}, we have
\begin{equation}\label{eq:chaos-share}
\E[(L_b^*)^2]=\sum_{G:\,b(G)=b}(s_G^*)^2.
\end{equation}
Because a graph is disconnected exactly when $b(G)\ge2$, \eqref{eq:disc-share-limit} and \eqref{eq:chaos-share} imply
\[
\sum_{b\ge2}\E[(L_b^*)^2]\ge\delta.
\]
Thus, the largest $B$ with $L_B^*\ne0$ satisfies $B\ge2$.

\item[\textbf{(iv)}] By Lemma~\ref{lem:chaosnotgauss}, applied with this $B$, the limit $F$ in (iii) is not Gaussian. On the other hand, (a), applied to the sequence $(N_k,P_k)$, implies convergence of the same standardized statistics to $\Nor(0,1)$. Since weak limits are unique, this is a contradiction.
\end{enumerate}

\textbf{Proof of (b)$\Rightarrow$(a).} We assume (b) and show that a failure of (a) leads to a contradiction. Suppose (a) fails. Then, there are $\varepsilon>0$, $N_k\to\infty$ and $P_k\in\Pcal$ with a Kolmogorov distance at least $\varepsilon$ between the law of $(U_{N_k}-\theta_{P_k})/\sigma_{N_k,P_k}$ and $\Phi$.
\begin{itemize}
\item[\textbf{(i)}] First, pass to a subsequence on which $s_{k,G}\to s_G^*$ for every $G$. By (b), $s_G^*=0$ for every disconnected $G$, so $\sum_{G\,\conn}(s_G^*)^2=1$.
\item[\textbf{(ii)}] Thus, every class with $s_G^*>0$ is connected, and (K) holds for these classes. Lemma~\ref{lem:S-replacement} therefore gives, along a further subsequence, fixed unit-norm limits $f_G^*\in\mathcal M_G^{\rm sym}$ and \eqref{eq:fixed-shape-replacement}, with the sum restricted to connected $G$ with $s_G^*>0$.
\item[\textbf{(iii)}] By Lemma~\ref{lem:wick}, $(\widetilde Z_{k,G})_{G\,\conn,\,s_G^*>0}$ converges jointly to first-chaos elements $W_G$ of the Gaussian process in \eqref{eq:wick-gaussian-space}. These limits are jointly Gaussian, pairwise orthogonal and of unit variance. Thus, they are independent standard normal random variables, and $\sum_{G:\,s_G^*>0}s_G^*W_G\sim\Nor(0,1)$ by (i).
\item[\textbf{(iv)}] Step~(iii) and \eqref{eq:fixed-shape-replacement} show that the standardized statistics converge in law to $\Nor(0,1)$ along the subsequence. Since $\Phi$ is continuous, P\'olya's theorem gives convergence in Kolmogorov distance, contradicting the choice of $(N_k,P_k)$.
\end{itemize}
\end{proof}

\subsection{Joint limit of the statistic and variance estimator}\label{sec:S-beyond}\label{sec:S-joint-notation}

We use the common Gaussian space of Section~\ref{app:wick}, taking $\Gamma_0=\Gamma$. We consider $F=\sum_{b=1}^BI_b(f_b)$, with $f_b\in\mathfrak h^{\odot b}$, so that $\E[F^2]=\sum_bb!\|f_b\|^2$, and we write $D(F)=\sum_{b\ge2}b!\|f_b\|^2$ for the variance of its components of order at least two. Define
\begin{equation}\label{eq:joint-limit-variances}
\mathcal V(F)=\sum_{b,b'=1}^B\sum_{t=1}^{b\wedge b'}
t!\binom bt\binom{b'}t I_{b+b'-2t}(f_b\widetilde\otimes_t f_{b'}),
\end{equation}
where $\widetilde\otimes$ denotes the symmetrized contraction and $I_0$ is the identity on constants.

\begin{theorem}\label{thm:S-joint}
Let $\Pcal$ be a class of laws with $\sigma_{N,P}>0$ for every $P\in\Pcal$. Assume (M4) and (K) on $\Pcal$, and let $P_N\in\Pcal$. Every subsequence has a further subsequence along which $D_N\to D_\infty$ and
\[
\left(\frac{U_N-\theta_N}{\sigma_N},\frac{\Vu}{\sigma_N^2}\right)\dto\bigl(F,\mathcal V(F)\bigr),
\]
where $F=\sum_{b=1}^{p+q}I_b(f_b)$ with $\E[F^2]=1$ and $D(F)=D_\infty$, for an isonormal Gaussian process on a space $\mathfrak h$ as above. 
\end{theorem}

To prove Theorem~\ref{thm:S-joint}, we use Lemma~\ref{lem:auto-decomp} to write
\[
\frac{\Vu}{\sigma_N^2}=\frac{J}{\sigma_N^2}+o_{L^2}(1).
\]
It therefore suffices to find the joint limit of $(U_N-\theta_N)/\sigma_N$ and $J/\sigma_N^2$. The following lemma bounds the change in $J$ when we approximate the graph projections.

Write $|A|=|A_1|+|A_2|$. For a family $a=(a_G)_G$ of ordinary graph projections $a_G\in \mathcal M_G^{\rm sym}$, with the normalization of Lemma~\ref{lem:decomp}, write $X_A(a)$ for the sum of the projections in $a$ on the vertex set $A$, with $c_A$ as in \eqref{eq:grouped-projections}. Write $\sigma(a)^2=\sum_Ac_A^2\E[X_A(a)^2]$ for the variance of $\sum_Ac_AX_A(a)$. For two families, define
\[
J(a,b)=\sum_{A\sim A'}c_Ac_{A'}X_A(a)X_{A'}(b).
\]
For the ordinary projections $q=(q_G)_G$ of the kernel, $J(q,q)=J$.

\begin{lemma}\label{lem:S-bilin}
For all families $a,b$ as above, $\E[|J(a,b)|]\le2^{p+q}\sigma(a)\sigma(b)$.
\end{lemma}

\begin{proof}[Proof of Lemma~\ref{lem:S-bilin}]
For a nonempty pair of row and column sets $T$, put $S_T(a)=\sum_{A\supseteq T}c_AX_A(a)$, and define $S_T(b)$ similarly. Inclusion--exclusion gives
\[
\one\{A\sim A'\}=\sum_{\varnothing\ne T\subseteq A\cap A'}(-1)^{|T|+1},
\qquad J(a,b)=\sum_{T\ne\varnothing}(-1)^{|T|+1}S_T(a)S_T(b).
\]
Since the grouped projections on distinct vertex sets are orthogonal,
\[
\sum_{T\ne\varnothing}\E[S_T(a)^2]
=\sum_A(2^{|A|}-1)c_A^2\E[X_A(a)^2]
\le2^{p+q}\sigma(a)^2,
\]
and the same bound holds for $b$. The triangle inequality and the Cauchy--Schwarz inequality give
\[
\begin{aligned}
\E[|J(a,b)|]&\le\sum_{T\ne\varnothing}\E[|S_T(a)S_T(b)|]\\
&\le\left(\sum_{T\ne\varnothing}\E[S_T(a)^2]\sum_{T\ne\varnothing}\E[S_T(b)^2]\right)^{1/2}
\le2^{p+q}\sigma(a)\sigma(b).
\end{aligned}
\]
\end{proof}

The following lemma gives the joint limit for the finite tensor approximations used in the proof of Theorem~\ref{thm:S-joint}.

\begin{lemma}\label{lem:S-joint-finite}
Let $g=(g_G)_{G\in\Gamma}$, where each $g_G\in\mathcal M_G^{\rm sym}$ is a fixed finite symmetrized tensor expansion of bounded functions on the connected components, as in Lemma~\ref{lem:wick-tensor}(i). Define
\[
q_{N,G}^g=\frac{g_G}{\sqrt{w_G}(p)_{r(G)}(q)_{c(G)}},
\qquad q_N^g=(q_{N,G}^g)_{G\in\Gamma},
\]
and
\[
Z_N(g)=\sum_Ac_AX_A(q_N^g),
\qquad J_N(g)=J(q_N^g,q_N^g).
\]
Then
\[
(Z_N(g),J_N(g))\dto(F_g,\mathcal V(F_g)),
\qquad
F_g=\sum_{G\in\Gamma}\frac{\mathcal W_G(g_G)}{\sqrt{\mu_G}},
\]
with $\mu_G$ and $\mathcal W_G$ defined before Lemma~\ref{lem:wick-tensor}.
\end{lemma}

\begin{proof}[Proof of Lemma~\ref{lem:S-joint-finite}]
First, we express $Z_N(g)$ through the sums over vertex-disjoint connected copies in Lemma~\ref{lem:wick-products}. By \eqref{eq:copy-normalization} and \eqref{eq:grouped-projections},
\[
Z_N(g)=\sum_{G\in\Gamma}\frac{1}{\sqrt{M_N(G)}}
\sum_{\Lambda\text{ copy of }G}g_G(\Lambda).
\]
For each $G$, Lemma~\ref{lem:wick-tensor}(ii) gives
\[
\frac{\sum_{\Lambda\text{ copy of }G}g_G(\Lambda)}{\sqrt{M_N(G)}}
=\frac{a_{N,G}}{\mu_G\sqrt{M_N(G)}}\frac{S_N(g_G)}{a_{N,G}}
=\frac{1+O(N^{-1})}{\sqrt{\mu_G}}\frac{S_N(g_G)}{a_{N,G}}.
\]
Replacing these coefficients by their limits changes $Z_N(g)$ by $o_{L^2}(1)$ and, by Lemma~\ref{lem:S-bilin}, $J_N(g)$ by $o_{L^1}(1)$.

Let $(\phi_k)$ be the finite orthonormal family of bounded functions used in the expansions $g_G$, with $\phi_k\in\mathcal M_{H_k}^{\rm sym}$ for a connected class $H_k$. For $\mathbf k=(k_1,\dots,k_b)$, write $\Xi_N(\mathbf k)$ for the normalized sum over ordered, pairwise vertex-disjoint copies from Lemma~\ref{lem:wick-products}. We have
\[
\begin{aligned}
Z_N(g)&=\sum_{b=1}^{p+q}\sum_{\mathbf k}\beta_b(\mathbf k)\Xi_N(\mathbf k)+o_{L^2}(1),\\
F_g&=\sum_{b=1}^{p+q}I_b(f_b),
\qquad f_b=\sum_{\mathbf k}\beta_b(\mathbf k)\phi_{k_1}\otimes\cdots\otimes\phi_{k_b},
\end{aligned}
\]
with symmetric coefficients $\beta_b$ independent of $N$. Lemma~\ref{lem:clt} and \eqref{eq:wick-disjoint} give $Z_N(g)\dto F_g$ on the common Gaussian space.

Next, we calculate the joint limit with $J_N(g)$. A pair of summands in $J_N(g)$ is indexed by tuples $(\mathbf\Lambda,\mathbf\Lambda')$ that share a row or a column. Each tuple consists of pairwise vertex-disjoint copies. If the concatenated tuple is regular in the sense of Section~\ref{app:wick-products}, every shared copy joins opposite tuples, and the shared vertices are those of the matched copies. In contrast, the nonregular pairs form a subset of the tuples defining $B_N$ for the concatenated indices. By Lemma~\ref{lem:wick-products}(i), their entire normalized sum is $O_{L^2}(N^{-1/2})$.

For $\mathbf k=(k_1,\dots,k_b)$ and $\mathbf k'=(k'_1,\dots,k'_{b'})$, let $\mathfrak M(\mathbf k,\mathbf k')$ be the set of nonempty matchings between the two position sets with $H_{k_i}=H_{k'_j}$ for every matched pair $(i,j)$. Write $(\mathbf k,\mathbf k')\setminus V(\mathfrak m)$ for the list of unmatched indices in their original order. Grouping the regular pairs as in \eqref{eq:matching-triangular} gives
\[
\begin{aligned}
J_N(g)
={}&\sum_{b,b'=1}^{p+q}\sum_{\mathbf k,\mathbf k'}
\beta_b(\mathbf k)\beta_{b'}(\mathbf k')
\sum_{\mathfrak m\in\mathfrak M(\mathbf k,\mathbf k')}\\
&\quad\times\Bigl(\prod_{(i,j)\in\mathfrak m}\widehat\chi_{N,k_i k'_j}\Bigr)
\Xi_N\bigl((\mathbf k,\mathbf k')\setminus V(\mathfrak m)\bigr)
+o_{L^1}(1),
\end{aligned}
\]
where $\widehat\chi_{N,kk'}$ is the average defined before Lemma~\ref{lem:wick-products}, with $\phi_k,\phi_{k'}$ in place of $\phi_r,\phi_s$. Expanding each product also includes terms in which a matched copy equals or meets another copy. These are nonregular, so their total is $o_{L^2}(1)$ by the same subset bound. There are finitely many matchings. By Lemma~\ref{lem:wick-products}(ii), $\widehat\chi_{N,kk'}\to\langle\phi_k,\phi_{k'}\rangle$ in probability. The $o_P(1)$ expansion in \eqref{eq:wick-disjoint} and Slutsky's lemma therefore give joint convergence with $Z_N(g)$, with each $\Xi_N$ replaced by the Wick product of its remaining coordinates.

Finally, for fixed $b,b'$ and $1\le t\le b\wedge b'$, there are $t!\binom bt\binom{b'}t$ ways to choose $t$ positions on each side and match them. Since $\beta_b$ and $\beta_{b'}$ are symmetric, their contributions are equal after summing over the coordinates. The matched coordinates contract $f_b$ and $f_{b'}$ over $\mathfrak h^{\otimes t}$, while the Wick product of the unmatched coordinates gives the multiple integral \citep[Thm.~7.26]{janson1997gaussian}. Thus, their contribution is
\[
t!\binom bt\binom{b'}t
I_{b+b'-2t}\bigl(f_b\widetilde\otimes_tf_{b'}\bigr).
\]
Summing over $b,b',t$ gives $\mathcal V(F_g)$ by \eqref{eq:joint-limit-variances}, proving the claimed joint convergence.
\end{proof}

\begin{proof}[Proof of Theorem~\ref{thm:S-joint}]
There are two approximations in this proof. First, the laws may vary with $N$, so we use (K) to replace their normalized projections by fixed $L^2$ limits. Then, we approximate those fixed functions by finite sums of products of bounded functions, for which Lemma~\ref{lem:S-joint-finite} gives the joint limit. Finally, we verify the conditions of Theorem~3.2 of \citet{billingsley1999convergence} and apply it to obtain the claimed convergence. Throughout, $C$ depends only on $(p,q,C_{\mathrm{M4}})$. 

\textbf{(i)} First, we pass to a subsequence along which all coefficients $s_{N,G}$ converge, and write $s_G^*=\lim_Ns_{N,G}$. For the classes with $s_G^*>0$, Lemma~\ref{lem:S-replacement} gives, along a further subsequence, normalized projections converging in $L^2$ to limit functions $f_G^*$. Recall that the ordinary projection is $q_G=p_{N,G}/\{(p)_r(q)_c\}$ for $r=r(G)$ and $c=c(G)$. Define
\[
\widetilde q_{N,G}
=\frac{\sigma_Ns_G^*}{\sqrt{w_G}(p)_r(q)_c}\,f_G^*
\quad\text{if }s_G^*>0,
\qquad \widetilde q_{N,G}=0\quad\text{otherwise}.
\]
Since projections on distinct copies are orthogonal, the squared relative approximation error $\sigma(q-\widetilde q)^2/\sigma_N^2$ is the sum over $G$ of the following terms. For $p_{N,G}\ne0$, the term of $G$ is
\[
\frac{w_G(p)_r^2(q)_c^2}{\sigma_N^2}
\|q_G-\widetilde q_{N,G}\|_2^2
=\left\|s_{N,G}\frac{p_{N,G}}{\|p_{N,G}\|_2}-s_G^*f_G^*\right\|_2^2,
\]
where the second term on the right is zero when $s_G^*=0$. The terms with $s_G^*>0$ tend to zero, because $s_{N,G}\to s_G^*$ and $p_{N,G}/\|p_{N,G}\|_2\to f_G^*$ in $L^2$ along this subsequence. All other terms are $s_{N,G}^2\to0$. Summing over the finite graph family gives $\sigma(q-\widetilde q)/\sigma_N\to0$ and $\sigma(\widetilde q)/\sigma_N\to1$.

Since $J$ is bilinear, Lemma~\ref{lem:S-bilin} now gives
\[
\frac{\E[|J(q,q)-J(\widetilde q,\widetilde q)|]}{\sigma_N^2}
\le C\frac{\sigma(q-\widetilde q)}{\sigma_N}
       \frac{\sigma(q)+\sigma(\widetilde q)}{\sigma_N}\longrightarrow0.
\]
Moreover, by \eqref{eq:fixed-shape-replacement}, the same replacement changes $(U_N-\theta_N)/\sigma_N$ by a quantity tending to zero in $L^2$.

\textbf{(ii)} To apply Lemma~\ref{lem:S-joint-finite}, we approximate the fixed functions $f_G^*$ by finite sums of products of bounded functions. For each $G$ with $s_G^*>0$, Lemma~\ref{lem:wick-tensor}(i) gives finite symmetrized tensor expansions $g_G^{(M)}$ such that $\|g_G^{(M)}-f_G^*\|_2\to0$ as $M\to\infty$. Let $q^{(M)}$ be the family obtained by replacing $f_G^*$ by $g_G^{(M)}$ in $\widetilde q$, and define
\[
Z_{N,M}=\frac{1}{\sigma_N}\sum_Ac_AX_A(q^{(M)}),
\qquad
J_{N,M}=\frac{J(q^{(M)},q^{(M)})}{\sigma_N^2}.
\]
By the variance identity in Lemma~\ref{lem:wick-tensor}(ii), Lemma~\ref{lem:S-bilin} and the first replacement above,
\begin{equation}\label{eq:S-joint-approximation}
\lim_{M\to\infty}\limsup_{N\to\infty}
\left\{
\left\|\frac{U_N-\theta_N}{\sigma_N}-Z_{N,M}\right\|_2
+\E\left[\left|\frac{J}{\sigma_N^2}-J_{N,M}\right|\right]
\right\}=0.
\end{equation}

For fixed $M$, we apply Lemma~\ref{lem:S-joint-finite} with $g_G=s_G^*g_G^{(M)}$ when $s_G^*>0$ and $g_G=0$ otherwise. Since $q_N^g=q^{(M)}/\sigma_N$, the lemma gives
\begin{equation}\label{eq:S-joint-finite-limit}
\begin{aligned}
(Z_{N,M},J_{N,M})&\dto(F_M,\mathcal V(F_M))
\qquad\text{as }N\to\infty,\\
F_M&=\sum_{G:\,s_G^*>0}\frac{s_G^*}{\sqrt{\mu_G}}\mathcal W_G(g_G^{(M)}).
\end{aligned}
\end{equation}

\textbf{(iii)} Now, we let $M\to\infty$ and put $F=\sum_{G:\,s_G^*>0}s_G^*W_G=\sum_bI_b(f_b)$, where $W_G$ is the limit associated with $f_G^*$ in Lemma~\ref{lem:wick}. Since the limits for distinct graphs are orthogonal, \eqref{eq:wick-isometry} gives
\[
\|F_M-F\|_2^2
=\sum_{G:\,s_G^*>0}(s_G^*)^2\|g_G^{(M)}-f_G^*\|_2^2
\longrightarrow0.
\]
In addition, expanding the differences of contractions in \eqref{eq:joint-limit-variances} and applying \eqref{eq:multiple-integral-isometry} and the Cauchy--Schwarz inequality, we obtain
\[
\E[|\mathcal V(F_M)-\mathcal V(F)|]
\le C\|F_M-F\|_2\bigl(\|F_M\|_2+\|F\|_2\bigr)
\longrightarrow0.
\]
Hence, 
\begin{equation}\label{eq:S-joint-approximation-limit}
(F_M,\mathcal V(F_M))
\xrightarrow[M\to\infty]{d}(F,\mathcal V(F)).
\end{equation}

\textbf{(iv)} Therefore, we can apply Theorem~3.2 of \citet{billingsley1999convergence} to \eqref{eq:S-joint-approximation},~\eqref{eq:S-joint-finite-limit} and~\eqref{eq:S-joint-approximation-limit}, and then replace $J$ by $\Vu$ with Slutsky's lemma. This gives
\[
\left(\frac{U_N-\theta_N}{\sigma_N},\frac{\Vu}{\sigma_N^2}\right)
\dto\bigl(F,\mathcal V(F)\bigr).
\]

Finally, by Lemma~\ref{lem:wick}, the $W_G$ are orthogonal, have unit variance, and lie in the chaoses of orders $b(G)$. Thus,
\[
\E[F^2]=\sum_G(s_G^*)^2=1,
\qquad
D(F)=\sum_{G\,\disc}(s_G^*)^2=D_\infty.
\]

\end{proof}

\subsection{Proof of Corollary~\ref{cor:S-consistency}}\label{sec:S-cor-consistency}

\begin{proof}[Proof of Corollary~\ref{cor:S-consistency}]
If $D_N\to0$ uniformly, bound~\eqref{eq:variance-D-bound} and Chebyshev's inequality give $\Vu/\sigma_N^2\to1$ uniformly in probability. To prove necessity, suppose that $D_N\ge\delta>0$ along a sequence. By Theorem~\ref{thm:S-joint}, along a subsequence, $\Vu/\sigma_N^2\dto\mathcal V(F)$ with $D(F)>0$. It suffices to show that this limit is not constant. 

We look at its highest possible chaos order, since no term of a different order can cancel it. Let $b^*$ be the largest order with $f_{b^*}\ne0$, so that $b^*\ge2$. In the expansion of $\mathcal V(F)$ in~\eqref{eq:joint-limit-variances}, the only term in the chaos of order $2b^*-2$ is $(b^*)^2I_{2b^*-2}(f_{b^*}\widetilde\otimes_1f_{b^*})$, from $b=b'=b^*$ and $t=1$.

By polarization, the tensors $h^{\otimes b^*}$, $h\in\mathfrak h$, span a dense subspace of $\mathfrak h^{\odot b^*}$. Since $f_{b^*}\ne0$, choose $h$ with $\langle f_{b^*},h^{\otimes b^*}\rangle\ne0$. By the definition of contraction and the Cauchy--Schwarz inequality,
\[
\begin{aligned}
\langle f_{b^*}\widetilde\otimes_1f_{b^*},h^{\otimes(2b^*-2)}\rangle
&=\|f_{b^*}(h,\dots,h,\cdot)\|^2\\
&\ge\frac{|\langle f_{b^*},h^{\otimes b^*}\rangle|^2}{\|h\|^2}>0.
\end{aligned}
\]
Thus, the chaos term is nonzero and $\mathcal V(F)$ is not constant. Along the selected subsequence, $\Vu/\sigma_N^2$ therefore cannot converge in probability to $1$, which rules out uniform ratio-consistency.
\end{proof}

\section{Conservative procedures}\label{sec:S-part-procedures}

In this section, we prove Proposition~\ref{prop:block-safe}, Proposition~\ref{prop:ident}, Lemma~\ref{lem:S-filler-span}, Proposition~\ref{prop:S-block-hybrid} and Theorem~\ref{thm:S-check} (Sec.~\ref{app:protected}, \ref{sec:S-ident} and~\ref{sec:S-procedure}), and we give an alternative diagnostic for the block procedure (Sec.~\ref{sec:S-chaos-procedure}). We also give details on the block average test (Sec.~\ref{sec:S-block-average}).

\subsection{The block test}\label{app:protected}\label{sec:S-block-cost}

We prove Proposition~\ref{prop:block-safe}, then explain why the block test has an increased variance.

\begin{proof}[Proof of Proposition~\ref{prop:block-safe}]
Recall the block construction of Section~\ref{sec:S-block-hybrid}: $m_b=\lfloor m/K_N\rfloor$ and $n_b=\lfloor n/K_N\rfloor$ are the block dimensions, $S_{K_N}^2$ is the sample variance in \eqref{eq:block-statistics}, and $\tau_N^2=\var(U_{N,1})$ is the variance of one block statistic.

\textbf{(i)} First, we show that $\tau_N^2/K_N$ is at least the variance $\sigma_N^2$ of the complete statistic. The graph projections are unchanged on a block. Write $\varsigma_G^2(m,n)=\var(T_G)$ for the variance of the graph contribution at these dimensions. For each class $G$ with $r$ rows and $c$ columns, the variance formula~\eqref{eq:graph-variance-weight} gives
\begin{equation}\label{eq:S-block-graph-variance}
\varsigma_G^2(m_b,n_b)=\frac{(m)_r(n)_c}{(m_b)_r(n_b)_c}\,\varsigma_G^2(m,n).
\end{equation}
Eventually, $m_b\ge2p$ and $n_b\ge2q$. For $0\le j<r$, we have $m-j-K_N(m_b-j)=m-K_Nm_b+(K_N-1)j\ge0$, and similarly for the columns. Thus,
\[
\frac{(m)_r(n)_c}{(m_b)_r(n_b)_c}
=\prod_{j=0}^{r-1}\frac{m-j}{m_b-j}
 \prod_{j=0}^{c-1}\frac{n-j}{n_b-j}
\ge K_N^{r+c}.
\]
Since every nonempty graph has $r+c\ge1$, summing gives
\begin{equation}\label{eq:S-block-upper}
\tau_N^2=\sum_G\varsigma_G^2(m_b,n_b)
\ge K_N\sum_G\varsigma_G^2(m,n)
=K_N\sigma_N^2>0.
\end{equation}

\textbf{(ii)} To estimate $\tau_N^2$ by $S_{K_N}^2$, we need a uniform fourth-moment bound for the block statistics. For all sufficiently large $N$, each block statistic is the complete statistic, with the same kernel and law, at dimensions $m_b\ge2p$ and $n_b\ge2q$. Lemma~\ref{lem:S-four} at these dimensions gives, with a constant $C_*$ depending only on $(p,q,C_{\mathrm{M4}})$,
\begin{equation}\label{eq:S-block-moment}
\E[(U_{N,1}-\theta_N)^4]\le C_*\tau_N^4.
\end{equation}

\textbf{(iii)} Now, we bound the rejection probability using $S_{K_N}^2$. The $K_N$ block statistics use disjoint sets of rows and columns. They are therefore independent by dissociation and identically distributed by exchangeability. The exact i.i.d. sample-variance identities and \eqref{eq:S-block-moment} give
\begin{equation}\label{eq:S-sample-variance}
\begin{aligned}
\E[S_{K_N}^2]&=\tau_N^2,\\
\E[(S_{K_N}^2/\tau_N^2-1)^2]
&=\frac1{K_N}\left\{\frac{\E[(U_{N,1}-\theta_N)^4]}{\tau_N^4}
-\frac{K_N-3}{K_N-1}\right\}
\le\frac{C_*+1}{K_N}.
\end{aligned}
\end{equation}
For $0<\delta<1$, on the event $S_{K_N}^2\ge(1-\delta)\tau_N^2$, the definition of $B_N$ and \eqref{eq:S-block-upper} give
\[
B_N\ge\frac{S_{K_N}^2}{K_N}
\ge\frac{(1-\delta)\tau_N^2}{K_N}
\ge(1-\delta)\sigma_N^2.
\]
Under any $P\in\Pcal_0$, we have $\theta_N=\theta_0$. Splitting according to this event and applying Chebyshev's inequality to both probabilities gives
\begin{equation}\label{eq:S-block-safe-bound}
\begin{aligned}
&\P_P\{|U_N-\theta_0|\ge\sqrt{B_N/\alpha},\ B_N>0\}\\
&\quad\le\P_P\{S_{K_N}^2<(1-\delta)\tau_N^2\} +\P_P\{|U_N-\theta_0|\ge\sigma_N\sqrt{(1-\delta)/\alpha}\}\\
&\quad\le\frac{C_*+1}{K_N\delta^2}+\frac{\alpha}{1-\delta}.
\end{aligned}
\end{equation}
This bound is uniform over $\Pcal_0$. Taking $\delta=\delta_N\to0$ with $K_N\delta_N^2\to\infty$ proves the two-sided assertion.

\textbf{(iv)} For the upper one-sided test, the second probability in \eqref{eq:S-block-safe-bound} is replaced by
\[
\P_P\{U_N-\theta_0\ge\sigma_N\sqrt{(1-\alpha)(1-\delta)/\alpha}\}
\le\frac{\alpha}{1-(1-\alpha)\delta},
\]
by Cantelli's inequality. This replaces the last term in \eqref{eq:S-block-safe-bound} and gives the one-sided assertion by the same choice of $\delta_N$.
\end{proof}

Recall the variance shares $s_{N,G}^2$ in \eqref{eq:variance-proportions}. Summing \eqref{eq:S-block-graph-variance} over $G$ and dividing by $K_N\sigma_N^2$ gives
\begin{equation}\label{eq:S-block-inflation}
 \frac{\tau_N^2}{K_N\sigma_N^2}
 =\sum_Gs_{N,G}^2\frac{(m)_{r(G)}(n)_{c(G)}}
 {K_N(m_b)_{r(G)}(n_b)_{c(G)}}
 =\{1+O(K_N/L_N)\}\sum_Gs_{N,G}^2K_N^{d(G)-1}.
\end{equation}
The multiplicative error is uniform over the finite graph family. At a fixed law with principal degree $d$, the ratio $\tau_N^2/(K_N\sigma_N^2)$ is asymptotic to $K_N^{d-1}$. Indeed, the full-sample variance shares of higher-degree graphs are $O(L_N^{d-d(G)})$, and their contributions after division by $K_N^{d-1}$ are $O\{(K_N/L_N)^{d(G)-d}\}$.

\subsection{A diagnostic}\label{sec:S-ident}

First, we show that finitely many expected products of kernel values determine a quantity comparable to $\Upsilon_N$. Estimating these products gives the diagnostic for $\Upsilon_N\to0$.

\paragraph{Expected products and contractions.}
Use the index $\tau$, disjoint vertex sets $A_i,B_j$, and products $\mathcal D(A,B)$ and $\Pi_\tau$ defined in \eqref{eq:diagnostic-differences} and~\eqref{eq:diagnostic-product}. We choose all rows and columns used in the four factors, including the replacement and additional vertices. We use the ordinary graph projections $q_\Lambda$ of Section~\ref{app:projection-notation}, so that $H_W-\theta_N=\sum_{\varnothing\ne \Lambda\subseteq W}q_\Lambda$. For fixed vertex sets $A,B$, let $g_\tau$ be the sum of the projections $q_\Lambda$ whose graphs have vertex set exactly $A\cup B$ and no edge between $A$ and $B$. Viewed as an integral kernel, it maps $L^2(x_B)$ to $L^2(x_A)$. Here, $x_A$ contains the row and column variables in $A$, together with cell variables whose two endpoints are in this vertex set. Set
\[
Z_\tau=\|g_\tau\|_{S_4}^4=\|g_\tau\otimes_Bg_\tau\|_2^2,
\]
with $c_\tau$ as in \eqref{eq:diagnostic-statistic}.
Each such choice uses $4p+2r$ rows and $4q+2c$ columns.

\begin{proof}[Proof of Proposition~\ref{prop:ident}]
\textbf{(i)} First, we prove $\E[\Pi_\tau]=Z_\tau$. By orthogonality, conditioning $H_{W_T}$ on $x_A,x_B$ removes the ordinary projections involving a replacement vertex, an additional vertex, or an edge between $A$ and $B$. By \eqref{eq:diagnostic-differences},
\[
\E[\mathcal D(A,B)\mid x_A,x_B]
=\sum_{T\subseteq A\cup B}(-1)^{|T|}
  \left(\theta_N+\sum_{\Lambda:\,V(\Lambda)\subseteq(A\cup B)\setminus T}q_\Lambda\right),
\]
where $\Lambda$ ranges over nonempty graphs with vertices in $A\cup B$ and no edge between $A$ and $B$. Since $A\cup B\ne\varnothing$, the constant term cancels:
\[
\theta_N\sum_{T\subseteq A\cup B}(-1)^{|T|}=0.
\]
If $\Lambda$ omits a vertex $u\in A\cup B$, the terms indexed by $T$ and $T\cup\{u\}$ cancel in its coefficient. Otherwise, only $T=\varnothing$ contributes. Thus,
\[
\begin{aligned}
\E[\mathcal D(A,B)\mid x_A,x_B]
&=\sum_\Lambda q_\Lambda
  \sum_{T\subseteq(A\cup B)\setminus V(\Lambda)}(-1)^{|T|}\\
&=\sum_{\Lambda:\,V(\Lambda)=A\cup B}q_\Lambda
=g_\tau(x_A,x_B).
\end{aligned}
\]

Now, all other variables in each factor of $\Pi_\tau$ occur only in that factor and are independent of $(x_{A_1},x_{A_2},x_{B_1},x_{B_2})$. These are the cell variables between $A_i$ and $B_j$, the separate replacement and additional vertex variables, and their cell variables. Therefore,
\[
\E[\Pi_\tau\mid x_{A_1},x_{A_2},x_{B_1},x_{B_2}]
=\prod_{i,j}\E[\mathcal D(A_i,B_j)\mid x_{A_i},x_{B_j}]
=\prod_{i,j}g_\tau(x_{A_i},x_{B_j}).
\]
Integrating over the four independent groups gives
\[
\E[\Pi_\tau]
=\int\!\!\int\left\{\int g_\tau(x,y)g_\tau(x',y)\,\d y\right\}^2\,\d x\,\d x'
=Z_\tau.
\]
Since $\widehat Z_\tau$ averages copies of $\Pi_\tau$, we also have $\E[\widehat Z_\tau]=Z_\tau$.

\textbf{(ii)} Next, we compare the weighted sum of these expectations with $\Upsilon_N$. All comparison constants below depend only on $(p,q)$. Decompose $g_\tau=\sum_S\varphi_S$ into labeled projections. Their number is bounded in terms of $(p,q)$, so the kernel bounds in Lemma~\ref{lem:S-blocks} give
\[
Z_\tau=\|g_\tau\|_{S_4}^4
\asymp\sum_S\|\varphi_S\|_{S_4}^4
=\sum_S\|\varphi_S\otimes_B\varphi_S\|_2^2.
\]

Fix a labeled vertex set $V$ with $r$ rows and $c$ columns. For an ordered split $V=A\sqcup B$ into nonempty sets, let $Z(A,B)$ be the corresponding contraction $Z_\tau$. There are $\binom r{r_A}\binom c{c_A}$ splits with $r_A$ rows and $c_A$ columns in $A$, all with the same contraction by exchangeability. A graph with no edge between $A$ and $B$ is disconnected, and $B$ is a nonempty proper union of its connected components. Thus,
\[
\begin{aligned}
\sum_{\tau:\,r_A+r_B=r,\,c_A+c_B=c}c_\tau Z_\tau
&=\{(m)_r(n)_c\}^{-2}
  \sum_{\substack{A\sqcup B=V\\A,B\ne\varnothing}}Z(A,B)\\
&\asymp\{(m)_r(n)_c\}^{-2}
  \sum_\Lambda\sum_S\sum_B\|\varphi_S\otimes_B\varphi_S\|_2^2,
\end{aligned}
\]
where, in the second line, $\Lambda$ ranges over disconnected labeled graphs on $V$, $S$ over the labeled projections of $q_\Lambda$, and $B$ over the nonempty proper unions of connected components of $\Lambda$. We identify each such union with the set of components used to index the contractions in \eqref{eq:gamma}.

Now, for $\ell(G)=(r,c)$, the normalization $p_G=(p)_r(q)_c q_G$ gives $f_S=(p)_r(q)_c\varphi_S$. Applying the same kernel bounds to $p_G=\sum_Sf_S$ gives
\[
\sum_S\|\varphi_S\otimes_B\varphi_S\|_2^2
\asymp\{(p)_r(q)_c\}^{-4}\|p_G\otimes_Bp_G\|_2^2.
\]
There are $r!c!/|\mathrm{Aut}(G)|$ labeled graphs of class $G$ on $V$, with equal contraction sums by exchangeability. Their coefficient satisfies
\[
\frac{r!c!}{|\mathrm{Aut}(G)|\{(p)_r(q)_c\}^4\{(m)_r(n)_c\}^2}
=\frac{r!c!|\mathrm{Aut}(G)|}{\{(p)_r(q)_c\}^4}\,w_G^2
\asymp w_G^2,
\]
because $1\le|\mathrm{Aut}(G)|\le r!c!$ and $1\le(p)_r(q)_c\le p!q!$. Grouping by $G$ and summing over $(r,c)$ therefore gives, by \eqref{eq:gamma},
\[
\begin{aligned}
\sum_\tau c_\tau\E[\widehat Z_\tau]
&=\sum_\tau c_\tau Z_\tau\\
&\asymp\sum_{G\ \disc}\sum_Bw_G^2\|p_G\otimes_Bp_G\|_2^2
=\sigma_N^4\Upsilon_N.
\end{aligned}
\]
\end{proof}

\paragraph{A sufficient condition for (S2).}
\begin{proof}[Proof of Lemma~\ref{lem:S-filler-span}]
Choose an orthonormal basis $(g_k)$ for the span in $L^2(U)$, and include the basis-change coefficients in the Hilbert--Schmidt kernels $f_k$, viewed as integral operators. Let $\mathcal T:L^2(x_B)\to L^2(x_{A\cup U})$ be the integral operator with kernel $F$. Since the $g_k$ are orthonormal, we have
\[
\mathcal T^*\mathcal T=\sum_kf_k^*f_k,\qquad
\|\mathcal T\|_{S_4}^4=\operatorname{tr}\!\left[\left(\sum_kf_k^*f_k\right)^2\right]
\ge\sum_k\|f_k\|_{S_4}^4,
\]
where the inequality follows by expanding the square and using $\operatorname{tr}(f_k^*f_kf_l^*f_l)\ge0$. Minkowski's inequality in $L^2(U;S_4)$ and H\"older's inequality give
\[
\left(\E[\left\|\sum_kg_kf_k\right\|_{S_4}^2]\right)^2
\le\left(\sum_k\|f_k\|_{S_4}\right)^4
\le M^3\sum_k\|f_k\|_{S_4}^4
\le M^3\|\mathcal T\|_{S_4}^4.
\]
\end{proof}

\paragraph{Estimation error and diagnostic guarantees.}\label{def:S-diagnostic-quantities}
Here, we prove Theorem~\ref{thm:S-check}. Write $x_\tau=c_\tau Z_\tau/\sigma_N^4$ and $\widehat x_\tau=c_\tau\widehat Z_\tau/\sigma_N^4$. Then $\sum_\tau x_\tau\asymp\Upsilon_N$. Define $X_\tau=\max_{\tau'\succ\tau}x_{\tau'}$, taking the maximum over feasible indices with strictly more vertices. By convention, this maximum is zero when there are none. The constants below depend only on $(p,q,C_{\mathrm{M8}},C_{\mathrm{S2}})$.

Let $\Omega_\tau$ be the set of choices of rows and columns, including replacement and additional vertices, used to compute $\widehat Z_\tau$. We call each choice a placement $\omega$ and write $\Pi_\omega$ for its product \eqref{eq:diagnostic-product}. Since $\E[\Pi_\omega]=Z_\tau$,
\begin{equation}\label{eq:S-diagnostic-covariance}
\E[(\widehat x_\tau-x_\tau)^2]
=\frac{c_\tau^2}{\sigma_N^8|\Omega_\tau|^2}
  \sum_{\omega,\omega'\in\Omega_\tau}\cov(\Pi_\omega,\Pi_{\omega'}).
\end{equation}
We bound this sum over two sets of ordered pairs. Let $\mathcal C_{\tau,1}$ be the set of pairs in $\Omega_\tau^2$ sharing a vertex that is original or a replacement in at least one placement. Let $\mathcal C_{\tau,2}$ be the set of pairs sharing at least one vertex, with every shared vertex additional in both placements. The following lemma uses (M8) and vertex counting to bound the sum over $\mathcal C_{\tau,1}$.

\begin{lemma}\label{lem:S-original-covariance}
Under the assumptions of Theorem~\ref{thm:S-check}, for all sufficiently large $N$ and every feasible $\tau$,
\[
\frac{c_\tau^2}{\sigma_N^8|\Omega_\tau|^2}
\sum_{(\omega,\omega')\in\mathcal C_{\tau,1}}
|\cov(\Pi_\omega,\Pi_{\omega'})|
\le CL_N^{-1/2}.
\]
\end{lemma}

\begin{proof}[Proof of Lemma~\ref{lem:S-original-covariance}]
To bound the covariances, we expand each alternating sum in $\Pi_\omega$ into labeled Hoeffding projections of the $q_\Lambda$, using \eqref{eq:diagnostic-differences}. A term remains precisely when it contains exactly one of $u$ and its replacement $u''$ for each original position $u$. If neither occurs, the two choices cancel. Each replacement changes the sign. Additional vertices may also occur. By (M8), each labeled projection has $L^8$ norm at most $C\|q_G\|_2$, since it is a finite signed sum of conditional expectations of $q_G$.

Write $r=r_A+r_B$ and $c=c_A+c_B$. For a product $P$ of four such projections, let $s$ count the replacement rows used, $z$ the original rows replaced in both factors containing them, and $t$ the additional rows used. The row count with multiplicity and the number of distinct rows are
\[
r_{\mathrm{tot}}(P)=4r+t,\qquad v_{P,r}=2r-z+s+t.
\]
Consequently,
\[
2r+v_{P,r}-r_{\mathrm{tot}}(P)=s-z=:\varphi_r(P)\ge0.
\]
Each original row replaced in at least one factor contributes one to $\varphi_r(P)$. Define $c_{\mathrm{tot}}(P)$, $v_{P,c}$ and $\varphi_c(P)$ in the same way for columns. The numbers of possible graphs, signs and assignments depend only on $(p,q)$.

Now, we fix a pair $(\omega,\omega')\in\mathcal C_{\tau,1}$, products $P,P'$ from their expansions and a pattern of shared vertices. Let $v_r$ count distinct rows used by the two products, and let $\chi_r$ count shared rows absent from at least one product, including rows absent from both. Define $v_c,\chi_c$ similarly. The proportion of placement pairs with this pattern is at most
\[
C m^{v_r-v_{P,r}-v_{P',r}-\chi_r}
  n^{v_c-v_{P,c}-v_{P',c}-\chi_c}.
\]
Indeed, each identification removes a free row or column index. There are $v_{P,r}+v_{P',r}-v_r$ shared rows used by both products and $\chi_r$ other shared rows, and similarly for columns.

Moreover, a nonzero expectation requires every independent AHK variable to occur at least twice. Each vertex either supplies its own variable or is an endpoint of an edge variable, so it also occurs at least twice. This applies to the expectation of all eight projections and to each of the two four-factor expectations subtracted in the covariance. Thus, a nonzero covariance requires
\[
\ell_r=r_{\mathrm{tot}}(P)+r_{\mathrm{tot}}(P')-2v_r\ge0,
\qquad
\ell_c=c_{\mathrm{tot}}(P)+c_{\mathrm{tot}}(P')-2v_c\ge0.
\]
Let $G_1,\ldots,G_8$ be the graph classes of the factors in $P,P'$. By \eqref{eq:copy-scale},
\[
\|q_G\|_2\le C\varsigma_Gm^{r(G)/2}n^{c(G)/2}.
\]
H\"older's inequality, the preceding $L^8$ bound, and $c_\tau\le Cm^{-2r}n^{-2c}$ give
\[
\begin{aligned}
\frac{c_\tau^2}{\sigma_N^8}|\cov(P,P')|
&\le\frac{Cc_\tau^2}{\sigma_N^8}\prod_{j=1}^8\|q_{G_j}\|_2\\
&\le C m^{\{r_{\mathrm{tot}}(P)+r_{\mathrm{tot}}(P')\}/2-4r}
 n^{\{c_{\mathrm{tot}}(P)+c_{\mathrm{tot}}(P')\}/2-4c}
 \prod_{j=1}^8\frac{\varsigma_{G_j}}{\sigma_N}.
\end{aligned}
\]
When we multiply by the proportion of pairs, the row exponent becomes
\[
\begin{aligned}
&v_r-v_{P,r}-v_{P',r}-\chi_r
+\frac{r_{\mathrm{tot}}(P)+r_{\mathrm{tot}}(P')}2-4r\\
&\qquad=-\frac{\ell_r}2-\varphi_r(P)-\varphi_r(P')-\chi_r.
\end{aligned}
\]
The same calculation for columns bounds the contribution of this pattern by
\[
C m^{-\ell_r/2-\varphi_r(P)-\varphi_r(P')-\chi_r}
 n^{-\ell_c/2-\varphi_c(P)-\varphi_c(P')-\chi_c}
 \prod_{j=1}^8\frac{\varsigma_{G_j}}{\sigma_N}.
\]

Finally, we consider a shared row that is original or a replacement in at least one placement. If it is absent from either product, then $\chi_r\ge1$. If it is used as a replacement, or an occurrence of it as an original row is replaced, then $\varphi_r(P)+\varphi_r(P')\ge1$. Otherwise, the row appears twice in one product and at least once in the other, so $\ell_r\ge1$. For a shared column, the same argument uses $\chi_c$, $\varphi_c(P)+\varphi_c(P')$ and $\ell_c$. Thus, at least one exponent is at most $-1/2$. Summing over the finitely many patterns and graph classes and applying Lemma~\ref{lem:auto-pattern}(ii), we obtain
\[
\begin{aligned}
&\frac{c_\tau^2}{\sigma_N^8|\Omega_\tau|^2}
\sum_{(\omega,\omega')\in\mathcal C_{\tau,1}}
|\cov(\Pi_\omega,\Pi_{\omega'})|\\
&\qquad\le CL_N^{-1/2}\left(\frac{\sum_G\varsigma_G}{\sigma_N}\right)^8
\le CL_N^{-1/2}.
\end{aligned}
\]
\end{proof}

The next lemma uses (S2) to bound the sum over $\mathcal C_{\tau,2}$ in terms of $x_\tau$ and $X_\tau$.

\begin{lemma}\label{lem:S-additional-covariance}
Under the assumptions of Theorem~\ref{thm:S-check}, for all sufficiently large $N$ and every feasible $\tau$,
\[
\frac{c_\tau^2}{\sigma_N^8|\Omega_\tau|^2}
\sum_{(\omega,\omega')\in\mathcal C_{\tau,2}}
|\cov(\Pi_\omega,\Pi_{\omega'})|
\le C\bigl(x_\tau^{3/2}X_\tau^{1/2}+X_\tau^2\bigr).
\]
\end{lemma}

\begin{proof}[Proof of Lemma~\ref{lem:S-additional-covariance}]
Fix $(\omega,\omega')\in\mathcal C_{\tau,2}$. To bound $\cov(\Pi_\omega,\Pi_{\omega'})$, we condition on their shared additional vertices. Let $U$ be this shared vertex set and $\mathcal U$ the $\sigma$-field of its row, column and cell variables. All remaining variables are independent between the two products. With $G_\omega=\E[\Pi_\omega\mid\mathcal U]$, we have
\[
\begin{aligned}
\cov(\Pi_\omega,\Pi_{\omega'})
&=\cov\big(\E[\Pi_\omega\mid\mathcal U],\E[\Pi_{\omega'}\mid\mathcal U]\big)\\
&=\E[(G_\omega-Z_\tau)(G_{\omega'}-Z_\tau)].
\end{aligned}
\]
For factor $(i,j)$, let $W_{ij}$ be its additional vertex set. For $\varnothing\ne U'\subseteq U\cap W_{ij}$, let $F^{U'}_{ij}$ sum the ordinary projections on exactly $A_i\cup B_j\cup U'$ with no edge between the three sets, and put $F^{\varnothing}_{ij}=g_\tau$. Write
\[
R(K_1,K_2,K_3,K_4)=\operatorname{tr}(K_1K_2^*K_4K_3^*)
\]
for the four-factor integral. After integrating over replacement vertices, unshared additional vertices and cell variables used by only one factor, we expand the four factors. The term with all $U'_{ij}=\varnothing$ is $R(g_\tau,g_\tau,g_\tau,g_\tau)=Z_\tau$, so
\[
G_\omega-Z_\tau
=\sum_{\substack{U'_{ij}\subseteq U\cap W_{ij}\\\bigcup_{i,j}U'_{ij}\ne\varnothing}}
R(F^{U'_{11}}_{11},F^{U'_{21}}_{21},F^{U'_{12}}_{12},F^{U'_{22}}_{22}).
\]
The four sets $W_{ij}$ are disjoint, and no factor uses cells between different sets. Thus, the four kernels in each term depend on independent groups of additional variables. The Schatten--H\"older inequality gives
\[
\|R(K_1,K_2,K_3,K_4)\|_{L^2(U)}
\le\prod_{j=1}^4\left(\E_U[\|K_j\|_{S_4}^2]\right)^{1/2},
\]
where each $\|K_j\|_{S_4}$ is computed with the additional variables fixed.
We number the four kernels of a term of $G_\omega-Z_\tau$ as $F_1,\dots,F_4$, and those of a term of $G_{\omega'}-Z_\tau$ as $F_5,\dots,F_8$, with additional vertex sets $U'_1,\dots,U'_8$. The triangle inequality and the Cauchy--Schwarz inequality therefore give
\begin{equation}\label{eq:S-additional-covariance-bound}
|\cov(\Pi_\omega,\Pi_{\omega'})|
\le\sum\prod_{j=1}^8\bigl(\E_U[\|F_j\|_{S_4}^2]\bigr)^{1/2},
\end{equation}
where the sum is over the finitely many pairs of terms, one from each expansion, so that at least one of $U'_1,\dots,U'_4$ and at least one of $U'_5,\dots,U'_8$ is nonempty.

To bound a factor $F=F^{U'}_{ij}$ with $U'\ne\varnothing$ in terms of $X_\tau$, we apply (S2) with $A=A_i$, $B=B_j$ and $U=U'$. For a split $(U_1,U_2)\in\mathcal A(F)$, let $\tau'=\tau'(U_1,U_2)$ be the diagnostic index of $A\cup U_1,B\cup U_2$. We view $F$ as an operator from $L^2(x_{B\cup U_2})$ to $L^2(x_{A\cup U_1})$. Let $P_0,Q_0$ be the orthogonal projections obtained by integrating the cell variables joining $A$ to $U_1$ and $B$ to $U_2$, respectively. No nonzero graph term of $F$ has an edge between $U_1$ and $U_2$. By graph orthogonality, it follows that
\begin{equation}\label{eq:S-operator-bound}
F=P_0g_{\tau'}Q_0,\qquad
\|F\|_{S_4}\le\|g_{\tau'}\|_{S_4}=Z_{\tau'}^{1/4}.
\end{equation}

Next, the proportion of pairs sharing $u_r$ rows and $u_c$ columns is at most $Cm^{-u_r}n^{-u_c}$. In a product of eight factors in \eqref{eq:S-additional-covariance-bound}, write $U'_{j,r},U'_{j,c}$ for the rows and columns of $U'_j$. Every shared vertex can be selected at most once in each placement, so
\[
\sum_{j=1}^8|U'_{j,r}|\le2u_r,\qquad
\sum_{j=1}^8|U'_{j,c}|\le2u_c.
\]
Consequently, we can split $m^{-u_r}n^{-u_c}$ among the eight factors, using
\[
m^{-u_r}n^{-u_c}
\le\prod_{j=1}^8m^{-|U'_{j,r}|/2}n^{-|U'_{j,c}|/2}.
\]
Since $c_\tau^2/\sigma_N^8=(c_\tau^{1/4}/\sigma_N)^8$, each factor also receives $c_\tau^{1/4}/\sigma_N$. For a nonempty $U'$, the index $\tau'$ has $|U'_r|$ more rows and $|U'_c|$ more columns, so \eqref{eq:diagnostic-statistic} gives
\[
c_\tau m^{-2|U'_r|}n^{-2|U'_c|}\asymp c_{\tau'}.
\]
Combining this comparison with (S2) and~\eqref{eq:S-operator-bound} gives
\[
\begin{aligned}
&\frac{c_\tau^{1/4}m^{-|U'_r|/2}n^{-|U'_c|/2}}{\sigma_N}
\left(\E_{x_{U'}}[\|F(\cdot,\cdot,x_{U'})\|_{S_4}^2]\right)^{1/2}\\
&\quad\le\frac{C}{\sigma_N}
\max_{(U_1,U_2)\in\mathcal A(F)}
\left(c_\tau m^{-2|U'_r|}n^{-2|U'_c|}Z_{\tau'(U_1,U_2)}\right)^{1/4}\\
&\quad\le C\max_{(U_1,U_2)\in\mathcal A(F)}
\left(\frac{c_{\tau'(U_1,U_2)}Z_{\tau'(U_1,U_2)}}{\sigma_N^4}\right)^{1/4}
\le CX_\tau^{1/4}.
\end{aligned}
\]
The last inequality holds because every $\tau'(U_1,U_2)$ has more vertices than $\tau$. In contrast, a factor without additional vertices contributes $c_\tau^{1/4}Z_\tau^{1/4}/\sigma_N=x_\tau^{1/4}$. Each of the two expansions contributes at least one nonempty $U'_j$, so the number $k$ of factors with nonempty $U'_j$ satisfies $2\le k\le8$. Their normalized contribution is at most
\[
Cx_\tau^{2-k/4}X_\tau^{k/4}
\le C\begin{cases}
x_\tau^{3/2}X_\tau^{1/2},&X_\tau\le x_\tau,\\
X_\tau^2,&X_\tau>x_\tau.
\end{cases}
\]
Summing over the finitely many terms and patterns proves the stated bound.
\end{proof}

Combining these bounds gives the following estimate for the proof of Theorem~\ref{thm:S-check}.

\begin{lemma}\label{lem:S-diagnostic-error}
Under the assumptions of Theorem~\ref{thm:S-check}, for all sufficiently large $N$ and every feasible $\tau$, the normalized estimates defined above satisfy
\begin{equation}\label{eq:S-noise}
\E[(\widehat x_\tau-x_\tau)^2]\le C\big(L_N^{-1/2}+x_\tau^{3/2}X_\tau^{1/2}+X_\tau^2\big).
\end{equation}
Here $C$ depends only on $(p,q,C_{\mathrm{M8}},C_{\mathrm{S2}})$.
\end{lemma}

\begin{proof}[Proof of Lemma~\ref{lem:S-diagnostic-error}]
In \eqref{eq:S-diagnostic-covariance}, the sum over $\Omega_\tau^2$ can be decomposed into one over $\mathcal C_{\tau,1}$ and one over $\mathcal C_{\tau,2}$. The remaining pairs share no vertices, so their products are independent by dissociation and their covariances vanish. Lemmas~\ref{lem:S-original-covariance} and~\ref{lem:S-additional-covariance} give
\[
\begin{aligned}
\E[(\widehat x_\tau-x_\tau)^2]
&\le\frac{c_\tau^2}{\sigma_N^8|\Omega_\tau|^2}
\sum_{j=1}^2\sum_{(\omega,\omega')\in\mathcal C_{\tau,j}}
|\cov(\Pi_\omega,\Pi_{\omega'})|\\
&\le C\bigl(L_N^{-1/2}+x_\tau^{3/2}X_\tau^{1/2}+X_\tau^2\bigr).
\end{aligned}
\]
\end{proof}

\begin{proof}[Proof of Theorem~\ref{thm:S-check}]
We use Lemma~\ref{lem:S-diagnostic-error} in two ways. When $\Upsilon_N$ is small, all products and their estimation errors are small relative to the selection threshold $\varepsilon_N$. When $\Upsilon_N$ stays positive, we must instead find at least one product whose estimate stays positive and is large relative to that threshold. We prove these claims in this order.

\textbf{Proof of (ii).} We show that the estimated numerators $\widehat x_\tau$ are negligible relative to $\varepsilon_N$ and that the estimated variance $\Vu$ is ratio-consistent. Since $x_\tau,X_\tau\le C\Upsilon_N$, applying \eqref{eq:S-noise} on $\Upsilon_N\le a_N\varepsilon_N$ gives
\[
\frac{\E[(\widehat x_\tau-x_\tau)^2]}{\varepsilon_N^2}
\le C\left(\frac{L_N^{-1/2}}{\varepsilon_N^2}+a_N^2\right)\longrightarrow0.
\]
Here $L_N^{-1/4}=o(\varepsilon_N)$ by assumption. Also, $x_\tau/\varepsilon_N\le Ca_N\to0$. Chebyshev's inequality and a union bound over the finite set of indices show that all positive numerators $\widehat x_\tau^+$ in the diagnostic are $o_P(\varepsilon_N)$, uniformly on these classes. Meanwhile, $\Vu/\sigma_N^2\to1$ by Theorem~\ref{thm:contraction}. Thus, the diagnostic selects the Gaussian test with the probability stated in (ii).

\textbf{Proof of (i).} If $\Upsilon_N\ge\delta$, Proposition~\ref{prop:ident} gives $\max_\tau x_\tau\ge a\delta$, where $a>0$ depends only on $(p,q)$. For $x_\tau>0$, \eqref{eq:S-noise} gives
\[
\frac{\E[(\widehat x_\tau-x_\tau)^2]}{x_\tau^2}
\le C\left\{\frac{L_N^{-1/2}}{x_\tau^2}
+\left(\frac{X_\tau}{x_\tau}\right)^{1/2}
+\left(\frac{X_\tau}{x_\tau}\right)^2\right\}.
\]
Put $\kappa_N=\varepsilon_N^{1/\{2(p+q)\}}$. Start from a maximizing index. Whenever $X_\tau>\kappa_Nx_\tau$, move to an index with more vertices attaining $X_\tau$. Each move increases the number of vertices and replaces $x_\tau$ by a value greater than $\kappa_Nx_\tau$. There can be at most $p+q$ moves. We therefore reach $\tau^*$ with
\[
x_{\tau^*}\ge a\delta\kappa_N^{p+q}=a\delta\sqrt{\varepsilon_N},\qquad
X_{\tau^*}\le\kappa_Nx_{\tau^*}.
\]
By \eqref{eq:diagnostic-statistic},
\[
\{\widehat\Upsilon_N\le\varepsilon_N\}
\subseteq\{\widehat x_{\tau^*}<x_{\tau^*}/2\}
\cup\{(\Vu/\sigma_N^2)^2\ge x_{\tau^*}/(2\varepsilon_N)\}.
\]
Theorem~\ref{thm:contraction} and $\Upsilon_N\le CD_N^2\le C$ give $\E[(\Vu/\sigma_N^2)^2]\le C$. Chebyshev's inequality and Markov's inequality therefore yield
\[
\begin{aligned}
\P(\widehat\Upsilon_N\le\varepsilon_N)
&\le\frac{4\E[(\widehat x_{\tau^*}-x_{\tau^*})^2]}{x_{\tau^*}^2}
+\frac{2\varepsilon_N\E[(\Vu/\sigma_N^2)^2]}{x_{\tau^*}}\\
&\le C_\delta\left\{\frac{L_N^{-1/2}}{\varepsilon_N}
+\varepsilon_N^{1/\{4(p+q)\}}+\varepsilon_N^{1/(p+q)}
+\sqrt{\varepsilon_N}\right\}\longrightarrow0.
\end{aligned}
\]
The constants are uniform over $\Pcal$, which proves (i).
\end{proof}

\subsection{The procedure}\label{sec:S-procedure}

To combine a conservative bound with the Gaussian test, we use the inclusion of the rejection events and the diagnostic guarantees. As in \eqref{eq:hybrid-definition}, $C_N$ is the event on which the Gaussian test is chosen, $R_{G,N}$ and $R_{S,N}$ are the Gaussian and conservative rejection events, and $R_{H,N}$ is their combined rejection event.

\begin{lemma}\label{lem:hybrid-switch}
Let $\alpha\in(0,1)$, and let $\Pcal_0$ be a null class of laws. Let $t_N:\Pcal_0\to[0,M]$ be deterministic, with $M<\infty$. Suppose that
\begin{enumerate}[(a)]
\item $R_{S,N}\subseteq R_{G,N}$ for every dataset, and $\limsup_N\sup_{P\in\Pcal_0}\P_P(R_{S,N})\le\alpha$;
\item along every sequence of increasing indices and laws in $\Pcal_0$ with $t_N(P)\to0$, $\P_P(R_{G,N})\to\alpha$;
\item for every $\delta>0$, $\sup_{P\in\Pcal_0:t_N(P)\ge\delta}\P_P(C_N)\to0$.
\end{enumerate}
Then $\limsup_N\sup_{P\in\Pcal_0}\P_P(R_{H,N})\le\alpha$. For specified nonempty subclasses $\mathcal Q_N\subseteq\Pcal_0$, if
\[
\sup_{P\in\mathcal Q_N}t_N(P)\to0,
\qquad \sup_{P\in\mathcal Q_N}\P_P(C_N^c)\to0,
\]
then $\sup_{P\in\mathcal Q_N}|\P_P(R_{H,N})-\alpha|\to0$.
\end{lemma}
\begin{proof}[Proof of Lemma~\ref{lem:hybrid-switch}]
For contradiction, suppose that the first conclusion of the lemma fails. Choose a sequence of indices and null laws along which the rejection probability $\P(R_{H,N})$ of the procedure exceeds $\alpha$ by a fixed positive amount. Since $t_N$ is bounded, we can pass to a further subsequence on which $t_N\to t$.

If $t=0$, the inclusion in (a) gives $R_{H,N}\subseteq R_{G,N}$, so the rejection probability of the procedure is at most that of the Gaussian test, $\P(R_{G,N})$, which tends to $\alpha$ by (b). If instead $t>0$, the combined rejection probability $\P(R_{H,N})$ is at most $\P(R_{S,N})+\P(C_N)$, whose upper limit is at most $\alpha$ by (a) and (c). Both cases contradict the choice of the sequence, so this proves the first conclusion.

Finally, on $\mathcal Q_N$, the same argument applied to (b) shows that the Gaussian rejection probability $\P(R_{G,N})$ converges uniformly to $\alpha$. The decisions of the procedure and the Gaussian test differ only on $C_N^c$, whose probability tends uniformly to zero.
\end{proof}

\begin{proof}[Proof of Proposition~\ref{prop:S-block-hybrid}]
Proposition~\ref{prop:block-safe} gives the conservative rejection bound required in Lemma~\ref{lem:hybrid-switch}(a). Since $B_N\ge A_N$ and $\alpha^{-1/2}\ge z_{1-\alpha/2}$ by Chebyshev's inequality applied to a standard normal variable, $R_{S,N}\subseteq R_{G,N}$ for every dataset. Both tests require $A_N>0$. For the upper one-sided tests, Cantelli's inequality applied to a standard normal variable gives \(z_{1-\alpha}\le\sqrt{(1-\alpha)/\alpha}\). With \(B_N\ge A_N\), this again yields \(R_{S,N}\subseteq R_{G,N}\).

Under (M8) and (S2), Theorem~\ref{thm:S-check} shows that the diagnostic is unlikely to choose the Gaussian test when $\Upsilon_N$ stays positive: for every fixed $\delta>0$, $\sup_{P:\Upsilon_N(P)\ge\delta}\P_P(C_N)\to0$. In addition, Theorem~\ref{thm:contraction} gives $\widehat Z_N\dto\Nor(0,1)$ along every sequence with $\Upsilon_N\to0$. Finally, $0\le\Upsilon_N\le CD_N^2\le C$ by Proposition~\ref{prop:gammaD}.

Lemma~\ref{lem:hybrid-switch} with $t_N=\Upsilon_N$ therefore gives the uniform upper bound on the rejection probability.

Moreover, on the subclasses where $\Upsilon_N\le a_N\varepsilon_N$, Theorem~\ref{thm:S-check}(ii) gives $\sup_P\P_P(C_N^c)\to0$. Also, $\sup_P\Upsilon_N\to0$ there. The last assertion of Lemma~\ref{lem:hybrid-switch} therefore gives convergence of the rejection probability to $\alpha$, uniformly on these shrinking classes.
\end{proof}

\subsection{An alternative diagnostic under additional assumptions}\label{sec:S-chaos-procedure}

\paragraph{Selected variance components.} We use assumption (L) from Section~\ref{sec:alternative-procedure} of the paper. Recall that $\ell(G)=(r(G),c(G))$, and that $\widehat V_\ell$ estimates the variance contributed by graphs with these vertex counts (Sec.~\ref{sec:S-variance} and Lem.~\ref{lem:overlap}). Let $\mathcal L_{\mathrm{disc}}=\{\ell(G):G\in\Gamma_{\mathrm{disc}}\}$, and define
\begin{equation}\label{eq:chaos-variance-adjustment}
F_N=\sum_{\ell\in\mathcal L_{\mathrm{disc}}}\widehat V_\ell^+ .
\end{equation}
The following lemma bounds its estimation error.

\begin{lemma}\label{thm:S-test}

Let $\Pcal$ be a class of laws with $\sigma_{N,P}>0$ for every $P\in\Pcal$, and let $m\ge2p$, $n\ge2q$. Assume (M4) and (L) on $\Pcal$. Put
\[
\rho_N=\sum_{\substack{G\ \disc\\G\notin\Gamma_{\mathrm{disc}}}}s_{N,G}^2,
\qquad v_N=\sum_{\ell\in\mathcal L_{\mathrm{disc}}}\frac{V_\ell}{\sigma_N^2}.
\]
Then for every $P\in\Pcal$
\begin{equation}\label{eq:S-selected-variance-bound}
\E_P\!\left[\left(\frac{F_N}{\sigma_N^2}-v_N\right)^2\right]
\le C(L_N^{-1/2}+\rho_N),\qquad D_N-\rho_N\le v_N\le1,
\end{equation}
where $C$ depends only on $(p,q,C_{\mathrm{M4}})$. Consequently, $F_N/\sigma_N^2\le1+o_P(1)$ uniformly over $\Pcal$. Moreover, if $\Upsilon_N\to0$ along a sequence of laws in $\Pcal$, then $\max\{\Vu,F_N\}/\sigma_N^2\to1$ in probability and $(U_N-\theta_N)/\max\{\Vu,F_N\}^{1/2}\dto\Nor(0,1)$.
\end{lemma}
\begin{proof}[Proof of Lemma~\ref{thm:S-test}]

\textbf{(i)} Fix $\ell\in\mathcal L_{\mathrm{disc}}$.  By orthogonality of the grouped projections $X_A$, the definition of $J_\ell$ in~\eqref{eq:component-J} gives
\[
\E[J_\ell]=\sum_{\ell(A)=\ell}c_A^2\E[X_A^2]=V_\ell.
\]
To bound $\var(J_\ell)$, we use the covariance expansion in the proof of Lemma~\ref{lem:auto-decomp}(ii), restricted to
\[
\ell\bigl(N(\Lambda_1)\cap N(\Lambda_2)\bigr)
=\ell\bigl(N(\Lambda_3)\cap N(\Lambda_4)\bigr)=\ell.
\]
The paired configurations with nonzero covariance have $d_{\mathrm{exc}}\ge2$, as in that proof. Lemma~\ref{lem:S-cover}(iii) and Lemma~\ref{lem:auto-pattern}(ii) bound their contribution by $CL_N^{-1}\sigma_N^4$ and that of the unpaired configurations with $d_{\mathrm{exc}}\ge1$ by $CL_N^{-1/2}\sigma_N^4$.

It remains to consider an unpaired configuration with $d_{\mathrm{exc}}=0$ and nonzero covariance. Every vertex belongs to two copies, and Lemma~\ref{lem:S-atoms} shows that connected components are shared whole. Thus, in each pair $(1,2)$ and $(3,4)$, at least one copy shares a nonempty proper subcollection of its components with its partner. Otherwise, the two copies in that pair would coincide,  and be disjoint from the other pair, giving zero covariance. Now, in the pair $(1,2)$, let $\Lambda_k$ be such a copy, with class $G_k$. The components it shares with its partner cover $N(\Lambda_1)\cap N(\Lambda_2)$, so their total numbers of rows and columns are given by $\ell$. Since $\ell\in\mathcal L_{\mathrm{disc}}$, we have $\ell=\ell(G')$ for some $G'\in\Gamma_{\mathrm{disc}}$. Therefore, by (L)(b), $G_k\notin\Gamma_{\mathrm{disc}}$. Moreover, since $G_k$ also has a component that is not shared, it is disconnected. The same argument applies to the pair $(3,4)$. Finally, Lemma~\ref{lem:S-cover}(iii) bounds the sum over labelings of each pattern by $C\prod_{k=1}^4\varsigma_{G_k}$. Summing over the finitely many patterns and graph classes gives
\[
C\left(\sum_{\substack{G\ \disc\\G\notin\Gamma_{\mathrm{disc}}}}\varsigma_G\right)^2
\left(\sum_G\varsigma_G\right)^2
\le C\rho_N\sigma_N^4
\]
by the Cauchy--Schwarz inequality. Thus,
\[
\frac{\var(J_\ell)}{\sigma_N^4}\le C(L_N^{-1/2}+\rho_N).
\]
Since $V_\ell\ge0$, we have $|\widehat V_\ell^+-V_\ell|\le|\widehat V_\ell-V_\ell|$. Lemma~\ref{lem:S-level} now gives
\[
\begin{aligned}
\E\!\left[\left(\frac{F_N}{\sigma_N^2}-v_N\right)^2\right]
&\le C\sum_{\ell\in\mathcal L_{\mathrm{disc}}}\frac{\E[(\widehat V_\ell-V_\ell)^2]}{\sigma_N^4}\\
&\le C\sum_{\ell\in\mathcal L_{\mathrm{disc}}}
\frac{\E[(\widehat V_\ell-J_\ell)^2]+\var(J_\ell)}{\sigma_N^4}
\le C(L_N^{-1/2}+\rho_N),
\end{aligned}
\]
which proves the mean-square bound in \eqref{eq:S-selected-variance-bound}.

\textbf{(ii)} Since $V_\ell$ is the variance contributed by the graphs with row and column counts $\ell$, we have $v_N=\sum_{G:\,\ell(G)\in\mathcal L_{\mathrm{disc}}}s_{N,G}^2$. Every $G\in\Gamma_{\mathrm{disc}}$ has $\ell(G)\in\mathcal L_{\mathrm{disc}}$. Thus,
\[
v_N\ge\sum_{G\in\Gamma_{\mathrm{disc}}}s_{N,G}^2=D_N-\rho_N,
\qquad v_N\le\sum_Gs_{N,G}^2=1.
\] 
By (L)(a) and the finiteness of $\Gamma$, $\rho_N\to0$ uniformly. Consequently,~\eqref{eq:S-selected-variance-bound} and Chebyshev's inequality give
\[
\frac{F_N}{\sigma_N^2}=v_N+o_P(1)\le1+o_P(1)
\]
uniformly over $\Pcal$.

\textbf{(iii)} Moreover, if $\Upsilon_N\to0$ along a sequence of laws in $\Pcal$, then Theorem~\ref{thm:contraction} gives both $\Vu/\sigma_N^2\to1$ in probability and $(U_N-\theta_N)/\sigma_N\dto\Nor(0,1)$. Together with (ii), this gives $\max\{\Vu,F_N\}/\sigma_N^2 \to 1$ in probability and by Slutsky's lemma, $(U_N-\theta_N)/\max\{\Vu,F_N\}^{1/2} \dto \Nor(0,1)$.
\end{proof}

\paragraph{The alternative block procedure.}
Let $0<\alpha<1$ and let $\Pcal_0$ be a null class satisfying the assumptions of Proposition~\ref{prop:block-safe} and (L). We can use an alternative test, alongside a simpler diagnostic. For the test, define
\[
\begin{aligned}
A_N&=\max\{\Vu,F_N\},\qquad B_N=\max\{A_N,S_{K_N}^2/K_N\},\\
R_{G,N}&=\{A_N>0,\ |U_N-\theta_0|\ge z_{1-\alpha/2}\sqrt{A_N}\},\\
R_{S,N}&=\{A_N>0,\ |U_N-\theta_0|\ge\sqrt{B_N/\alpha}\}.
\end{aligned}
\]
For a diagnostic event $C_N$, the combined procedure rejects on $(C_N\cap R_{G,N})\cup(C_N^c\cap R_{S,N})$. Now, for the diagnostic, we consider two possibilities: the complete diagnostic gives the event
\[
C_N^{\rm ctr}=\{\widehat\Upsilon_N\le\varepsilon_N\},\qquad
\varepsilon_N\to0,\quad \varepsilon_NL_N^{1/4}\to\infty,
\]
and the simpler diagnostic uses $F_N$ in \eqref{eq:chaos-variance-adjustment}:
\begin{equation}\label{eq:S-floor-gate}
C_N^{\rm fl}=\{\Vu>0,\ F_N\le\varepsilon_N\Vu\},\qquad
\varepsilon_N>0,\quad\varepsilon_N\to0,
\end{equation}
where each diagnostic has its own $\varepsilon_N$ sequence, with the corresponding specified conditions.

\begin{proposition}\label{prop:S-alt-procedure}
Let $\alpha$, $\Pcal_0$ and the combined procedure be as above.
\begin{enumerate}[(i)]
\item With $C_N=C_N^{\rm fl}$,
\[
\limsup_N\sup_{P\in\Pcal_0}\P_P\{(C_N\cap R_{G,N})\cup(C_N^c\cap R_{S,N})\}\le\alpha.
\]
Moreover, if, along a sequence of laws in $\Pcal_0$,
\[
D_N\to0,\qquad \frac{F_N}{\varepsilon_N\sigma_N^2}\to0
\quad\text{in probability},
\]
then the rejection probability tends to $\alpha$. If both convergences hold uniformly on a subclass, then the rejection probability tends to $\alpha$ uniformly on that subclass.
\item Assume in addition (M8) and (S2) on $\Pcal_0$, with the sample-size requirements of Theorem~\ref{thm:S-check} for all sufficiently large $N$. With $C_N=C_N^{\rm ctr}$, the same bound holds, and the rejection probability tends to $\alpha$ uniformly on $\{P\in\Pcal_0:\Upsilon_N(P)\le a_N\varepsilon_N\}$ for every deterministic $a_N\to0$.
\end{enumerate}
\end{proposition}

\begin{proof}[Proof of Proposition~\ref{prop:S-alt-procedure}]

Since $F_N\ge0$, we have $A_N\ge0$, so Proposition~\ref{prop:block-safe} gives the conservative rejection bound. Also, $B_N\ge A_N$ and $\alpha^{-1/2}\ge z_{1-\alpha/2}$, so $R_{S,N}\subseteq R_{G,N}$.

\textbf{Proof of (i).} To verify condition (c) of Lemma~\ref{lem:hybrid-switch}, fix $\delta>0$. By (L)(a), $\rho_N\le\delta/4$ uniformly for all sufficiently large $N$. If $D_N\ge\delta$, then $v_N\ge3\delta/4$ by~\eqref{eq:S-selected-variance-bound}. Furthermore, Theorem~\ref{thm:contraction} and Proposition~\ref{prop:gammaD} give $\E[(\Vu^+/\sigma_N^2)^2]\le C$. By \eqref{eq:S-selected-variance-bound} and Chebyshev's inequality,
\[
\begin{aligned}
\P_P(C_N^{\rm fl})
&\le\P_P\!\left(\left|\frac{F_N}{\sigma_N^2}-v_N\right|>\frac\delta4\right)
+\P_P\!\left(\frac{\varepsilon_N\Vu^+}{\sigma_N^2}\ge\frac\delta2\right)\\
&\le C_\delta(L_N^{-1/2}+\rho_N+\varepsilon_N^2).
\end{aligned}
\]
Thus,
\begin{equation}\label{eq:S-hybrid-gate}
\sup_{P\in\Pcal_0:D_N(P)\ge\delta}\P_P(C_N^{\rm fl})\to0.
\end{equation}
If $D_N\to0$, then $\Upsilon_N\to0$ by Proposition~\ref{prop:gammaD}, and Lemma~\ref{thm:S-test} gives $(U_N-\theta_0)/A_N^{1/2}\dto\Nor(0,1)$. Using \eqref{eq:S-hybrid-gate} and the opening paragraph, we can therefore apply Lemma~\ref{lem:hybrid-switch} with $t_N=D_N$ to obtain the uniform rejection bound. Moreover, under the additional rate condition in (i), Theorem~\ref{thm:contraction} gives $\Vu/\sigma_N^2\to1$. Since $F_N/(\varepsilon_N\sigma_N^2)\to0$, we also have $\P_P(C_N^{\rm fl})\to1$. The last assertion of Lemma~\ref{lem:hybrid-switch} gives convergence of the rejection probability to $\alpha$, with the stated uniformity.

\textbf{Proof of (ii).} We apply the proof of Proposition~\ref{prop:S-block-hybrid} with $A_N=\max\{\Vu,F_N\}$ in place of $\Vu^+$. The conservative rejection bound and the inclusion $R_{S,N}\subseteq R_{G,N}$ hold as above. In addition, Lemma~\ref{thm:S-test} gives $(U_N-\theta_0)/A_N^{1/2}\dto\Nor(0,1)$ whenever $\Upsilon_N\to0$. The diagnostic is unchanged, so Theorem~\ref{thm:S-check} and Lemma~\ref{lem:hybrid-switch} with $t_N=\Upsilon_N$ give both conclusions.
\end{proof}

\subsection{The block-average test}\label{sec:S-block-average}
Alternatively, we can average the block statistics and avoid the diagnostic. The blocks use disjoint rows and columns, so their statistics are independent. An increasing number of blocks can therefore give a Gaussian average even when a single block has a non-Gaussian limit.

\begin{corollary}\label{cor:S-block-average}
Under the assumptions of Proposition~\ref{prop:block-safe},
\[
 \frac{\sqrt{K_N}(\overline U_{\rm block}-\theta_N)}{S_{K_N}}
 \dto\Nor(0,1)
\]
uniformly over the class.
\end{corollary}
\begin{proof}[Proof of Corollary~\ref{cor:S-block-average}]
Recall that $\tau_N^2=\var(U_{N,1})$. By \eqref{eq:S-block-moment}, the Lyapunov ratio of the independent centered block statistics satisfies
\[
\frac{\sum_{j=1}^{K_N}\E[(U_{N,j}-\theta_N)^4]}
     {(K_N\tau_N^2)^2}
\le\frac{C_*}{K_N}\to0.
\]
Also, the sample-variance bound~\eqref{eq:S-sample-variance} gives $S_{K_N}^2/\tau_N^2\to1$ uniformly. Lyapunov's central limit theorem and Slutsky's lemma therefore give the stated limit along every sequence of laws in the class, which proves uniformity.
\end{proof}

However, the variance of the block average $\overline U_{\rm block}$ is increased by \eqref{eq:S-block-inflation} throughout the class. Averaging an increasing number of blocks, and the resulting loss of the faster rate under degeneracy, are established properties of incomplete statistics \citep[see][Sec.~4.1.2]{leminh2023ustatisticsrowcolumn}.

\section{Proofs for the product-form applications}\label{sec:S-applications}\label{sec:S-product}
This section supplies the linearization and plug-in results used in Section~\ref{sec:product}, followed by the graph projections and application proofs.

\subsection{Preliminaries: Linearization and studentized limits}\label{sec:S-smooth}

Let $A_N,B_N$ be fixed-order RCE $U$-statistics with means $a_N,b_N$ and standard deviations $\tau_A,\tau_B$, defined by kernels $h_A,h_B$ on the same $p\times q$ subarray, in the asymptotic regime of Section~\ref{sec:framework} of the paper, with $m\ge2p$, $n\ge2q$ eventually and $\min(m,n)\to\infty$. Let $\psi$ be twice continuously differentiable on an open interval $I$ containing a compact set $K$, with $b_N\in K$ for all $N$. We study the statistic $A_N\psi(B_N)$ as an estimator of $a_N\psi(b_N)$. Its linear approximation $\ell_N$ is the first-order Taylor expansion of $A_N\psi(B_N)-a_N\psi(b_N)$ at $(a_N,b_N)$:
\[
\begin{gathered}c_1=\psi(b_N),\qquad c_2=a_N\psi'(b_N),\\ \ell_N=c_1(A_N-a_N)+c_2(B_N-b_N),\qquad \sigma_N^2=\var(\ell_N)>0.\end{gathered}
\]
Thus, $\ell_N$ is the centered $U$-statistic with the \emph{linearized kernel} $c_1h_A+c_2h_B$. To define $A_N\psi(B_N)$ and the plug-in coefficients below for every sample, we extend $\psi$ to a $C^2$ function on $\mathbb R$ that agrees with it on a fixed neighborhood of $K$. A smooth cutoff inside $I$ gives such an extension.

To estimate $\sigma_N^2$, we replace $c_1,c_2$ by $\widehat c_1=\psi(B_N)$, $\widehat c_2=A_N\psi'(B_N)$, and let $\widehat V=\Vu(\widehat c_1h_A+\widehat c_2h_B)$, computed with the same estimated coefficients in every subarray evaluation. We establish the following key lemma.

\begin{lemma}\label{lem:ratiofull}
Let $\Pcal$ be a class of laws and $P_N\in\Pcal$. Assume $\tau_B\to0$ and
\begin{equation}\label{eq:ratio-remainder}
\tau_A\tau_B=o(\sigma_N),\qquad |a_N|\tau_B^2=o(\sigma_N).
\end{equation}
\begin{enumerate}
\item[(i)] $A_N\psi(B_N)-a_N\psi(b_N)=\ell_N+o_P(\sigma_N)$.
\item[(ii)] Assume in addition that $\inf|\psi|>0$ on a neighborhood of $K$, that the linearized kernel $c_1h_A+c_2h_B$ satisfies (M4) on $\Pcal$, that $\sup_N\E[h_B^4]<\infty$, and that
\begin{equation}\label{eq:ratio-plugin}
|a_N|\,L_N^{-1}=o(\sigma_N).
\end{equation}
If $c_1h_A+c_2h_B$ satisfies $D_N\to0$, then $\widehat V/\sigma_N^2\to1$ in probability.
\item[(iii)] If the assumptions of (ii) hold and the linearized kernel $c_1h_A+c_2h_B$ satisfies $D_N\to0$, then $\widehat Z_N=\{A_N\psi(B_N)-a_N\psi(b_N)\}/(\widehat V^+)^{1/2}$, with $\widehat V^+=\max(\widehat V,0)$ and $\widehat Z_N=0$ when $\widehat V^+=0$, satisfies $\widehat Z_N\dto\Nor(0,1)$.
\end{enumerate}
All conclusions hold uniformly over $\Pcal$ if the assumptions hold uniformly over $\Pcal$, that is, along every sequence of laws in $\Pcal$.
\end{lemma}

Part (ii) of Lemma~\ref{lem:ratiofull} uses the following moment bound to replace population coefficients by their estimates. Here $\Vu(f,k)$ is the symmetric bilinear form associated with the variance estimator: $\Vu(f,k)=\{\Vu(f+k)-\Vu(f)-\Vu(k)\}/2$.

\begin{lemma}\label{lem:bilinear}
Let $f$ and $k$ be two kernels on the subarray $K_{p,q}$, with $m\ge2p$ and $n\ge2q$. Let $\Pcal$ be a class of laws on which the graph projections of $f$ satisfy (M4) with constant $C_{\mathrm{M4}}$ and $\E[k^4]\le C_{k,4}$ for every $P\in\Pcal$. Then for every $P\in\Pcal$
\[
\E[\Vu(f,k)^2]\le C\,\var(U_N(f))\,L_N^{-1},\qquad \E[\Vu(k)^2]\le C\,L_N^{-2},
\]
where $C$ depends only on $(p,q,C_{\mathrm{M4}},C_{k,4})$.
\end{lemma}

\begin{proof}[Proof of Lemma~\ref{lem:bilinear}]
By translation invariance, we may center both kernels. Put $Q(A,A')=\sum_\ell b_\ell(A,A')$, where $b_\ell$ are the coefficients in Lemma~\ref{lem:S-coefficients}. These coefficients depend on $m$, $n$, $p$ and $q$, but not on the kernel. For every kernel, Proposition~\ref{prop:pathwise}(ii) and Lemma~\ref{lem:S-coefficients} give $\Vu=\sum_{A,A'}Q(A,A')X_AX_{A'}$, where $X_A$ are the grouped projections of that kernel. Moreover, swapping $W$ and $W'$ in the definition of $N_{a,b}(A,A')$ in the proof of Lemma~\ref{lem:S-coefficients} gives $N_{a,b}(A,A')=N_{a,b}(A',A)$, so $Q$ is symmetric. Thus, polarization gives
\[
\Vu(f,k)=\sum_{A,A'}Q(A,A')X_A^fX_{A'}^k,
\qquad |Q(A,A')|\le Cc_Ac_{A'},
\]
where $X_A^f,X_A^k$ are the grouped projections in~\eqref{eq:grouped-projections}.  By \eqref{eq:S-level-coefficients}, $|b_\ell(A,A')|\le(1+CL_N^{-1})c_Ac_{A'}$ for each $\ell$. Summing over the finitely many $\ell$ gives the coefficient bound. 

Next, expanding these projections and the square gives
\[
\E[\Vu(f,k)^2]
\le C\sum_{\Lambda_1,\dots,\Lambda_4}
\Big(\prod_{j=1}^4c_{G_j}\Big)
\left|\E[q^f_{\Lambda_1}q^f_{\Lambda_2}q^k_{\Lambda_3}q^k_{\Lambda_4}]\right|,
\]
where the sum is over ordered quadruples of nonempty graph copies, and $G_j$ is the class of $\Lambda_j$.

Let $\varsigma_G^f$ be the standard deviation of the graph contribution $T_G$ for the kernel $f$. By (M4) and~\eqref{eq:copy-scale},
\[
c_G\|q_G^f\|_4\le C\frac{\varsigma_G^f}{\sqrt{M_N(G)}}.
\]
Now, each projection of $k$ is a finite signed sum of conditional expectations, so $\|q_G^k\|_4\le C$. Since $c_G\asymp M_N(G)^{-1}$ and $M_N(G)^{-1}\le CL_N^{-1}$ for every nonempty $G$, we have
\[
c_G\|q_G^k\|_4\le\frac{C}{M_N(G)}
\le\frac{CL_N^{-1/2}}{\sqrt{M_N(G)}}.
\]

Finally, fix a pattern among copies of $G_1,\dots,G_4$, with $v_r,v_c,r_{\mathrm{tot}},c_{\mathrm{tot}}$ defined as in Lemma~\ref{lem:auto-pattern}(i). If the mixed expectation is nonzero, Lemma~\ref{lem:config} gives $v_r\le r_{\mathrm{tot}}/2$ and $v_c\le c_{\mathrm{tot}}/2$. Thus, H\"older's inequality and the counting argument of Lemma~\ref{lem:auto-pattern}(i), applied with the two displayed norm bounds, bound its contribution by
\[
CL_N^{-1}m^{v_r-r_{\mathrm{tot}}/2}n^{v_c-c_{\mathrm{tot}}/2}
\varsigma_{G_1}^f\varsigma_{G_2}^f
\le CL_N^{-1}\varsigma_{G_1}^f\varsigma_{G_2}^f.
\]
Summing over the finitely many patterns and graph classes, and using Lemma~\ref{lem:auto-pattern}(ii), gives
\[
\E[\Vu(f,k)^2]
\le CL_N^{-1}\Big(\sum_G\varsigma_G^f\Big)^2
\le C\,\var(U_N(f))\,L_N^{-1}.
\]
The same calculation with four projections of $k$ gives $\E[\Vu(k)^2]\le CL_N^{-2}$.
\end{proof}

\begin{proof}[Proof of Lemma~\ref{lem:ratiofull}]
We suppress $N$ in $A$, $B$, $a$, $b$, $\sigma$ and $\ell$.

\textbf{Proof of (i).} Chebyshev's inequality gives $A-a=O_P(\tau_A)$ and $B-b=O_P(\tau_B)$. Since $\tau_B\to0$, $B$ lies in a fixed compact neighborhood of $K$ inside $I$ with probability tending to one. On this neighborhood, the first two derivatives of $\psi$ are bounded. Taylor's formula and~\eqref{eq:ratio-remainder} therefore give
\[
\begin{aligned}
A\psi(B)-a\psi(b)-\ell
&=(A-a)\{\psi(B)-\psi(b)\} +a\{\psi(B)-\psi(b)-\psi'(b)(B-b)\}\\
&=O_P(\tau_A\tau_B+|a|\tau_B^2)=o_P(\sigma).
\end{aligned}
\]

\textbf{Proof of (ii).} Put $h_\ell=c_1h_A+c_2h_B$, so $\var(U_N(h_\ell))=\sigma^2$. By (M4), $D_N\to0$ and~\eqref{eq:variance-D-bound}, $\Vu(h_\ell)/\sigma^2\to1$. To compare $\widehat V$ with $\Vu(h_\ell)$, write
\[
\widehat c_1h_A+\widehat c_2h_B=rh_\ell+sh_B,
\qquad r=\frac{\widehat c_1}{c_1},\quad s=\widehat c_2-rc_2.
\]
Since $|c_1|=|\psi(b)|$ is bounded away from zero and $\ell=O_P(\sigma)$, the local derivative bounds give
\begin{align*}
r-1&=\frac{\psi(B)-\psi(b)}{\psi(b)}=O_P(\tau_B),\\
A-a&=\frac{\ell-c_2(B-b)}{c_1}=O_P(\sigma+|a|\tau_B).
\end{align*}
Using the boundedness of $\psi''$, we obtain
\begin{align*}
s&=(A-a)\psi'(B)+a\{\psi'(B)-\psi'(b)\}-a(r-1)\psi'(b)\\
&=O_P(\sigma+|a|\tau_B).
\end{align*}
Under (ii), $\E[h_B^2]$ is bounded, so $\tau_B\le CL_N^{-1/2}$, and condition~\eqref{eq:ratio-plugin} gives
\[
sL_N^{-1/2}
=O_P(\sigma L_N^{-1/2}+|a|L_N^{-1})
=o_P(\sigma).
\]
Lemma~\ref{lem:bilinear} gives $\Vu(h_\ell,h_B)=O_P(\sigma L_N^{-1/2})$ and $\Vu(h_B)=O_P(L_N^{-1})$. Thus, by bilinearity,
\begin{align*}
\widehat V-\Vu(h_\ell)
&=(r^2-1)\Vu(h_\ell)+2rs\Vu(h_\ell,h_B)+s^2\Vu(h_B)\\
&=O_P(\tau_B\sigma^2)+o_P(\sigma^2)
=o_P(\sigma^2),
\end{align*}
which proves (ii).

\textbf{Proof of (iii).} By (ii), $\widehat V^+/\sigma^2\to1$ and $\P\{\widehat V^+=0\}\to0$. Moreover, Proposition~\ref{prop:gammaD} and Theorem~\ref{thm:contraction} give $\ell/\sigma\dto\Nor(0,1)$. Using (i), Slutsky's lemma then gives the claimed limit. All estimates use only the common bounds and rates in the assumptions, so the same argument gives the uniform conclusions.
\end{proof}

\subsection{Poisson product-form projections and inference}\label{sec:S-poisson-projections}
This section treats the Poisson graphon model of Section~\ref{sec:dw}. We bound the graph projections of $h_\gamma$, then prove Corollary~\ref{cor:dwnull} on the null limit and Corollary~\ref{cor:dwuniform} on uniform inference for $\gamma$. The same bounds give the variance comparisons used for the test in Section~\ref{sec:dw-uniform}.

Recall the linear approximation $\ell_{N,P}$ to $\widehat\gamma_U-\gamma_P$ in \eqref{eq:gamma-linearization}, its variance $\sigma_{N,P}^2=\var_P(\ell_{N,P})$ and the corresponding unnormalized kernel $\widetilde h_{\ell,P}=h_\gamma-4\gamma_P\lambda^3Y_{11}$, with graph projections $p_G^{\widetilde\ell}$ as in Corollary~\ref{cor:dwuniform}. We retain $\Delta=\wbar-f\otimes g$ and $\gamma=\|\Delta\|_2^2$ from Section~\ref{sec:dw}.

\paragraph{\texorpdfstring{Projections of $h_\gamma$}{Projections of h-gamma}.}\label{sec:S-ratio-projections}

We express the edgeless projections of $h_\gamma$ as derivatives of the kernel expectation. Propositions~\ref{prop:latent} and~\ref{prop:S-two-group} both use this representation. For row and column laws $\mu,\nu$ on $[0,1]$, put
\[
T=\iint\wbar\,\d\mu\,\d\nu,\qquad
f_\nu(x)=\int\wbar(x,y)\,\d\nu(y),\qquad
g_\mu(y)=\int\wbar(x,y)\,\d\mu(x).
\]
Using conditional independence and $\E[Y_{11}(Y_{11}-1)\mid\xi_1,\eta_1]=\lambda^2\wbar(\xi_1,\eta_1)^2$, we obtain from~\eqref{eq:gamma-kernel}
\begin{align*}
\lambda^{-4}\E_{\mu,\nu}[h_\gamma]
=\mathcal E(\mu,\nu)
&:=T^2\iint\wbar^2-2T\iint\wbar\,f_\nu\otimes g_\mu
+\int f_\nu^2\,\d\mu\int g_\mu^2\,\d\nu\\
&=\|T\wbar-f_\nu\otimes g_\mu\|_{L^2(\mu\otimes\nu)}^2,
\end{align*}
where the double integrals are with respect to $\mu\otimes\nu$.

Let $G$ be edgeless with $r$ rows and $c$ columns, and fix $\mathbf x=(x_1,\dots,x_r)$ and $\mathbf y=(y_1,\dots,y_c)$. Set
\[
\mu_s=\mathrm{Leb}+\sum_{a=1}^rs_a(\delta_{x_a}-\mathrm{Leb}),
\qquad
\nu_t=\mathrm{Leb}+\sum_{b=1}^ct_b(\delta_{y_b}-\mathrm{Leb}),
\]
where $\mathrm{Leb}$ is the Lebesgue measure on $[0,1]$, the law $\mathcal U[0,1]$ of the AHK variables. The resulting expression for $\mathcal E(\mu_s,\nu_t)$ is polynomial in $(s,t)$. Differentiating once in each $s_a$ and $t_b$ replaces the measures at distinct row and column positions by the corresponding centered directions. Expanding these directions gives the ordinary Hoeffding projection for each of the $(3)_r(3)_c$ placements. Thus, the normalization in Lemma~\ref{lem:decomp} gives
\begin{equation}\label{eq:S-latent-projection}
\left.\partial_{s_1}\cdots\partial_{s_r}\partial_{t_1}\cdots\partial_{t_c}
\mathcal E(\mu_s,\nu_t)\right|_{s=t=0}
=\lambda^{-4}p_G(\mathbf x,\mathbf y).
\end{equation}

\begin{proposition}\label{prop:latent}
Let $\wbar$ be bounded by $M$. For every \emph{edgeless} graph $G$, the projection $p_G$ of the kernel $h_\gamma$ in~\eqref{eq:gamma-kernel} is a finite sum of multilinear integrals, each containing at least two factors $\Delta$. Therefore, $\E[p_G^2]\le C(M,\lambda)\,\|\Delta\|_\infty^2\,\gamma(\wbar)$.
\end{proposition}

\begin{proof}[Proof of Proposition~\ref{prop:latent}]
We expand the squared expression defining $\mathcal E$ in powers of $\Delta$. Write $w_0=f\otimes g$, keeping $f,g$ fixed as $\mu,\nu$ vary. Denote the preceding integrals with $\wbar$ replaced by $w$ by $T(w)$, $f_{\nu,w}$ and $g_{\mu,w}$. Since $T(w_0)w_0=f_{\nu,w_0}\otimes g_{\mu,w_0}$, substituting $\wbar=w_0+\Delta$ gives
\begin{align*}
T\wbar-f_\nu\otimes g_\mu
&=T(w_0)\Delta+T(\Delta)w_0\\
&\quad-f_{\nu,w_0}\otimes g_{\mu,\Delta}
      -f_{\nu,\Delta}\otimes g_{\mu,w_0}\\
&\quad+T(\Delta)\Delta-f_{\nu,\Delta}\otimes g_{\mu,\Delta}.
\end{align*}
Every term contains at least one factor $\Delta$, possibly integrated. Squaring gives at least two such factors, and differentiating in $(s,t)$ changes only the integration measures. Hence, by~\eqref{eq:S-latent-projection}, $p_G$ is a finite sum of terms of the form
\[
\mathcal I(z)=\int B(z,u)\prod_{\ell=1}^k
\Delta(x_{a_\ell},y_{b_\ell})\,\d u,
\qquad 2\le k\le4,
\]
where $z$ contains the fixed variables, $u$ the integrated variables, and $\|B\|_\infty\le C(M,\lambda)$. All remaining integrations are with respect to product Lebesgue measures. Jensen's inequality and Fubini's theorem give
\begin{align*}
\|\mathcal I\|_2^2
&\le\int |B(z,u)|^2\prod_{\ell=1}^k
       |\Delta(x_{a_\ell},y_{b_\ell})|^2\,\d u\,\d z\\
&\le C(M,\lambda)\|\Delta\|_\infty^{2(k-1)}
       \int|\Delta(x_{a_1},y_{b_1})|^2\,\d u\,\d z\\
&=C(M,\lambda)\|\Delta\|_\infty^{2(k-1)}\|\Delta\|_2^2\\
&\le C(M,\lambda)\|\Delta\|_\infty^2\gamma(\wbar).
\end{align*}
For the equality, the variables absent from the retained $\Delta$ factor integrate to one. The last inequality uses $2\le k\le4$ and $\|\Delta\|_\infty\le M+M^2$. The triangle inequality for the finitely many terms gives
\[
\|p_G\|_2\le C(M,\lambda)\|\Delta\|_\infty\|\Delta\|_2.
\]
Squaring proves the bound.
\end{proof}

Together with the preceding proposition, the following bound identifies the principal support graph at product form and gives the variance comparisons used in Corollary~\ref{cor:dwuniform}.

\begin{proposition}\label{prop:edge}
Write $\lambda_{11}=\lambda\wbar(\xi_1,\eta_1)$ and $e_{11}=Y_{11}-\lambda_{11}$, and let
\[
a(x)=\int\Delta(x,y)g(y)\,\d y,\quad b(y)=\int\Delta(x,y)f(x)\,\d x,\quad c_0=\langle\Delta,f\otimes g\rangle .
\]
Set $B_\gamma(x,y)=\Delta(x,y)-a(x)-b(y)+c_0+\gamma$. The $K_{1,1}$ projection of $h_\gamma$ is
\[
p_{K_{1,1}}=\lambda^2(e_{11}^2-e_{11}-\lambda_{11})+2\lambda^3e_{11}\,B_\gamma(\xi_1,\eta_1).
\]
Its squared norm satisfies
\[
\E[p_{K_{1,1}}^2]=2\lambda^6\iint\wbar^2+4\lambda^7\E[\wbar\,B_\gamma^2]\ \ge\ 2\lambda^6 .
\]
For $\widetilde h_{\ell,P}$, both formulas hold with $B_\gamma$ replaced by $B_\gamma-2\gamma$, and so does the lower bound.
\end{proposition}

\begin{proof}[Proof of Proposition~\ref{prop:edge}]
Put $x=\xi_1$ and $y=\eta_1$. To obtain $p_{K_{1,1}}$, we integrate out the other variables and subtract the conditional mean given $x,y$, then sum over the nine possible positions of the conditioned cell, as in Section~\ref{app:projection-notation}. Terms not containing this cell vanish. Set $A=e_{11}^2-e_{11}-\lambda_{11}$. Then,
\begin{align*}
\E[Y_{11}(Y_{11}-1)\mid x,y]&=\lambda_{11}^2,\\
Y_{11}(Y_{11}-1)-\lambda_{11}^2&=A+2\lambda_{11} e_{11}.
\end{align*}
Write
\[
\begin{gathered}
S=\iint\wbar^2,\qquad \bar c_0=\langle\wbar,f\otimes g\rangle,\\
\bar a(x)=\int\wbar(x,y)g(y)\,\d y,\qquad
\bar b(y)=\int\wbar(x,y)f(x)\,\d x.
\end{gathered}
\]
Since $\wbar=f\otimes g+\Delta$, we have
\[
\bar a=G_2f+a,\qquad \bar b=F_2g+b,\qquad
S-\bar c_0=\langle\wbar,\Delta\rangle=c_0+\gamma,
\]
with $F_2,G_2$ as in Section~\ref{sec:dw}.

Next, in the first monomial of \eqref{eq:gamma-kernel}, the conditioned cell can occur in the factorial term or in either of the two linear factors. The second monomial has four linear factors. The third has two cells in row $1$ and two in column $3$. Summing their contributions and using $\lambda_{11}=\lambda\{f(x)g(y)+\Delta(x,y)\}$ gives
\begin{align*}
p_{K_{1,1}}
&=\lambda^2(A+2\lambda_{11} e_{11})+2\lambda^3e_{11}S\\
&\quad-2\lambda^3e_{11}\{\bar a(x)+f(x)g(y)+\bar b(y)+\bar c_0\}\\
&\quad+2\lambda^3e_{11}\{G_2f(x)+F_2g(y)\}\\
&=\lambda^2A+2\lambda^3e_{11}B_\gamma(x,y).
\end{align*}

Now, conditional on $x,y$, the centered Poisson moments are
\[
\begin{gathered}
\E[e_{11}\mid x,y]=0,\qquad
\E[e_{11}^2\mid x,y]=\E[e_{11}^3\mid x,y]=\lambda_{11},\\
\E[e_{11}^4\mid x,y]=\lambda_{11}+3\lambda_{11}^2.
\end{gathered}
\]
Thus,
\begin{align*}
\E[Ae_{11}\mid x,y]&=\lambda_{11}-\lambda_{11}=0,\\
\E[A^2\mid x,y]
&=(\lambda_{11}+3\lambda_{11}^2)-2\lambda_{11}+(1-2\lambda_{11})\lambda_{11}+\lambda_{11}^2=2\lambda_{11}^2.
\end{align*}
Expanding the square and conditioning on $x,y$ now gives
\begin{align*}
\E[p_{K_{1,1}}^2]
&=2\lambda^4\E[\lambda_{11}^2]+4\lambda^6\E[\lambda_{11} B_\gamma(x,y)^2]\\
&=2\lambda^6\iint\wbar^2+4\lambda^7\E[\wbar B_\gamma^2]
\ge 2\lambda^6,
\end{align*}
where Jensen's inequality gives $\iint\wbar^2\ge(\iint\wbar)^2=1$.

Finally, subtracting the linearization term $4\gamma\lambda^3Y_{11}$ changes the projection by $-4\gamma\lambda^3e_{11}$. It therefore replaces $B_\gamma$ by $B_\gamma-2\gamma$ in the projection and in the same moment calculation, with the same lower bound.
\end{proof}

The next lemma combines Proposition~\ref{prop:latent} with the variance weights to bound the disconnected variance of the linear approximation, used to bound $D_N$ in Corollary~\ref{cor:dwuniform} and in the test of Section~\ref{sec:dw-uniform}.

\begin{lemma}\label{lem:S-poisson-disconnected}
Let $\Pcal_N$ be classes of Poisson graphon laws with fixed $\lambda>0$ and $\|\wbar_P\|_\infty\le M$ for every $N$ and $P\in\Pcal_N$, and let $m,n\asymp N$. Then the disconnected variance of $\ell_{N,P}$ in \eqref{eq:gamma-linearization} satisfies
\[
\sigma_{\mathrm{disc},N,P}^2\le C(\gamma_P N^{-2}+N^{-3}),
\]
where $C$ depends only on $\lambda$, $M$ and the constants in $m,n\asymp N$.
\end{lemma}

\begin{proof}[Proof of Lemma~\ref{lem:S-poisson-disconnected}]
Fix $N$ and $P\in\Pcal_N$, and omit the subscript $P$. Proposition~\ref{prop:latent} bounds the squared norms of the edgeless projections of $h_\gamma$ by $C\|\Delta\|_\infty^2\gamma\le C\gamma$, since $\|\Delta\|_\infty\le C$ under the bound on $\wbar$. Moreover, the linearization term $-4\gamma\lambda^3Y_{11}$ has squared projection norms bounded by $C\gamma^2$. Since every disconnected graph with two vertices is edgeless and $\gamma$ is bounded,
\[
\E[(p_G^{\widetilde\ell})^2]\le C(\gamma+\gamma^2)\le C\gamma
\qquad(G\ \disc,\ d(G)=2).
\]
Because the mean graphon $\lambda\wbar$ is uniformly bounded, we also have $\E[\widetilde h_\ell^2]\le C$, so all remaining squared projection norms are bounded. There are finitely many graph classes, and their weights in \eqref{eq:graph-variance-weight} satisfy $w_G\le CN^{-d(G)}$. Thus,
\[
\begin{aligned}
\sigma_{\mathrm{disc}}^2
&=\lambda^{-8}\left\{
\sum_{\substack{G\ \disc\\d(G)=2}}w_G\E[(p_G^{\widetilde\ell})^2]
+\sum_{\substack{G\ \disc\\d(G)\ge3}}w_G\E[(p_G^{\widetilde\ell})^2]\right\}\\
&\le C\gamma N^{-2}+CN^{-3}.
\end{aligned}
\]
\end{proof}

\paragraph{The null limit.} We now prove Corollary~\ref{cor:dwnull} from Propositions~\ref{prop:latent} and~\ref{prop:edge}.

\begin{proof}[Proof of Corollary~\ref{cor:dwnull}]
At product form \(\Delta=0\), every edgeless projection vanishes by Proposition~\ref{prop:latent}, at every order. In particular, all degree-one projections vanish. At degree two, $K_{1,1}$ is the only graph containing an edge. Proposition~\ref{prop:edge} gives
\[
p_{K_{1,1}}
=\lambda^2(e_{11}^2-e_{11}-\lambda_{11})
=\lambda^2(e_{11}^2-Y_{11}),
\]
which is nonzero because its conditional variance is positive whenever \(\lambda_{11}>0\). Thus, the principal degree is two, with this graph as the unique principal support graph.

For this graph, \(r=c=1\), \(|\operatorname{Aut}(G)|=1\), and~\eqref{eq:graph-variance-weight} gives
\[
\var(T_{K_{1,1}})=\frac{\E[p_{K_{1,1}}^2]}{mn}.
\]
Under product form, \(\wbar=f\otimes g\) and the second term in Proposition~\ref{prop:edge} vanishes, so
\[
\E[p_{K_{1,1}}^2]
=2\lambda^6\iint f(x)^2g(y)^2\,\d x\,\d y
=2\lambda^6F_2G_2.
\]
The sum of the higher-degree contributions is \(o((mn)^{-1})\), giving
\[
\var(U)\sim\frac{2\lambda^6F_2G_2}{mn}.
\]
The unique principal support graph is connected, so Lemma~\ref{lem:fixedclt} gives the Gaussian limit.
\end{proof}

\paragraph{\texorpdfstring{Uniform inference for $\gamma$}{Uniform inference for gamma}.} We now prove Corollary~\ref{cor:dwuniform} by combining Proposition~\ref{prop:edge} and Lemma~\ref{lem:S-poisson-disconnected} with Lemma~\ref{lem:ratiofull}.

\begin{proof}[Proof of Corollary~\ref{cor:dwuniform}]
We verify the remainder and variance-estimation conditions of Lemma~\ref{lem:ratiofull}. The $K_{1,1}$ lower bound applies in both regimes. We then check the remaining conditions separately under (a) and (b).

Fix $N$ and $P\in\Pcal_N$, and omit the subscript $P$. All constants and rates below are uniform over $P\in\Pcal_N$. In Lemma~\ref{lem:ratiofull}, take $A_N=U$, $B_N=\bar Y$ and $\psi(b)=b^{-4}$, so $a_N=\lambda^4\gamma$, the linear term is $\ell_N$ in \eqref{eq:gamma-linearization}, and $\widehat V=\Vu^{(\mathrm{lin})}$, since $\widehat c_1=\bar Y^{-4}$ and $\widehat c_2=-4U\bar Y^{-5}$.

Since the mean graphon $\lambda\wbar$ is uniformly bounded, we have $\E[h_A^2]\le C$ and $\E[h_B^4]\le C$ for the kernels of $U$ and $\bar Y$. By the variance decomposition of $\bar Y$ in Example~\ref{ex:mean},
\[
\tau_B^2\le C\{m^{-1}+n^{-1}+(mn)^{-1}\}\le CN^{-1}.
\]
We set the ratio $\widehat\gamma_U$ and the estimated coefficients $\widehat c_1,\widehat c_2$ to zero when $\bar Y=0$. Chebyshev's inequality gives $\P(\bar Y=0)\le\tau_B^2/\lambda^2=O(N^{-1})$.

By \eqref{eq:gamma-linearization} and $L_N\asymp N$,
\begin{align*}
\tau_A&\le\lambda^4\sigma_N+4\gamma\lambda^3\tau_B,\\
\frac{\tau_A\tau_B+|a_N|\tau_B^2+|a_N|L_N^{-1}}{\sigma_N}
&\le C\left(N^{-1/2}+\frac{\gamma}{N\sigma_N}\right).
\end{align*}
Thus, both \eqref{eq:ratio-remainder} and \eqref{eq:ratio-plugin} follow if
\[
\frac{\gamma}{N\sigma_N}\longrightarrow0.
\]
Proposition~\ref{prop:edge} gives
\begin{equation}\label{eq:C2-edge-floor}
\sigma_N^2\ge\frac{\lambda^{-8}}{mn}\E[(p^{\widetilde\ell}_{K_{1,1}})^2]\ge cN^{-2}.
\end{equation}

\textbf{Proof under (a).} By \eqref{eq:C2-edge-floor} and Lemma~\ref{lem:S-poisson-disconnected},
\[
\frac{\gamma}{N\sigma_N}\le C\sup_{P\in\Pcal_N}\gamma_P\to0,
\qquad
D_N\le C(\sup_{P\in\Pcal_N}\gamma_P+N^{-1})\to0.
\]

\textbf{Proof under (b).} The first-order projections and the $K_{1,1}$ projection are orthogonal. Thus, (b) and Proposition~\ref{prop:edge} give
\[
\sigma_N^2\ge c'\left(N^{-2}+\frac{\gamma^2}{N}\right).
\]
Combining this with Lemma~\ref{lem:S-poisson-disconnected}, we obtain
\begin{align*}
\frac{\gamma}{N\sigma_N}
&\le\frac{C\gamma}{\sqrt{1+N\gamma^2}}\le CN^{-1/2}\to0,\\
D_N&\le\frac{C(\gamma+N^{-1})}{1+N\gamma^2}
\le C\left(\frac1{2\sqrt N}+\frac1N\right)\to0,
\end{align*}
where $2\sqrt N\gamma\le1+N\gamma^2$ gives the last bound.

In both regimes, these bounds and the assumed (M4) verify all conditions of Lemma~\ref{lem:ratiofull}. Parts (ii) and (iii) give the stated ratio-consistency and studentized Gaussian limit uniformly over $P\in\Pcal_N$.
\end{proof}

\subsection{The two-group model}\label{sec:S-two-group}
For the two-group model of Section~\ref{sec:planted}, define the degree-two variance coefficients by
\[
A_G:=\lambda^{-8}\frac{\E[(p_G^{\widetilde\ell})^2]}{|\operatorname{Aut}(G)|},
\qquad
\var(T_G^\ell)=\frac{A_G}{(m)_{r(G)}(n)_{c(G)}}.
\]

\begin{proposition}\label{prop:S-two-group}
Consider the two-group model of Section~\ref{sec:planted} at a fixed $\varepsilon\in(0,1]$.
\begin{enumerate}[(i)]
\item The first-order projections of $\bar Y$ and $U$ vanish.
\item The degree-two variance coefficients are
\[
\begin{gathered}
A_{K_{2,0}}=A_{K_{0,2}}=2\varepsilon^4,\qquad A_{K_{1,0}\sqcup K_{0,1}}=16\varepsilon^6,\\
A_{K_{1,1}}=\frac{2(1+\varepsilon^2)}{\lambda^2}+\frac{4\varepsilon^2(1-\varepsilon^2)}{\lambda}.
\end{gathered}
\]
\item The statistics $\ell_N/\sigma_N$ and $(\widehat\gamma_U-\varepsilon^2)/\sigma_N$ converge to the same non-Gaussian limit.
\end{enumerate}
\end{proposition}

\begin{proof}[Proof of Proposition~\ref{prop:S-two-group}]
\textbf{Proof of (i).} To compute the edgeless projections, we use $\mathcal E(\mu,\nu)$ defined before Proposition~\ref{prop:latent}. For a row law $\mu$ and a column law $\nu$, put $u=\int s\,\d\mu$, $v=\int s\,\d\nu$. Then,
\[
T=1+\varepsilon uv,\qquad
f_\nu(x)=1+\varepsilon v s(x),\qquad
g_\mu(y)=1+\varepsilon u s(y),
\]
and, since $s^2=1$,
\begin{align*}
T\wbar-f_\nu\otimes g_\mu
&=\varepsilon(s(x)-u)(s(y)-v),\\
\mathcal E(\mu,\nu)
&=\varepsilon^2\int(s-u)^2\,\d\mu\int(s-v)^2\,\d\nu\\
&=\varepsilon^2(1-u^2)(1-v^2).
\end{align*}
At the uniform laws, $u=v=0$. Neither $T$ nor $\mathcal E$ has a linear term there, so the first-order projections of $\bar Y$ and $U$ vanish.

\textbf{Proof of (ii).} Since $\gamma=\varepsilon^2$, the linearized kernel satisfies
\[
\lambda^{-4}\E_{\mu,\nu}[\widetilde h_\ell]
=\mathcal E(\mu,\nu)-4\varepsilon^2T.
\]
For the perturbations in \eqref{eq:S-latent-projection}, $\partial_{s_a}u=s(x_a)$ and $\partial_{t_b}v=s(y_b)$. Differentiating at $u=v=0$ therefore gives
\begin{align*}
p_{K_{2,0}}^{\widetilde\ell}(x_1,x_2)&=-2\lambda^4\varepsilon^2s(x_1)s(x_2),\\
p_{K_{0,2}}^{\widetilde\ell}(y_1,y_2)&=-2\lambda^4\varepsilon^2s(y_1)s(y_2),\\
p_{K_{1,0}\sqcup K_{0,1}}^{\widetilde\ell}(x,y)&=-4\lambda^4\varepsilon^3s(x)s(y).
\end{align*}
Their automorphism factors are $2$, $2$ and $1$, respectively. Squaring and using $s^2=1$ gives the first three coefficients in (ii).

Finally, for $K_{1,1}$, we have $f=g=1$, $a=b=c_0=0$ and $B_\gamma(x,y)-2\gamma=\varepsilon s(x)s(y)-\varepsilon^2$. Since $s(\xi_1)s(\eta_1)$ takes the values $1$ and $-1$ with equal probabilities, Proposition~\ref{prop:edge} gives
\begin{align*}
A_{K_{1,1}}
&=\frac{2(1+\varepsilon^2)}{\lambda^2} +\frac{2\varepsilon^2}{\lambda}\{(1+\varepsilon)(1-\varepsilon)^2
+(1-\varepsilon)(1+\varepsilon)^2\}\\
&=\frac{2(1+\varepsilon^2)}{\lambda^2}
+\frac{4\varepsilon^2(1-\varepsilon^2)}{\lambda}.
\end{align*}

\textbf{Proof of (iii).} By (i) and (ii), $\tau_A,\tau_B=O(N^{-1})$ and $\sigma_N\asymp N^{-1}$. Thus, $\tau_A\tau_B+|a_N|\tau_B^2=O(N^{-2})=o(\sigma_N)$, and Lemma~\ref{lem:ratiofull}(i), with $A_N=U$, $B_N=\bar Y$ and $\psi(b)=b^{-4}$, gives
\[
\widehat\gamma_U-\varepsilon^2
=\ell_N+o_P(\sigma_N).
\]
By Theorem~\ref{thm:chaos} and Lemma~\ref{lem:chaosnotgauss}, the positive coefficients $A_G$ of the disconnected graphs give a non-Gaussian limit for $\ell_N/\sigma_N$. Consequently, the same limit holds for $(\widehat\gamma_U-\varepsilon^2)/\sigma_N$.
\end{proof}

\subsection{The degree-four result for general cell distributions}\label{sec:S-degree-four}

\begin{figure}[!ht]
\centering
\begin{tikzpicture}[x=0.6cm, y=-0.6cm,
  dot/.style={circle, fill=black, minimum size=3.5pt, inner sep=0pt},
  path/.style={thick}]
\begin{scope}
\foreach \c/\r in {1/1,2/1,1/2,2/2,3/3,4/3,3/4,4/4} \fill[black!18] (\c-0.5,\r-0.5) rectangle ++(1,1);
\foreach \i in {1.5,2.5,3.5} {\draw[black!45] (\i,0.5) -- (\i,4.5); \draw[black!45] (0.5,\i) -- (4.5,\i);}
\draw (0.5,0.5) rectangle (4.5,4.5);
\draw[path] (1,1) -- (2,1) -- (2,2) -- (1,2) -- cycle;
\draw[path] (3,3) -- (4,3) -- (4,4) -- (3,4) -- cycle;
\foreach \c/\r in {1/1,2/1,1/2,2/2,3/3,4/3,3/4,4/4} \node[dot] at (\c,\r) {};
\node[anchor=base] at (2.5,5.5) {(a)};
\end{scope}
\begin{scope}[xshift=4.2cm]
\foreach \c/\r in {1/1,1/2,2/2,2/3,3/3,3/4,4/4,4/1} \fill[black!18] (\c-0.5,\r-0.5) rectangle ++(1,1);
\foreach \i in {1.5,2.5,3.5} {\draw[black!45] (\i,0.5) -- (\i,4.5); \draw[black!45] (0.5,\i) -- (4.5,\i);}
\draw (0.5,0.5) rectangle (4.5,4.5);
\draw[path] (1,1) -- (1,2) -- (2,2) -- (2,3) -- (3,3) -- (3,4) -- (4,4) -- (4,1) -- cycle;
\foreach \c/\r in {1/1,1/2,2/2,2/3,3/3,3/4,4/4,4/1} \node[dot] at (\c,\r) {};
\node[anchor=base] at (2.5,5.5) {(b)};
\end{scope}
\end{tikzpicture}
\caption{Products of cells in the kernel $h_\Psi$ in a subarray of four rows and four columns. The cells in each product are shaded and joined in cycle order. (a) A copy of $C_4\sqcup C_4$, one of the $18$ possible edge sets $P$. (b) A copy of $C_8$ through all four rows and columns, one of the $72$ possible edge sets $Q$. Each product uses exactly two cells in each row and each column. The $2\times2$ blocks of (a) never occur in a copy of $C_8$.}\label{fig:psi-patterns}
\end{figure}
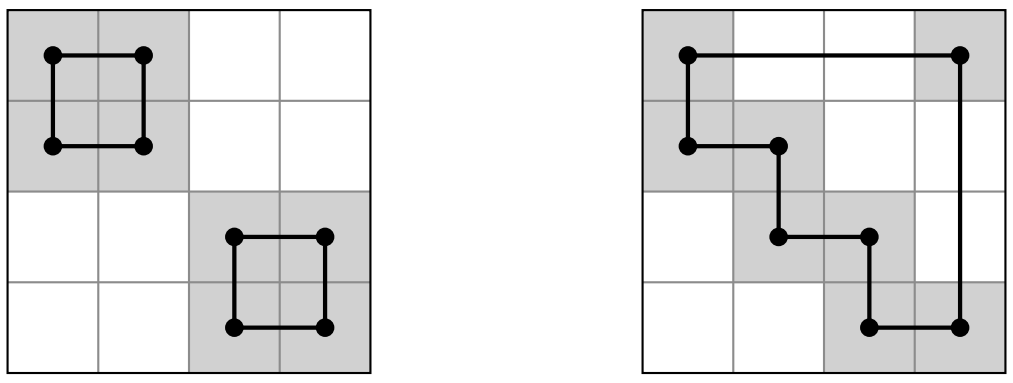

The next proposition identifies the graph projections of $h_\Psi$ used in Corollary~\ref{cor:d3}.

\begin{proposition}\label{prop:d3}
Assume that the cells have conditional means $\E[Y_{ij}\mid\xi,\eta]=\lambda f(\xi_i)g(\eta_j)$ and arbitrary conditional distributions. Write $e_{ij}=Y_{ij}-\lambda f(\xi_i)g(\eta_j)$. Then every projection $p_G$ of $h_\Psi$ with $d(G)\le3$ vanishes, and among the graphs with $d(G)=4$, the only projection that may be nonzero is that on $C_4$,
\[
p_{C_4}=c\,\lambda^4F_2^2G_2^2\prod_{e\in C_4}e_e,\qquad c>0,
\]
which is nonzero whenever $\lambda F_2G_2\ne0$ and $\E[\prod_{e\in C_4}\var(Y_e\mid\xi,\eta)]>0$.
\end{proposition}

\begin{proof}[Proof of Proposition~\ref{prop:d3}]
Every product in $h_\Psi$ uses two cells in each row and each column. Substituting $Y_{ij}=\lambda f(\xi_i)g(\eta_j)+e_{ij}$ therefore gives
\begin{align*}
h_\Psi
&=\sum_A c(A)\lambda^{8-|A|}\prod_{e\in A}e_e \times\prod_{i\in V_1(K_{4,4})}f(\xi_i)^{2-\deg_A(i)}
\prod_{j\in V_2(K_{4,4})}g(\eta_j)^{2-\deg_A(j)},\\
c(A)&=\Pr\{A\subseteq P\}-\Pr\{A\subseteq Q\}.
\end{align*}
Here, $\Pr\{A\subseteq P\}$ denotes the proportion of the $18$ edge sets $P$ in the definition of $h_\Psi$ that contain $A$, and $\Pr\{A\subseteq Q\}$ the proportion of the $72$ edge sets $Q$ that contain $A$. The sum is over sets of cells $A$ with $\deg_A(v)\le2$ at every vertex, where $\deg_A(v)$ counts the cells of $A$ incident to $v$.

Conditional independence and $\E[e_e\mid\xi,\eta]=0$ eliminate terms with $A\setminus E(G)\ne\varnothing$ when projecting onto $G$. If $e\in E(G)\setminus A$, the term does not depend on $\zeta_e$, so centering it in $\zeta_e$ also gives zero. Thus, only $A=E(G)$ can contribute, and $c(E(G))=0$ implies $p_G=0$.

By row and column symmetry, each proportion is the number of copies of $A$ in one member of the family divided by the number of copies in $K_{4,4}$. Since every member has eight edges and degree two at each vertex, these ratios are
\[
\begin{array}{lc}
\text{Shape of }A&\Pr\{A\subseteq P\}=\Pr\{A\subseteq Q\}\\ \hline
\varnothing&1\\
\text{one cell}&8/16\\
\text{two cells sharing a vertex}&4/24\\
\text{two disjoint cells}&\bigl(\binom82-8\bigr)/72=20/72\\
\text{a path of three cells}&8/(16\cdot3^2)=8/144\\
\text{three cells sharing a vertex}&0.
\end{array}
\]
For the path count, each edge is the middle edge of exactly one such path. These are all edge sets on at most four vertices except $C_4$, so the corresponding coefficients $c(A)$ vanish.

In contrast, a fixed $C_4$ belongs to exactly one of the $18$ sets $P$ and to no set $Q$, so $c(C_4)=1/18$. Its vertices have degree two. Integrating over the two unused rows and the two unused columns gives $F_2^2G_2^2$. Thus, the normalization in Lemma~\ref{lem:decomp} gives
\[
p_{C_4}=\frac{(4)_2^2}{18}\lambda^4F_2^2G_2^2\prod_{e\in C_4}e_e
=8\lambda^4F_2^2G_2^2\prod_{e\in C_4}e_e.
\]
By conditional independence,
\[
\E[p_{C_4}^2]
=64\lambda^8F_2^4G_2^4\,
\E\!\left[\prod_{e\in C_4}\var(Y_e\mid\xi,\eta)\right],
\]
which is positive under the stated condition.
\end{proof}

Now, we prove Corollary~\ref{cor:d3}.

\begin{proof}[Proof of Corollary~\ref{cor:d3}]
We first identify the principal support graph, then compare the remaining graph variances with its variance. This gives a uniform bound on $D_N$ and allows us to apply Theorem~\ref{thm:contraction}. Proposition~\ref{prop:d3} and the assumption $\inf_{\Pcal_0}\|p_{C_4}\|_2>0$ show that the principal degree is four and that the unique principal support graph is the connected graph \(C_4\).

Every graph projection is a signed sum of conditional expectations of $h_\Psi$, so $\E[p_G^2]\le C\E[h_\Psi^2]\le C$ uniformly over $\Pcal_0$. By $\var(T_G)=w_G\E[p_G^2]$ with $w_G\asymp N^{-d(G)}$, the \(C_4\) class contributes at least $cN^{-4}$, with $c>0$ from $\inf\|p_{C_4}\|_2>0$, while every class of degree at least five contributes at most \(CN^{-5}\). The only degree-four projection with nonzero variance is indexed by a connected graph. Each disconnected graph therefore contributes variance of order at most \(N^{-5}\). Hence,
\[
D_N=\frac{\sum_{G\,\disc}\var(T_G)}{\sigma_N^2}\le CN^{-1}
\]
uniformly over $\Pcal_0$. Proposition~\ref{prop:gammaD} gives $\Upsilon_N\le CD_N^2\le CN^{-2}$, so Theorem~\ref{thm:contraction} gives $\widehat Z_N\dto\Nor(0,1)$ uniformly over $\Pcal_0$.
\end{proof}

\section{Additional results for Section~\ref{sec:fallback-experiment}}\label{sec:S-numerics}\label{sec:S-block-experiments}

In this section, we compute the null variance in the simulation model of Section~\ref{sec:fallback-experiment} and verify that the block procedure has asymptotic level at most $\alpha$ (Sec.~\ref{sec:S-model-hybrid}). We then derive the limits of $D_N$ in the three settings (Sec.~\ref{sec:S-model-limits}) and give the figures for every array size (Sec.~\ref{sec:S-figures}).

\subsection{Verifications for the procedure}\label{sec:S-model-hybrid}

Under the null hypothesis $a=0$ in Section~\ref{sec:fallback-experiment}, the simulated array is
\begin{equation}\label{eq:S-model-array}
Y_{ij}=\frac{t}{\sqrt r}\sum_{k=1}^rU_{ik}V_{jk}+e_{ij},
\qquad m=n\ge4,\quad r\in\{1,4\},\quad 0\le t\le2,
\end{equation}
with independent standard normal coordinates. We use the kernel $(Y_{ij}-Y_{kj})(Y_{il}-Y_{kl})$ and the complete statistic of Section~\ref{sec:fallback-experiment}. Put $\sigma_0^2=\var(U_N)$ for the null variance. The bounds below are uniform over these null laws.

The alternative diagnostic of Section~\ref{sec:S-chaos-procedure} uses an estimate of the variance contributed by $K_{0,2}$. With the complete averages $\Thetah(a,b)$ from \eqref{eq:kernel-product-averages}, define
\begin{equation}\label{eq:S-model-floor}
\widehat V_{0,2}
=\frac{\Thetah(0,2)-2\Thetah(0,1)+\Thetah(0,0)}{\binom n2},
\qquad F_N=(\widehat V_{0,2})^+.
\end{equation}

\begin{lemma}\label{lem:S-model-variance}\label{lem:S-model-floor-error}
For the null model \eqref{eq:S-model-array}, the following hold.
\begin{enumerate}
\item[(i)] The variance $V_{0,2}$ contributed by $K_{0,2}$ and the disconnected share satisfy
\begin{equation}\label{eq:S-model-variance}
\begin{gathered}
\sigma_0^2=\frac{8(nt^4/r+(1+t^2)^2)}{n(n-1)^2},\qquad
V_{0,2}=\frac{8t^4}{rn(n-1)},\\
D_N=1-\frac{n-1}{n(nt^4/r+(1+t^2)^2)}.
\end{gathered}
\end{equation}
Consequently,
\begin{equation}\label{eq:S-model-share-gap}
0\le D_N-\frac{V_{0,2}}{\sigma_0^2}
\le\frac{r+2}{n}+\sqrt{\frac rn}.
\end{equation}
\item[(ii)] Assumption (L) holds when the selected graph is $K_{0,2}$.
\item[(iii)] Assumption (M4) holds on this class of null laws.
\item[(iv)] The estimator $\widehat V_{0,2}$ is unbiased, and there is a constant $C$ independent of $n,r,t$ such that
\begin{equation}\label{eq:S-model-floor-error}
\frac{\E[(\widehat V_{0,2}-V_{0,2})^2]}{\sigma_0^4}\le\frac Cn.
\end{equation}
\end{enumerate}
\end{lemma}

Lemma~\ref{lem:S-model-variance}, proved below, verifies (M4) and (L), and gives $\sigma_0^2>0$. The selected graph is $K_{0,2}$, so $F_N=(\widehat V_{0,2})^+$. The procedure in Section~\ref{sec:fallback-experiment} uses $A_N=\max\{\Vu,F_N\}$ and the diagnostic \eqref{eq:S-floor-gate}. Since $m=n=L$, its choices satisfy
\[
K_N=\lfloor n^{1/3}\rfloor\to\infty,\qquad K_N=o(L),
\qquad \varepsilon_N=1/\lceil n^{1/4}\rceil\to0.
\]
Thus, Proposition~\ref{prop:S-alt-procedure}(i) gives an asymptotic rejection probability at most $\alpha$, uniformly over this null class.

Moreover, if $D_N=o(\varepsilon_N)$ uniformly on a null subclass, then the rejection probability tends to $\alpha$ uniformly on that subclass. Indeed, since $V_{0,2}\ge0$, \eqref{eq:S-model-floor-error} and \eqref{eq:S-model-share-gap} give
\[
\left|\frac{F_N}{\sigma_0^2}-D_N\right|
\le\frac{|\widehat V_{0,2}-V_{0,2}|}{\sigma_0^2}
+\left|\frac{V_{0,2}}{\sigma_0^2}-D_N\right|
=O_P(n^{-1/2}),
\]
uniformly on this subclass.
Consequently,
\[
\frac{F_N}{\varepsilon_N\sigma_0^2}
=\frac{D_N}{\varepsilon_N}
+O_P(n^{-1/4})
\longrightarrow0
\quad\text{in probability, uniformly}.
\]
Proposition~\ref{prop:S-alt-procedure}(i) therefore gives the claimed convergence to $\alpha$, which concludes our verification.

\begin{proof}[Proof of Lemma~\ref{lem:S-model-variance}]
\textbf{Proof of (i).} First, we compute $\sigma_0^2$. Choose an $n\times(n-1)$ matrix $B$ whose columns form an orthonormal basis of $\{x\in\mathbb R^n:\sum_i x_i=0\}$. Put $X_j=B^\top(Y_{1j},\ldots,Y_{nj})^\top$ and $Z=B^\top(U_{ik})_{i,k}$. Since $BB^\top=I_n-\mathbf1\mathbf1^\top/n$,
\begin{equation}\label{eq:S-model-row-pair}
\begin{aligned}
\sum_{i<i'}(Y_{ij}-Y_{i'j})(Y_{il}-Y_{i'l})
&=n\sum_iY_{ij}Y_{il}-\Bigl(\sum_iY_{ij}\Bigr)\Bigl(\sum_iY_{il}\Bigr)\\
&=nX_j^\top X_l,\\
U_N&=\frac{2}{(n-1)\binom n2}\sum_{j<l}X_j^\top X_l.
\end{aligned}
\end{equation}
Conditional on $(U_{ik})_{i,k}$, the vectors $X_j$ are independent centered Gaussian vectors with covariance
\begin{equation}\label{eq:S-model-covariance}
\Sigma=(t^2/r)ZZ^\top+I_{n-1}.
\end{equation}
Distinct column pairs have zero conditional covariance, and $\E[(X_j^\top X_l)^2\mid (U_{ik})_{i,k}]=\operatorname{tr}(\Sigma^2)$. Thus,
\[
\sigma_0^2=\frac{4\E[\operatorname{tr}(\Sigma^2)]}{(n-1)^2\binom n2}.
\]
The entries of $Z$ are independent standard normal variables, so
\[
\E\!\left[\Bigl(\sum_{k=1}^rZ_{ak}Z_{bk}\Bigr)^2\right]
=\begin{cases}r(r+2),&a=b,\\r,&a\ne b.\end{cases}
\]
Consequently,
\[
\begin{aligned}
\E[\operatorname{tr}(\Sigma^2)]
&=(n-1)(1+2t^2) +\frac{t^4}{r^2}\bigl\{(n-1)r(r+2)+(n-1)(n-2)r\bigr\}\\
&=(n-1)(nt^4/r+(1+t^2)^2),\\
\sigma_0^2&=\frac{8(nt^4/r+(1+t^2)^2)}{n(n-1)^2}.
\end{aligned}
\]

Next, we compute $V_{0,2}$ and $D_N$. For $W=(\{i,i'\},\{j,l\})$, put $\delta_j=(t/\sqrt r)\sum_{k=1}^r(U_{ik}-U_{i'k})V_{jk}$. Expanding the kernel gives
\begin{equation}\label{eq:S-model-kernel-expansion}
\begin{aligned}
H_W&=(\delta_j+e_{ij}-e_{i'j})(\delta_l+e_{il}-e_{i'l})\\
&=\delta_j\delta_l+\delta_j(e_{il}-e_{i'l})+\delta_l(e_{ij}-e_{i'j}) +e_{ij}e_{il}+e_{i'j}e_{i'l}-e_{ij}e_{i'l}-e_{i'j}e_{il}.
\end{aligned}
\end{equation}
Write $V_j=(V_{j1},\ldots,V_{jr})$ and $V_l=(V_{l1},\ldots,V_{lr})$. We have
\[
\begin{aligned}
\E[H_W\mid V_j,V_l]
&=\E[\delta_j\delta_l\mid V_j,V_l]\\
&=\frac{t^2}{r}\sum_{k,k'=1}^r V_{jk}V_{lk'}
\E[(U_{ik}-U_{i'k})(U_{ik'}-U_{i'k'})].
\end{aligned}
\]
Since $\E[H_W\mid V_j]=\E[H_W\mid V_l]=\E[H_W]=0$, this conditional expectation is the term $q_\Lambda$ of $H_W$ for the copy $\Lambda$ of $K_{0,2}$ on the columns $j$ and $l$. Thus, the formula for $V_{0,2}$ in Lemma~\ref{lem:overlap} gives $V_{0,2}=\E[q_\Lambda^2]/\binom n2$. Moreover, we have
\[
\E[(U_{ik}-U_{i'k})(U_{ik'}-U_{i'k'})]=2\one\{k=k'\}.
\]
Thus,
\[
V_{0,2}=\frac{1}{\binom n2}\E\!\left[\Bigl(\frac{2t^2}{r}\sum_{k=1}^rV_{jk}V_{lk}\Bigr)^2\right]
=\frac{4t^4}{r\binom n2}.
\]

Now, in~\eqref{eq:S-model-kernel-expansion}, only the products $e_{ij}e_{il}$ and $e_{i'j}e_{i'l}$ represent connected graphs. Each row belongs to $n-1$ row pairs, so their complete contribution and its variance are
\[
\begin{aligned}
U_{\mathrm{conn}}
&=\frac{n-1}{\binom n2^2}\sum_i\sum_{j<l}e_{ij}e_{il}
=\frac{2}{n\binom n2}\sum_i\sum_{j<l}e_{ij}e_{il},\\
\var(U_{\mathrm{conn}})&=\frac{4}{n\binom n2},\\
D_N&=1-\frac{\var(U_{\mathrm{conn}})}{\sigma_0^2}
=1-\frac{n-1}{n(nt^4/r+(1+t^2)^2)}.
\end{aligned}
\]

Finally, put $x=t^2\sqrt{n/r}$, so that $nt^4/r=x^2$ and $2t^2=\sqrt{r/n}\,2x$. Using $(1+t^2)^2\ge1$, we obtain
\[
\begin{aligned}
0\le D_N-\frac{V_{0,2}}{\sigma_0^2}
&=\frac{n(nt^4/r+(1+t^2)^2)-(n-1)-n(n-1)t^4/r}{n(nt^4/r+(1+t^2)^2)}\\
&=\frac{(r+1)t^4/r+2t^2+n^{-1}}{nt^4/r+(1+t^2)^2}\\
&\le\frac{r+1}{n}
+\sqrt{\frac rn}\frac{2x}{x^2+1}+\frac1n
\le\frac{r+2}{n}+\sqrt{\frac rn}.
\end{aligned}
\]

\textbf{Proof of (ii).} Take $\Gamma_{\mathrm{disc}}=\{K_{0,2}\}$, and let $\rho_N$ be as in Lemma~\ref{thm:S-test}. Then $\rho_N=D_N-V_{0,2}/\sigma_0^2$, so \eqref{eq:S-model-share-gap} gives $\rho_N\le(r+2)/n+\sqrt{r/n}$. Since $r\le4$, this bound tends to zero uniformly over these null laws, which proves (L)(a). Moreover, the connected components of $K_{0,2}$ are single column vertices, so (L)(b) holds.

\textbf{Proof of (iii).} To verify (M4), we use equivalence of norms on finite-dimensional polynomial spaces. We regard $H_W$ as a polynomial in the $2r$ row coordinates $U_{ik},U_{i'k}$, the $2r$ column coordinates $V_{jk},V_{lk}$, and the four cell noises $e_{ij},e_{il},e_{i'j},e_{i'l}$. These are $4r+4\le20$ independent standard normal variables. In~\eqref{eq:S-model-array}, each term has total degree at most two in these Gaussian variables, and their product of differences $H_W$ in~\eqref{eq:S-model-kernel-expansion} has total degree at most four. The same degree bound holds for each projection $p_G$, since it is a linear combination of conditional expectations of $H_W$ and integrating out independent coordinates replaces their powers by Gaussian moments. Thus, for each $r\in\{1,4\}$, the projections belong to the fixed finite-dimensional space of polynomials of total degree at most four in $4r+4$ standard normal variables. By equivalence of the $L^2$ and $L^4$ norms on this space, we have
\[
\E[p_G^4]\le C\bigl(\E[p_G^2]\bigr)^2,
\]
where $C$ is independent of $n,r,t$. This proves (M4).

\textbf{Proof of (iv).} Each $\Thetah(a,b)$ is unbiased for $\Theta(a,b)$, and Lemma~\ref{lem:overlap} gives $\Theta(0,2)-2\Theta(0,1)+\Theta(0,0)=\omega(0,2)=\binom n2V_{0,2}$. Thus, $\widehat V_{0,2}$ is unbiased. By \eqref{eq:S-model-floor} and \eqref{eq:kernel-product-averages},
\[
\begin{gathered}
Q:=\binom n2\widehat V_{0,2}
=\sum_{b=0}^2\frac{(-1)^b\binom2b}{|\Omega_b|}
\sum_{(W,W')\in\Omega_b}H_WH_{W'},\\
\Omega_b=\{(W,W'):o(W,W')=(0,b)\}.
\end{gathered}
\]
To bound $\var(Q)$, we expand each product $H_WH_{W'}$ with~\eqref{eq:S-model-kernel-expansion}. Up to sign, each resulting term is a product of four factors, each equal to some $\delta_j$ or some cell noise $e_{ij}$. Fix one such term with $k$ cell noises, and write $T_\omega$ for its value at $\omega=(W,W')\in\Omega_b$. Its squared $L^2$ norm is at most $Ct^{8-2k}$, since $r\le4$ and the Gaussian moments are bounded. Its average has variance
\[
\var\!\left(\frac1{|\Omega_b|}\sum_{\omega\in\Omega_b}T_\omega\right)
=\frac1{|\Omega_b|^2}\sum_{\omega,\omega'\in\Omega_b}
\cov(T_\omega,T_{\omega'}).
\]
The row pairs in $W,W'$ are disjoint, and each kernel uses distinct columns. Thus, every cell noise occurs once, and at most two cell noises share a row or a column. Consequently, a nonzero covariance requires that $T_\omega$ and $T_{\omega'}$ contain the same cell noises. When $k=0$, however, it requires at least a common row or a common column. Each placement uses four rows and $4-b$ columns, so $|\Omega_b|\asymp n^{8-b}$. Fixing $s$ row and column indices of $\omega'$ leaves at most $Cn^{8-b-s}$ choices. Therefore, the necessary matches give
\[
\begin{array}{clc}
k&\text{Required common vertices}&\text{Fraction of pairs}\\ \hline
0&\text{a row or a column}&O(n^{-1})\\
1&\text{one row and one column}&O(n^{-2})\\
2&\text{at least three vertices}&O(n^{-3})\\
3,4&\text{at least two rows and two columns}&O(n^{-4}).
\end{array}
\]
By the Cauchy--Schwarz inequality, each covariance is at most $Ct^{8-2k}$ in absolute value. Summing the finitely many terms in $Q$ gives
\[
\var(Q)\le C\left\{\frac{t^8}{n}+\frac{t^6}{n^2}
+\frac{t^4}{n^3}+\frac{t^2+1}{n^4}\right\}
\le\frac Cn(t^4+n^{-1})^2,
\]
where the last bound uses $t\le2$ and $t^6/n^2\le(t^8/n+t^4/n^3)/2$. Finally, part (i) and $r\le4$ give
\[
\begin{gathered}
\binom n2\sigma_0^2=\frac{4(nt^4/r+(1+t^2)^2)}{n-1}\ge t^4+n^{-1},\\
\frac{\E[(\widehat V_{0,2}-V_{0,2})^2]}{\sigma_0^4}
=\frac{\var(Q)}{[\binom n2\sigma_0^2]^2}\le\frac Cn.
\end{gathered}
\]
\end{proof}

\subsection{Limits of $D_N$}\label{sec:S-model-limits}

\begin{lemma}\label{lem:S-model-alternative}
Let $Y_{ij}=aA_i+Y_{ij}^0$, where $(Y_{ij}^0)$ follows the null model \eqref{eq:S-model-array} and $A_1,\dots,A_n$ are independent standard normal variables, independent of $(Y_{ij}^0)$. Put $\theta=2a^2=\vartheta\sigma_0$ with $\vartheta\ge0$ fixed, where $\sigma_0^2$ is as in Lemma~\ref{lem:S-model-variance}. In the settings (i), (ii) and (iii) of Section~\ref{sec:fallback-experiment}, the disconnected share of this law satisfies $D_N\to0$, $D_N\to4/5$ and $D_N\to1$, respectively, as $N\to\infty$.
\end{lemma}

\begin{proof}[Proof of Lemma~\ref{lem:S-model-alternative}]
Write $D_N(\theta)$ for the disconnected share at $\theta=2a^2$. We prove that $D_N(\vartheta\sigma_0)-D_N(0)\to0$. Use the matrix $B$ from the proof of Lemma~\ref{lem:S-model-variance}, and put $\alpha=B^\top(A_1,\ldots,A_n)^\top$, $X_j=B^\top(Y_{1j}^0,\ldots,Y_{nj}^0)^\top$ and $S=\sum_jX_j$. The row-pair identity~\eqref{eq:S-model-row-pair} gives
\[
\begin{aligned}
U_N
&=\frac{2}{(n-1)\binom n2}\sum_{j<l}(a\alpha+X_j)^\top(a\alpha+X_l)\\
&=\underbrace{\frac{2a^2\|\alpha\|^2}{n-1}}_{R}
+\underbrace{\frac{4a\,\alpha^\top S}{n(n-1)}}_{Q}+U_N^0,
\end{aligned}
\]
where $U_N^0$ is the statistic computed from the null array. Here, each column belongs to $n-1$ column pairs. Since $\alpha$ is standard normal in $\mathbb R^{n-1}$ and independent of the null array,
\[
\E[R]=\theta,\qquad
\var(R)=\frac{8a^4}{n-1}=\frac{2\theta^2}{n-1},\qquad
\E[Q\mid (Y_{ij}^0)]=0.
\]

Next, conditional on $(U_{ik})$, the vectors $X_j$ are independent and centered, with covariance $\Sigma$ given by~\eqref{eq:S-model-covariance}. This gives
\[
\begin{aligned}
\E\|S\|^2&=n\E[\operatorname{tr}\Sigma]=n(n-1)(1+t^2),\\
\var(Q)&=\frac{16a^2\E\|S\|^2}{n^2(n-1)^2}
=\frac{8(1+t^2)\theta}{n(n-1)}.
\end{aligned}
\]
Moreover, the variables $R$ and $U_N^0$ are independent, and $Q$ is odd in $\alpha$, whereas $R$ is even and $U_N^0$ does not depend on $\alpha$. Thus, the three terms are pairwise uncorrelated, and
\[
v:=\var(U_N)-\sigma_0^2
=\var(R+Q-\theta)
=\frac{2\theta^2}{n-1}+\frac{8(1+t^2)\theta}{n(n-1)}.
\]

Now, integrating out the independent variables $A_i$ gives
\[
\E[R+Q-\theta\mid (U_{ik}),(V_{jk}),(e_{ij})]=0.
\]
Graph projections are sums of conditional expectations over subsets of the underlying variables. These expectations commute with integration over the variables $A_i$, so the graph contributions $T_G$ and $T_G^0$ of the two arrays satisfy
\[
\E[T_G-T_G^0\mid (U_{ik}),(V_{jk}),(e_{ij})]=0,
\qquad \E[(T_G-T_G^0)T_G^0]=0.
\]
By orthogonality across graph classes,
\[
\begin{aligned}
\sum_G\var(T_G-T_G^0)&=v,\\
v_{\mathrm{disc}}:=\sum_{G\,\disc}\var(T_G-T_G^0)&\in[0,v],\\
D_N(\theta)-D_N(0)
&=\frac{D_N(0)\sigma_0^2+v_{\mathrm{disc}}}{\sigma_0^2+v}-D_N(0)
=\frac{v_{\mathrm{disc}}-D_N(0)v}{\sigma_0^2+v}.
\end{aligned}
\]
Since $0\le D_N(0)\le1$ and~\eqref{eq:S-model-variance} gives $\sigma_0\ge2\sqrt2\,(1+t^2)/\{\sqrt n(n-1)\}$, substituting $\theta=\vartheta\sigma_0$ yields
\[
|D_N(\vartheta\sigma_0)-D_N(0)|
\le\frac{v}{\sigma_0^2}
=\frac{2\vartheta^2}{n-1}+\frac{8\vartheta(1+t^2)}{n(n-1)\sigma_0}
\le\frac{2\vartheta^2}{n-1}+\frac{2\sqrt2\,\vartheta}{\sqrt n}\longrightarrow0.
\]
Finally, \eqref{eq:S-model-variance} gives
\[
D_N(0)=1-\frac{n-1}{n\{nt^4/r+(1+t^2)^2\}}
=\begin{cases}
1/n\to0,&\text{in (i)},\\
1-\dfrac{n-1}{n\{4+(1+4n^{-1/2})^2\}}\to4/5,&\text{in (ii)},\\
1-\dfrac{n-1}{n(n+4)}\to1,&\text{in (iii)}.
\end{cases}
\]
Thus, the same limits hold for $D_N(\vartheta\sigma_0)$.
\end{proof}

\subsection{Additional figures}\label{sec:S-figures}

Figure~\ref{fig:S-block-coverage} shows the coverage frequencies as $\vartheta$ increases, for all settings and sizes considered. It shows undercoverage at smaller array sizes. Figure~\ref{fig:S-block-power} shows the rejection frequencies as $\vartheta$ increases, for all settings and sizes considered. Without interaction, the block procedure and the Gaussian test exceed the target rejection probability at $\vartheta=0$ for $L=64,128$. Figure~\ref{fig:S-block-widths} shows the median interval half-widths as $\vartheta$ increases, for all settings and sizes considered. These can drop sharply when more replicates select the Gaussian test.

\begin{figure}[p]
\centering
\includegraphics[width=\textwidth,height=0.76\textheight,keepaspectratio]{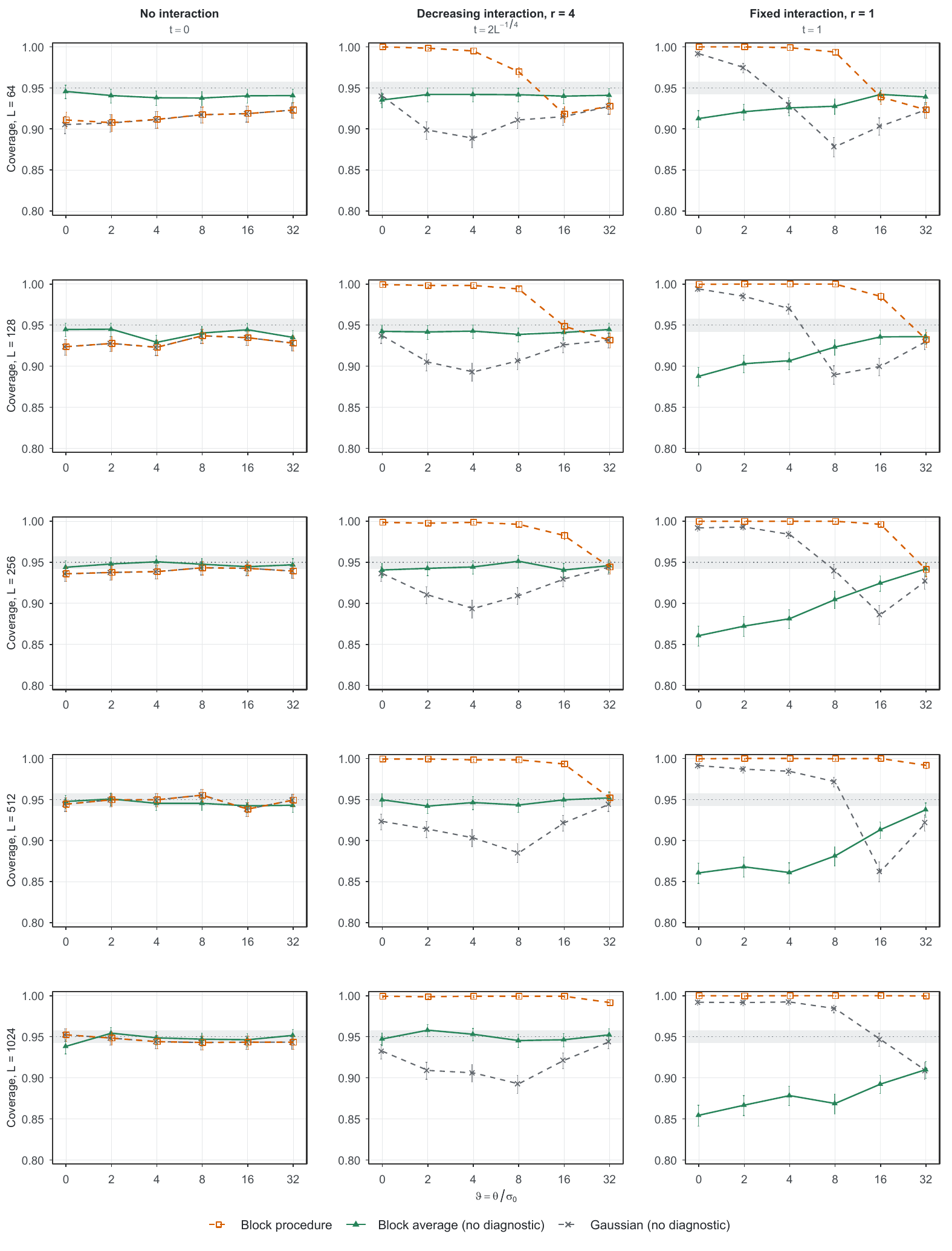}
\caption{Coverage of $95\%$ confidence intervals for the settings of
Figure~\ref{fig:block-chaos} (columns), the different sizes considered (rows), at different values of $\vartheta$ (x-axis, equally spaced). The bars are pointwise $95\%$ Wilson intervals. The dotted grey line is the target $95\%$ coverage. The shaded bands are the central $95\%$ ranges of $\operatorname{Bin}(3000,0.95)/3000$.}\label{fig:S-block-coverage}
\end{figure}

\begin{figure}[p]
\centering
\includegraphics[width=\textwidth,height=0.76\textheight,keepaspectratio]{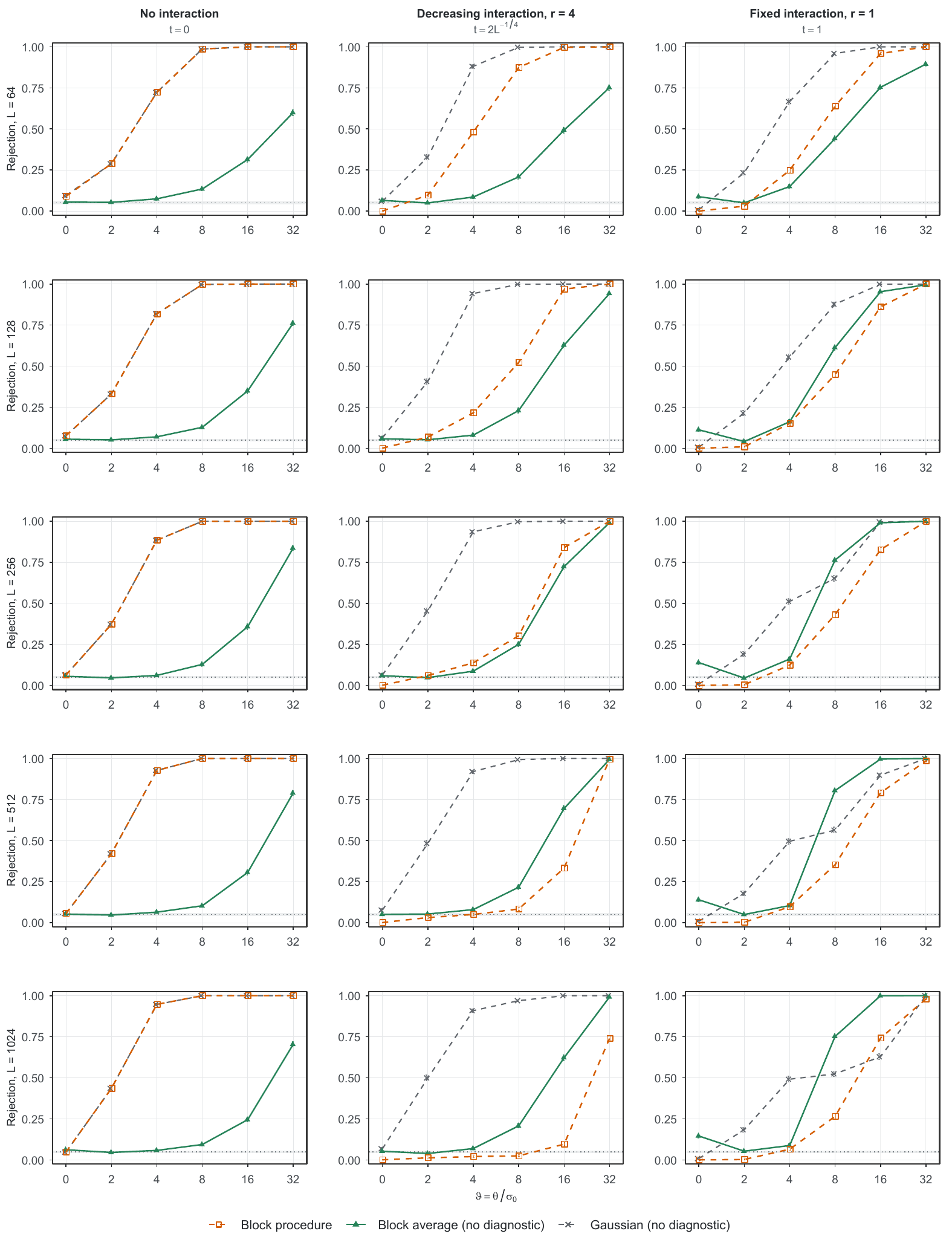}
\caption{Rejection frequencies for tests of zero for the settings of
Figure~\ref{fig:block-chaos} (columns), the different sizes considered (rows), at different values of $\vartheta$ (x-axis, equally spaced). The bars are pointwise $95\%$ Wilson intervals. The dotted grey line is the target $5\%$ rejection probability. The shaded bands are the central $95\%$ ranges of $\operatorname{Bin}(3000,0.05)/3000$.}\label{fig:S-block-power}
\end{figure}

\begin{figure}[p]
\centering
\includegraphics[width=\textwidth,height=0.76\textheight,keepaspectratio]{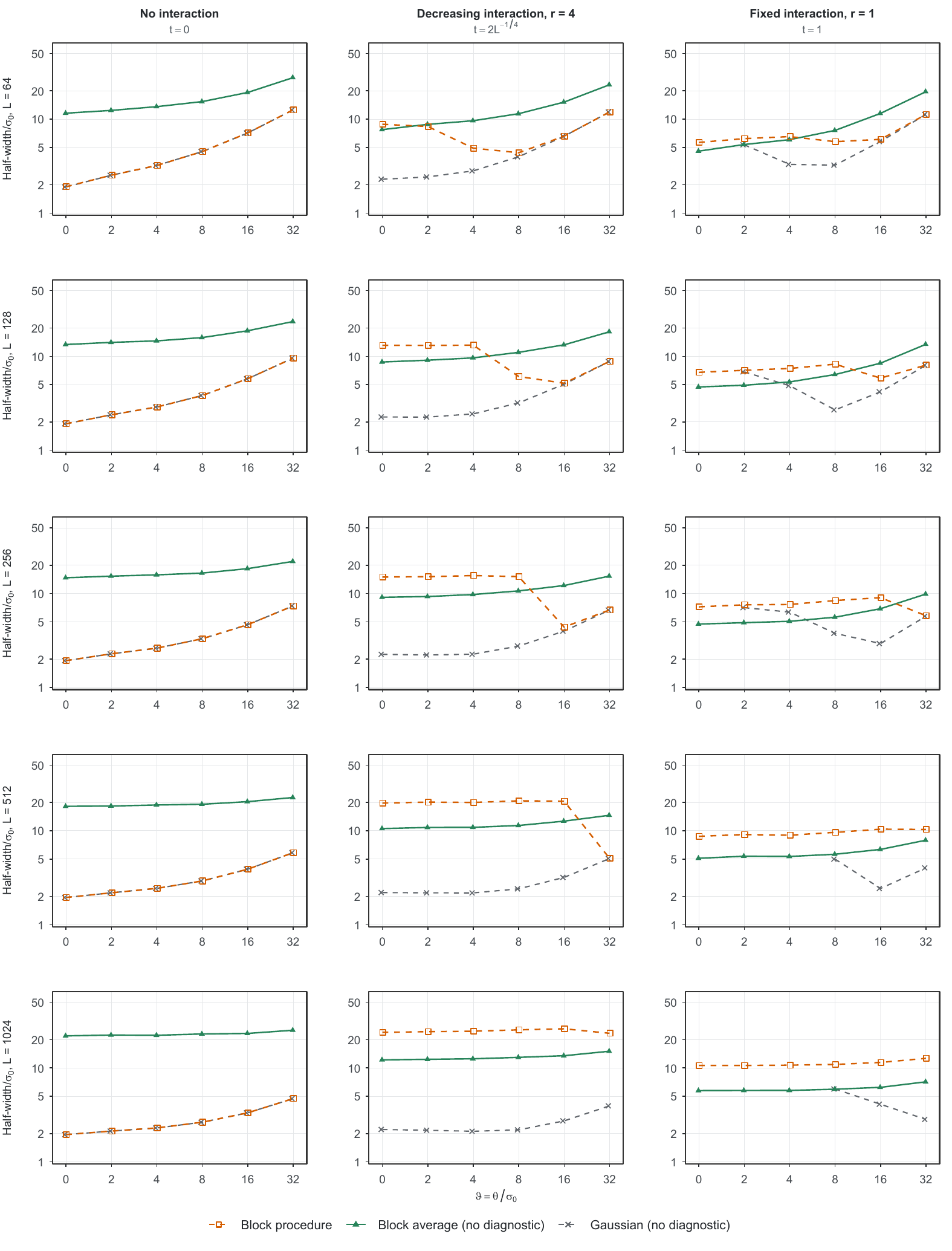}
\caption{Median confidence-interval half-widths divided by the null standard deviation $\sigma_0$ (log-scale), for the settings of Figure~\ref{fig:block-chaos} (columns), the different sizes considered (rows), at different values of $\vartheta$ (x-axis, equally spaced).}\label{fig:S-block-widths}
\end{figure}

\clearpage
\fi

\end{document}